\documentclass[a4paper,14pt]{extarticle}
\newcommand{\mbf}[1]{\protect\text{\boldmath$#1$}}
\newcommand{\comment}[1]{\begin{small}\noindent\texttt{// #1}\end{small}}

\usepackage[left=3cm,right=1.5cm,top=2cm,bottom=2cm,bindingoffset=0cm]{geometry}
\usepackage{setspace} 

\usepackage{cmap}					
\usepackage[T2A]{fontenc}			
\usepackage[utf8]{inputenc}			
\usepackage[english,russian]{babel}	

\usepackage{geometry} 
\usepackage{mathtext} 
\usepackage{mathrsfs} 
\usepackage{graphicx}
\usepackage{amsmath,amsthm,amsfonts,amssymb} 
\usepackage{icomma} 

\usepackage{array,tabularx,tabulary,booktabs} 
\usepackage{longtable}  
\usepackage{multirow} 
\usepackage{hhline}

\usepackage{graphicx}  
\graphicspath{{images/}}  
\usepackage{wrapfig} 
\usepackage{epstopdf} 

\usepackage{lscape} 

\usepackage[usenames]{color}
\usepackage{color}
\usepackage{colortbl}

\makeatletter
\renewcommand{\@biblabel}[1]{#1.}
\makeatother

\usepackage{tocloft}

\usepackage{titlesec} 
\titlelabel{\thetitle.\quad} 

\usepackage{indentfirst}

\numberwithin{equation}{section}
\numberwithin{figure}{section}
\numberwithin{table}{section}

\usepackage{soul}
\usepackage{csquotes}

\usepackage{pdfpages}

\begin{document}

\includepdf[pages=-]{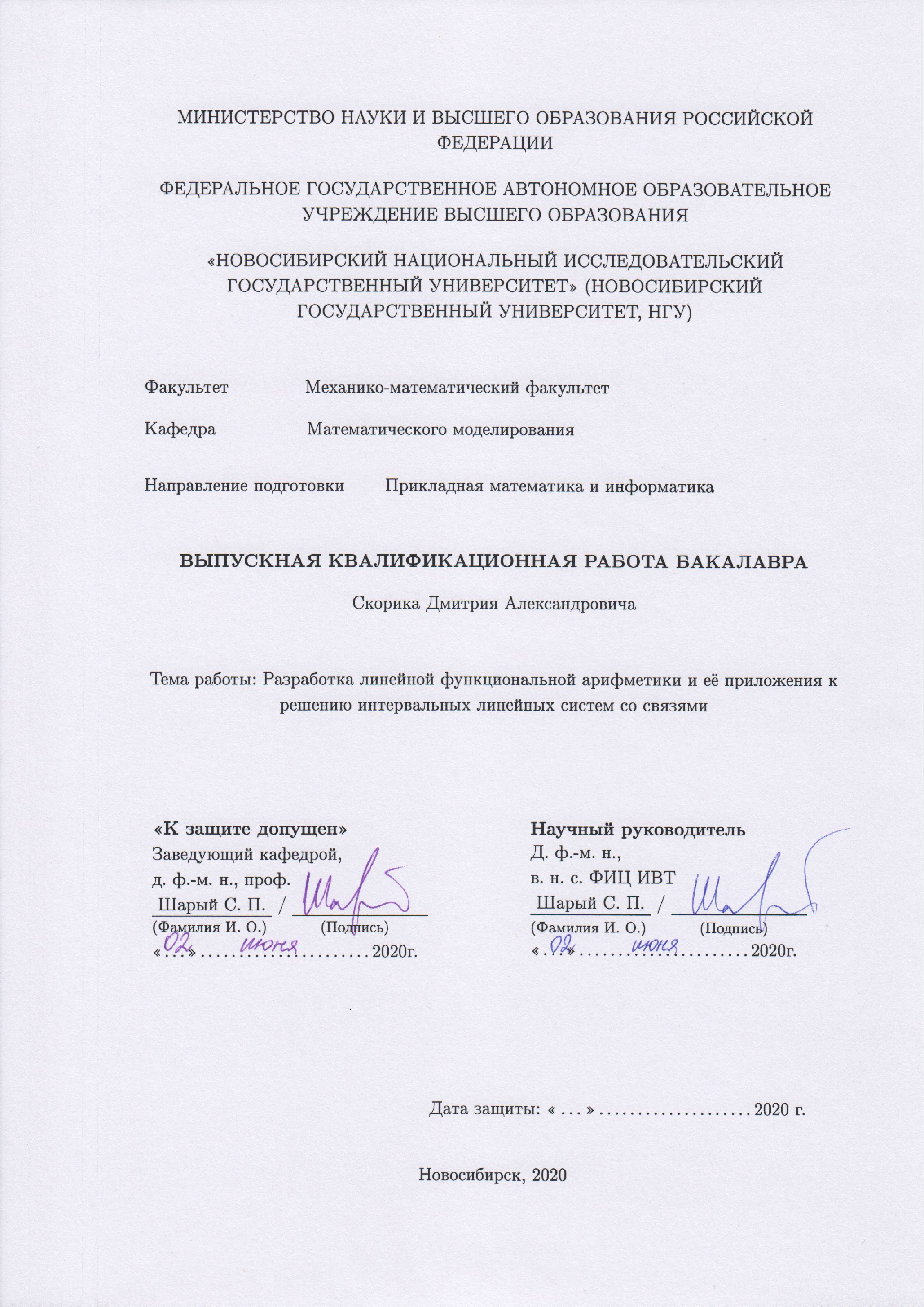}
\clearpage

\begin{titlepage}
\end{titlepage}

\begin{doublespace}

\setcounter{page}{2}
\clearpage
\tableofcontents

\clearpage
\section* {Реферат}
\addcontentsline{toc}{section}{Реферат}

\textbf{Тема работы:}
Разработка линейной функциональной арифметики и её приложения к решению интервальных линейных систем со связями.

\textbf{Объём работы} составляет 115 страниц, список использованной литературы включает 13 источников, в работе приводятся 18 таблиц и 65 рисунков.

\textbf{Ключевые слова:} интервал, аффинная арифметика, классическая интервальная арифметика, арифметика, СЛАУ, объединённое множество решений, интервальный анализ, ИСЛАУ.

Работа посвящена построению новой интервальной арифметики, которая сочетала бы в себе алгоритмическую эффективность и высокое качество оценивания областей значений выражений.

Разработана интервальная арифметика, основная идея которой близка к идее известной аффинной арифметики и заключается в приближении области совместных значений операндов. При этом сконструированная арифметика имеет особенности, повышающие качества получаемых интервальных оценок.

Проведены сравнительные эксперименты, которые показали эффективность построенной арифметики. Она была применена для решения задачи нахождения оценок объединённого множества решений ИСЛАУ, в том числе и со связанными параметрами.

\clearpage
\section {Введение}

Предметом текущей работы является задача нахождения оптимальной внешней оценки для решения интервальных систем линейных алгебраических уравнений (ИСЛАУ) вида:
\begin{equation}
\label{eq:fullIntervalMatrix}
    \left\{ 
    \begin{array}{ccccccccc}
        \mbf{a}_{11} x_{1} & + & \mbf{a}_{12} x_{2} & + & \dots & + & \mbf{a}_{1n} x_{n} & = & \mbf{b}_{1},\\
        
        \mbf{a}_{21} x_{1} & + & \mbf{a}_{22} x_{2} & + & \dots & + & \mbf{a}_{2n} x_{n} & = & \mbf{b}_{2},\\
        
        \vdots &  & \vdots & & \ddots &  & \vdots &  & \vdots \: \, \, ,\\
        
        \mbf{a}_{n1} x_{1} & + & \mbf{a}_{n2} x_{2} & + & \dots & + & \mbf{a}_{nn} x_{n} & = & \mbf{b}_{n}.
    \end{array} \right.
\end{equation}

\begin{center}
или
\end{center}
\begin{equation}
\label{eq:intervalMatrix}
    \mbf{A}x = \mbf{b}.
\end{equation}

Здесь $\mbf{A} = (\mbf{a}_{ij})$ --- интервальная $n \times n$ - матрица, а $\mbf{b} = (\mbf{b}_{i})$ --- $n$-вектор. Будем считать, что системы \eqref{eq:fullIntervalMatrix},  \eqref{eq:intervalMatrix} --- семейства классических точечных СЛАУ $Ax = b$ с элементами $a_{ij} \in \mbf{a}_{ij}$, $b_{i} \in \mbf{b}_{i}$. 

Также будем полагать, что матрица $\mbf{A}$ невырожденная, то есть неособенны все точечные матрицы $A$ с элементами $a_{ij} \in \mbf{a}_{ij}$. Таким образом, система $Ax = b$ имеет ограниченное решение независимо от выбора элементов $a_{ij}$ в заданных интервалах $\mbf{a}_{ij}$.

\subsection{Определения и обозначения}

Приведём определения из книги \cite{SharyIntervalBook}, используемые далее в работе.

\textit{Интервалом} будем называть замкнутое связное ограниченное подможество вещественной оси $\mathbb{R}$, то есть множество вида $\{ \, x \in \mathbb{R} \, | \, a \leq x \leq b \, \}$, где $a, \: b \in \mathbb{R}$ --- \textit{концы интервала}.

Согласно \cite{IntervalDefinitions} интервал будем обозначать латинскими буквами жирного шрифта $(\mbf{A}$, $\mbf{B}$, $\dots$, $\mbf{Y}$, $\mbf{Z}$, $\mbf{a}$, $\mbf{b}$, $\dots$, $\mbf{y}$, $\mbf{z})$. Левый (нижний) конец интервала будем обозначать нижним подчёркиванием символа интервала, а правый (верхний) --- верхним подчёркиванием. Семейство всех вещественных интервалов будем обозначать $\mathbb{I}\mathbb{R}$.

\textit{Интервальной величиной} или \textit{интервальным параметром} будем называть упорядоченную пару $(a, \: \mbf{a})$, где $a$ --- некоторая переменная, а $\mbf{a}$ --- интервал её изменения. Далее в работе интервальный параметр будет обозначаться как $a \in \mbf{a}$.

\begin{figure}
    \centering
        \includegraphics[width = 0.6 \linewidth]{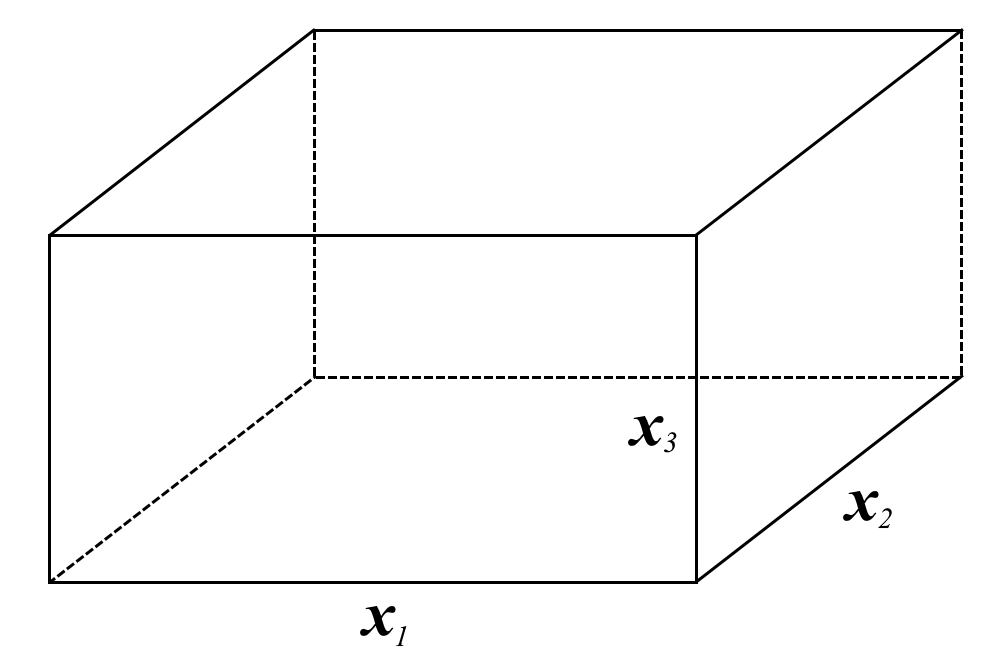}
        \caption{Геометрическая интерпретация интервального вектора --- брус}
    \label{fig:brus}
\end{figure}

\textit{Интервальным вектором} размера $n$ будем называть упорядоченный кортеж из $n$ интервалов, расположенный вертикально. Семейство интервальных векторов будем обозначать $\mathbb{I}\mathbb{R}^{n}$. В силу того, что геометрической интерпретацией интервального вектора является прямоугольный параллелепипед (рис. \ref{fig:brus}), будем также называть интервальный вектор \textit{брусом}.

\textit{Интервальной матрицей} будем называть прямоугольную таблицу, составленную из интервалов. 

Интервалы, входящие в интервальные векторы и матрицы, будем называть \textit{компонентами интервального вектора} и \textit{компонентами интервальной матрицы} соотвественно.

Будем говорить, что интервальная функция $\mbf{f} : \mathbb{I}\mathbb{R}^{n} \rightarrow \mathbb{I}\mathbb{R}^{m}$ является \textit{интервальным продолжением точечной функции} $f : \mathbb{R}^{n} \rightarrow \mathbb{R}^{m}$ на множестве $D \subseteq \mathbb{R}^{n}$, если $\mbf{f}(x) = f(x)$ для всех точечных аргументов $x \in D$.

Интервальная функция $\mbf{f} : \mathbb{I}\mathbb{R}^{n} \rightarrow \mathbb{I}\mathbb{R} ^ {m}$ называется \textit{интервальным расширением точечной функции} $f : \mathbb{R}^{n} \rightarrow \mathbb{R}^{m}$ \textit{на} $D \subseteq \mathbb{R}^{n}$, если:
\begin{enumerate}

	\item является интервальным продолжением $f$ на $D$
	
	\item из $\mbf{x} \subseteq \mbf{y}$ следует, что $\mbf{f}(\mbf{x}) \subseteq \mbf{f}(\mbf{y})$ для любых $\mbf{x}, \: \mbf{y} \in \mathbb{I}D$.

\end{enumerate} 

\textit{Элементарными функциональными выражениями} будем называть аналитические выражения, которые составлены из символов переменных, констант, четырёх арифметических операций — сложения, вычитания, умножения и деления — и элементарных функций.

\textit{Естественным интервальным
расширением} будем называть интервальное расширение элементарного функционального выражения, которое получается в результате замены его аргументов на интервалы их изменения, а арифметических операций и элементарных функций — на
их интервальные аналоги и расширения.

\subsection{Характеристики интервалов}

Напомним основные характеристики интервалов.

Пусть в рассмотрении имеется интервал $\mbf{a} =  \big[ \, \underline{\mbf{a}}, \: \overline{\mbf{a}} \, \big]$. Тогда:

\begin{enumerate}

    \item Ширина интервала:
    
    $\text{wid} \, \mbf{a} = \overline{\mbf{a}} - \underline{\mbf{a}}$
    
    \item Радиус интервала:
    
    $\text{rad} \, \mbf{a} = \frac{1}{2} \cdot \text{wid} \, \mbf{a} $
    
    \item Середина интервала:
    
    $\text{mid} \, \mbf{a} = \frac{1}{2} \cdot (\underline{\mbf{a}} + \overline{\mbf{a}})$
    
    \item Мигнитуда интервала:
    
    $\langle \mbf{a} \rangle = \left\{
    \begin{array}{rl}
    
    	\text{min} \big\{ \, |\underline{\mbf{a}}|, \: |\overline{\mbf{a}}| \, \big\}, & \text{если }0 \notin \mbf{a}, \\
    	
    	0, & \text{иначе.} 
    	
    \end{array}\right.$   
    
    \item Магнитуда интервала:
    
    $|\mbf{a}| = \text{max} \big\{ \, |\underline{\mbf{a}}|, \: |\overline{\mbf{a}}| \, \big\}$
    
    \item Отклонение интервала:
    
    $\text{dev} \, \mbf{a} = \left\{ 
    \begin{array}{rl}
    
        \underline{\mbf{a}} \text{,} & \text{если } |\underline{\mbf{a}}| \geq |\overline{\mbf{a}}| \text{,}\\
        
        \overline{\mbf{a}} \text{,} & \text{иначе.}
        
    \end{array} \right.$
    
\end{enumerate}

Исходя из постановок практических задач, исследованию подвергаются различные множества решений ИСЛАУ. Самыми распространнёными являются:

\begin{enumerate}

    \item Объединённое множество решений:
    
    \begin{center}
		$\varXi_{uni} ( \mbf{A}, \: \mbf{b} ) = \big\{ \,  x \in \mathbb{R}^{n} \, | \, ( \exists A \in \mbf{A} ) \: ( \exists b \in \mbf{b} ) \: ( Ax = b ) \, \big\},$
	\end{center}

    \item Допусковое множество решений:
    
    \begin{center}
		$\varXi_{tol} ( \mbf{A}, \: \mbf{b} ) = \big\{ \, x \in \mathbb{R}^{n} \, | \, ( \forall A \in \mbf{A} ) \: ( \exists b \in \mbf{b} ) \: ( Ax = b ) \, \big\},$
	\end{center}

    \item Управляемое множество решений:
    
    \begin{center}
		$\varXi_{ctl} ( \mbf{A}, \: \mbf{b} ) = \big\{ \, x \in \mathbb{R}^{n} \, | \, ( \forall b \in \mbf{b} ) \: ( \exists A \in \mbf{A} ) \: ( Ax = b ) \, \big\}.$
	\end{center}

\end{enumerate}

\subsection[Проблема оценки множеств решений ИСЛАУ]{Проблема оценки множеств решений ИСЛАУ}

Одна из популярных задач в интервальном анализе --- нахождение внешних оценок для множества решений ИСЛАУ. Рассмотрим в контексте данной задачи объединённое множество решений.

Если необходимо найти координатные границы точного бруса, описанного около этого множества, то возникает проблема, описанная в хрестоматийной книге по трудоёмкости интервальных вычислений \cite{NPHard}. 

Дело в том, что нахождение точных границ объединённого множества решений является $NP$-трудной задачей. Следовательно, мы неизбежно столкнёмся с экспоненциальным ростом алгоритмической сложности.

Так, независимо друг от друга Бекком \cite{Beeck} и Никелем \cite{Nickel} был рассмотрен алгоритм нахождения точных координатных границ бруса. Трудоёмкость такой процедуры растёт как $O \big( 2 ^ {n ^ {2} + n} \big)$, так как задача сводится к перебору и решению всех <<точечных>> СЛАУ, элементы которых являются границами компонентов ИСЛАУ.

Таким образом, уже при небольшом для практики размере матрицы $n = 10$, получаем характерное число операций для нахождения граней бруса $S = 2 ^ {10 ^ 2 + 10} = 2 ^ {110} \approx 1.3 \, \cdot \, 10^{33}$. Или, по-другому, если ЭВМ выполняет $10 ^ {9}$ операций в секунду, то время работы алгоритма будет составлять примерно $4.2 \, \cdot \, 10 ^ {16}$ лет, что превышает время жизни Вселленной на несколько порядков.

\begin{figure}
    \centering
        \includegraphics[width = 0.6 \linewidth]{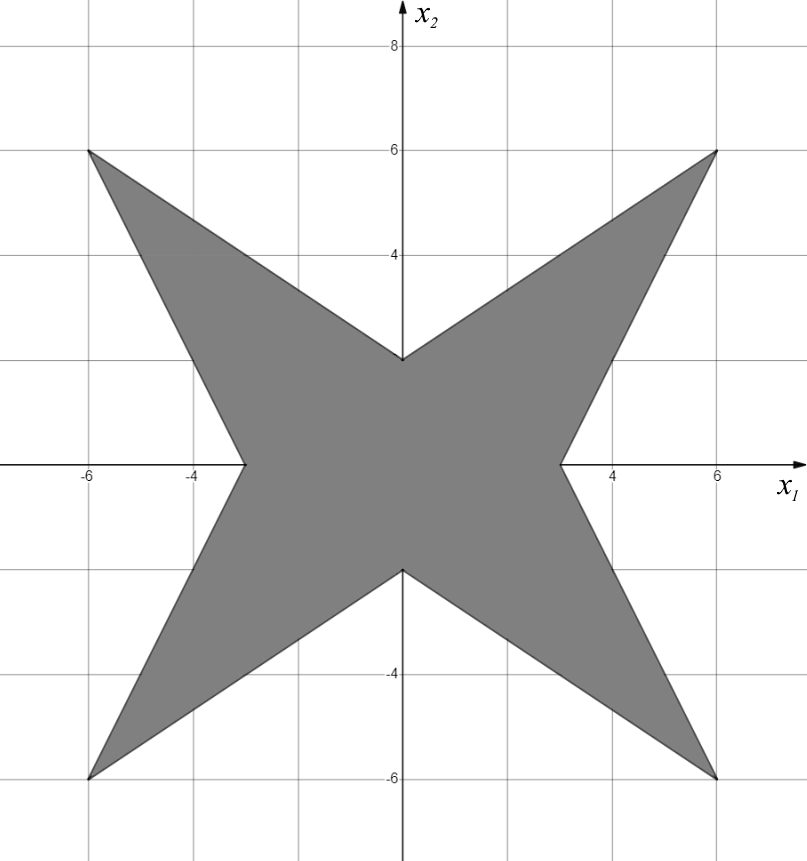}
        \caption{Геометрическая интерпретация объединённого множества решений ИСЛАУ $2 \times 2$ --- полиэдр, выпуклый в каждом ортанте пространства.}
    \label{fig:union_set}
\end{figure}

Первым шагом снижения трудоёмкости стало использование независимости интервальных величин, входящих в ИСЛАУ. В книге \cite{SharyIntervalBook} доказывается, что в таком случае объединённое множество решений представляет собой полиэдр, выпуклый в каждом ортанте пространства (рис. \ref{fig:union_set}). 

Данный факт позволяет получить асимптотику решения $P(n) \cdot O(2^{n})$, где $P(n)$ --- полином от переменной $n$. Это достигается за счёт того, что исходная задача разбивается на $2 ^ {n}$ (по числу ортантов) независимых задач линейного программирования.

Следующим шагом снижения трудоёмкости стало ослабление постановки задачи с нахождения точных границ объединённого множества решений на нахождение лишь его приближённой оценки. Это позволило преодолеть ограничение на экспоненциальный рост трудоёмкости и создать алгоритмы, работающие за полиномиальное время.

В интервальном анализе были разработаны различные методы получения оценок объединённого множества с такой асимптотикой:

\begin{enumerate}

    \item Интервальный алгоритм Гаусса.
    
    \item Интервальный алгоритм Холесского.
    
    \item Интервальный алгоритм Гаусса-Зейделя.
    
    \item \dots
    
\end{enumerate}

Однако при использовании данных методов возникают проблемы. В частности, ни один метод не может похвастаться своей универсальностью --- для каждого из приведённых подходов существуют примеры таких \mbox{ИСЛАУ}, в которых он превосходит остальные в качестве оценивания, а также --- которые делают его неприменимым для решения. 

Также все методы, использующие для получения результата интервальную арифметику, имеют общий недостаток --- при прочих равных условиях, чем больше размер рассматриваемой ИСЛАУ, тем грубее будет получаться оценка для вектора неизвестных. Данное явление связано с эффектом зависимости интервальных величин.

\subsection[Зависимые интервальные величины]{Зависимые интервальные величины}

Назовём интервальные величины $a_{1} \in \mbf{a}_{1}, \, \dots, \, a_{n} \in \mbf{a}_{n}$ \textit{независимыми}, если упорядоченный набор $(a_{1}, \, \dots, \, a_{n})$ принимает все значения из декартова произведения $\mbf{a}_{1} \times \dots \times \mbf{a}_{n}$. В противном случае будем говорить, что интервальные величины \textit{зависимые}.

\textit{Множеством совместных значений} интервальных параметров $a_{1} \in \mbf{a}_{1}$, $\dots$, $a_{n} \in \mbf{a}_{n}$ будем называть множество всех кортежей $(x_{1}, \dots, x_{n})$, где 
\begin{center}

	$x_{1} = a_{1}, \quad \dots, \quad x_{n} = a_{n}$.

\end{center}

Будем говорить, что на интервальные параметры $a_{1} \in \mbf{a}_{1}$, $\dots$, $a_{n} \in \mbf{a}_{n}$ \textit{наложены связи}, если имеются какие-либо соотношения между $a_{1}$, $\dots$, $a_{n}$, которые ограничивают множество их совместных значений.

Пусть в рассмотрении имеются интервальные параметры $a_{1} \in \mbf{a}$, $\dots$, $a_{n} \in \mbf{a}_{n}$. Тогда \textit{диаграммой зависимости} будем называть упорядоченную пару $(A, \: B),$ где $A$ --- множество совместных значений $a_{1} \in \mbf{a}_{1}$, $\dots$, $a_{n} \in \mbf{a}_{n}$, а $B$ --- декартово произведение $\mbf{a}_{1} \times \dots \times\mbf{a}_{n}$.

Таким образом, множество совместных значений диаграммы зависимости интервальных параметров вложено в декартово произведение интервалов изменений этих параметров. Причём равенство достигается в случае, когда эти интервальные параметры независимы.

В качестве примера зависимых интервальных величин рассмотрим кривошипно-шатунный механизм колеса паровоза. Поршень толкает шатун на величину $\Delta x$, при этом происходит движение кривошипа на величину $\Delta y$. В силу жесткой сцепки данных элементов будет естественно ожидать, что $\Delta y = f(\Delta x)$, где $f$ --- некоторая возрастающая или убывающая функция (в зависимости от выбора $\Delta x$). Таким образом, значение одной физической величины является аргументом функции для нахождения другой.

Пусть изначально система находилась в состоянии
\begin{center}

	$x = x_{0}$, \qquad $y = f(x_{0})$.

\end{center}
После смещения шатуна на величину $\Delta x$ система перешла в состояние
\begin{center}

	$x = x_{0} + \Delta x$, \qquad $y = f(x_{0} + \Delta x).$

\end{center} 

Таким образом, величины
\begin{center}

	$\mbf{x} = \big[ \, x_{0}, \: x_{0} + \Delta x \, \big]$ \: и \: $\mbf{y} = \big[ \, f(x_{0}), \: f(x_{0} + \Delta x) \, \big]$ 
	
\end{center}
зависимые, поскольку в декартовом произведении $\mbf{x} \times \mbf{y}$ не достигаются точки
\begin{center}

	$\big( x_{0}, \:f (x_{0} + \Delta x) \big)$ \: и \: $\big( x_{0} + \Delta x, \: f(x_{0}) \big)$.

\end{center}

\begin{figure}
    \centering
        \includegraphics[width = 0.5 \linewidth]{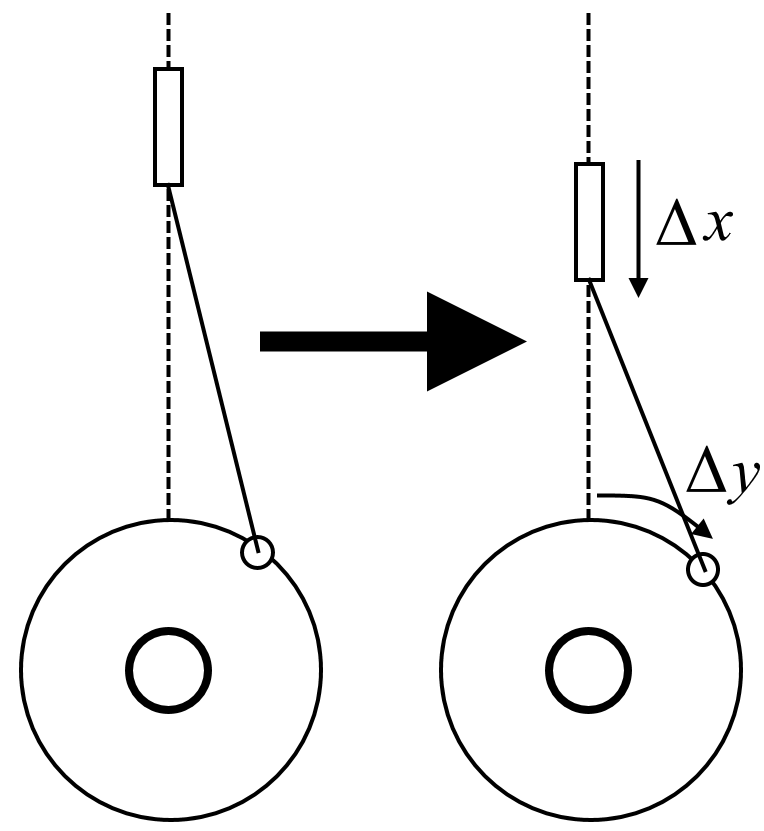}
        \caption{Схема кривошипно-шатунного механизма}
    \label{fig:shatun}
\end{figure}

Приведём пример более сложной зависимости. Пусть имеется световой источник, подвешенный на высоте $h$. Также известна мощность светового потока $n$ и масса $m$ этого источника. Рассмотрим величину интенсивности светового потока от источника $I = f_{1}(h, \: n)$, а также величину потенциальной энергии $U = f_{2}(h, \: m)$, где $f_{1}$, $f_{2}$ --- некоторые функции.

Зафиксируем значения параметров $m = m^{*}$ и $n = n^{*}$. Также заметим, что функция $I = f_{1}(h)$ --- убывающая, а функция $U = f_{2}(h)$ --- возрастающая. Тогда величины $I$ и $H$ --- зависимые, поскольку точка, где одновременно достигается максимум этих двух величин, не входит в декартово произведение $I(\mbf{h}) \times U(\mbf{h})$.

Рассмотрим подробнее свойства зависимых интервальных величин. Итак, пусть в рассмотрении имеются интервалы $\mbf{a} \in \mathbb{I}\mathbb{R}$ и $\mbf{b} \in \mathbb{I}\mathbb{R}$, тогда:

\begin{enumerate}

    \item $\text{rad}(\mbf{a} + \mbf{b}) = \text{rad} \, \mbf{a} + \text{rad} \, \mbf{b}.$
    
    \item $\text{rad} (\mbf{a} - \mbf{b}) = \text{rad} \, \mbf{a} + \text{rad} \, \mbf{b}.$
    
    \item $\text{rad} (\mbf{a} \cdot \mbf{b}) \leq |\mbf{a}| \cdot \text{rad} \, \mbf{b} + \text{rad} \, \mbf{a} \cdot |\mbf{b}|.$
    
    \item $\text{rad} (\mbf{a} \cdot \mbf{b}) \geq \text{max} \big\{ \, |\mbf{a}| \cdot \text{rad} \, \mbf{b}, \: \text{rad} \, \mbf{a} \cdot |\mbf{b}| \, \big\}.$
    
    \item $\text{rad} \big( 1 \, / \, \mbf{a} \big) = \frac{\mbox{rad} \, \mbox{\mbf{a}}}{\langle \mbox{\mbf{a}} \rangle \, \cdot \, |\mbox{\mbf{a}}|}$, если $0 \notin \mbf{a}.$
    
\end{enumerate}

Первые два свойства радиусов говорят о том, что при сложении или вычитании интервалов получается интервал, ширина которого равна сумме ширин операндов. Однако то, что результирующий интервал в общем случае расширяется в случаях умножения и деления, не очевидно.

Для демонстрации данного факта будем оперировать не абсолютной характеристикой интервала --- его шириной, а воспользуемся величиной, характеризующую относительную ширину интервала  --- так называемым функционалом Рачека. Данный функционал введён в работе \cite{Rachek} следующим образом:

\begin{center}
    $\chi(\mbf{a}) := 
    \left\{\begin{array}{cl}
    
        \underline{\mbf{a}} \, / \, \overline{\mbf{a}}, & \text{если } |\overline{\mbf{a}}| \geq |\underline{\mbf{a}}|,\\
        
        \overline{\mbf{a}} \, / \, \underline{\mbf{a}}, & \text{иначе.}
        
    \end{array}\right.$
\end{center}

С помощью данного функционала можно представить любой интервал $\mbf{a} \in \mathbb{I}\mathbb{R}$ в виде:
\begin{center}

	$\mbf{a} = \text{dev}(\mbf{a}) \cdot \big[ \, \chi(\mbf{a}), \: 1 \, \big].$
	
\end{center}

В работе Ирины Шарой \cite{IrinaSharaya} доказывается следующая лемма.

\noindent\textbf{Лемма.}
Пусть $\mbf{a} \in \mathbb{I}\mathbb{R}$ и $\mbf{b} \in \mathbb{I}\mathbb{R}$, тогда
\begin{center}

   $\chi(\mbf{a} \cdot \mbf{b}) = \text{min} \big\{ \, \chi(\mbf{a}), \: \chi(\mbf{b}), \: \chi(\mbf{a}) \cdot \chi(\mbf{b}) \, \big\}.$
   
\end{center}

\noindent\textbf{Доказательство}

Пусть $\mbf{a}$ и $\mbf{b}$ --- ненулевые интервалы.

Тогда $\mbf{a} \cdot \mbf{b}$ --- также ненулевой интервал и верно

\begin{center}
	$\begin{array}{lll}
	
		\mbf{a} \cdot \mbf{b} & = & \text{dev}(\mbf{a}) \cdot \big[ \, \chi(\mbf{a}), \: 1 \, \big] \cdot \text{dev}(\mbf{b}) \cdot \big[ \, \chi(\mbf{b}), \: 1 \, \big] = \\
	
		& = & \text{dev}(\mbf{a}) \cdot \text{dev}(\mbf{b}) \cdot \big[ \, \chi(\mbf{a}), \: 1 \big] \cdot \big[ \, \chi(\mbf{b}), \: 1 \, \big] = \\ 
	
		& = & \text{dev}(\mbf{a}) \cdot \text{dev}(\mbf{b}) \cdot \Big[ \, \text{min} \big\{ \, \chi(\mbf{a}), \: \chi(\mbf{b}), \: \chi(\mbf{a}) \cdot \chi(\mbf{b}), \: 1 \, \big\}, \\
	
		& & \hspace{4cm} \text{max} \big\{ \, \chi(\mbf{a}), \: \chi(\mbf{b}), \: \chi(\mbf{a}) \cdot \chi(\mbf{b}), \: 1 \, \big\} \, \Big] = \\
	
		& = & \text{dev}(\mbf{a}) \cdot \text{dev}(\mbf{b}) \cdot \Big[ \, \text{min} \big\{ \, \chi(\mbf{a}), \: \chi(\mbf{b}), \: \chi(\mbf{a}) \cdot \chi(\mbf{b}) \, \big\}, \: 1 \, \Big].
	
	\end{array}$
\end{center}

Поскольку 
\begin{center}
	$\Big| \text{min} \big\{ \, \chi(\mbf{a}), \: \chi(\mbf{b}), \: \chi(\mbf{a}) \cdot \chi(\mbf{b}) \, \big\} \Big| \leq 1$
\end{center}
\noindentи любой интервал $\mbf{z} \in \mathbb{I}\mathbb{R}$ представляется в виде
\begin{center}

	$\mbf{z} = \text{dev}(\mbf{z}) \cdot \big[ \, \chi(\mbf{z}), \: 1 \, \big],$
	
\end{center}
\noindentто
\begin{center}

	$\chi(\mbf{a} \cdot \mbf{b}) = \text{min} \big\{ \, \chi(\mbf{a}), \: \chi(\mbf{b}), \: \chi(\mbf{a}) \cdot \chi(\mbf{b})\, \big\}.$
	
\end{center}

$\hfill \blacksquare$

В работе $\cite{IrinaSharaya}$ также доказано свойство

\begin{center}
	
	$\text{wid}(\mbf{a} \cdot \mbf{b}) = |\mbf{a}| \cdot |\mbf{b}| \cdot \big( 1 - \chi(\mbf{a} \cdot \mbf{b}) \big).$	
	
\end{center}

$\hfill \blacksquare$

Таким образом, в силу доказанной леммы и свойства получаем, что
\begin{center}

	$\text{wid}(\mbf{a} \cdot \mbf{b}) = |\mbf{a}| \cdot |\mbf{b}| \cdot \Big( 1 - \text{min} \big\{ \, \chi(\mbf{a}), \, \chi(\mbf{b}), \, \chi(\mbf{a}) \cdot \chi(\mbf{b}) \, \big\} \Big).$
	
\end{center}

Используем данный результат. Зафиксируем интервал $\mbf{a}$. Рассмотрим, при каких условиях после умножения на интервал $\mbf{b}$ ширина результирующего интервала не будет увеличиваться.

Это возможно при одновременном выполнении двух условий:
\begin{enumerate}
	
	\item $|\mbf{b}| \leq 1$.
	
	То есть, $\mbf{b} \subseteq [ \, -1, \: 1 \, ]$.
		
	\item $1 - \text{min} \big\{ \, \chi(\mbf{a}), \, \chi(\mbf{b}), \, \chi(\mbf{a}) \cdot \chi(\mbf{b}) \, \big\} \leq 1.$
	
	Перепишем данное неравенство в виде	
	\begin{center}
	
		$ \text{min} \big\{ \, \chi(\mbf{a}), \, \chi(\mbf{b}), \, \chi(\mbf{a}) \cdot \chi(\mbf{b}) \, \big\} \geq 0$.
		
	\end{center}
	
	\begin{enumerate}
	
		\item Если $\chi(\mbf{a}) < 0$, то поставленное неравенство не может быть удовлетворено.
	
		\item Если $\chi(\mbf{a}) \geq 0$, то необходимо, чтобы $\chi(\mbf{b}) \geq 0$.
		
	\end{enumerate}
		
	Теперь рассмотрим введённый ранее функционал Рачека для интервала $\mbf{b}$:
	
	\begin{enumerate}
	
    	\item Если $\chi(\mbf{b}) = 0$, то $\underline{\mbf{b}} = 0$ или $\overline{\mbf{b}} = 0.$
    
    	\item Если $\chi(\mbf{b}) > 0$, то $0 \notin \mbf{b}.$
    
    \end{enumerate}	
    
    То есть получаем,что $\mbf{b} \subseteq ( \, -\infty, \: 0 \, ]$ или $\mbf{b} \subseteq [ \, 0, \: +\infty \, )$.
	
\end{enumerate}

Таким образом, чтобы
\begin{center}

	$\text{wid}(\mbf{a}) \geq \text{wid}(\mbf{a} \cdot \mbf{b})$ 
	
\end{center}
необходимо
\begin{center}

	$\mbf{b} \subseteq [ \, -1, \: 0 \,] \quad \text{или} \quad \mbf{b} \subseteq [ \, 0, \: 1 \, ].$
	
\end{center}

Получаем, что в общем случае, когда умножение происходит на интервалы, которые не удовлетворяют указанным условиям, ширина результирующего интервала будет увеличиваться.

Для демонстрации эффекта непременного расширения интервалов при вычислениях проведём следующий численный эксперимент. Допустим, стоит задача, где при вычислениях используются величины, значения которых получены с некоторой погрешностью. 

В условиях данного эксперимента будем считать, что неточность каждого измерения (ширина интервала, в котором гарантированно содержится измеряемая величина) --- это случайная величина, которая распределена, как модуль стандартной нормальной случайной величины, а каждое измерение (середина интервала) --- случайная величина, равномерно распределённая на некотором интервале.

Проделаем следующее:

\begin{enumerate}

    \item Создадим файл, куда будем записывать результаты эксперимента.
    
    \item Выберем длину цепочки вычислений $S \in \mathbb{Z}_{+}$.
    
    \item Выберем количество интервальных величин, которые будут участвовать в расчетах $n \in \mathbb{N}_{+}$ .
    
    \item Сгенерируем <<измеренные в эксперименте величины>> --- множество интервалов
    \begin{center}
    
    	$M = \big\{ \, \mbf{a}_{1}, \: \dots, \: \mbf{a}_{S} \, \big\}.$
    	
    \end{center}
    
    Генерацию будем производить случайным образом. Для получения ширины интервалов будем использовать нормальное распределение с параметрами
    \begin{center}
    
    	$\alpha = 0, \qquad \sigma ^ {2} = 0.01,$
    	
    \end{center}
    а для получения середины интервала --- равномерное распределение на интервале
    \begin{center}
    
    	$[ \, -5, \: 5 \, ].$
    	
    \end{center}
    
    При этом будем помещать в множество $M$ только те интервалы, которые заведомо не содержат в себе нуль. Это обусловлено тем, что мы не хотим исключать из рассмотрения операцию деления на интервал. Данное требование не противоречит практике, поскольку величины, малоотличимые от нуля в характерных для задачи единицах измерения, зачастую не имеют физического смысла.
    
    \item Инициализируем вычисляемый интервал из множества $M$, выбор будем производить равновероятно случайно.
    
    \item $S$ раз преобразуем вычисляемый интервал следующим образом:
    \begin{enumerate}
    
        \item Выберем арифметическую операцию из множества 
        \begin{center}
        
        	$\big\{ \, +, \: -, \: \cdot, \: / \, \big\},$
        	
        \end{center}
        выбор будем производить случайно равновероятно.
        
        \item Выберем интервал из множества $M$, выбор будем производить случайно равновероятно. 
        
        Выбранный интервал будет выступать как операнд нашей арифметической операции, которая будет далее применяться к вычисляемому интервалу.
        
        \item Применим к вычисляемому интервалу арифметическую операцию, выбранную в пункте (a) с операндом, выбранным в пункте (b).
        
    \end{enumerate}
    
    \item После получения результирующего интервала, запишем в файл (см. п. 1) его ширину.
    
\end{enumerate}

Для выбранных чисел $S$ и $n$ будем запускать данный алгоритм большое число раз, чтобы получить репрезентативную выборку для анализа ширины итогового интервала.

В ходе эксперимента величины $S$ и $n$ варьировались в следующих интервалах:
\begin{center}

	$S \in \big\{ \, 1, \: \dots \: , \: 6 \, \big\}, \qquad n \in \big\{ \, 0, \: \dots \: , \: 20 \, \big\}.$
	
\end{center}

Количество экспериментов для каждой упорядоченной пары $\big( S, \: n \big)$ было равно $5 \: \cdot \: 10^{6}$.

Далее (таб. \ref{tab:s1}, \ref{tab:s2}, \ref{tab:s3}, \ref{tab:s4}, \ref{tab:s5}, \ref{tab:s6}) приведены обработанные результаты проведённых экспериментов.

В первой колонке таблице приведены интервалы, в которых может находиться ширина итогового интервала, а в последующих --- относительная частота попадания ширины итогового интервала в данный интервал.

\begin{table}	
	\centering
	\caption{\label{tab:s1}Количество измерений $S = 1$}
	\vspace{3mm}
	$\begin{array}{c||c|c|c|c|c|c|c}	
	
		\hline \hline \rule[-1mm]{0mm}{6mm}

		\text{wid}(\mbf{r}) & n = 0 & n = 1 & n = 2 & n = 3 & n = 4 & n = 5 & n = 6 \\ \hline \hline \rule[-1mm]{0mm}{8mm}
		
		\big[ \, 0, \: \frac{1}{2} \, \big) & 1.000 & 0.986 & 0.936 & 0.885 & 0.839 & 0.797 & 0.758 \\ \rule[-1mm]{0mm}{6mm}
		
		\big[ \, \frac{1}{2}, \: 1 \, \big) & 0.000 & 0.013 & 0.032 & 0.050 & 0.059 & 0.068 & 0.074 \\ \rule[-1mm]{0mm}{6mm}
		
		[ \, 1, \: 2 \, ) & 0.000 & 0.001 & 0.014 & 0.026 & 0.035 & 0.044 & 0.050 \\ \rule[-1mm]{0mm}{6mm}
		
		[ \, 2, \: 3 \, ) & 0.000 & 0.000 & 0.006 & 0.010 & 0.015 & 0.018 & 0.022 \\ \rule[-1mm]{0mm}{6mm}
		
		[ \, 3, \: 4 \, ) & 0.000 & 0.000 & 0.004 & 0.007 & 0.009 & 0.011 & 0.013 \\ \rule[-1mm]{0mm}{6mm}
		
		[ \, 4, \: 5 \, ) & 0.000 & 0.000 & 0.003 & 0.005 & 0.007 & 0.008 & 0.010 \\ \rule[-1mm]{0mm}{6mm}
		
		[ \, 5, \: 10 \, ) & 0.000 & 0.000 & 0.005 & 0.009 & 0.015 & 0.020 & 0.025 \\ \rule[-1mm]{0mm}{6mm}
		
		[ \, 10, \: 100 \, ) & 0.000 & 0.000 & 0.001 & 0.008 & 0.017 & 0.027 & 0.037 \\ \rule[-1mm]{0mm}{6mm}
		
		[ \, 100, \: +\infty \, ) & 0.000 & 0.000 & 0.000 & 0.001 & 0.003 & 0.006 & 0.011 \\ \hline \hline \rule[-1mm]{0mm}{6mm} 
		
		\text{wid}(\mbf{r}) & n = 7 & n = 8 & n = 9 & n = 10 & n = 11 & n = 12 & n = 13 \\ \hline \hline \rule[-1mm]{0mm}{8mm}
		
		\big[ \, 0, \: \frac{1}{2} \, \big) & 0.722 & 0.689 & 0.659 & 0.631 & 0.606 & 0.583 & 0.561 \\ \rule[-1mm]{0mm}{6mm}
		
		\big[ \, \frac{1}{2}, \: 1 \, \big) & 0.079 & 0.083 & 0.086 & 0.089 & 0.090 & 0.092 & 0.093 \\ \rule[-1mm]{0mm}{6mm}
		
		[ \, 1, \: 2 \, ) & 0.056 & 0.060 & 0.064 & 0.067 & 0.069 & 0.071 & 0.073 \\ \rule[-1mm]{0mm}{6mm}
		
		[ \, 2, \: 3 \, ) & 0.024 & 0.027 & 0.029 & 0.031 & 0.032 & 0.033 & 0.035 \\ \rule[-1mm]{0mm}{6mm}
		
		[ \, 3, \: 4 \, ) & 0.015 & 0.017 & 0.018 & 0.019 & 0.020 & 0.021 & 0.022 \\ \rule[-1mm]{0mm}{6mm}
		
		[ \, 4, \: 5 \, ) & 0.011 & 0.012 & 0.013 & 0.014 & 0.015 & 0.015 & 0.015 \\ \rule[-1mm]{0mm}{6mm}
		
		[ \, 5, \: 10 \, ) & 0.029 & 0.032 & 0.035 & 0.038 & 0.040 & 0.042 & 0.043 \\ \rule[-1mm]{0mm}{6mm}
		
		[ \, 10, \: 100 \, ) & 0.047 & 0.056 & 0.064 & 0.071 & 0.078 & 0.085 & 0.091 \\ \rule[-1mm]{0mm}{6mm}
		
		[ \, 100, \: +\infty \, ) & 0.018 & 0.025 & 0.032 & 0.041 & 0.050 & 0.058 & 0.067 \\ \hline \hline \rule[-1mm]{0mm}{6mm}
		
		\text{wid}(\mbf{r}) & n = 14 & n = 15 & n = 16 & n = 17 & n = 18 & n = 19 & n = 20 \\ \hline \hline \rule[-1mm]{0mm}{8mm}
		
		\big[ \, 0, \: \frac{1}{2} \, \big) & 0.540 & 0.522 & 0.504 & 0.488 & 0.473 & 0.458 & 0.445 \\ \rule[-1mm]{0mm}{6mm}
		
		\big[ \, \frac{1}{2}, \: 1 \, \big) & 0.095 & 0.096 & 0.097 & 0.096 & 0.096 & 0.096 & 0.095 \\ \rule[-1mm]{0mm}{6mm}
		
		[ \, 1, \: 2 \, ) & 0.074 & 0.076 & 0.076 & 0.077 & 0.079 & 0.080 & 0.080 \\ \rule[-1mm]{0mm}{6mm}
		
		[ \, 2, \: 3 \, ) & 0.036 & 0.036 & 0.037 & 0.038 & 0.039 & 0.039 & 0.039 \\ \rule[-1mm]{0mm}{6mm}
		
		[ \, 3, \: 4 \, ) & 0.022 & 0.023 & 0.023 & 0.024 & 0.024 & 0.024 & 0.025 \\ \rule[-1mm]{0mm}{6mm}
		
		[ \, 4, \: 5 \, ) & 0.016 & 0.016 & 0.017 & 0.017 & 0.018 & 0.018 & 0.018 \\ \rule[-1mm]{0mm}{6mm}
		
		[ \, 5, \: 10 \, ) & 0.044 & 0.046 & 0.047 & 0.048 & 0.048 & 0.050 & 0.050 \\ \rule[-1mm]{0mm}{6mm}
		
		[ \, 10, \: 100 \, ) & 0.096 & 0.100 & 0.105 & 0.108 & 0.112 & 0.115 & 0.117 \\ \rule[-1mm]{0mm}{6mm}
		
		[ \, 100, \: +\infty \, ) & 0.076 & 0.085 & 0.094 & 0.103 & 0.112 & 0.121 & 0.129

	\end{array}$
\end{table}

\begin{table}	
	\centering
	\caption{\label{tab:s2}Количество измерений $S = 2$}
	\vspace{3mm}
	$\begin{array}{c||c|c|c|c|c|c|c}

		\hline \hline \rule[-1mm]{0mm}{6mm}

		\text{wid}(\mbf{r}) & n = 0 & n = 1 & n = 2 & n = 3 & n = 4 & n = 5 & n = 6 \\ \hline \hline \rule[-1mm]{0mm}{8mm}
		
		\big[ \, 0, \: \frac{1}{2} \, \big) & 1.000 & 0.987 & 0.926 & 0.863 & 0.804 & 0.749 & 0.702 \\ \rule[-1mm]{0mm}{6mm}
		
		\big[ \, \frac{1}{2}, \: 1 \, \big) & 0.000 & 0.010 & 0.035 & 0.058 & 0.072 & 0.082 & 0.088 \\ \rule[-1mm]{0mm}{6mm}
		
		[ \, 1, \: 2 \, ) & 0.000 & 0.001 & 0.018 & 0.032 & 0.045 & 0.055 & 0.063 \\ \rule[-1mm]{0mm}{6mm}
		
		[ \, 2, \: 3 \, ) & 0.000 & 0.000 & 0.008 & 0.014 & 0.020 & 0.025 & 0.029 \\ \rule[-1mm]{0mm}{6mm}
		
		[ \, 3, \: 4 \, ) & 0.000 & 0.000 & 0.004 & 0.008 & 0.012 & 0.015 & 0.018 \\ \rule[-1mm]{0mm}{6mm}
		
		[ \, 4, \: 5 \, ) & 0.000 & 0.000 & 0.003 & 0.005 & 0.008 & 0.010 & 0.013 \\ \rule[-1mm]{0mm}{6mm}
		
		[ \, 5, \: 10 \, ) & 0.000 & 0.000 & 0.004 & 0.010 & 0.017 & 0.024 & 0.031 \\ \rule[-1mm]{0mm}{6mm}
		
		[ \, 10, \: 100 \, ) & 0.000 & 0.001 & 0.002 & 0.009 & 0.020 & 0.032 & 0.044 \\ \rule[-1mm]{0mm}{6mm}
		
		[ \, 100, \: +\infty \, ) & 0.000 & 0.000 & 0.001 & 0.001 & 0.003 & 0.007 & 0.013 \\ \hline \hline \rule[-1mm]{0mm}{6mm}
		
		\text{wid}(\mbf{r}) & n = 7 & n = 8 & n = 9 & n = 10 & n = 11 & n = 12 & n = 13 \\ \hline \hline \rule[-1mm]{0mm}{8mm}
		
		\big[ \, 0, \: \frac{1}{2} \, \big) & 0.660 & 0.623 & 0.590 & 0.560 & 0.533 & 0.509 & 0.487 \\ \rule[-1mm]{0mm}{6mm}
		
		\big[ \, \frac{1}{2}, \: 1 \, \big) & 0.092 & 0.094 & 0.095 & 0.096 & 0.096 & 0.095 & 0.094 \\ \rule[-1mm]{0mm}{6mm}
		
		[ \, 1, \: 2 \, ) & 0.069 & 0.073 & 0.077 & 0.078 & 0.080 & 0.081 & 0.082 \\ \rule[-1mm]{0mm}{6mm}
		
		[ \, 2, \: 3 \, ) & 0.032 & 0.035 & 0.037 & 0.039 & 0.040 & 0.041 & 0.042 \\ \rule[-1mm]{0mm}{6mm}
		
		[ \, 3, \: 4 \, ) & 0.020 & 0.022 & 0.024 & 0.025 & 0.026 & 0.026 & 0.027 \\ \rule[-1mm]{0mm}{6mm}
		
		[ \, 4, \: 5 \, ) & 0.014 & 0.016 & 0.017 & 0.018 & 0.019 & 0.019 & 0.020 \\ \rule[-1mm]{0mm}{6mm}
		
		[ \, 5, \: 10 \, ) & 0.036 & 0.041 & 0.045 & 0.048 & 0.051 & 0.053 & 0.055 \\ \rule[-1mm]{0mm}{6mm}
		
		[ \, 10, \: 100 \, ) & 0.057 & 0.069 & 0.080 & 0.091 & 0.100 & 0.108 & 0.116 \\ \rule[-1mm]{0mm}{6mm}
		
		[ \, 100, \: +\infty \, ) & 0.019 & 0.027 & 0.036 & 0.046 & 0.056 & 0.066 & 0.077 \\ \hline \hline \rule[-1mm]{0mm}{6mm}
		
		\text{wid}(\mbf{r}) & n = 14 & n = 15 & n = 16 & n = 17 & n = 18 & n = 19 & n = 20 \\ \hline \hline \rule[-1mm]{0mm}{8mm}
		
		\big[ \, 0, \: \frac{1}{2} \, \big) & 0.468 & 0.450 & 0.434 & 0.418 & 0.404 & 0.391 & 0.379 \\ \rule[-1mm]{0mm}{6mm}
		
		\big[ \, \frac{1}{2}, \: 1 \, \big) & 0.093 & 0.092 & 0.091 & 0.089 & 0.088 & 0.087 & 0.086 \\ \rule[-1mm]{0mm}{6mm}
		
		[ \, 1, \: 2 \, ) & 0.082 & 0.082 & 0.082 & 0.082 & 0.081 & 0.080 & 0.079 \\ \rule[-1mm]{0mm}{6mm}
		
		[ \, 2, \: 3 \, ) & 0.042 & 0.043 & 0.043 & 0.043 & 0.043 & 0.043 & 0.043 \\ \rule[-1mm]{0mm}{6mm}
		
		[ \, 3, \: 4 \, ) & 0.027 & 0.028 & 0.028 & 0.028 & 0.029 & 0.029 & 0.029 \\ \rule[-1mm]{0mm}{6mm}
		
		[ \, 4, \: 5 \, ) & 0.020 & 0.021 & 0.021 & 0.021 & 0.021 & 0.021 & 0.021 \\ \rule[-1mm]{0mm}{6mm}
		
		[ \, 5, \: 10 \, ) & 0.056 & 0.058 & 0.059 & 0.060 & 0.060 & 0.061 & 0.061 \\ \rule[-1mm]{0mm}{6mm}
		
		[ \, 10, \: 100 \, ) & 0.123 & 0.128 & 0.133 & 0.137 & 0.141 & 0.145 & 0.148 \\ \rule[-1mm]{0mm}{6mm}
		
		[ \, 100, \: +\infty \, ) & 0.088 & 0.099 & 0.110 & 0.121 & 0.132 & 0.143 & 0.153

	\end{array}$
\end{table}

\begin{table}	
	\centering
	\caption{\label{tab:s3}Количество измерений $S = 3$}
	\vspace{3mm}
	$\begin{array}{c||c|c|c|c|c|c|c}

		\hline \hline \rule[-1mm]{0mm}{6mm}

		\text{wid}(\mbf{r}) & n = 0 & n = 1 & n = 2 & n = 3 & n = 4 & n = 5 & n = 6 \\ \hline \hline \rule[-1mm]{0mm}{8mm}
		
		\big[ \, 0, \: \frac{1}{2} \, \big) & 1.000 & 0.987 & 0.923 & 0.856 & 0.792 & 0.734 & 0.683 \\ \rule[-1mm]{0mm}{6mm}
		
		\big[ \, \frac{1}{2}, \: 1 \, \big) & 0.000 & 0.009 & 0.036 & 0.060 & 0.075 & 0.086 & 0.092 \\ \rule[-1mm]{0mm}{6mm}
		
		[ \, 1, \: 2 \, ) & 0.000 & 0.001 & 0.019 & 0.034 & 0.048 & 0.058 & 0.067 \\ \rule[-1mm]{0mm}{6mm}
		
		[ \, 2, \: 3 \, ) & 0.000 & 0.001 & 0.008 & 0.015 & 0.022 & 0.027 & 0.031 \\ \rule[-1mm]{0mm}{6mm}
		
		[ \, 3, \: 4 \, ) & 0.000 & 0.000 & 0.004 & 0.009 & 0.013 & 0.017 & 0.019 \\ \rule[-1mm]{0mm}{6mm}
		
		[ \, 4, \: 5 \, ) & 0.000 & 0.000 & 0.002 & 0.005 & 0.009 & 0.011 & 0.014 \\ \rule[-1mm]{0mm}{6mm}
		
		[ \, 5, \: 10 \, ) & 0.000 & 0.001 & 0.003 & 0.010 & 0.018 & 0.025 & 0.033 \\ \rule[-1mm]{0mm}{6mm}
		
		[ \, 10, \: 100 \, ) & 0.000 & 0.001 & 0.002 & 0.010 & 0.020 & 0.034 & 0.047 \\ \rule[-1mm]{0mm}{6mm}
		
		[ \, 100, \: +\infty \, ) & 0.000 & 0.000 & 0.001 & 0.002 & 0.004 & 0.007 & 0.013 \\ \hline \hline \rule[-1mm]{0mm}{6mm}
		
		\text{wid}(\mbf{r}) & n = 7 & n = 8 & n = 9 & n = 10 & n = 11 & n = 12 & n = 13 \\ \hline \hline \rule[-1mm]{0mm}{8mm}
		
		\big[ \, 0, \: \frac{1}{2} \, \big) & 0.640 & 0.601 & 0.566 & 0.536 & 0.509 & 0.485 & 0.464 \\ \rule[-1mm]{0mm}{6mm}
		
		\big[ \, \frac{1}{2}, \: 1 \, \big) & 0.095 & 0.097 & 0.097 & 0.098 & 0.096 & 0.095 & 0.094 \\ \rule[-1mm]{0mm}{6mm}
		
		[ \, 1, \: 2 \, ) & 0.073 & 0.077 & 0.080 & 0.082 & 0.083 & 0.083 & 0.083 \\ \rule[-1mm]{0mm}{6mm}
		
		[ \, 2, \: 3 \, ) & 0.035 & 0.037 & 0.039 & 0.041 & 0.043 & 0.043 & 0.043 \\ \rule[-1mm]{0mm}{6mm}
		
		[ \, 3, \: 4 \, ) & 0.022 & 0.024 & 0.025 & 0.026 & 0.028 & 0.028 & 0.029 \\ \rule[-1mm]{0mm}{6mm}
		
		[ \, 4, \: 5 \, ) & 0.016 & 0.017 & 0.018 & 0.019 & 0.020 & 0.021 & 0.021 \\ \rule[-1mm]{0mm}{6mm}
		
		[ \, 5, \: 10 \, ) & 0.039 & 0.044 & 0.048 & 0.052 & 0.055 & 0.057 & 0.058 \\ \rule[-1mm]{0mm}{6mm}
		
		[ \, 10, \: 100 \, ) & 0.062 & 0.075 & 0.087 & 0.098 & 0.108 & 0.117 & 0.125 \\ \rule[-1mm]{0mm}{6mm}
		
		[ \, 100, \: +\infty \, ) & 0.020 & 0.028 & 0.038 & 0.048 & 0.059 & 0.070 & 0.082 \\ \hline \hline \rule[-1mm]{0mm}{6mm}
		
		\text{wid}(\mbf{r}) & n = 14 & n = 15 & n = 16 & n = 17 & n = 18 & n = 19 & n = 20 \\ \hline \hline \rule[-1mm]{0mm}{8mm}
		
		\big[ \, 0, \: \frac{1}{2} \, \big) & 0.445 & 0.427 & 0.411 & 0.396 & 0.382 & 0.370 & 0.359 \\ \rule[-1mm]{0mm}{6mm}
		
		\big[ \, \frac{1}{2}, \: 1 \, \big) & 0.092 & 0.090 & 0.089 & 0.087 & 0.086 & 0.084 & 0.082 \\ \rule[-1mm]{0mm}{6mm}
		
		[ \, 1, \: 2 \, ) & 0.083 & 0.082 & 0.082 & 0.080 & 0.080 & 0.079 & 0.078 \\ \rule[-1mm]{0mm}{6mm}
		
		[ \, 2, \: 3 \, ) & 0.044 & 0.044 & 0.044 & 0.044 & 0.044 & 0.043 & 0.043 \\ \rule[-1mm]{0mm}{6mm}
		
		[ \, 3, \: 4 \, ) & 0.029 & 0.029 & 0.030 & 0.030 & 0.030 & 0.029 & 0.029 \\ \rule[-1mm]{0mm}{6mm}
		
		[ \, 4, \: 5 \, ) & 0.022 & 0.022 & 0.022 & 0.022 & 0.022 & 0.022 & 0.022 \\ \rule[-1mm]{0mm}{6mm}
		
		[ \, 5, \: 10 \, ) & 0.060 & 0.061 & 0.062 & 0.063 & 0.063 & 0.064 & 0.064 \\ \rule[-1mm]{0mm}{6mm}
		
		[ \, 10, \: 100 \, ) & 0.133 & 0.139 & 0.143 & 0.148 & 0.152 & 0.156 & 0.158 \\ \rule[-1mm]{0mm}{6mm}
		
		[ \, 100, \: +\infty \, ) & 0.094 & 0.106 & 0.118 & 0.130 & 0.141 & 0.152 & 0.164

	\end{array}$
\end{table}

\begin{table}	
	\centering
	\caption{\label{tab:s4}Количество измерений $S = 4$}
	\vspace{3mm}
	$\begin{array}{c||c|c|c|c|c|c|c}

		\hline \hline \rule[-1mm]{0mm}{6mm}

		\text{wid}(\mbf{r}) & n = 0 & n = 1 & n = 2 & n = 3 & n = 4 & n = 5 & n = 6 \\ \hline \hline \rule[-1mm]{0mm}{8mm}
		
		\big[ \, 0, \: \frac{1}{2} \, \big) & 1.000 & 0.987 & 0.922 & 0.852 & 0.785 & 0.726 & 0.674 \\ \rule[-1mm]{0mm}{6mm}
		
		\big[ \, \frac{1}{2}, \: 1 \, \big) & 0.000 & 0.008 & 0.037 & 0.062 & 0.078 & 0.088 & 0.094 \\ \rule[-1mm]{0mm}{6mm}
		
		[ \, 1, \: 2 \, ) & 0.000 & 0.002 & 0.020 & 0.035 & 0.049 & 0.060 & 0.068 \\ \rule[-1mm]{0mm}{6mm}
		
		[ \, 2, \: 3 \, ) & 0.000 & 0.001 & 0.009 & 0.016 & 0.023 & 0.028 & 0.033 \\ \rule[-1mm]{0mm}{6mm}
		
		[ \, 3, \: 4 \, ) & 0.000 & 0.000 & 0.004 & 0.009 & 0.013 & 0.017 & 0.020 \\ \rule[-1mm]{0mm}{6mm}
		
		[ \, 4, \: 5 \, ) & 0.000 & 0.000 & 0.002 & 0.005 & 0.009 & 0.012 & 0.014 \\ \rule[-1mm]{0mm}{6mm}
		
		[ \, 5, \: 10 \, ) & 0.000 & 0.001 & 0.003 & 0.010 & 0.018 & 0.026 & 0.034 \\ \rule[-1mm]{0mm}{6mm}
		
		[ \, 10, \: 100 \, ) & 0.000 & 0.001 & 0.003 & 0.010 & 0.021 & 0.035 & 0.050 \\ \rule[-1mm]{0mm}{6mm}
		
		[ \, 100, \: +\infty \, ) & 0.000 & 0.000 & 0.001 & 0.002 & 0.004 & 0.008 & 0.013 \\ \hline \hline \rule[-1mm]{0mm}{6mm}
		
		\text{wid}(\mbf{r}) & n = 7 & n = 8 & n = 9 & n = 10 & n = 11 & n = 12 & n = 13 \\ \hline \hline \rule[-1mm]{0mm}{8mm}
		
		\big[ \, 0, \: \frac{1}{2} \, \big) & 0.628 & 0.59 & 0.555 & 0.525 & 0.498 & 0.474 & 0.452 \\ \rule[-1mm]{0mm}{6mm}
		
		\big[ \, \frac{1}{2}, \: 1 \, \big) & 0.097 & 0.098 & 0.098 & 0.098 & 0.096 & 0.095 & 0.093 \\ \rule[-1mm]{0mm}{6mm}
		
		[ \, 1, \: 2 \, ) & 0.075 & 0.078 & 0.081 & 0.083 & 0.083 & 0.083 & 0.083 \\ \rule[-1mm]{0mm}{6mm}
		
		[ \, 2, \: 3 \, ) & 0.036 & 0.039 & 0.041 & 0.042 & 0.043 & 0.044 & 0.044 \\ \rule[-1mm]{0mm}{6mm}
		
		[ \, 3, \: 4 \, ) & 0.023 & 0.025 & 0.027 & 0.027 & 0.028 & 0.029 & 0.029 \\ \rule[-1mm]{0mm}{6mm}
		
		[ \, 4, \: 5 \, ) & 0.016 & 0.018 & 0.019 & 0.020 & 0.021 & 0.021 & 0.022 \\ \rule[-1mm]{0mm}{6mm}
		
		[ \, 5, \: 10 \, ) & 0.040 & 0.046 & 0.050 & 0.053 & 0.057 & 0.059 & 0.061 \\ \rule[-1mm]{0mm}{6mm}
		
		[ \, 10, \: 100 \, ) & 0.064 & 0.078 & 0.090 & 0.103 & 0.113 & 0.122 & 0.130 \\ \rule[-1mm]{0mm}{6mm}
		
		[ \, 100, \: +\infty \, ) & 0.021 & 0.029 & 0.039 & 0.050 & 0.061 & 0.073 & 0.085 \\ \hline \hline \rule[-1mm]{0mm}{6mm}
		
		\text{wid}(\mbf{r}) & n = 14 & n = 15 & n = 16 & n = 17 & n = 18 & n = 19 & n = 20 \\ \hline \hline \rule[-1mm]{0mm}{8mm}
		
		\big[ \, 0, \: \frac{1}{2} \, \big) & 0.433 & 0.415 & 0.400 & 0.386 & 0.372 & 0.360 & 0.349 \\ \rule[-1mm]{0mm}{6mm}
		
		\big[ \, \frac{1}{2}, \: 1 \, \big) & 0.091 & 0.089 & 0.087 & 0.086 & 0.084 & 0.083 & 0.081 \\ \rule[-1mm]{0mm}{6mm}
		
		[ \, 1, \: 2 \, ) & 0.083 & 0.082 & 0.082 & 0.080 & 0.080 & 0.078 & 0.077 \\ \rule[-1mm]{0mm}{6mm}
		
		[ \, 2, \: 3 \, ) & 0.045 & 0.045 & 0.044 & 0.044 & 0.044 & 0.044 & 0.043 \\ \rule[-1mm]{0mm}{6mm}
		
		[ \, 3, \: 4 \, ) & 0.030 & 0.030 & 0.030 & 0.030 & 0.030 & 0.030 & 0.030 \\ \rule[-1mm]{0mm}{6mm}
		
		[ \, 4, \: 5 \, ) & 0.022 & 0.022 & 0.023 & 0.022 & 0.022 & 0.022 & 0.022 \\ \rule[-1mm]{0mm}{6mm}
		
		[ \, 5, \: 10 \, ) & 0.062 & 0.063 & 0.063 & 0.064 & 0.064 & 0.065 & 0.065 \\ \rule[-1mm]{0mm}{6mm}
		
		[ \, 10, \: 100 \, ) & 0.137 & 0.143 & 0.149 & 0.154 & 0.158 & 0.161 & 0.163 \\ \rule[-1mm]{0mm}{6mm}
		
		[ \, 100, \: +\infty \, ) & 0.097 & 0.110 & 0.122 & 0.134 & 0.146 & 0.158 & 0.170

	\end{array}$
\end{table}

\begin{table}	
	\centering
	\caption{\label{tab:s5}Количество измерений $S = 5$}
	\vspace{3mm}
	$\begin{array}{c||c|c|c|c|c|c|c}

		\hline \hline \rule[-1mm]{0mm}{6mm}

		\text{wid}(\mbf{r}) & n = 0 & n = 1 & n = 2 & n = 3 & n = 4 & n = 5 & n = 6 \\ \hline \hline \rule[-1mm]{0mm}{8mm}
		
		\big[ \, 0, \: \frac{1}{2} \, \big) & 1.000 & 0.987 & 0.921 & 0.849 & 0.781 & 0.721 & 0.669 \\ \rule[-1mm]{0mm}{6mm}
		
		\big[ \, \frac{1}{2}, \: 1 \, \big) & 0.000 & 0.008 & 0.037 & 0.062 & 0.079 & 0.089 & 0.095 \\ \rule[-1mm]{0mm}{6mm}
		
		[ \, 1, \: 2 \, ) & 0.000 & 0.002 & 0.020 & 0.036 & 0.050 & 0.061 & 0.070 \\ \rule[-1mm]{0mm}{6mm}
		
		[ \, 2, \: 3 \, ) & 0.000 & 0.001 & 0.009 & 0.016 & 0.023 & 0.029 & 0.033 \\ \rule[-1mm]{0mm}{6mm}
		
		[ \, 3, \: 4 \, ) & 0.000 & 0.000 & 0.004 & 0.009 & 0.014 & 0.018 & 0.021 \\ \rule[-1mm]{0mm}{6mm}
		
		[ \, 4, \: 5 \, ) & 0.000 & 0.000 & 0.002 & 0.005 & 0.009 & 0.012 & 0.014 \\ \rule[-1mm]{0mm}{6mm}
		
		[ \, 5, \: 10 \, ) & 0.000 & 0.001 & 0.003 & 0.010 & 0.018 & 0.027 & 0.034 \\ \rule[-1mm]{0mm}{6mm}
		
		[ \, 10, \: 100 \, ) & 0.000 & 0.001 & 0.003 & 0.010 & 0.021 & 0.035 & 0.051 \\ \rule[-1mm]{0mm}{6mm}
		
		[ \, 100, \: +\infty \, ) & 0.000 & 0.000 & 0.001 & 0.002 & 0.004 & 0.008 & 0.014 \\ \hline \hline \rule[-1mm]{0mm}{6mm}
		
		\text{wid}(\mbf{r}) & n = 7 & n = 8 & n = 9 & n = 10 & n = 11 & n = 12 & n = 13 \\ \hline \hline \rule[-1mm]{0mm}{8mm}
		
		\big[ \, 0, \: \frac{1}{2} \, \big) & 0.623 & 0.584 & 0.549 & 0.518 & 0.492 & 0.468 & 0.447 \\ \rule[-1mm]{0mm}{6mm}
		
		\big[ \, \frac{1}{2}, \: 1 \, \big) & 0.098 & 0.099 & 0.099 & 0.098 & 0.096 & 0.095 & 0.093 \\ \rule[-1mm]{0mm}{6mm}
		
		[ \, 1, \: 2 \, ) & 0.075 & 0.079 & 0.082 & 0.083 & 0.083 & 0.084 & 0.084 \\ \rule[-1mm]{0mm}{6mm}
		
		[ \, 2, \: 3 \, ) & 0.036 & 0.039 & 0.041 & 0.043 & 0.044 & 0.044 & 0.044 \\ \rule[-1mm]{0mm}{6mm}
		
		[ \, 3, \: 4 \, ) & 0.023 & 0.025 & 0.027 & 0.028 & 0.029 & 0.029 & 0.030 \\ \rule[-1mm]{0mm}{6mm}
		
		[ \, 4, \: 5 \, ) & 0.016 & 0.018 & 0.019 & 0.020 & 0.021 & 0.022 & 0.022 \\ \rule[-1mm]{0mm}{6mm}
		
		[ \, 5, \: 10 \, ) & 0.041 & 0.046 & 0.051 & 0.054 & 0.057 & 0.059 & 0.061 \\ \rule[-1mm]{0mm}{6mm}
		
		[ \, 10, \: 100 \, ) & 0.066 & 0.080 & 0.093 & 0.105 & 0.115 & 0.125 & 0.133 \\ \rule[-1mm]{0mm}{6mm}
		
		[ \, 100, \: +\infty \, ) & 0.021 & 0.030 & 0.040 & 0.051 & 0.062 & 0.074 & 0.087 \\ \hline \hline \rule[-1mm]{0mm}{6mm}
		
		\text{wid}(\mbf{r}) & n = 14 & n = 15 & n = 16 & n = 17 & n = 18 & n = 19 & n = 20 \\ \hline \hline \rule[-1mm]{0mm}{8mm}
		
		\big[ \, 0, \: \frac{1}{2} \, \big) & 0.427 & 0.410 & 0.394 & 0.380 & 0.367 & 0.355 & 0.344 \\ \rule[-1mm]{0mm}{6mm}
		
		\big[ \, \frac{1}{2}, \: 1 \, \big) & 0.091 & 0.089 & 0.087 & 0.085 & 0.083 & 0.081 & 0.080 \\ \rule[-1mm]{0mm}{6mm}
		
		[ \, 1, \: 2 \, ) & 0.083 & 0.082 & 0.082 & 0.080 & 0.079 & 0.078 & 0.077 \\ \rule[-1mm]{0mm}{6mm}
		
		[ \, 2, \: 3 \, ) & 0.045 & 0.045 & 0.044 & 0.044 & 0.044 & 0.043 & 0.043 \\ \rule[-1mm]{0mm}{6mm}
		
		[ \, 3, \: 4 \, ) & 0.030 & 0.030 & 0.030 & 0.030 & 0.030 & 0.030 & 0.030 \\ \rule[-1mm]{0mm}{6mm}
		
		[ \, 4, \: 5 \, ) & 0.022 & 0.023 & 0.023 & 0.023 & 0.023 & 0.023 & 0.022 \\ \rule[-1mm]{0mm}{6mm}
		
		[ \, 5, \: 10 \, ) & 0.062 & 0.063 & 0.065 & 0.065 & 0.066 & 0.065 & 0.065 \\ \rule[-1mm]{0mm}{6mm}
		
		[ \, 10, \: 100 \, ) & 0.140 & 0.146 & 0.151 & 0.155 & 0.159 & 0.163 & 0.166 \\ \rule[-1mm]{0mm}{6mm}
		
		[ \, 100, \: +\infty \, ) & 0.099 & 0.112 & 0.125 & 0.137 & 0.149 & 0.161 & 0.173

	\end{array}$
\end{table}

\begin{table}	
	\centering
	\caption{\label{tab:s6}Количество измерений $S = 6$}
	\vspace{3mm}
	$\begin{array}{c||c|c|c|c|c|c|c}

		\hline \hline \rule[-1mm]{0mm}{6mm}

		\text{wid}(\mbf{r}) & n = 0 & n = 1 & n = 2 & n = 3 & n = 4 & n = 5 & n = 6 \\ \hline \hline \rule[-1mm]{0mm}{8mm}
		
		\big[ \, 0, \: \frac{1}{2} \, \big) & 1.000 & 0.987 & 0.920 & 0.847 & 0.779 & 0.718 & 0.665 \\ \rule[-1mm]{0mm}{6mm}
		
		\big[ \, \frac{1}{2}, \: 1 \, \big) & 0.000 & 0.008 & 0.037 & 0.062 & 0.079 & 0.090 & 0.095 \\ \rule[-1mm]{0mm}{6mm}
		
		[ \, 1, \: 2 \, ) & 0.000 & 0.002 & 0.020 & 0.037 & 0.051 & 0.062 & 0.070 \\ \rule[-1mm]{0mm}{6mm}
		
		[ \, 2, \: 3 \, ) & 0.000 & 0.001 & 0.009 & 0.017 & 0.024 & 0.029 & 0.034 \\ \rule[-1mm]{0mm}{6mm}
		
		[ \, 3, \: 4 \, ) & 0.000 & 0.000 & 0.004 & 0.009 & 0.014 & 0.018 & 0.021 \\ \rule[-1mm]{0mm}{6mm}
		
		[ \, 4, \: 5 \, ) & 0.000 & 0.000 & 0.002 & 0.005 & 0.009 & 0.012 & 0.015 \\ \rule[-1mm]{0mm}{6mm}
		
		[ \, 5, \: 10 \, ) & 0.000 & 0.001 & 0.003 & 0.010 & 0.018 & 0.027 & 0.035 \\ \rule[-1mm]{0mm}{6mm}
		
		[ \, 10, \: 100 \, ) & 0.000 & 0.001 & 0.003 & 0.010 & 0.022 & 0.036 & 0.051 \\ \rule[-1mm]{0mm}{6mm}
		
		[ \, 100, \: +\infty \, ) & 0.000 & 0.000 & 0.001 & 0.002 & 0.004 & 0.008 & 0.014 \\ \hline \hline \rule[-1mm]{0mm}{6mm}
		
		\text{wid}(\mbf{r}) & n = 7 & n = 8 & n = 9 & n = 10 & n = 11 & n = 12 & n = 13 \\ \hline \hline \rule[-1mm]{0mm}{8mm}
		
		\big[ \, 0, \: \frac{1}{2} \, \big) & 0.619 & 0.578 & 0.544 & 0.514 & 0.488 & 0.464 & 0.443 \\ \rule[-1mm]{0mm}{6mm}
		
		\big[ \, \frac{1}{2}, \: 1 \, \big) & 0.098 & 0.099 & 0.099 & 0.097 & 0.095 & 0.093 & 0.092 \\ \rule[-1mm]{0mm}{6mm}
		
		[ \, 1, \: 2 \, ) & 0.076 & 0.080 & 0.082 & 0.084 & 0.084 & 0.084 & 0.084 \\ \rule[-1mm]{0mm}{6mm}
		
		[ \, 2, \: 3 \, ) & 0.037 & 0.040 & 0.041 & 0.043 & 0.044 & 0.044 & 0.045 \\ \rule[-1mm]{0mm}{6mm}
		
		[ \, 3, \: 4 \, ) & 0.024 & 0.026 & 0.027 & 0.028 & 0.029 & 0.030 & 0.030 \\ \rule[-1mm]{0mm}{6mm}
		
		[ \, 4, \: 5 \, ) & 0.017 & 0.019 & 0.020 & 0.021 & 0.021 & 0.022 & 0.022 \\ \rule[-1mm]{0mm}{6mm}
		
		[ \, 5, \: 10 \, ) & 0.041 & 0.047 & 0.052 & 0.055 & 0.058 & 0.060 & 0.062 \\ \rule[-1mm]{0mm}{6mm}
		
		[ \, 10, \: 100 \, ) & 0.067 & 0.082 & 0.095 & 0.107 & 0.118 & 0.127 & 0.135 \\ \rule[-1mm]{0mm}{6mm}
		
		[ \, 100, \: +\infty \, ) & 0.021 & 0.030 & 0.040 & 0.051 & 0.063 & 0.076 & 0.089 \\ \hline \hline \rule[-1mm]{0mm}{6mm}
		
		\text{wid}(\mbf{r}) & n = 14 & n = 15 & n = 16 & n = 17 & n = 18 & n = 19 & n = 20 \\ \hline \hline \rule[-1mm]{0mm}{8mm}
		
		\big[ \, 0, \: \frac{1}{2} \, \big) & 0.424 & 0.407 & 0.392 & 0.377 & 0.364 & 0.352 & 0.341 \\ \rule[-1mm]{0mm}{6mm}
		
		\big[ \, \frac{1}{2}, \: 1 \, \big) & 0.089 & 0.088 & 0.086 & 0.084 & 0.083 & 0.081 & 0.079 \\ \rule[-1mm]{0mm}{6mm}
		
		[ \, 1, \: 2 \, ) & 0.083 & 0.082 & 0.081 & 0.080 & 0.079 & 0.078 & 0.076 \\ \rule[-1mm]{0mm}{6mm}
		
		[ \, 2, \: 3 \, ) & 0.044 & 0.044 & 0.044 & 0.044 & 0.044 & 0.043 & 0.043 \\ \rule[-1mm]{0mm}{6mm}
		
		[ \, 3, \: 4 \, ) & 0.030 & 0.030 & 0.030 & 0.030 & 0.030 & 0.030 & 0.030 \\ \rule[-1mm]{0mm}{6mm}
		
		[ \, 4, \: 5 \, ) & 0.022 & 0.022 & 0.022 & 0.022 & 0.022 & 0.022 & 0.022 \\ \rule[-1mm]{0mm}{6mm}
		
		[ \, 5, \: 10 \, ) & 0.063 & 0.064 & 0.065 & 0.065 & 0.065 & 0.065 & 0.065 \\ \rule[-1mm]{0mm}{6mm}
		
		[ \, 10, \: 100 \, ) & 0.142 & 0.148 & 0.153 & 0.157 & 0.161 & 0.164 & 0.167 \\ \rule[-1mm]{0mm}{6mm}
		
		[ \, 100, \: +\infty \, ) & 0.102 & 0.114 & 0.127 & 0.140 & 0.153 & 0.165 & 0.177

	\end{array}$
\end{table}

Видно, что для каждого значения $S$ выявляется одинаковая закономерность --- чем большее число раз вовлекаются в вычисления интервалы, тем шире становится итоговый интервал. При рассмотрении случаев $n > S$ можно увидеть влияние зависимых интервальных параметров на ширину итогового интервала при вычислениях с ними.

Таким образом, необходимо учитывать различные связи между интервальными параметрами в вычисляемом выражении, чтобы получать после вычислений как можно более узкий интервал.

\subsection[Случай матрицы со связями. Постановка задачи]{Случай матрицы со связями.\\Постановка задачи} 

Ранее, при постановке задачи нахождения оценки объединённого множества решений ИСЛАУ, мы полагали, что элементы, входящие в неё, независимые. Однако на практике это зачастую не так.

В случае, когда элементы системы имеют связи сложного вида, не имеется алгоритмов для нахождения точных оценок объединённого множества решений. Задача в такой постановке является более общей, и, тем самым, более сложной для решения. 

Примером может служить интервальная линейная задача наименьших квадратов. Она была рассмотрена в работе Дэвида Гея \cite{Gay}. Специфика данной задачи состоит в том, что её можно свести к решению ИСЛАУ нормального вида. При этом, согласно работе \cite{SharyMoradi}, можно иметь различные представления этого нормального вида:

\begin{center}

    $\mbf{A}^{\top} (\mbf{b} - \mbf{A} x) = 0$, \qquad $\left\{\begin{array}{rll}
     
		y + \mbf{A}x & = & \mbf{b},\\
        
        \mbf{A}^{\top} y & = & 0.
		
	\end{array}\right.$
    
\end{center}

Итак, стоит задача нахождения оценки границ объединённого множества решений данной ИСЛАУ. Матрица системы является симметричной, то есть её элементы меняются в заданных пределах, но при этом любая <<точечная>> СЛАУ является симметричной.

Оценки объединённого множества ИСЛАУ с зависимыми элементами, которые могут быть получены методами интервального анализа, в общем случае, являются ещё более грубыми, чем с независимыми.

Таким образом, поставленная задача позволяет в полной мере продемонстрировать негативное влияние эффекта зависимости на получение ответа. Также она позволит качественно сравнить различные интервальные арифметики для её решения.

\clearpage
\section [Численные погрешности. Использование интервального анализа] {Численные погрешности.\\ Использование интервального анализа}

При решении практических задач часто возникают численные неопредёленности, имеющие разную природу. Одни, например, происходят из неточности измерительных приборов или человеческого фактора. Другие --- из-за использования на ЭВМ арифметики с ограниченной точностью (или арифметики с плавающей точкой).

Эти неточности могут достаточно сильно отклонить итоговый ответ от идеально точного, в особенности, когда он является результатом длинной цепочки вычислений.

Более того, существуют погрешности, которые мы не можем нивелировать. Примерами могут служить атомные массы химических элементов, для значений которых известны лишь доверительные интервалы, в которых они гарантированно содержатся. 

\begin{figure}
    \centering
        \includegraphics[width=0.8\linewidth]{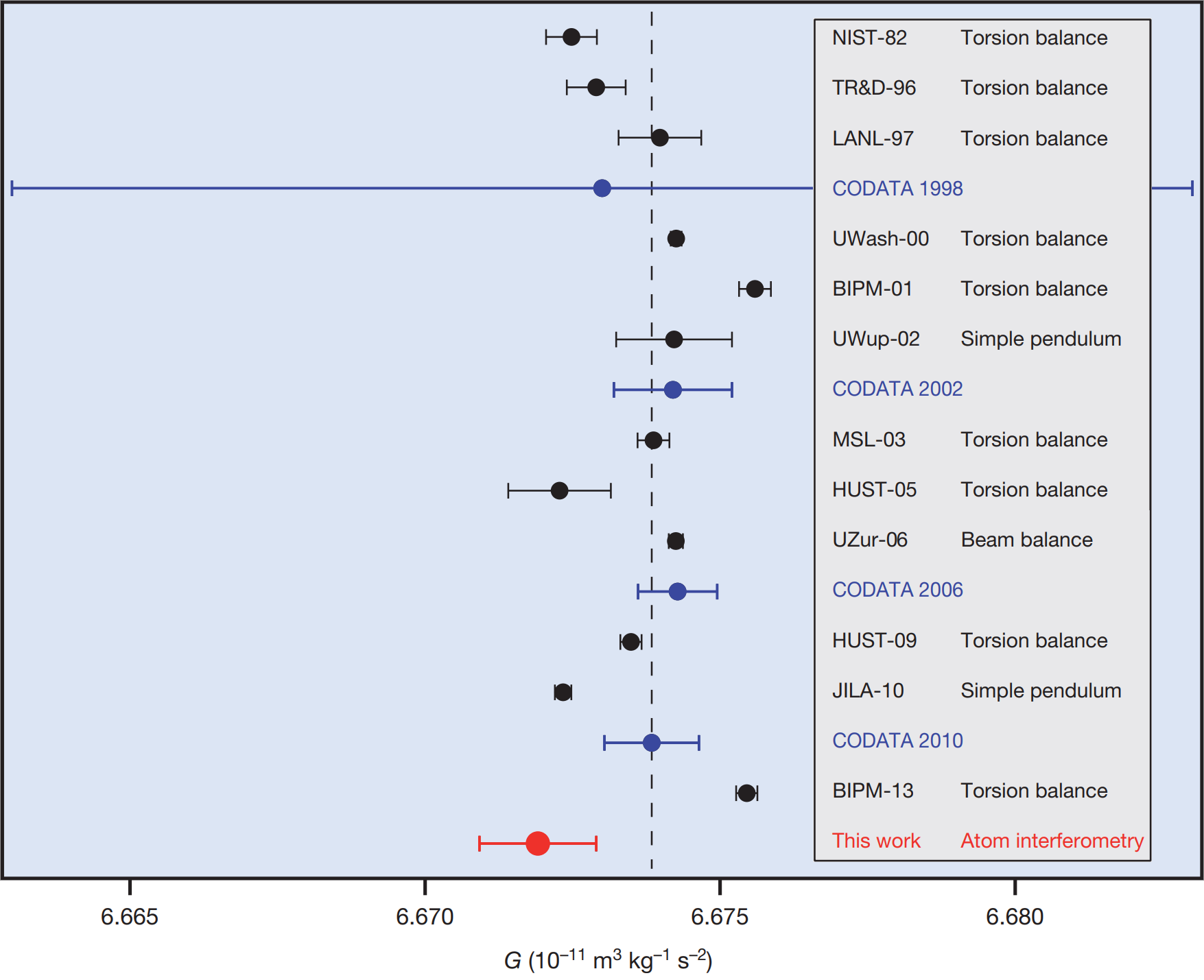}
        \caption{История опубликованных результатов экспериментов по измерению константы $G$}
    \label{fig:G}
\end{figure}

Только некоторые теоретические константы могут быть точно записаны в числовом или аналитическом виде. В качестве примера можно рассмотреть историю получения значений гравитационной постоянной $G$. На (рис. \ref{fig:G}) приведена иллюстрация из статьи журнала <<Nature>>, где показана история опубликованных результатов экспериментов по измерению данной константы.

Видно, что полученные в различных экспериментах доверительные интервалы для значения гравитационной постоянной не пересекаются даже в одной точке. Сложность сопоставления результатов состоит в том, что невозможно в той или иной мере учесть все искажающие факторы экспериментальной среды. 

Таким образом, встаёт вопрос: какой из полученных доверительных интервалов считать наиболее достоверным для проведения расчетов? Ответ на данный вопрос важен, поскольку при наличии такого интервала можно использовать его среднюю точку как наиболее представительную <<точечную>> величину, характеризующую константу $G$.

Одним из решений данной проблемы является выбор представительной не <<точечной>>, а интервальной величины. Например, можно взять интервал, полученный в результате объединения всех (или только некоторых) имеющихся интервалов и далее использовать его в дальнейших вычислениях.

Однако тогда встает другой вопрос: каким образом необходимо решать задачу с представительной интервальной величиной, чтобы полученный в результате вычислений ответ был допустим для дальнейшего использования и анализа?

Ответом на этот вопрос может быть использование интервального анализа. Его особенностью является то, что он использует в качестве базисного объекта интервал. Эта характерная черта данной математической дисциплины позволяет в полной мере учитывать неопределенности, встречающиеся в задачах.

В примере с представительной интервальной величиной для константы $G$ можно проводить все вычисления, используя классическую интервальную арифметику. Тем самым можно гарантировать, что полученный в результате интервал будет содержать в себе идеально точное искомое решение.

\subsection[Пример Мюллера]{Пример Мюллера}

Пусть стоит задача найти тридцатый элемент последовательности

\begin{center}

	$\begin{array}{l}
	
		x_{0} = 4,\\
		
		x_{1} = 4.25,\\
		
		x_{i} = 108 - \frac{815 - 1500 / x_{i - 2}}{x_{i - 1}}.
		
	\end{array}$
	
\end{center}

Как обычно принято на практике, для решения данной задачи будем использовать машинную арифметику двойной точности (тип <<\texttt{double}>> согласно стандарту \cite{IEEE754}). Также решим эту задачу, учитывая ошибки округления при оперировании числами с плавающей точкой в процессе вычислений. Для этого будем использовать классическую интервальную арифметику. Полученные результаты представлены в таблице (\ref{tab}).

\begin{table}
	\centering
	\caption{\label{tab}Результаты вычислительных экспериментов}
	\vspace{3mm}
	\begin{small}
	\begin{tabular}{r||c|c}
	
		 & \hspace{3mm} <<Точечная>> арифметика \hspace{3mm} & Классическая интервальная \\ 	
	
		$i$ & с плавающей точкой & арифметика \\ \hline \hline \rule[-1mm]{0mm}{6mm}
		
		$1$ & $4$ & $[ \, 4.000000000000, \: 4.000000000000 \, ]$\\ \hline \rule[-1mm]{0mm}{6mm}
		
		$2$ & $4.25$ & $[ \, 4.250000000000, \: 4.250000000000 \, ]$ \\ \hline \rule[-1mm]{0mm}{6mm}
		
		$3$ & $4.470588235294116$ & $[ \, 4.470588235294, \: 4.470588235294 \, ]$ \\ \hline \rule[-1mm]{0mm}{6mm}
		
		$4$ & $4.6447368421052317$ & $[ \, 4.644736842105, \: 4.644736842106 \, ]$ \\ \hline \rule[-1mm]{0mm}{6mm}
		
		$5$ & $4.7705382436253814$ & $[ \, 4.770538243625, \: 4.770538243633 \, ]$ \\ \hline \rule[-1mm]{0mm}{6mm}
		
		$6$ & $4.8557007125748015$ & $[ \, 4.855700712558, \: 4.855700712739 \, ]$ \\ \hline \rule[-1mm]{0mm}{6mm}
		
		$7$ & $4.9108474987890247$ & $[ \, 4.910847498332, \: 4.910847502291 \, ]$ \\ \hline \rule[-1mm]{0mm}{6mm}
		
		$8$ & $4.9455373981440545$ & $[ \, 4.945537386426,\: 4.945537471849 \, ]$ \\ \hline \rule[-1mm]{0mm}{6mm}
		
		$9$ & $4.9669624608756351$ & $[ \, 4.966962172699, \: 4.966964002158 \, ]$ \\ \hline \rule[-1mm]{0mm}{6mm}
		
		$10$ & $4.98004326786787$ & $[ \, 4.980036380732, \: 4.980075380346 \, ]$\\ \hline \rule[-1mm]{0mm}{6mm}
		
		$11$ & $4.9879305877668543$ & $[ \, 4.987769309963, \: 4.988598347566 \, ]$ \\ \hline \rule[-1mm]{0mm}{6mm}
		
		$12$ & $4.9917907586369807$ & $[ \, 4.988070629252, \: 5.005662648814 \, ]$ \\ \hline \rule[-1mm]{0mm}{6mm}
		
		$13$ & $4.9760336813273796$ & $[ \, 4.891127743800, \: 5.263480342516 \, ]$ \\ \hline \rule[-1mm]{0mm}{6mm}
		
		$14$ & $4.6030641192723891$ & $[ \, 2.637923691460, \: 10.292306138550 \, ]$ \\ \hline \rule[-1mm]{0mm}{6mm}
		
		$15$ & $-3.568095025888212$ & $[ \, -92.922210598236, \: 58.611432105461 \, ]$ \\ \hline \rule[-1mm]{0mm}{6mm}
	
		$16$ & $245.084384261717$ & $[ \, -\infty, \: +\infty \, ]$ \\ \hline \rule[-1mm]{0mm}{6mm}
	
		$17$ & $102.95931821632739$ & $[ \, -\infty, \: +\infty \, ]$ \\ \hline \rule[-1mm]{0mm}{6mm}
		
		$18$ & $100.14369643238381$ & $[ \, -\infty, \: +\infty \, ]$ \\ \hline \rule[-1mm]{0mm}{6mm}
		
		$19$ & $100.00717401619831$ & $[ \, -\infty, \: +\infty \, ]$ \\ \hline\rule[-1mm]{0mm}{6mm}
		
		$20$ & $100.00035866020249$ & $[ \, -\infty, \: +\infty \, ]$ \\ \hline \rule[-1mm]{0mm}{6mm}
		
		$21$ & $100.00001793249953$ & $[ \, -\infty, \: +\infty \, ]$ \\ \hline \rule[-1mm]{0mm}{6mm}
		
		$22$ & $100.00000089661144$ & $[ \, -\infty, \: +\infty \, ]$ \\ \hline \rule[-1mm]{0mm}{6mm}
		
		$23$ & $100.00000004483017$ & $[ \, -\infty, \: +\infty \, ]$ \\ \hline \rule[-1mm]{0mm}{6mm}
		
		$24$ & $100.00000000224151$ & $[ \, -\infty, \: +\infty \, ]$ \\ \hline \rule[-1mm]{0mm}{6mm}
		
		$25$ & $100.00000000011208$ & $[ \, -\infty, \: +\infty \, ]$ \\ \hline \rule[-1mm]{0mm}{6mm}
		
		$26$ & $100.00000000000561$ & $[ \, -\infty, \: +\infty \, ]$ \\ \hline \rule[-1mm]{0mm}{6mm}
		
		$27$ & $100.00000000000028$ & $[ \, -\infty, \: +\infty \, ]$ \\ \hline \rule[-1mm]{0mm}{6mm}
		
		$28$ & $100.00000000000001$ & $[ \, -\infty, \: +\infty \, ]$ \\ \hline \rule[-1mm]{0mm}{6mm}
		
		$29$ & $100$ & $[ \, -\infty, \: +\infty \, ]$ \\ \hline \rule[-1mm]{0mm}{6mm}
		
		$30$ & $100$ & $[ \, -\infty, \: +\infty \, ]$
		
	\end{tabular}
	\end{small}
\end{table}

Проведём анализ полученных результатов. Рассмотрим последовательность, полученную вычислением арифметикой с плавающей точкой. Видно, что последовательность сошлась к $100$.

Теперь обратимся к результатам, полученным использованием интервальной классической арифметики. Начиная с $14$-ой итерации, ширина интервала начинает расти. Заметим, что это происходит в тот момент, когда в <<точечной>> последовательности начались отклонения от числа $5$.

Далее ширина интервала, содержащего значение $x_{i}$, продолжает расти и уже для $i = 15$ мы получаем интервал, ширина которого превышает $100$. Так же он содержит в себе $0$, что ставит под сомнение возможность деления на $x_{15}$ при вычислении $x_{16}$.

Принимая во внимание данные наблюдения, мы должны поставить вопрос: может ли идти речь о численной сходимости к какому-либо числу в принципе.

Рассмотренное рекуррентное соотношение называется \textit{примером Мюллера}. Аналитически данная последовательность сходится к числу $5$, однако при вычислениях с помощью арифметики с плавающей точкой, даже более точной, чем <<\texttt{double}>> (например <<длинная>> арифметика в языке программирования \texttt{Python}), последовательность, начиная с некоторого номера, по-прежнему начинает сходиться к числу $100$. 

Данный эффект объясняется тем, что дискретный процесс, который задается приведённым рекуррентным соотношением, имеет неустойчивую точку $x = 5$ и устойчивую точку $x = 100$. В процессе вычислений с помощью <<точечной>> арифметики неизбежно появляются ошибки округлений. Они и провоцируют переход из окрестности неустойчивой точки в окрестность устойчивой.

Пример Мюллера наглядно демонстрирует пользу интервальной арифметики для проведения практических расчетов. Она позволяет учесть все вычислительные эффекты, которые могут исказить результат, а также даёт наглядное представление <<точности>> результата в виде ширины интервала. 

Это свойство очень полезно на практике. Представим, что перед нами стоит задача по вычислению последовательности Мюллера, при этом мы не имеем никаких априорных представлений о её свойствах. 

Тогда в процессе вычислений, после $14$-ой итерации, мы могли бы сделать вывод о том, что проводимые расчёты не имеют смысла, поскольку последовательность ведет себя неустойчиво и, соответственно, результат не пригоден для дальнейшего анализа.

\clearpage
\section [Обзор интервальных арифметик] {Обзор интервальных арифметик}

Указанные ранее проблемы, связанные с погрешностями представления чисел, а также с округлением при вычислениях, мотивировали создание математического аппарата, который был бы построен на следующем принципе:
\vspace{-5mm}
\begin{equation}
	\label{eq:BaseIntervalPrincipe}
	\begin{array}{c}

		\text{Пусть } \mbf{a} \text{ и } \mbf{b} \text{ --- интервалы, тогда} \\
	
		\mbf{a} \, * \, \mbf{b} = \{ \, a \, * \, b \mid a \in \mbf{a}, \: b \in \mbf{b} \, \},\\
		
		\text{где <<} * \text{>> --- арифметическая операция между интервалами.}
	
	\end{array}
\end{equation} 

Иными словами, мы требуем, чтобы результирующий интервал содержал в себе все возможные результаты применения арифметического оператора для любой упорядоченной пары $(a, \: b)$, где $a \in \mbf{a}$, $b \in \mbf{b}$.

В интервальном анализе было сконструировано несколько видов интервальных арифметик, основанных на этом принципе.

\subsection{Классическая интервальная арифметика}

Классическая интервальная арифметика определена на множестве $\mathbb{I}\mathbb{R}$. Операции, которые мы определим далее, являются прямой реализацией базового принципа интервального анализа (\ref{eq:BaseIntervalPrincipe}):

\begin{enumerate}
    \item Сложение:
    
    \begin{center}
		$\begin{array}{rcl}
	
			\mbf{a} + \mbf{b} & = & \big[ \, \underline{\mbf{a}}, \: \overline{\mbf{a}} \, \big] + \big[ \, \underline{\mbf{b}}, \: \overline{\mbf{b}} \, \big] = \\
	
			& = & \big\{ \, a + b \, \big| \, a \in \mbf{a}, \: b \in \mbf{b} \, \big\} = \\
		
			& = & \big\{ \, x \, \big| \, \underline{\mbf{a}} + \underline{\mbf{b}} \leq x \leq \overline{\mbf{a}} + \overline{\mbf{b}} \, \big\} = \\
	
			& = & \big[ \, \underline{\mbf{a}} + \underline{\mbf{b}}, \: \overline{\mbf{a}} + \overline{\mbf{b}} \, \big].
	
		\end{array}$
	\end{center}    
    
    \item Вычитание:
    
    \begin{center}
    	$\begin{array}{rcl}
    
			\mbf{a} - \mbf{b} & = & \big[ \, \underline{\mbf{a}}, \: \overline{\mbf{a}} \, \big] - \big[ \, \underline{\mbf{b}}, \: \overline{\mbf{b}} \, \big] = \\
    
    		& = & \big\{ \, a - b \, \big| \, a \in \mbf{a}, \: b \in \mbf{b} \, \big\} = \\
    
			& = & \big\{ \, x \, \big| \, \underline{\mbf{a}} - \overline{\mbf{b}} \leq x \leq \overline{\mbf{a}} - \underline{\mbf{b}} \, \big\} = \\    
    
    		& = & \big[ \, \underline{\mbf{a}} - \overline{\mbf{b}}, \: \overline{\mbf{a}} - \underline{\mbf{b}} \, \big].
    
    	\end{array}$
    \end{center}
    
    \item Умножение:  
    
	Рассмотрим функцию 
	\begin{center}
		
		$f(x, \: y) = x \, \cdot \, y.$
		
	\end{center}
	
	Найдём её стационарные точки:
	
	\begin{center}
		\begin{tabular}{l}
		
			$\left \{ \begin{array}{rcl}
		
				f_{x}^{'}(x, \: y) & = & 0, \\
			
				f_{y}^{'}(x, \: y) & = & 0.			
			
			\end{array}\right.$ \\
		
			$\left \{ \begin{array}{rcl}
		
				y & = & 0, \\
			
				x & = & 0.		
		
			\end{array}\right.$
			
		\end{tabular}
	\end{center}
	
	Единственная стационарная точка данной функции --- $(0, \: 0)$. Она является седловой, поэтому при рассмотрении $x \in \mbf{x}$, $y \in \mbf{y}$ получаем, что экстремумы достигаются на границах интервалов $\mbf{x}$ и $\mbf{y}$.
	
	Используя это, получаем:    
    
    \begin{center}
    	$\begin{array}{rcl}
    
        	\mbf{a} \cdot \mbf{b} & = & \big[ \, \underline{\mbf{a}}, \: \overline{\mbf{a}} \, \big] \cdot \big[ \, \underline{\mbf{b}}, \: \overline{\mbf{b}} \, \big] = \\
        	
        	& = & \big\{ \, a \, \cdot \, b \mid a \in \mbf{a}, \: b \in \mbf{b} \, \big\} = \\
        	
        	& = & \Big\{ \, x \, \big| \, \text{min} \big\{ \, \underline{\mbf{a}} \cdot \underline{\mbf{b}}, \: \underline{\mbf{a}} \cdot \overline{\mbf{b}}, \: \overline{\mbf{a}} \cdot \underline{\mbf{b}}, \: \overline{\mbf{a}} \cdot \overline{\mbf{b}} \, \big\} \leq x \: \bigcap \\
        	
        	& & \hspace{2cm} \bigcap \: x \leq \, \text{max} \big\{ \, \underline{\mbf{a}} \cdot \underline{\mbf{b}}, \: \underline{\mbf{a}} \cdot \overline{\mbf{b}}, \: \overline{\mbf{a}} \cdot \underline{\mbf{b}}, \: \overline{\mbf{a}} \cdot \overline{\mbf{b}} \, \big\} \Big\} = \\
        
         	& = & \Big[ \, \text{min} \big\{ \, \underline{\mbf{a}} \cdot \underline{\mbf{b}}, \: \underline{\mbf{a}} \cdot \overline{\mbf{b}}, \: \overline{\mbf{a}} \cdot \underline{\mbf{b}}, \: \overline{\mbf{a}} \cdot \overline{\mbf{b}} \, \big\}, \\
         
        	& & \hspace{3cm} \text{max} \big\{ \, \underline{\mbf{a}} \cdot \underline{\mbf{b}}, \: \underline{\mbf{a}} \cdot \overline{\mbf{b}}, \: \overline{\mbf{a}} \cdot \underline{\mbf{b}}, \: \overline{\mbf{a}} \cdot \overline{\mbf{b}} \, \big\} \, \Big].
         
    	\end{array}$
    \end{center}
    
    \item Деление: 
    
	Рассмотрим функцию
	\begin{center}
	
		$f(x, \: y) = x \, / \, y.$
		
	\end{center}
	
	Найдём её стационарные точки:
	\begin{center}
	
	\begin{tabular}{l}
		$\left \{ \begin{array}{rcl}
		
			f_{x}^{'}(x, \: y) & = & 0, \\
			
			f_{y}^{'}(x, \: y) & = & 0.
		
		\end{array} \right.$ \\
		
		$\left \{ \begin{array}{rcl}
		
			1 \, / \, y & = & 0, \\
			
			- x \, / \, y ^ {2} & = & 0.		
		
		\end{array} \right.$ \\
		
		$\left \{ \begin{array}{rcl}
		
			y & \in & \varnothing, \\
			
			x & = & 0.
		
		\end{array} \right.$
	\end{tabular}
	
	\end{center}
	
	Данная функция не имеет стационарных точек, поэтому при рассмотрении $x \in \mbf{x}$, $y \in \mbf{y}$ получаем, что экстремумы достигаются на границах интервалов.
    
    Используя это, получаем, что если $0 \notin \mbf{b}$, то:
    
    \begin{center}
    	$\begin{array}{rcl}
    
        	\mbf{a} \, / \, \mbf{b} & = & \big[ \, \underline{\mbf{a}}, \: \overline{\mbf{a}} \, \big] \, / \, \big[ \, \underline{\mbf{b}}, \: \overline{\mbf{b}} \, \big] = \\
        	
        	& = & \big\{ \, a \, / \, b \, \big| \, a \in \mbf{a}, \: b \in \mbf{b} \, \big\} = \\
        	
        	& = & \Big\{ \, x \mid \text{min} \big\{ \, \underline{\mbf{a}} \, / \, \underline{\mbf{b}}, \: \underline{\mbf{a}} \, / \, \overline{\mbf{b}}, \: \overline{\mbf{a}} \, / \, \underline{\mbf{b}}, \: \overline{\mbf{a}} \, / \, \overline{\mbf{b}} \, \big\} \leq x \: \bigcap \\
        	
        	& & \hspace{2cm} \bigcap \: x \leq \text{max} \big\{ \, \underline{\mbf{a}} \, / \, \underline{\mbf{b}}, \: \underline{\mbf{a}} \, / \, \overline{\mbf{b}}, \: \overline{\mbf{a}} \, / \, \underline{\mbf{b}}, \: \overline{\mbf{a}} \, / \, \overline{\mbf{b}} \, \big\} \, \Big\} = \\ 
        
        	& = & \Big[ \, \text{min} \big\{ \, \underline{\mbf{a}} \, / \, \underline{\mbf{b}}, \: \underline{\mbf{a}} \, / \, \overline{\mbf{b}}, \: \overline{\mbf{a}} \, / \, \underline{\mbf{b}}, \: \overline{\mbf{a}} \, / \, \overline{\mbf{b}} \, \big\}, \\
        
        	& & \hspace{3cm} \text{max} \big\{ \, \underline{\mbf{a}} \, / \, \underline{\mbf{b}}, \: \underline{\mbf{a}} \, / \, \overline{\mbf{b}}, \: \overline{\mbf{a}} \, / \, \underline{\mbf{b}}, \: \overline{\mbf{a}} \, / \, \overline{\mbf{b}} \, \big\} \, \Big].
        
    	\end{array}$
    \end{center}
\end{enumerate}
    
Главным достоинством классической интервальной арифметики можно назвать простоту её реализации для использования на ЭВМ --- результирующие интервалы получаются из концов интервалов-операндов, что облегчает написание и тестирование программ. К недостаткам этой арифметики можно отнести:
\begin{enumerate}
    \item \textit{Примитивность вида интервальных границ.}
        
		В классической интервальной арифметике множество совместных значений диаграммы зависимости описывается покоординатно. Поэтому зачастую это происходит неэффективно. Подтверждением данного утверждения является то, что на практике при использовании данной арифметики часто обнаруживает себя <<эффект обёртывания>>.
         
    \item \textit{Эффект зависимости.}
        
    Классическая интервальная арифметика не позволяет учитывать какие-либо связи между интервальными параметрами в вычисляемом выражении. Для демонстрации этого факта достаточно взять произвольный интервал ширины, отличной от нуля
    \begin{center}
    
    	$\mbf{x} \in \mathbb{I}\mathbb{R} \setminus \{ \mbf{y} \, | \, \text{wid} \, \mbf{y} \, \neq 0 \}$
    	
    \end{center}
        
    и рассмотреть функцию 
    \begin{center}
    
    	$f(x) = x - x$, где $x \in \mbf{x}.$
    	
    \end{center}
    Тогда при нахождении её естественного интервального расширения получим, что   
    \begin{center}
    
    	$\begin{array}{rcl}
    	
    	f(x) & \subseteq & \mbf{f}(\mbf{x}) = \mbf{x} - \mbf{x} = \\
    	
    	& = & \big[ \, \underline{\mbf{x}}, \: \overline{\mbf{x}} \, \big] - \big[ \, \underline{\mbf{x}}, \: \overline{\mbf{x}} \, \big] = \big[ \, \underline{\mbf{x}} - \overline{\mbf{x}}, \: \overline{\mbf{x}} - \underline{\mbf{x}} \, \big] \neq [ \, 0, \: 0 \, ].
    	
    	\end{array}$
    	
    \end{center}
    
	Отдельно отметим, что     
    \begin{center}
    
    	$\text{wid}(\mbf{x} - \mbf{x}) = 2 \cdot \text{wid}(\mbf{x}).$
    	
	\end{center}
	
\end{enumerate}

\subsection{Аффинная арифметика}

\begin{figure}
\centering
    \includegraphics[width=0.65\linewidth]{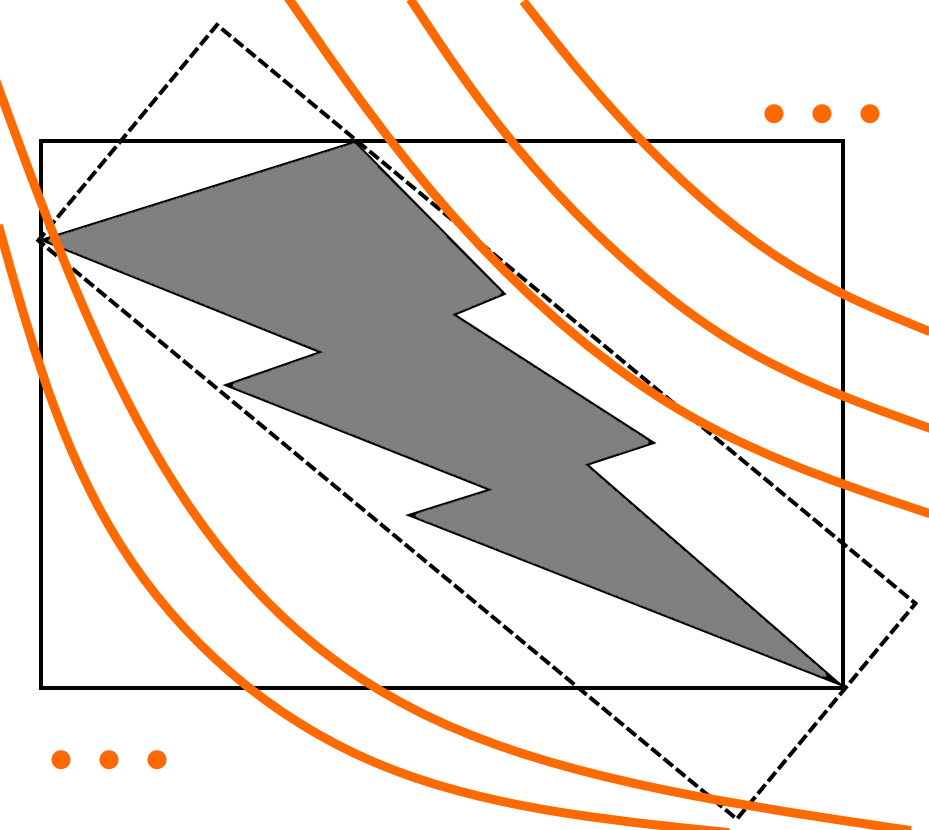}
    \caption{Пример, когда множество совместных значений диаграммы зависимости можно описать эффективнее (то есть <<ближе>> к линиям уровня оцениваемых функций), чем покоординатными границами (оранжевыми линиями обозначены линии уровня для операции умножения).}
\label{fig:affine_example_fill_set}
\end{figure}

На практике при вычислениях часто встречаются зависимые интервальные величины. Таким образом, нередка ситуация, когда множество совместных значений можно описать лучше, чем покоординатными границами (рис. \ref{fig:affine_example_fill_set}). Эта идея легла в основу построения аффинной арифметики.

В данной арифметике рассматриваются аффинные формы

\begin{center}

	$\mbf{x} = x_{0} + x_{1} \cdot \varepsilon_{1} + \dots + x_{n} \cdot \varepsilon_{n}$, 
	
	где $\varepsilon_{i} \in [ \, -1, \: 1 \, ]$, $x_{i} \in \mathbb{R}.$
	
\end{center}

\noindent$n$ будем называть \textit{размером аффинной формы} $\mbf{x}$, $\varepsilon_{i}$ --- \textit{символами шума}, а $x_{i}$ --- \textit{коэффициентами при соответствующих символах шума}. Семейство всех аффинных форм обозначим $\mathbb{A}\mathbb{R}$.

Рассмотрим две аффинных формы размера $n$:

\begin{center}

	$\mbf{x} = x_{0} + x_{1} \cdot \varepsilon_{1} + \dots + x_{n} \cdot \varepsilon_{n}$,

	$\mbf{y} = y_{0} + y_{1} \cdot \varepsilon_{1} + \dots + y_{n} \cdot \varepsilon_{n}$.

\end{center}

Рассмотрим область, образуемую множеством всех точек $(x, \: y)$ таких, что:
\begin{center}

	$x = x_{0} + x_{1} \cdot \varepsilon_{1} + \dots + x_{n} \cdot \varepsilon_{n}$,
	
	$y = y_{0} + x_{1} \cdot \varepsilon_{1} + \dots + y_{n} \cdot \varepsilon_{n}$,
	
	$\varepsilon_{i} \in [ \, -1, \: 1 \, ], \qquad i = 1, \dots, n$.

\end{center}

\begin{figure}
\centering
    \includegraphics[width=0.8\linewidth]{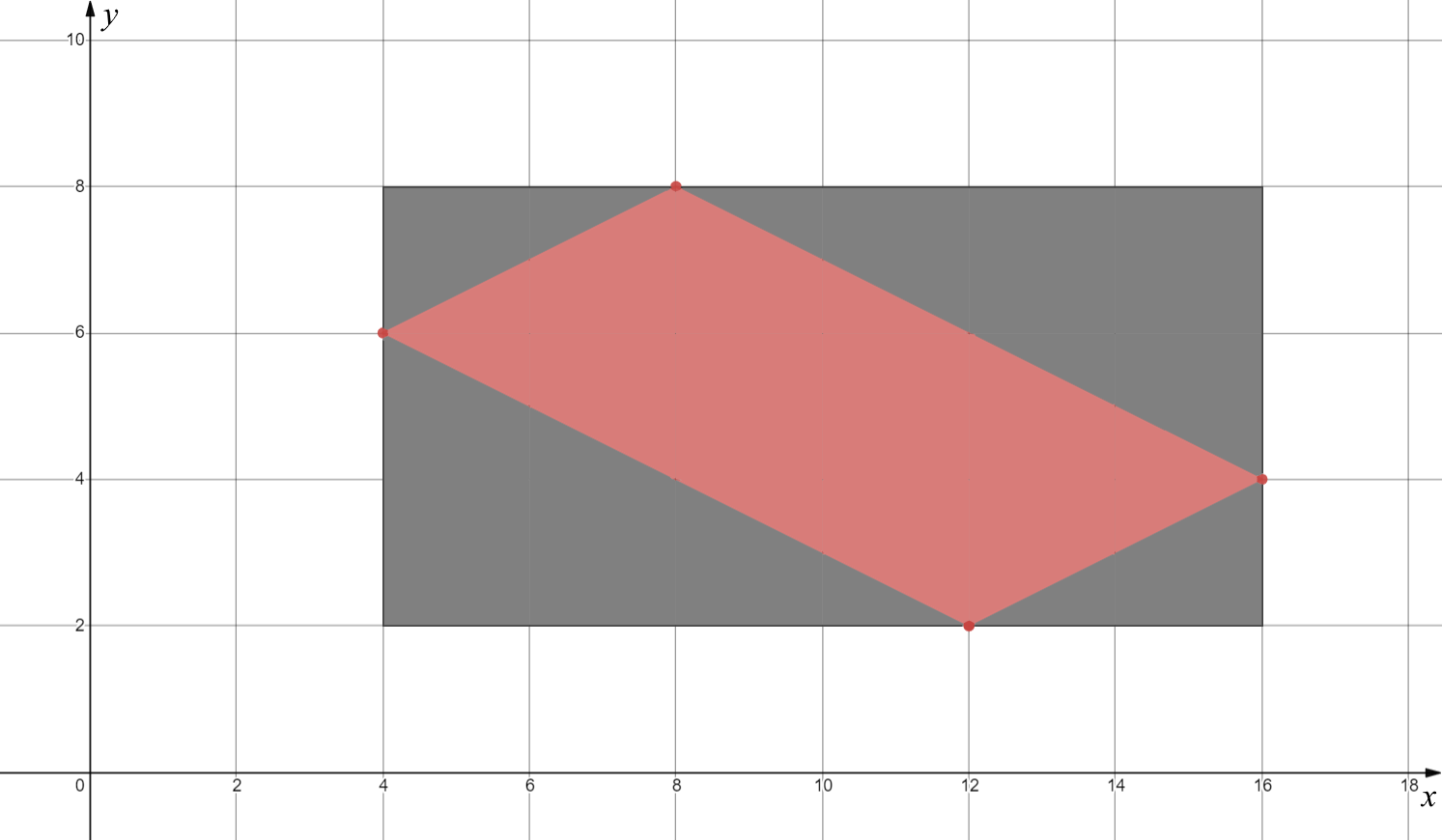}
    \caption{Диаграмма зависимости упорядоченной пары $(\mbf{x}, \: \mbf{y})$, где $\mbf{x} = 10 - 4 \cdot \varepsilon_{1} + 2 \cdot \varepsilon_{2}$, $\mbf{y} = 5 + 2 \cdot \varepsilon_{1} + 1 \cdot \varepsilon_{2}$. Розовым цветом показано множество совместных значений этой пары --- зонотоп.}
\label{fig:affine_zonotope}
\end{figure}

В результате получим центрально-симметричный параллелепипед, называемый \textit{зонотопом} (рис. \ref{fig:affine_zonotope}), имеющий $2n$ сторон.

Отметим, что аффинная форма является обобщением классического интервала (то есть $\mathbb{I}\mathbb{R} \subseteq \mathbb{A}\mathbb{R}$), поскольку любому классическому интервалу
\begin{center}

	$\mbf{x} = \big[ \, \underline{\mbf{x}}, \: \overline{\mbf{x}} \, \big] \in \mathbb{I}\mathbb{R}$
	
\end{center}
можно сопоставить аффинную форму
\begin{center}

	$\mbf{x} = \text{mid}(\mbf{x}) + \text{rad}(\mbf{x}) \cdot \varepsilon_{1} \in \mathbb{A}\mathbb{R}$.
	
\end{center}

\subsubsection{Арифметические операции аффинной арифметики}

Предположим, что

\begin{center}
	$\begin{array}{c}
	
		\mbf{x} = x_{0} + x_{1} \cdot \varepsilon_{1} + \dots + x_{n} \cdot \varepsilon_{n} \in \mathbb{A}\mathbb{R}, \\
		
	\mbf{y} = y_{0} + y_{1} \cdot \varepsilon_{1} + \dots + y_{n} \cdot \varepsilon_{n} \in \mathbb{A}\mathbb{R}, \\
	
	\alpha, \: \beta, \: \gamma \in \mathbb{R}.
	
	\end{array}$
\end{center}

Тогда арифметические операции в аффинной арифметике вводятся следующим образом:

\begin{enumerate}

    \item Линейные операции:    
    \begin{center}
    	$\begin{array}{rcl}
    	
    	\alpha \mbf{x} + \beta \mbf{y} + \gamma & = & (\alpha x_{0} + \beta y_{0} + \gamma) + (\alpha x_{1} + \beta y_{1}) \cdot \varepsilon_{1} \: + \\
    	
    	& & \hspace{5cm} + \: \dots + (\alpha x_{n} + \beta y_{n}) \cdot \varepsilon_{n}.
    	
    	\end{array}$
    \end{center}
    
    \item Умножение аффинных форм:
    \begin{center}
    	$\begin{array}{rcl}
    
        	\alpha \mbf{x} \cdot \beta \mbf{y} & = & (\alpha x_{0} \cdot \beta y_{0}) + (\alpha x_{0} \cdot \beta y_{1} + \alpha x_{1} \cdot \beta y_{0}) \cdot \varepsilon_{1} + \dots + \\
        
        	& & \hspace{2mm} + \: (\alpha x_{n} \cdot \beta y_{0} + \alpha x_{0} \cdot \beta y_{n}) \cdot \varepsilon_{n} +\\
        
        	& & \hspace{4mm} + \: $\color{blue}\setlength{\fboxsep}{1.5mm}\fbox{\color{black}{$(x_{1} \cdot \varepsilon_{1} + \dots + x_{n} \cdot \varepsilon_{n}) \cdot (y_{1} \cdot \varepsilon_{1} + \dots + y_{n} \cdot \varepsilon_{n})$}} \color{black}$ = \\
        
        	& = & (\alpha x_{0} \cdot \beta y_{0}) + (\alpha x_{0} \cdot \beta y_{1} + \alpha x_{1} \cdot \beta y_{0}) \cdot \varepsilon_{1} + \dots + \\
        
        	& & \hspace{5mm} + \: (\alpha x_{n} \cdot \beta y_{0} + \alpha x_{0} \cdot \beta y_{n}) \cdot \varepsilon_{n} + \: $\color{blue}\setlength{\fboxsep}{1.5mm}\fbox{\color{black}{$\mbf{Q}$}} \color{black}$ = \\
        
        	& = & (\alpha x_{0} \cdot \beta y_{0} + \text{mid}(\mbf{Q})) + (\alpha x_{0} \cdot \beta y_{1} + \alpha x_{1} \cdot \beta y_{0}) \cdot \varepsilon_{1} \: + \\
        
        	& & \hspace{5mm} + \dots + \: (\alpha x_{n} \cdot \beta y_{0} + \alpha x_{0} \cdot \beta y_{n}) \cdot \varepsilon_{n} + \text{rad}(\mbf{Q}) \cdot \varepsilon_{extra}.
        
    	\end{array}$
    \end{center}
    
Здесь $\mbf{Q}$ --- интервальная оценка квадратичных членов в классической интервальной арифметике (в преобразованиях обведены прямоугольником), $\varepsilon_{extra}$ --- специальный фиктивный символ шума, не ассоциированный ни с каким интервальным параметром. Данный символ характеризует ширину интервала погрешности $\mbf{Q}$, которую нельзя представить, используя символы шума. 
     
\end{enumerate}

Мы рассматривали арифметические операции в предположении, что $\mbf{x}$ и $\mbf{y}$ --- аффинные формы размера $n$. Если же формы имеют разный размер, то необходимо проделать следующую процедуру уравнивания длин. 

Пусть $X$ --- множество символов шума аффинной формы $\mbf{x}$, $Y$ --- формы $\mbf{y}$ соответственно, а $Z = X \cap Y$. Тогда можно представить $\mbf{x}$ и $\mbf{y}$ в новом расширенном виде
\begin{center}
	$\begin{array}{l}
	
		\mbf{x} \leftarrow \sum_{x_{i} \in X} x_{i} \cdot \varepsilon_{i} + \sum_{x_{i} \in Y \setminus Z} 0 \cdot \varepsilon_{i}, \\
		
		\mbf{y} \leftarrow \sum_{x_{i} \in X \setminus Z} 0 \cdot \varepsilon_{i} + \sum_{x_{i} \in Y} x_{i} \cdot \varepsilon_{i}.
		
	\end{array}$
\end{center}

То есть в данном виде все коэффициенты при символах шума, которые не встречались в рассматриваемой форме, равны нулю.

\subsection[Другие интервальные арифметики]{Другие интервальные арифметики}

Существуют и другие виды интервальных арифметик. Они  не вошли в обзор, поскольку имеют особенности, которые по некоторым причинам не позволяют применять их для решения поставленных в работе задач:

\begin{enumerate}

    \item Арифметика Каухера --- интервальная арифметика, которая была разработана, как расширение классической интервальной арифметики для пополнения алгебраических свойств.
    
    \item Комплексные арифметики --- используются для оперирования интервальными объектами, которые содержат в себе комплексные числа.
    
    \item Твинная арифметика --- интервальная арифметика, которая оперирует интервалами с интервальными концами.
    
    \item Арифметика Кахана --- интервальная арифметика, которая доопределяет классическую интервальную арифметику операцией деления на нуль-содержащий интервал, а также определяет операции между интервалами, концы которых могут быть <<$-\infty$>> или <<$+\infty$>>.
    
    \item Мультиинтервальная арифметика --- интервальная арифметика, идея которой состоит в представлении оперируемого объекта в виде объединения классических интервалов.
    
\end{enumerate}

\clearpage
\section [Анализ аффинной арифметики] {Анализ аффинной арифметики}

Использование аффинных форм позволяет группировать члены при символах шума во время вычислений. Это свойство позволяет учитывать линейные связи между ними.

Так, например, рассмотрим функцию
\begin{center}

	$f(x) = x - x, \text{ где } x \in [ \, -1, \: 1 \, ].$
	
\end{center}

В аффинной арифметике любой классический интервал $\mbf{x} \in \mathbb{I}\mathbb{R}$ будет представляться в виде

\begin{center}
	$\mbf{x} = \text{mid}(\mbf{x}) + \text{rad}(\mbf{x}) \cdot \varepsilon_{1} \in \mathbb{A}\mathbb{R}$.
\end{center}

\noindentТогда

\begin{center}
	$\begin{array}{rcl}
	
f(x) \subseteq f(\mbf{x}) & = & \mbf{x} - \mbf{x} = \\

 & = & \big( \text{mid}(\mbf{x}) - \text{mid}(\mbf{x}) \big) + \big( \text{rad}(\mbf{x}) - \text{rad}(\mbf{x}) \big) \cdot \varepsilon_{1} = \\
 
 & = & 0 + 0 \cdot \varepsilon_{1} = [ \, 0, \: 0 \, ].
 
	\end{array}$
\end{center}

Таким образом, использование афинной арифметики позволяет учесть группировку по суммированию и вычитанию при символах шума, тем самым исправляя недостаток классической интервальной арифметики, в которой
\begin{center}

	$\forall \mbf{x} \in \mathbb{I}\mathbb{R} \setminus \{ \mbf{y} \, | \, \text{wid} \, \mbf{y} = 0 \} : \: \mbf{x} - \mbf{x} \neq [ \, 0, \: 0 \,].$
	
\end{center}
Однако данная арифметика имеет и недостатки:
\begin{enumerate}

    \item \textit{Реализация аффинной арифметики более трудоёмка по сравнению с классической интервальной арифметикой.}
    
    Вместо операций с парой вещественных чисел --- концами интервалов, необходимо реализовывать арифметические действия между массивами коэффициентов при соответствующих символах шума.
    
    \item \textit{Особенность, позволяющая учитывать группировку членов при символах шума в вычислениях влечёт за собой увеличение алгоритмической сложности и потребления памяти.} 
    
    Исходя из определения операций, сложение, вычитание и умножение аффинных форм размера $n$ требуют объёма памяти $O(n)$, и выполнения $O(n)$ операций.
    
    \item \textit{При обработке квадратичных членов, мы преобразуем их в классические интервалы, а учитываем только линейные члены при символах шума.}
    
    То есть, в аффинной арифметике верно
    \begin{center}
    
    	$\forall \mbf{x} \in \mathbb{I}\mathbb{R}: \mbf{x} - \mbf{x} = [ \, 0, \: 0 \, ],$
    	
    \end{center}
    но при этом
    \begin{center}
    
    	$\forall \mbf{x} \in \mathbb{I}\mathbb{R} \setminus \{ \, \mbf{y} \, | \, \text{wid} \, \mbf{y} = 0 \, \} : \mbf{x}^{2} - \mbf{x}^{2} \neq [ \, 0, \: 0 \, ].$
   
   \end{center}
    
    \item \textit{Существуют примеры, когда использование аффинной арифметики даёт результирующий интервал шире, чем при использовании классической интервальной арифметики.} 
    
    Рассмотрим функцию
    
    \begin{equation}
    	\label{eq:bad_affine_example}    	
    	f(x) = \big( x - \tfrac{1}{2} \big) ^ {2} \text{, где } x \in [ \, -1, \: 1 \, ].    
    \end{equation}
    
    Аффинная форма для интервального расширения этой функции будет иметь вид
\begin{center}

    $\begin{array}{rcl}
    
    	f(x) & \subseteq & f(\mbf{x}) = \big( \mbf{x} - \frac{1}{2} \big) ^ {2} =  \big( \mbf{x} - \frac{1}{2} \big) \cdot \big( \mbf{x} - \frac{1}{2} \big) =\\ 
    	
    	& = & \big( -\frac{1}{2} + 1 \cdot \varepsilon_{1} \big) \cdot \big( -\frac{1}{2} + 1 \cdot \varepsilon_{1} \big) = \\
    
    	& = & \frac{1}{4} - 1 \cdot \varepsilon_{1} + (1 \cdot \varepsilon_{1}) \cdot (1 \cdot \varepsilon_{1}) = \\
    	
    	& = & \frac{1}{4} - 1 \cdot \varepsilon_{1} + [ \, 0, \: 1 \, ] = \\
    	
    	& = & \frac{1}{4} - 1 \cdot \varepsilon_{1} + \frac{1}{2} + \frac{1}{2} \cdot \varepsilon_{\textit{extra}} = \\
    
    	& = & \frac{3}{4} - 1 \cdot \varepsilon_{1} + \frac{1}{2} \cdot \varepsilon_{\textit{extra}}.
    
    \end{array}$
\end{center}
    
    Преобразуем аффинную форму в классический интервал, получим
    \begin{center}
    
        $\frac{3}{4} - 1 \cdot [ \, -1, \: 1 \, ] + \big[ \, -\frac{1}{2}, \: \frac{1}{2} \, \big] = \big[ \, -\frac{3}{4}, \: \frac{9}{4} \, \big].$
        
    \end{center}
    
     Теперь используем классическую интервальную арифметику для нахождения естественного интервального расширения $\mbf{f}(\mbf{x})$, тогда получим   
     \begin{center}
     
         $\big( \mbf{x} - \frac{1}{2} \big) ^ {2} = \big[ \, -\frac{3}{2}, \: \frac{1}{2} \, \big]^{2} = \big[ \, 0, \: \frac{9}{4} \, \big].$
         
     \end{center}
    
    Данный пример демонстрирует случай, когда при использовании классической интервальной арифметики результирующий классический интервал получается уже, чем при использовании аффинной арифметики (рис. \ref{fig:n11}).
    
    Это получается из-за того, что результирующий интервал, который используется на практике, представляет собой классический интервал, а его конвертация в какой-либо иной формат и наоборот может сопровождаться погрешностью, которая будет больше, чем выигрыш от учёта группирования коэффициентов аффинной формы при соответствующих символах шума.
    
    \begin{figure}
    \centering
        \includegraphics[width=0.75\linewidth]{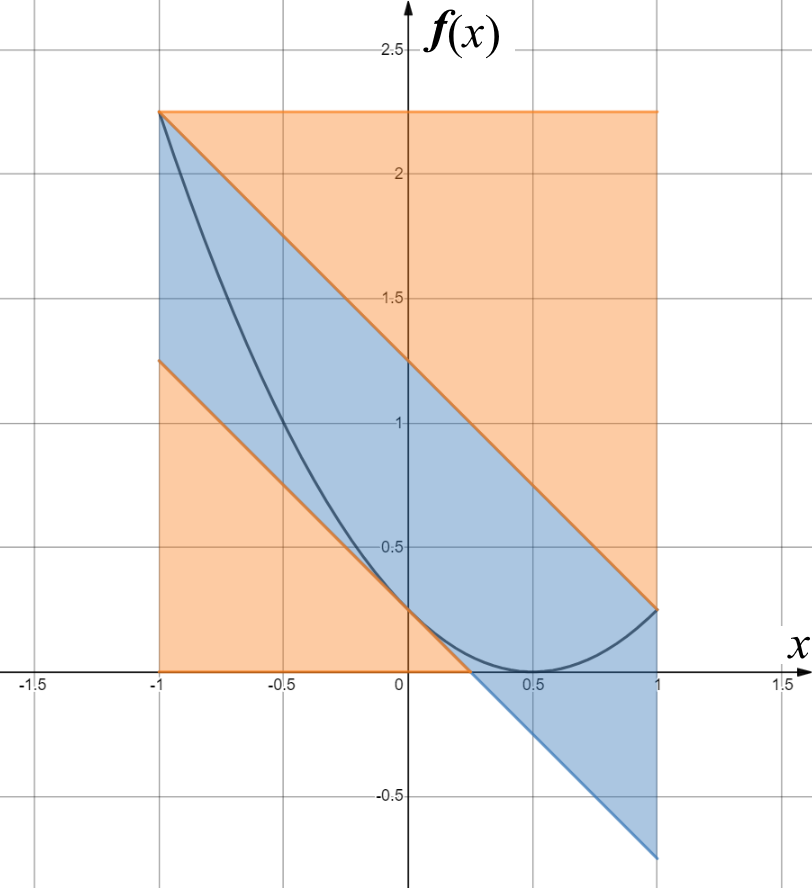}
        \caption{Пример, когда классическая интервальная арифметика даёт в результате вычислений классический интервал более узкий, чем при использовании аффинной арифметики. Чёрной линией показано множество совместных значений упорядоченной пары $\big( x, \: f(x) \big)$, где $f(x)$ задана в (\ref{eq:bad_affine_example}). Оранжевой областью показано его приближение в классической интервальной арифметике (интервал значений $\big[ \, 0, \: \frac{9}{4} \, \big]$), голубой областью --- в аффинной арифметике (интервал значений $\big[ \, -\frac{3}{4}, \: \frac{9}{4} \, \big]$).}
    	\label{fig:n11}
    \end{figure}
    
    Идея исправления данного недостатка была предложена Ахмеровым в работе \cite{AchmerovZonotopesQuestion}. Она состоит в использовании информации от вычислений в классической интервальной арифметике при каждом шаге обработки аффинных форм. Полученная модификация аффинной арифметики была названа \textit{аффинно-интервальной}.
       
    \item \textit{Аффинная арифметика позволяет относительно просто вводить только некоторые операции, а именно: сложение, вычитание, произведение, умножение на константу.}
    
    Для введения иных неаффинных операций между аффинными формами требуется создание специальных процедур, которые будут каким-либо оптимальным образом описывать результирующее множество зонотопом.
    
     В Работе Ахмерова \cite{AchmerovZonotopesQuestion} была упомянута процедура построения зонотопа, которая приблизительно наилучшим образом аппроксимирует диаграмму зависимости для результатов деления за время $O(n)$. Однако построение наилучшего приближения этого множества зонотопами остаётся открытым вопросом.
     
\end{enumerate}

\subsection[Пример неэффективности аффинной арифметики]{Пример неэффективности аффинной\\арифметики}

Чтобы продемонстрировать важность учета и вовлечения в вычисления квадратичных членов, рассмотрим функцию

\begin{center}
    $\begin{array}{c}
    
    f = (x - a) \cdot (x - b) \cdot (x - c), \\
    
    x \in [ \, -1, \: 1 \, ], \: a \in \big[ \, -\frac{1}{2}, \: \frac{1}{2} \, \big], \: b \in \big[ \, -\frac{1}{2}, \: \frac{1}{2} \, \big], \: c \in \big[ \, -\frac{1}{2}, \: \frac{1}{2} \, \big].
    
    \end{array}$
\end{center}

\noindentНайдём естественное интервальное расширение

\begin{center}
    $\begin{array}{rcl}
    
        f(x) & = & (x - a) \cdot (x - b) \cdot (x - c) \subseteq (\mbf{x} - \mbf{a}) \cdot (\mbf{x} - \mbf{b}) \cdot (\mbf{x} - \mbf{c}) = \\
        
        & = & \Big( \big[ \, -1, \: 1 \, \big] - \big[ \, -\frac{1}{2}, \: \frac{1}{2} \, \big] \Big) \: \cdot \\
        
        & & \hspace{2cm} \cdot \: \Big( \big[ \, -1, \: 1 \, \big] - \big[ \, -\frac{1}{2}, \: \frac{1}{2} \, \big] \Big) \: \cdot \\
        
        & & \hspace{4cm} \cdot \: \Big( \big[ \, -1, \: 1 \, \big] - \big[ \, -\frac{1}{2}, \: \frac{1}{2} \, \big] \Big) =\\
        
        & = &  \big[ \, -\frac{3}{2}, \: \frac{3}{2} \, \big] \cdot \big[ \, -\frac{3}{2}, \: \frac{3}{2} \, \big] \cdot \big[ \, -\frac{3}{2}, \: \frac{3}{2} \, \big] = \big[ \, -\frac{9}{4}, \: \frac{9}{4} \, \big] \cdot \big[ \, -\frac{3}{2}, \: \frac{3}{2} \, \big] = \\
        
        & = & \big[ \, -\frac{27}{8}, \: \frac{27}{8} \, \big]
        
    \end{array}$ 
\end{center}

\noindentВ аффинной арифметике

\begin{center}
	$\begin{array}{rcl}

    	\mbf{x} - \mbf{a} & = & 0 + 1 \cdot \varepsilon_{1} + \frac{1}{2} \cdot \varepsilon_{2}, \\
    
    	\mbf{x} - \mbf{b} & = & 0 + 1 \cdot \varepsilon_{1} + \frac{1}{2} \cdot \varepsilon_{3}, \\
    
    	\mbf{x} - \mbf{c} & = & 0 + 1 \cdot \varepsilon_{1} + \frac{1}{2} \cdot \varepsilon_{4}.
    
	\end{array}$
\end{center}

Итого имеем

\begin{center}
	$\begin{array}{rcl}
	
		f(x) & = & (x - a) \cdot (x - b) \cdot (x - c) \subseteq (\mbf{x} - \mbf{a}) \cdot (\mbf{x} - \mbf{b}) \cdot (\mbf{x} - \mbf{c}) =\\

		& = & \big( 0 + 1 \cdot \varepsilon_{1} + \frac{1}{2} \cdot \varepsilon_{2} \big) \cdot \big( 0 + 1 \cdot \varepsilon_{1} + \frac{1}{2} \cdot \varepsilon_{3} \big) \cdot \big( 0 + 1 \cdot \varepsilon_{1} + \frac{1}{2} \cdot \varepsilon_{4} \big) =\\

		& = & \big[ \, -\frac{3}{2}, \: \frac{3}{2} \, \big] \cdot \big[ \, -\frac{3}{2}, \: \frac{3}{2} \, \big] \cdot \big[ \, -\frac{3}{2}, \: \frac{3}{2} \, \big] = \big[ \, -\frac{9}{4}, \: \frac{9}{4} \, \big] \cdot \big[ \, -\frac{3}{2}, \: \frac{3}{2} \, \big] = \\

		& = & \big[ \, -\frac{27}{8}, \: \frac{27}{8} \, \big]. 

	\end{array}$
\end{center}

Получаем одинаковый результат в классической интервальной и аффинной арифметиках.

Теперь рассмотрим следующее выражение в аффинной арифметике, которое получается из предыдущего уточнением интервальных параметров $\mbf{a}, \: \mbf{b}, \: \mbf{c}$:
\begin{center}

    $f(x) = \big( x - \frac{1}{2} \big) \cdot \big( x - \frac{1}{2} \big) \cdot \big( x - \frac{1}{2} \big), \: x \in [ \, -1, \: 1 \, ]$
    
\end{center}

\noindentНайдём аффинную форму данного выражения:

\begin{center}
	$\begin{array}{rcl}

    	f(x) & \subseteq & f(\mbf{x}) = \big( \mbf{x} - \frac{1}{2} \big) \cdot \big( \mbf{x} - \frac{1}{2} \big) \cdot \big( \mbf{x} - \frac{1}{2} \big) = \\
    
    	& = & \big( -\frac{1}{2} + 1 \cdot \varepsilon_{1} \big) \cdot \big( -\frac{1}{2} + 1 \cdot \varepsilon_{1} \big) \cdot \big( -\frac{1}{2} + 1 \cdot \varepsilon_{1} \big) = \\
    
    	& = & \big( \frac{3}{4} - 1 \cdot \varepsilon_{1} + \frac{1}{2} \cdot \varepsilon_{\textit{extra}} \big) \cdot \big( -\frac{1}{2} + 1 \cdot \varepsilon_{1} \big) = \\ 
    
    	& = & \big( \frac{3}{4} - 1 \cdot \varepsilon_{1} + \frac{1}{2} \cdot \varepsilon_{\textit{extra}} \big) \cdot \big( -\frac{1}{2} \big) + \frac{3}{4} \cdot \varepsilon_{1} + \\
    	
    	& & \hspace{4cm} + \: \big( -1 \cdot \varepsilon_{1} + \frac{1}{2} \cdot \varepsilon_{\textit{extra}} \big) \cdot ( 1 \cdot \varepsilon_{1} ) = \\
    
    	& = & \big( -\frac{3}{8} + \frac{5}{4} \cdot \varepsilon_{1} \big) + \frac{1}{4} \cdot \varepsilon_{\textit{extra}} + \big( -1 \cdot \varepsilon_{1} + \frac{1}{2} \cdot \varepsilon_{\textit{extra}} \big) \cdot ( 1 \cdot \varepsilon_{1} ) = \\
    
    	& = & \big( -\frac{3}{8} + \frac{5}{4} \cdot \varepsilon_{1} \big) + \frac{1}{4} \cdot [ \, -1, \: 1 \, ] + \Big( [ \, -1, \: 1 \, ] + \big[ \, -\frac{1}{2}, \: \frac{1}{2} \, \big] \Big) \cdot [ \, -1, \: 1 \, ] = \\
    
    	& = & \big( -\frac{3}{8} + \frac{5}{4} \cdot \varepsilon_{1} \big) + \big[ \, -\frac{1}{4}, \: \frac{1}{4} \, \big] + \big[ \, -\frac{3}{2}, \: \frac{3}{2} \, \big] = \\
    
    	& = & \big( -\frac{3}{8} + \frac{5}{4} \cdot \varepsilon_{1} \big) + \frac{7}{4} \cdot \varepsilon_{\textit{extra}}
    
	\end{array}$
\end{center}

Проанализируем полученный результат. Посчитаем площадь $S_{\textit{effective}}$, заметаемую линейными членами интервальных параметров, которые мы удерживаем в выражении. Также посчитаем и площадь $S_{\textit{extra}}$ от неассоциированного с интервальными параметрами символа шума.

\begin{center}
	$\begin{array}{c}
	
    	S_{\textit{effective}} = \int_{-1}^{+1} \frac{5}{4} \: d\varepsilon_{1} = \frac{5}{2}, \\
    	
    	S_{\textit{extra}} = \int_{-1}^{+1} \frac{7}{4} \: d\varepsilon_{2} = \frac{7}{2}.

    \end{array}$
\end{center}

Приведённый пример показывает важность учета квадратичных членов при вычислении выражений. Был приведён пример, когда всего за $3$ операции умножения площадь от фиктивного символа шума погрешности $\varepsilon_{\textit{extra}}$ начинает превалировать над <<полезной площадью>> $S_{\textit{effective}}$.

\clearpage
\section [Построение новой интервальной арифметики] {Построение новой\\интервальной арифметики}

Пусть в рассмотрении имеется аффинная форма
\begin{center}

	$\mbf{x}(\varepsilon_{1}, \dots, \varepsilon_{n}) = x_{0} + x_{1} \cdot \varepsilon_{1} + \dots + \varepsilon_{n} + x_{\textit{extra}} \cdot \varepsilon_{\textit{extra}}.$

\end{center}

\begin{figure}
	\centering
    \includegraphics[width=0.6\linewidth]{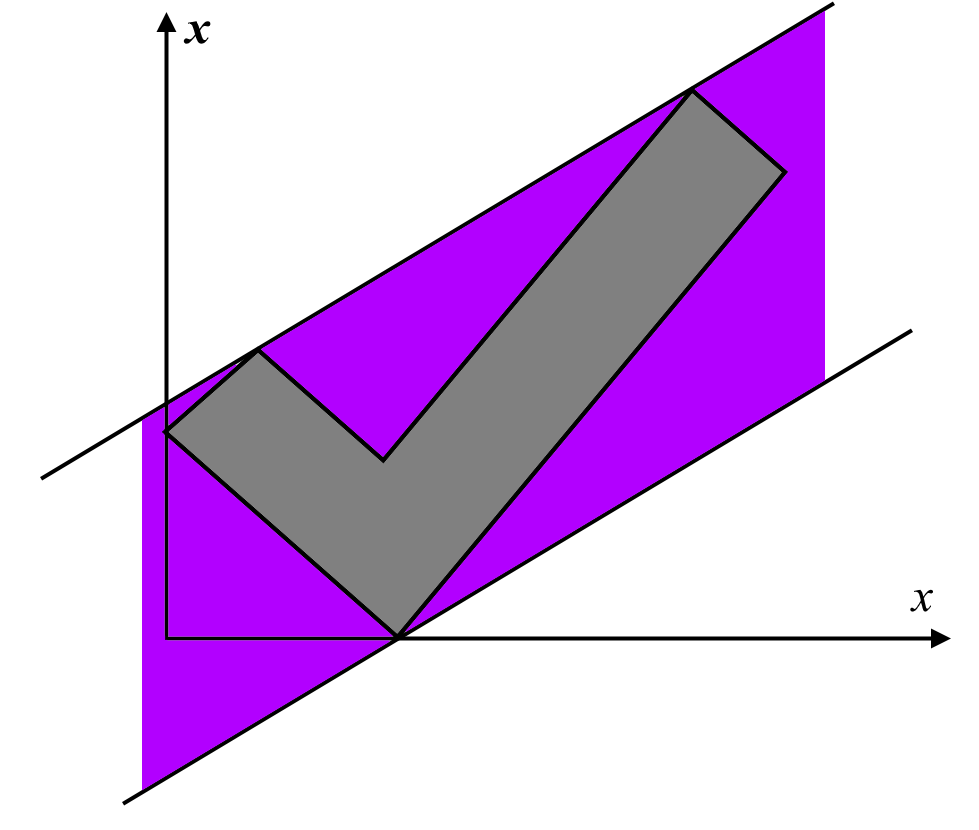}
    \caption{Множество совместных значений для упорядоченой пары $\big( x, \: \mbf{x}(\varepsilon_{1}) \big)$, где $\mbf{x}(\varepsilon_{1}) = x_{0} + x_{1} \cdot \varepsilon_{1} + x_{\textit{extra}} \cdot \varepsilon_{\textit{extra}}$.}
	\label{fig:galochka_affine}
\end{figure}

\noindentИдеей построения новой интервальной арифметики стало наблюдение, что аффинная арифметика позволяет описывать множество совместных значений упорядоченной пары 
\begin{center}

	$\big( (\varepsilon_{1}, \dots, \varepsilon_{n}),  \: \mbf{x}(\varepsilon_{1}, \dots, \varepsilon_{n}) \big)$

\end{center}
исключительно линейными членами при соответствующих символах шума. Тогда множества, представляющие собой, например, вогнутые фигуры, будут описываться с большой погрешностью (рис. \ref{fig:galochka_affine})

Также аффинная арифметика не предусматривает учёта квадратичных членов и в случае оперирования ими сводится к использованию классической интервальной арифметики. Значит, она мало пригодна к использованию на практике при большом числе умножений в процессе вычислений.

Автором данной работы была выдвинута идея, которая позволяет частично устранить недостатки, указанные выше. Она состоит в создании новой интервальной арифметики, в которой концы интервалов описываются при помощи функционалов границ.

Будем мыслить интервал, как параметрическое семейство классических интервалов. То есть $\mbf{x} = \mbf{x}(x_{1}, \dots, x_{n})$, где $x_{i} \in \mbf{x}_{i}$ интервальный параметры.

Функционал $L(l_{1}, \dots, l_{n_{l}})$ назовём \textit{функционалом нижней границы}, если $l_{i} \in \mbf{l}_{i} \: (i = 1, \dots, l_{n})$ --- интервальные параметры и
\begin{center}

	$L(l_{1}, \dots, l_{n_{l}}) \leq \underline{\mbf{x}(x_{1}, \dots, x_{n})}$,
	
	$\begin{array}{rl}
	
		l_{i} \in \mbf{l}_{i}, & i = 1, \dots, n_{l}, \\
		
		x_{k} \in \mbf{x}_{k}, & k = 1, \dots, n.
		
	\end{array}$

\end{center}

Функционал $U(u_{1}, \dots, u_{n_{u}})$ назовём \textit{функционалом верхней границы}, если $u_{i} \in \mbf{u}_{i} \: (i = 1, \dots, u_{n})$ --- интервальные параметры и
\begin{center}

	$U(u_{1}, \dots, u_{n_{u}}) \geq \overline{\mbf{x}(x_{1}, \dots, x_{n})},$
	
	$\begin{array}{rl}
	
		u_{i} \in \mbf{u}_{i}, & i = 1, \dots, n_{u}, \\
		
		x_{k} \in \mbf{x}_{k}, & k = 1, \dots, n.	
	
	\end{array}$

\end{center}

Интервал $\mbf{x}$, границы которого представляются функционалами нижней и верхней границы соответственно, будем называть \textit{функционально-граничным}, а арифметику, которая оперирует интервалами такого типа, будем называть \textit{функционально-граничной}. Семейство функционально-граничных интервалов будем обозначать $\mathbb{F}\mathbb{R}$.

\begin{figure}
\centering
    \includegraphics[width=1\linewidth]{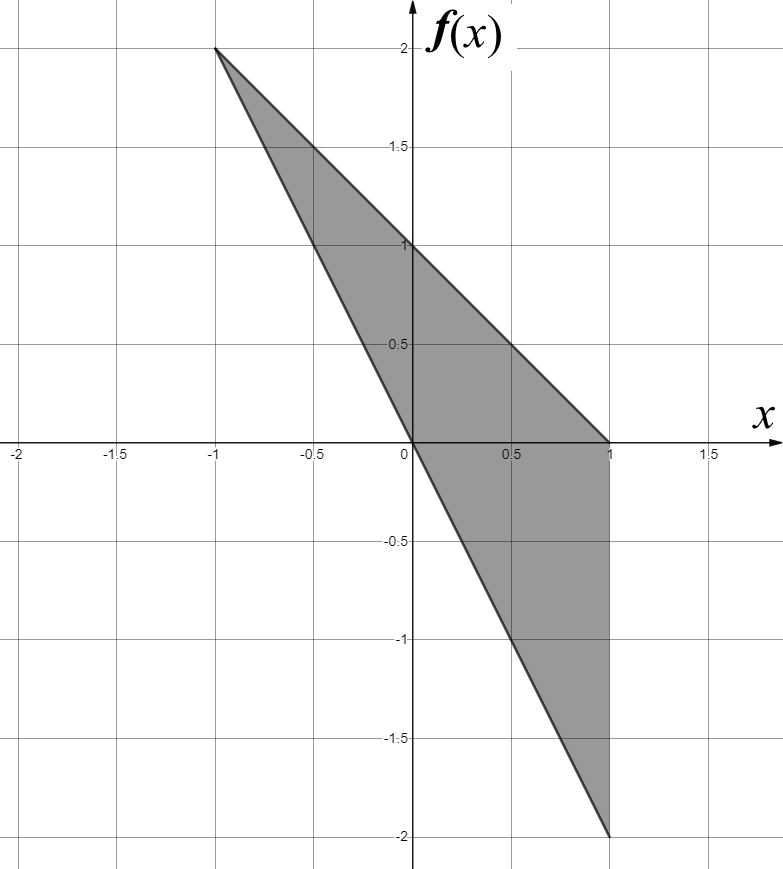}
    \caption{Множество совместных значений упорядоченной пары $\big( x, \: \mbf{f}(x) \big)$, где $\mbf{f}(x) = \big[ \, -2x, \: -x + 1 \, \big] \in \mathbb{F}\mathbb{R}, \: x \in [ \, -1, \: 1 \, ]$.}
\label{fig:example1}
\end{figure}

\begin{figure}
\centering
    \includegraphics[width=1\linewidth]{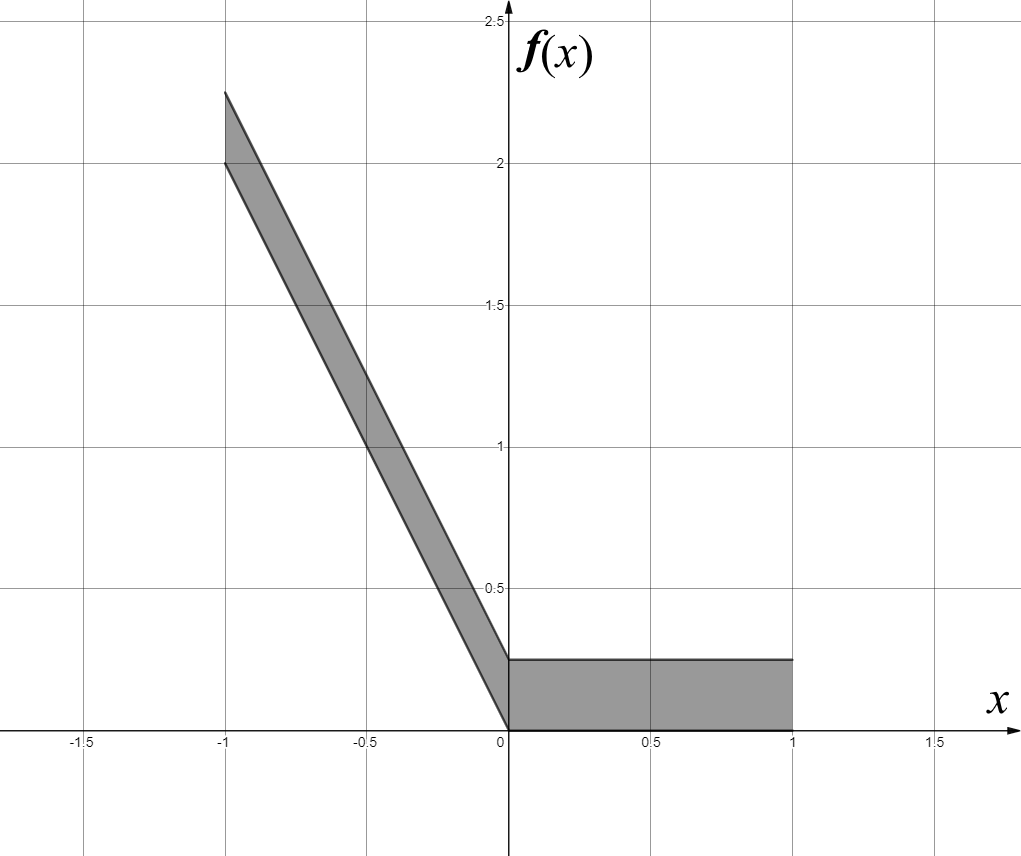}
    \caption{Множество совместных значений упорядоченной пары $\big( x, \: \mbf{f}(x) \big)$, где $\mbf{f}(x) = \big[ \, |x| - x, \: |x| - x + \frac{1}{4} \, \big] \in \mathbb{F}\mathbb{R}$.}
\label{fig:example2}
\end{figure}

Для записи функционально-граничных интервалов предлагается следующая нотация: вместо числовых границ интервалов будем писать на месте левого конца интервала функционал $L(l_{1}, \dots, l_{n_{l}})$, а вместо правого конца --- функционал $U(u_{1}, \dots, u_{n_{u}})$. Данные обозначения были использованы в описании рисунков (рис. \ref{fig:example1},  \ref{fig:example2}).

\textit{Интервальной оценкой} функции $f(x_{1}, \dots, x_{n})$ будем называть функционально-граничный интервал $\Big[ \, L(l_{1}, \dots, l_{n_{l}}), \: U(u_{1}, \dots, u_{n_{u}}) \, \Big]$, если
\begin{center}

	$L(l_{1}, \dots, l_{n_{l}}) \leq f(x_{1}, \dots, x_{n}),$
	
	$U(u_{1}, \dots, u_{n_{u}}) \geq f(x_{1}, \dots, x_{n}),$
	
	$\begin{array}{rl}
	
		l_{i} \in \mbf{l}_{i}, & i = 1, \dots, n_{l},	\\
	
		u_{j} \in \mbf{u}_{j}, & j = 1, \dots, n_{u}, \\
	
		x_{k} \in \mbf{x}_{k}, & k = 1, \dots, n.
	
	\end{array}$

\end{center}

Обозначать интервальную оценку функции будем жирным символом этой функции. Чтобы понимать, к какой интервальной оценке относится тот или иной функционал $L(l_{1}, \dots, l_{n_{l}})$ и $U(u_{1}, \dots, u_{n_{u}})$, будем писать обозначение этой интервальной оценки нижним индексом возле соответствующего функционала. В данных обозначениях, $\mbf{f}(x_{1}, \dots, x_{n})$ --- интервальная оценка для функции $f(x_{1}, \dots, x_{n})$, причём
\begin{center}

	$\mbf{f}(x_{1}, \dots, x_{n}) = \Big[ \, L_{\!\mbf{f}}(l_{1}, \dots, l_{n_{l}}), \: U_{\!\mbf{f}}(u_{1}, \dots, u_{n_{u}}) \, \Big]$.
	
\end{center}

Далее в работе будем полагать, что при нахождении интервальной оценки некоторой функции $f(x_{1}, \dots, x_{n})$ функционалы границ зависят только от интервальных параметров этой функции, то есть 
\begin{center}

	$L_{\!\mbf{f}} = L_{\!\mbf{f}}(x_{1}, \dots, x_{n})$ и $U_{\!\mbf{f}} = U_{\!\mbf{f}}(x_{1}, \dots, x_{n})$.
	
\end{center}
Также будем полагать, что
\begin{center}

	$x_{i} \in [ \, -1, \: 1 \, ], \qquad i = 1, \dots, n.$

\end{center}
В случае, если параметр функции изменяется в интервале $x^{*} \in \mbf{x}^{*} \neq [ \, -1, \: 1 \, ]$, то можно применить линейное преобразование:
\begin{center}

	$x^{*} \leftarrow \text{mid}(\mbf{x}^{*}) + x \cdot \text{rad}(\mbf{x}^{*}),$

\end{center}
которое позволяет свести интервал изменения переменной $x^{*}$ к интервалу $[ \, -1, \: 1 \, ]$.

Заметим, что в рамках данных соглашений, $\mathbb{I}\mathbb{R} \subseteq \mathbb{F}\mathbb{R}$, так как любой классический интервал
\begin{center}

	$\mbf{x} = \big[ \, \underline{\mbf{x}}, \: \overline{\mbf{x}} \, \big] \in \mathbb{I}\mathbb{R}$
	
\end{center}
представляет собой функционально-граничный интервал
\begin{center}

	$\mbf{x} \big( \underline{\mbf{x}}, \: \overline{\mbf{x}} \big) = \Big[ \, L_{\mbf{x}} \big( \underline{\mbf{x}}, \: \overline{\mbf{x}} \big), \: U_{\mbf{x}} \big( \underline{\mbf{x}}, \: \overline{\mbf{x}} \big) \, \Big] \in \mathbb{F}\mathbb{R}$,
	
	$L_{\mbf{x}} \big( \underline{\mbf{x}}, \: \overline{\mbf{x}} \big) = \underline{\mbf{x}}, \qquad U_{\mbf{x}} \big( \underline{\mbf{x}}, \: \overline{\mbf{x}} \big) = \overline{\mbf{x}}$.
	
\end{center}

Также $\mathbb{A}\mathbb{R} \subseteq \mathbb{F}\mathbb{R}$, так как любую аффинную форму
\begin{center}

	$\mbf{x} = x_{0} + x_{1} \cdot \varepsilon_{1} + \dots + x_{n} \cdot \varepsilon_{n} + x_{\textit{extra}} \cdot \varepsilon_{\textit{extra}} \in \mathbb{A}\mathbb{R}$
	
\end{center}
можно представить как функционально-граничный интервал
\begin{center}
	
	$\mbf{x} \big( \varepsilon_{1}, \dots, \varepsilon_{n} \big) = \Big[ \, L_{\mbf{x}} \big( \varepsilon_{1}, \dots, \varepsilon_{n} \big), \: U_{\mbf{x}} \big( \varepsilon_{1}, \dots, \varepsilon_{n} \big) \, \Big] \in \mathbb{F}\mathbb{R}$,

	$L_{\mbf{x}} \big( \varepsilon_{1}, \dots, \varepsilon_{n} \big) = x_{0} + x_{1} \cdot \varepsilon_{1} + \dots + x_{n} \cdot \varepsilon_{n} - |x_{\textit{extra}}|,$
	
	$U_{\mbf{x}} \big( \varepsilon_{1}, \dots, \varepsilon_{n} \big) = x_{0} + x_{1} \cdot \varepsilon_{1} + \dots + x_{n} \cdot \varepsilon_{n} + |x_{\textit{extra}}|.$
	
\end{center}

Новизной функционально-граничной интервальной арифметики является то, что в зависимости от выбранного функционального базиса, в котором будут описываться функционалы $L(x_{1}, \dots, x_{n})$ и $U(x_{1}, \dots, x_{n})$, она начинает обладать специфичными свойствами.

В данной работе будет выбран такой функциональный базис, который позволит описывать множества совместных значений упорядоченной пары $\big( (x_{1}, \dots, x_{n}), \: \mbf{x}(x_{1}, \dots, x_{n}) \big)$ лучше, чем в уже рассмотренных арифметиках. Примеры некоторых таких множеств показаны на рисунках (\ref{fig:example1}, \ref{fig:example2}).

Введём арифметические операции между интервалами в функционально-граничной арифметике. Пусть имеются \begin{center}

	$\mbf{a}(x_{1}, \dots, x_{n}) \in \mathbb{F}\mathbb{R},$
	
	$\mbf{b}(x_{1}, \dots, x_{n}) \in \mathbb{F}\mathbb{R},$
	
	$\alpha \in \mathbb{R}$. 
	
\end{center}

Положим, что $F_{\downarrow} \big\{ \, f_{1}(x_{1}, \dots, x_{n}), \dots, f_{m}(x_{1}, \dots, x_{n}) \, \big\}$ --- некоторая процедура, строящая функционал $L_{\textit{new}}(x_{1}, \dots, x_{n})$ с тем свойством, что:
\begin{center}
	$\begin{array}{cl}
	
		L_{\textit{new}}(x_{1}, \dots, x_{n}) \leq f_{j}(x_{1}, \dots, x_{n}), & j = 1, \dots, m, \\
		
		x_{i} \in [ \, -1, \: 1 \, ], & i = 1, \dots, n, 
		 
	\end{array}$
\end{center} 
а $F_{\uparrow} \big\{ \, f_{1}(x_{1}, \dots, x_{n}), \dots, f_{m}(x_{1}, \dots, x_{n}) \, \big\}$ --- функционал $U_{\textit{new}}(x_{1}, \dots, x_{n})$ со свойством:
\begin{center}
	$\begin{array}{cl}
	
		U_{\textit{new}}(x_{1}, \dots, x_{n}) \geq f_{j}(x_{1}, \dots, x_{n}), & j = 1, \dots, m, \\

		x_{i} \in [ \, -1, \: 1 \, ], & i = 1, \dots, n.
		
	\end{array}$
\end{center}

Тогда операции между $\mbf{a}(x_{1}, \dots, x_{n})$ и $\mbf{b}(x_{1}, \dots, x_{n})$ можно ввести следующим образом:

\begin{enumerate}

    \item Сложение:  
        \begin{center}
        	$\begin{array}{l}
        	
				\mbf{a}(x_{1}, \dots, x_{n}) + \mbf{b}(x_{1}, \dots, x_{n}) = \\
				
				\hspace{3cm} = \Big[ \, L_{\mbf{a}}(x_{1}, \dots, x_{n}), \: U_{\mbf{a}}(x_{1}, \dots, x_{n}) \, \Big] \: + \\
				
				\hspace{6cm} + \: \Big[ \, L_{\mbf{b}}(x_{1}, \dots, x_{n}), \: U_{\mbf{b}}(x_{1}, \dots, x_{n}) \, \Big] = \\
				
				\hspace{3cm} = \Big[ \, L_{\mbf{a}}(x_{1}, \dots, x_{n}) + L_{\mbf{b}}(x_{1}, \dots, x_{n}), \\
				
				\hspace{6cm} U_{\mbf{a}}(x_{1}, \dots, x_{n}) + U_{\mbf{b}}(x_{1}, \dots, x_{n}) \, \Big].
				 
    		\end{array}$ 
		\end{center} 
      
    \item Вычитание:  
    	\begin{center}
    		$\begin{array}{l}
    		
				\mbf{a}(x_{1}, \dots, x_{n}) - \mbf{b}(x_{1}, \dots, x_{n}) = \\
				
				\hspace{3cm} = \Big[ \, L_{\mbf{a}}(x_{1}, \dots, _{n}), \: U_{\mbf{a}}(x_{1}, \dots, x_{n}) \, \Big] \: - \\
				
				\hspace{6cm} - \: \Big[ \, L_{\mbf{b}}(x_{1}, \dots, x_{n}), \: U_{\mbf{b}}(x_{1}, \dots, x_{n}) \, \Big] = \\
				
				\hspace{3cm} = \Big[ \, L_{\mbf{a}}(x_{1}, \dots, x_{n}) - U_{\mbf{b}}(x_{1}, \dots, x_{n}), \\
				
				\hspace{6cm} U_{\mbf{a}}(x_{1}, \dots, x_{n}) - L_{\mbf{b}}(x_{1}, \dots, x_{n}) \, \Big].
				 
			\end{array}$
		\end{center}      
    \item Умножение:    
    
    \begin{center}
    $\begin{array}{l}
    
    	\mbf{a}(x_{1}, \dots, x_{n}) \cdot \mbf{b}(x_{1}, \dots, x_{n}) = \\
    
    	\hspace{1.5cm} \Big[ \, F_{\downarrow} \big\{ \, L_{\mbf{a}}(x_{1}, \dots, x_{n}) \cdot L_{\mbf{b}}(x_{1}, \dots, x_{n}), \\
     
     	\hspace{3cm} L_{\mbf{a}}(x_{1}, \dots, x_{n}) \cdot U_{\mbf{b}}(x_{1}, \dots, x_{n}), \\
     
     	\hspace{4.5cm} U_{\mbf{a}}(x_{1}, \dots, x_{n}) \cdot L_{\mbf{b}}(x_{1}, \dots, x_{n}), \\
     
     	\hspace{6cm} U_{\mbf{a}}(x_{1}, \dots, x_{n}) \cdot U_{\mbf{b}}(x_{1}, \dots, x_{n}) \big\}, \\
    
    	\hspace{2mm} \hspace{1.5cm} F_{\uparrow} \big\{ L_{\mbf{a}}(x_{1}, \dots, x_{n}) \cdot L_{\mbf{b}}(x_{1}, \dots, x_{n}), \\
    
    	\hspace{3cm} L_{\mbf{a}}(x_{1}, \dots, x_{n}) \cdot U_{\mbf{b}}(x_{1}, \dots, x_{n}), \\
    
    	\hspace{4.5cm} U_{\mbf{a}}(x_{1}, \dots, x_{n}) \cdot L_{\mbf{b}}(x_{1}, \dots, x_{n}), \\
    
    	\hspace{6cm} U_{\mbf{a}}(x_{1}, \dots, x_{n}) \cdot U_{\mbf{b}}(x_{1}, \dots, x_{n}) \big\} \, \Big].
    
    \end{array}$     
    \end{center}
    
    \item Частный случай предыдущего пункта --- умножение на константу:    
    
    \begin{center}
    
    $\mbf{a}(x_{1}, \dots, x_{n}) \cdot \alpha = 
    	\left\{ \begin{array}{l}
    
        \Big[ \, L_{\mbf{a}}(x_{1}, \dots, x_{n}) \cdot \alpha, \\
        
        \hspace{2cm} U_{\mbf{a}}(x_{1}, \dots, x_{n}) \cdot \alpha \, \Big], \text{ если } \alpha \geq 0, \\
        
        \Big[ \, U_{\mbf{a}}(x_{1}, \dots, x_{n}) \cdot \alpha, \\
        
        \hspace{2cm} L_{\mbf{a}}(x_{1}, \dots, x_{n}) \cdot \alpha \, \Big], \text{ иначе.} 
        
    \end{array}\right.$    
    
    \end{center}
    
    \item Деление. 
    
    Если $0 \notin \Big[ \, L_{\mbf{b}}(x_{1}, \dots, x_{n}), \: U_{\mbf{b}}(x_{1}, \dots, x_{n}) \, \Big]$, то:    
    
    \begin{center}
    
	$\begin{array}{l}
	
		\mbf{a}(x_{1}, \dots, x_{n}) \, / \, \mbf{b}(x_{1}, \dots, x_{n}) = \\
	
		\hspace{1.5cm} \Big[ \, F_{\downarrow} \big\{ \, L_{\mbf{a}}(x_{1}, \dots, x_{n}) \, / \, L_{\mbf{b}}(x_{1}, \dots, x_{n}), \\
	
		\hspace{3cm} L_{\mbf{a}}(x_{1}, \dots, x_{n}) \, / \, U_{\mbf{b}}(x_{1}, \dots, x_{n}), \\
	
		\hspace{4.5cm} U_{\mbf{a}}(x_{1}, \dots, x_{n}) \, / \, L_{\mbf{b}}(x_{1}, \dots, x_{n}), \\
	
		\hspace{6cm} U_{\mbf{a}}(x_{1}, \dots, x_{n}) \, / \, U_{\mbf{b}}(x_{1}, \dots, x_{n}) \, \big\}, \\
	
		\hspace{2mm} \hspace{1.5cm} F_{\uparrow} \big\{ \, L_{\mbf{a}}(x_{1}, \dots, x_{n}) \, / \, L_{\mbf{b}}(x_{1}, \dots, x_{n}), \\
	
		\hspace{3cm} L_{\mbf{a}}(x_{1}, \dots, x_{n}) \, / \, U_{\mbf{b}}(x_{1}, \dots, x_{n}), \\
	
		\hspace{4.5cm} U_{\mbf{a}}(x_{1}, \dots, x_{n}) \, / \, L_{\mbf{b}}(x_{1}, \dots, x_{n}), \\
	
		\hspace{6cm}U_{\mbf{a}}(x_{1}, \dots, x_{n}) \, / \, U_{\mbf{b}}(x_{1}, \dots, x_{n}) \, \big\} \, \Big].
	
	\end{array}$
	\end{center}
\end{enumerate}

Обоснованием того, что данные арифметические операции между интервалами функционально-граничной арифметики сохраняют базовое свойство интервальных арифметик (\ref{eq:BaseIntervalPrincipe}), является то, что упорядоченная пара функционалов границ $\big( L(x_{1}, \dots, x_{n}), \: U(x_{1}, \dots, x_{n}) \big)$ описывает параметрическое семейство классических интервалов, для которых эти операции уже определены в классической интервальной арифметике.

\subsection{Построение функционалов границ $L$ и $U$}

Ключевой проблемой функционально-граничной интервальной арифметики является введение правил построения функционалов границ $L$ и $U$, которыми будут описываться множества совместных значений диаграмм зависимости.

В предположении того, что $L = L(x_{1}, \dots, x_{n})$ и $U = U(x_{1}, \dots, x_{n})$, естественно будет использовать какой-либо функциональный базис, связанный с интервальными параметрами $x_{1} \in \mbf{x}_{1}$, $\dots$, $x_{n} \in \mbf{x}_{n}$.

Каким требованиям должен удовлетворять выбираемый функциональный базис? Ранее в работе были озвучены некоторые проблемы, присущие классической интервальной и аффинной арифметикам. На основе этого были сформулированы следующие требования к базису функционально-граничной арифметики:

\begin{enumerate}

	\item При вычислении должны учитываться квадратичные члены интервальных параметров.
	
	\item Можно несложным образом ввести процедуры $F_{\downarrow}$ и $F_{\uparrow}$.
	
	\item Процедуры $F_{\downarrow}$ и $F_{\uparrow}$ должны иметь алгоритмическую сложность не более $O(n)$.
	
\end{enumerate}

Пусть в рассмотрении имеется функция $f(x_{1}, \dots, x_{n})$. Для нахождения её интервальной оценки можно ипользовать различные базисы. При написании работы рассматривались следующие:
\begin{enumerate}

	\item $\big\{ \, 1$, $x_{1}$, $x_{1}^{2}$, $\dots$, $x_{n}$, $x_{n}^2 \, \big\}$.
	
	Данный базис оказался неподходящим, поскольку этим базисом сложно описывать множества совместных значений диаграмм связности, границы которых представляют собой, например, криволинейные ломаные.
	
	В качестве примера можно рассмотреть нахождение интервальной оценки для функции (рис. \ref{fig:example3}): 
	
	\begin{equation}	
	\label{eq:f_example}
		\begin{array}{c}
		
			f(x) = (x - a) \cdot (x - b), \\
		
			x \in \big[ \, -1, \: 1 \, \big], \: a \in \big[ \, -\frac{1}{2}, \: \frac{1}{2} \, \big], \: b \in \big[ \, -\frac{1}{2}, \: \frac{1}{2} \, \big].
		
		\end{array}		
	\end{equation}	
	
	Эффективному приближению множества совместных значений данной диаграммы зависимости непрерывными функциями мешает излом верхней границы семейства в точке $\big( 0, \: \frac{1}{4} \big)$. На момент написания работы остаётся открытым вопрос о существовании и целесообразности построения процедур $F_{\downarrow}$ и $F_{\uparrow}$ в гладком базисе.
	
	\begin{figure}
	\centering
    \includegraphics[width=0.7\linewidth]{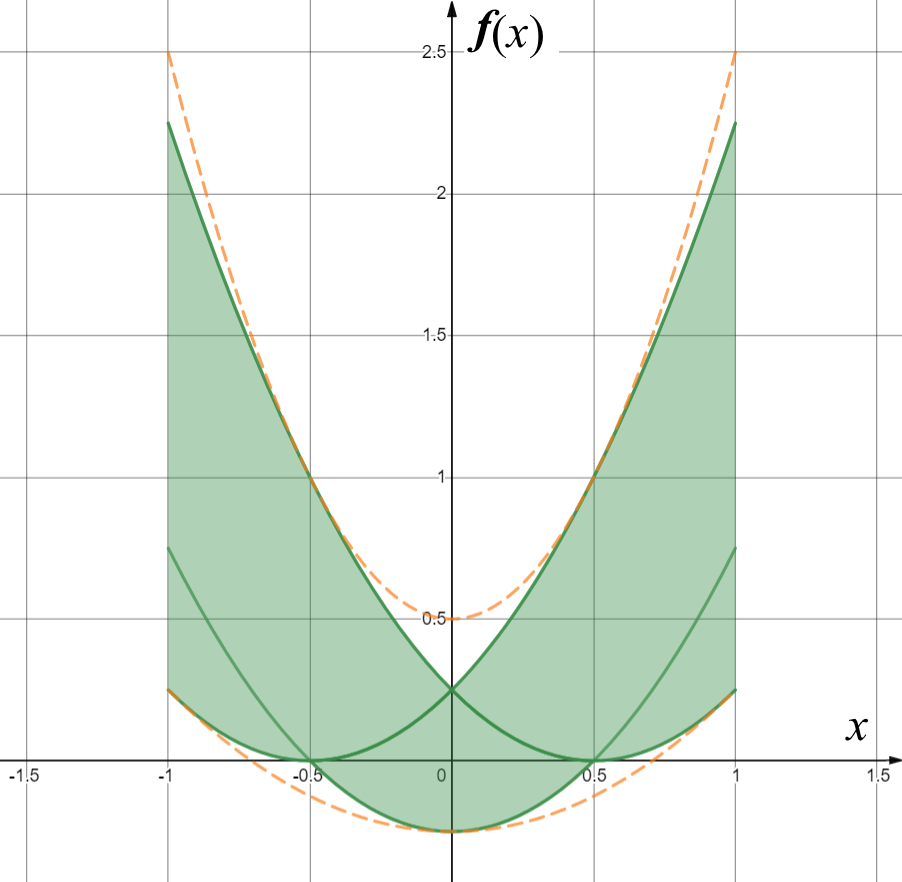}
    \caption{Зелёная область --- геометрическая интерпретация множества совместных значений $\big( x, \: \mbf{f}(x) \big)$, где $f(x)$ задаётся в (\ref{eq:f_example}). Оранжевая пунктирная линия --- пример приближения множества совместных значений в базисе $\big\{ \, 1$, $x$, $x^{2} \, \big\}$ сверху и снизу.}
	\label{fig:example3}
	\end{figure}

	\item $\begin{array}{l}
	
		\big\{ \, 1, \: x_{1}, \: \dots \: , \: x_{n}, \: x_{1} \cdot x_{1}, \: x_{1} \cdot x_{2}, \: \dots \: , \: x_{1} \cdot x_{n}, \\
	
		\hspace{7cm} x_{2} \cdot x_{2}, \: \dots \: , \: x_{2} \cdot x_{n}, \: \dots \:, \: x_{n - 1} \cdot x_{n} \, \big\}.
	
	\end{array}$
	
	Данный базис рассматривался, поскольку он позволяет точно описывать перекрёстные члены $x_{i} \cdot x_{j}$ при умножении. Но данный базис не подходит к рассмотрению в работе, поскольку:
	\begin{enumerate}
	
		\item Размер базиса растет как $n^{2}$, следовательно, асимптотика операции умножения составляет $O(n^{4})$, где $n$ --- число интервальных параметров.
		
		\item Остаётся открытым вопрос построения процедур $F_{\downarrow}$ и $F_{\uparrow}$.
	
	\end{enumerate}
	
	\item $\begin{array}{l}
	
	\big\{ \, 1, x_{1}, \dots, x_{n}, |x_{1} - x_{2}|, |x_{1} - x_{3}|, \dots, |x_{1} - x_{n}|, \\
	
	\hspace{5.5cm} |x_{2} - x_{3}|, \dots, |x_{2} - x_{n}|, \dots, |x_{n - 1} - x_{n}| \, \big\}.
	
	\end{array}$
	
	Данный базис также позволяет оптимально описывать перекрёстные члены $x_{i} \cdot x_{j}$ при умножении. Но при этом он имеет также имеет недостатки, подобно базису прошлого пункта. В теории он позволяет лучшим образом описывать криволинейные ломаные, но при этом остальные недостатки он не исправляет.
	
	\item $\big\{ \, 1$, $x_{1}$, $|x_{1}|$, $\dots$, $x_{n}$, $|x_{n}| \, \big\}$.
	
	Данный базис предназначен для эффективной работы с функциями одной переменной, он позволяет учитывать квадратичности в вычислениях.
	
	Также в данном базисе идёт оперирование линейными функциями, что позволяет достаточно просто ввести процедуры $F_{\downarrow}$ и $F_{\uparrow}$, причём, алгоритмическая сложность этих процедур не будет превышать $O(n)$.	
	
	Размер базиса растёт как $O(n)$, что позволяет не допускать экпоненциального роста при увеличении размерности задачи. Трудоёмкость операции сложения и вычитаний --- $O(n)$, умножения и деления --- $O(n ^ {2})$.
	
	В текущей работе именно этот базис был выбран для использования в вычислениях, поскольку удовлетворил всем выдвигаемым требованиям.
	
\end{enumerate}

Итак, был выбран базис для описания интервалов $\mathbb{F}\mathbb{R}$:
\begin{center}

	$\big\{ \, 1$, $x_{1}$, $|x_{1}|$, $\dots$, $x_{n}$, $|x_{n}| \, \big\}.$
	
\end{center}
Далее будем называть его \textit{центрально-ломаным}. 

Теперь любой функционально-граничный интервал $\mbf{x} \in \mathbb{F}\mathbb{R}$ будет представляться в виде:

\begin{center}

	$\mbf{x}(x_{1}, \dots, x_{n}) = \Big[ \, \sum_{i = 1}^{n} \underline{a}_{i} \cdot x_{i} + \sum_{i = 1}^{n} \underline{b}_{i} \cdot |x_{i}| + \underline{c},$
	
	$\hspace{7cm} \sum_{i = 1}^{n} \overline{a}_{i} \cdot x_{i} + \sum_{i = 1}^{n} \overline{b}_{i} \cdot |x_{i}| + \overline{c} \, \Big].$

\end{center}

Получим явные выражения для арифметических операций между интервалами в этом базисе.

\subsubsection{Сложение и вычитание}

Итак, пусть $\mbf{a}(x_{1}, \dots, x_{n}), \: \mbf{b}(x_{1}, \dots, x_{n}) \in \mathbb{F}\mathbb{R}$. Аналогично классической интервальной арифметике введём операции сложения и вычитания (далее опустим аргументы при $\mbf{a}(x_{1}, \dots, x_{n})$ и $\mbf{b}(x_{1}, \dots, x_{n})$).

Сложение:

\begin{center}	
	$\begin{array}{rcl}
	
		\mbf{a} + \mbf{b} & = & \Big[ \, \sum_{i = 1}^{n} \underline{a}_{a, i} \cdot x_{i} + \sum_{i = 1}^{n} \underline{b}_{a, i} \cdot |x_{i}| + \underline{c}_{a}, \\
	
		& & \hspace{4cm} \sum_{i = 1}^{n} \overline{a}_{a, i} \cdot x_{i} + \sum_{i = 1}^{n} \overline{b}_{a, i} \cdot |x_{i}| + \overline{c}_{a} \, \Big] \: + \\
	
		& + & \Big[ \, \sum_{i = 1}^{n} \underline{a}_{b, i} \cdot x_{i} + \sum_{i = 1}^{n} \underline{b}_{b, i} \cdot |x_{i}| + \underline{c}_{b}, \\
	
		& & \hspace{4cm} \sum_{i = 1}^{n} \overline{a}_{b, i} \cdot x_{i} + \sum_{i = 1}^{n} \overline{b}_{b, i} \cdot |x_{i}| + \overline{c}_{b} \, \Big] = \\
	
		& = & \Big[ \, \sum_{i = 1}^{n} \big( \underline{a}_{a, i} + \underline{a}_{b, i} \big) \cdot x_{i} + \sum_{i = 1}^{n} \big( \underline{b}_{a, i} + \underline{b}_{b, i} \big) \cdot |x_{i}| + \big( \underline{c}_{a} + \underline{c}_{b} \big), \\
	
		& & \hspace{6mm} \sum_{i = 1}^{n} \big( \overline{a}_{a, i} + \overline{a}_{b, i} \big) \cdot x_{i} + \sum_{i = 1}^{n} \big( \overline{b}_{a, i} + \overline{b}_{b, i} \big) \cdot |x_{i}| + \big( \overline{c}_{a} + \overline{c}_{b} \big) \, \Big].
	
	\end{array}$
\end{center}
	
Вычитание:
\begin{center}
	$\begin{array}{rcl}
	
		\mbf{a} - \mbf{b} & = & \Big[ \, \sum_{i = 1}^{n} \underline{a}_{a, i} \cdot x_{i} + \sum_{i = 1}^{n} \underline{b}_{a, i} \cdot |x_{i}| + \underline{c}_{a}, \\
	
		& & \hspace{4cm} \sum_{i = 1}^{n} \overline{a}_{a, i} \cdot x_{i} + \sum_{i = 1}^{n} \overline{b}_{a, i} \cdot |x_{i}| + \overline{c}_{a} \, \Big] \: - \\
	
		& - & \Big[ \, \sum_{i = 1}^{n} \underline{a}_{b, i} \cdot x_{i} + \sum_{i = 1}^{n} \underline{b}_{b, i} \cdot |x_{i}| + \underline{c}_{b}, \\
	
		& & \hspace{4cm} \sum_{i = 1}^{n} \overline{a}_{b, i} \cdot x_{i} + \sum_{i = 1}^{n} \overline{b}_{b, i} \cdot |x_{i}| + \overline{c}_{b} \, \Big] = \\
	
		& = & \Big[ \, \sum_{i = 1}^{n} \big( \underline{a}_{a, i} - \overline{a}_{b, i} \big) \cdot x_{i} + \sum_{i = 1}^{n} \big( \underline{b}_{a, i} - \overline{b}_{b, i} \big) \cdot |x_{i}| + \big( \underline{c}_{a} - \underline{c}_{b} \big), \\
	
		& & \hspace{6mm} \sum_{i = 1}^{n} \big( \overline{a}_{a, i} - \underline{a}_{b, i} \big) \cdot x_{i} + \sum_{i = 1}^{n} \big( \overline{b}_{a, i} - \underline{b}_{b, i} \big) \cdot |x_{i}| + \big( \overline{c}_{a} - \underline{c}_{b} \big) \, \Big].
	
	\end{array}$
\end{center}

\subsubsection{Умножение. Начало построения}

Введём операцию умножения между интервалами. По определению, данному ранее (далее опустим аргументы функционалов границ и интервалов):

\begin{center}
    $\begin{array}{rcl}
    
    	\mbf{a} \cdot \mbf{b} & = & \Big[ \, F_{\downarrow} \big\{ \, L_{\mbf{a}} \cdot L_{\mbf{b}}, \: L_{\mbf{a}} \cdot U_{\mbf{b}}, \: U_{\mbf{a}} \cdot L_{\mbf{b}}, \: U_{\mbf{a}} \cdot U_{\mbf{b}} \, \big\}, \\
    
    	& & \hspace{2cm} F_{\uparrow} \big\{ \, L_{\mbf{a}} \cdot L_{\mbf{b}}, \: L_{\mbf{a}} \cdot U_{\mbf{b}}, \: U_{\mbf{a}} \cdot L_{\mbf{b}}, \: U_{\mbf{a}} \cdot U_{\mbf{b}} \, \big\} \, \Big].
    
    \end{array}$ 
\end{center}    
    
Перепишем результат умножения в центрально-ломаном базисе (далее опустим аргументы функционалов границ):
\begin{center}
	\begin{longtable}{rcl}

		$L_{\mbf{a}} \cdot L_{\mbf{b}}$ & $ = $ & $\Big( \sum_{i = 1}^{n} \underline{a}_{a, i} \cdot x_{i} + \sum_{i = 1}^{n} \underline{b}_{a, i} \cdot |x_{i}| + \underline{c}_{a} \Big) \: \cdot$ \\
		
		& & $\hspace{6mm} \cdot \: \Big( \sum_{i = 1}^{n} \underline{a}_{b, i} \cdot x_{i} + \sum_{i = 1}^{n} \underline{b}_{b, i} \cdot |x_{i}| + \underline{c}_{b} \Big) = $ \\
		
		& $ = $ & $\sum_{i = 1}^{n} \big( \underline{a}_{a, i} \cdot \underline{c}_{b} + \underline{a}_{b, i} \cdot \underline{c}_{a} \big) \cdot x_{i} \: + $\\
		
		& & $ \hspace{3mm} + \: \sum_{i = 1}^{n} \big( \underline{b}_{a, i} \cdot \underline{c}_{b} + \underline{b}_{b, i} \cdot \underline{c}_{a} \big) \cdot |x_{i}| \: + $ \\
		
		& & $ \hspace{6mm} + \: \big( \underline{c}_{a} \cdot \underline{c}_{b} \big) + \sum_{i = 1}^{n} \sum_{j = 1}^{n} \big( \underline{a}_{a, i} \cdot \underline{a}_{b, j} \big) \cdot x_{i} \cdot x_{j} \: +$ \\
		
		& & $\hspace{9mm} + \sum_{i = 1}^{n} \sum_{j = 1}^{n} \big( \underline{b}_{a, i} \cdot \underline{b}_{b, j} \big) \cdot |x_{i}| \cdot |x_{j}| \: +$ \\
		
		& & $\hspace{12mm} + \: \sum_{i = 1}^{n} \sum_{j = 1}^{n} \big( \underline{a}_{a, i} \cdot \underline{b}_{b, j} + \underline{b}_{a, i} \cdot \underline{a}_{b, j} \big) \cdot |x_{i}| \cdot x_{j} =$ \\
		
		& $ = $ & $\sum_{i = 1}^{n} \big( \underline{a}_{a, i} \cdot \underline{c}_{b} + \underline{a}_{b, i} \cdot \underline{c}_{a} \big) \cdot x_{i} \: +$ \\
		
		& & $\hspace{3mm} + \: \sum_{i = 1}^{n} \big( \underline{b}_{a, i} \cdot \underline{c}_{b} + \underline{b}_{b, i} \cdot \underline{c}_{a} \big) \cdot |x_{i}| \: +$ \\
		
		& & $\hspace{6mm} + \: \sum_{i = 1}^{n} \sum_{j = 1, j \neq i}^{n} \big( \underline{a}_{a, i} \cdot \underline{a}_{b, j} \big) \cdot x_{i} \cdot x_{j} \: +$ \\
		
		& & $\hspace{9mm} + \: \sum_{i = 1}^{n} \big( \underline{a}_{a, i} \cdot \underline{a}_{b, i} \big) \cdot x_{i}^{2} + \big( \underline{c}_{a} \cdot \underline{c}_{b} \big) \: +$ \\
		
		& & $\hspace{12mm} + \: \sum_{i = 1}^{n} \sum_{j = 1, j \neq i}^{n} \big( \underline{b}_{a, i} \cdot \underline{b}_{b, j} \big) \cdot |x_{i}| \cdot |x_{j}| \: +$ \\
		
		& & $\hspace{15mm} + \: \sum_{i = 1}^{n} \big( \underline{b}_{a, i} \cdot \underline{b}_{b, i} \big) \cdot |x_{i}|^{2} \: +$ \\
		
		& & $\hspace{18mm} + \: \sum_{i = 1}^{n} \sum_{j = 1, j \neq i}^{n} \big(\underline{a}_{a, i} \cdot \underline{b}_{b, j} + \underline{b}_{a, i} \cdot \underline{a}_{b, j} \big) \cdot |x_{i}| \cdot x_{j} \: +$ \\
		
		& & $\hspace{21mm} + \: \sum_{i = 1}^{n} \big( \underline{a}_{a, i} \cdot \underline{b}_{b, i} + \underline{b}_{a, i} \cdot \underline{a}_{b, i} \big) \cdot |x_{i}| \cdot x_{i} =$ \\
		
		& $=$ & $\sum_{i = 1}^{n} \Big( \big( \underline{a}_{a, i} \cdot \underline{c}_{b} + \underline{a}_{b, i} \cdot \underline{c}_{a} \big) \cdot x_{i} \: + $ \\
		
		& & $\hspace{6mm} + \: \big( \underline{b}_{a, i} \cdot \underline{c}_{b} + \underline{b}_{b, i} \cdot \underline{c}_{a} \big) \cdot |x_{i}| \: +$ \\
		
		& & $ \hspace{12mm} + \: \big( \underline{a}_{a, i} \cdot \underline{a}_{b, i} + \underline{b}_{a, i} \cdot \underline{b}_{b, i} \big) \cdot x_{i}^{2} \Big) \: +$ \\
		
		& $+$ & $ \sum_{i = 1}^{n} \sum_{j = 1, j \neq i} \Big( \big( \underline{b}_{a, i} \cdot \underline{b}_{b, j}) \cdot |x_{i}| \cdot |x_{j}| +$ \\
		
		& & $\hspace{6mm} + \: (\underline{a}_{a, i} \cdot \underline{a}_{b, j}) \cdot x_{i} \cdot x_{j} +$  \\
		
		& & $\hspace{12mm} + \: \big( \underline{a}_{a, i} \cdot \underline{b}_{b, j} + \underline{a}_{b, i} \cdot \underline{b}_{a, j} \big) \cdot |x_{i}| \cdot x_{j} \Big)$.
				
	\end{longtable}
\end{center}

Проблемой является создание процедур $F_{\downarrow}$ и $F_{\uparrow}$. Так как выбранный базис содержит модули интервальных параметров, будем рассматривать для каждого интервального параметра случай $x_{i} \leq 0$ и $x_{i} \geq 0$.

Далее будем использовать чебышёвское приближение для каждого интервального параметра $x_{i}$ на интервалах $[ \, -1, \: 0 \, ]$ и $[ \, 0, \: 1 \, ]$, а затем <<склеивать>> результаты приближения до представимого в центрально-ломаном базисе.

\subsubsection[Схема приближения выпуклых или вогнутых\\ функций с помощью чебышёвского альтернанса]{Схема приближения выпуклых или вогнутых\\функций с помощью чебышёвского альтернанса}

Пусть в рассмотрении имеется некоторая выпуклая (вогнутая) функция $f(x)$ на интервале $[ \, x_{0}, \: x_{2} \, ]$. Тогда согласно книге \cite{Alternance}, пользуясь аппаратом чебышёвского альтернанса, можно найти линейную функцию 
\begin{center}
	$\begin{array}{c}
	
		g(x) = g_{1} \cdot x + g_{0}, \\
		g_{0}, \: g_{1} \in \mathbb{R},
		
	\end{array}$
\end{center}
\noindent которая приближает функцию $f(x)$ наилучшим образом в чебышёвской метрике
\begin{center}

	$\big\|f(x) - g(x) \big\|_{\infty} = \text{max}_{x \in [ \, x_{0}, \: x_{2} \, ]} \big\{ \, | f(x) - g(x) | \, \big\} \rightarrow \text{min}$.
	
\end{center}

Функцию $g(x)$ можно найти из решения системы

\begin{center}
	$\begin{array}{r|ccl}
	
		(1) & f(x_{0}) - g_{1} \cdot x_{0} - g_{0} & = & q, \\
		
		(2) & f(x_{1}) - g_{1} \cdot x_{1} - g_{0} & = & -q, \\

		(3) & f(x_{2}) - g_{1} \cdot x_{2} - g_{0} & = & q,

	\end{array}$	
	
	\vspace{2mm}
	
	$x_{1} \in [ \, x_{0}, \: x_{2} \, ].$
	
\end{center}

\noindentРазрешим данную систему.
\begin{enumerate}

	\item Вычтем из уравнения $(3)$ уравнение $(1)$. Тем самым найдём $g_{1}$:

	\begin{center}
		\color{blue}\setlength{\fboxsep}{3mm}\fbox{\color{black}{$g_{1} = \frac{f(x_{2}) - f(x_{0})}{x_{2} - x_{0}}$}}
	\end{center}

	\item Продифференцируем $(2)$ уравнение системы по $x$ в точке $x_{1}$, тогда сможем найти $x_{1}$:

	\begin{center}
		$\begin{array}{c}
	
			f^{'}_{x}(x_{1}) - g_{1} = 0 \\

			$\color{blue}\setlength{\fboxsep}{3mm}\fbox{\color{black}{$x_{1} = \big( f _{x}^{'} \big)^{-1}(g_{1})$}}$
		
		\end{array}$
	\end{center}

	\item Для нахождения $g_{0}$ можно сложить уравнения $(1)$ и $(2)$:

	\begin{center}

		\color{blue}\setlength{\fboxsep}{3mm}\fbox{\color{black}{$g_{0} = \tfrac{1}{2} \cdot \big( f(x_{0}) + f(x_{1}) - g_{1} \cdot (x_{0} + x_{1}) \big)$}}
	
	\end{center}

	Или сложить уравнения $(2)$ и $(3):$

	\begin{center}

		\color{blue}\setlength{\fboxsep}{3mm}\fbox{\color{black}{$g_{0} = \tfrac{1}{2} \cdot \big( f(x_{1}) + f(x_{2}) - g_{1} \cdot (x_{1} + x_{2}) \big)$}}

	\end{center}

	\item Значение $q$ можно найти тремя способами, подставив найденные величины в одно из исходных уравнений системы:

	\begin{center}
		$\begin{array}{rl}
	
	\vspace{1mm}
			(1): & $\color{blue}\setlength{\fboxsep}{3mm}\fbox{\color{black}{$q = f(x_{0}) - g_{1} \cdot x_{0} - g_{0}$}}$ \\ \vspace{1mm}

			(2): & $\color{blue}\setlength{\fboxsep}{3mm}\fbox{\color{black}{$q = -f(x_{1}) + g_{1} \cdot x_{1} + g_{0}$}}$ \\

			(3): & $\color{blue}\setlength{\fboxsep}{3mm}\fbox{\color{black}{$q = f(x_{2}) - g_{1} \cdot x_{2} - g_{0}$}}$
			
		\end{array}$
	\end{center}

\end{enumerate}

Итого, были найдены $g_{0}$, $g_{1}$, $q \in \mathbb{R}$ такие, что
\begin{center}

	$f(x) \subseteq \mbf{f}(x) = \big[ \, g_{1} \cdot x + g_{0} - |q|, \: g_{1} \cdot x + g_{0} + |q| \, \big].$
	
\end{center}

\noindentТаким образом, теперь известно, как можно найти интервальные оценки для любой выпуклой или вогнутой функции $f(x)$ на интервалах $[ \, -1, \: 0 \, ]$ и $[ \, 0, \: 1 \, ]$. Для представления объединения этих оценок в центрально-ломаном базисе обсудим вопрос их <<склейки>>.

\subsubsection{<<Склейка>> левой и правой интервальных оценок}

Пусть на интервале $[ \, -1, \: 0 \, ]$ получена интервальная оценка в линейном базисе:
\begin{center}

	$\mbf{f}_{\!\!\textit{left}}(x) = \big[ \, \underline{k}_{\textit{left}} \cdot x + \underline{b}_{\textit{left}}, \: \overline{k}_{\textit{left}} \cdot x + \overline{b}_{\textit{left}} \, \big],$
	
\end{center}
а на интервале $[ \, 0, \: 1 \, ]$:
\begin{center}

	$\mbf{f}_{\!\!\textit{right}}(x) = \big[ \, \underline{k}_{\textit{right}} \cdot x + \underline{b}_{\textit{right}}, \: \overline{k}_{\textit{right}} \cdot x + \overline{b}_{\textit{right}} \, \big].$
	
\end{center}

\noindentНеобходимо получить интервальную оценку в центрально-ломаном базисе:
\begin{center}

	$\mbf{f}(x) = \big[ \, \underline{a} \cdot x + \underline{b} \cdot |x| + \underline{c}, \: \overline{a} \cdot x + \overline{b} \cdot |x| + \overline{c} \, \big].$
	
\end{center}

Рассмотрим построение функционала границы $U_{\!\mbf{f}}(x)$ для
\begin{center}

		$U_{\!\mbf{f}_{\!\!\textit{left}}}(x) = \overline{k}_{\textit{left}} \cdot x + \overline{b}_{\textit{left}},  \qquad U_{\!\mbf{f}_{\!\!\textit{right}}}(x) = \overline{k}_{\textit{right}} \cdot x + \overline{b}_{\textit{right}}.$
		
\end{center}

\noindentВ точке $x = -1$
\begin{center}

	$U_{\!\mbf{f}}(-1) = U_{\!\mbf{f}_{\!\!\textit{left}}}(-1) = -\overline{k}_{\textit{left}} + \overline{b}_{\textit{left}},$
	
\end{center}
в точке $x = 1$
\begin{center}

	$U_{\!\mbf{f}}(1) = U_{\!\mbf{f}_{\!\!\textit{right}}}(1) = \overline{k}_{\textit{right}} + \overline{b}_{\textit{right}}.$

\end{center}

\noindentРассмотрим точку $x = 0$. Так как мы строим функционал верхней границы, то необходимо
\begin{center}
	$\begin{array}{rl}
	
		U_{\!\mbf{f}}(x) \geq U_{\!\mbf{f}_{\!\!\textit{left}}}(x), & x \in [ \, -1, \: 0 \, ], \\
		
		U_{\!\mbf{f}}(x) \geq U_{\!\mbf{f}_{\!\!\textit{right}}}(x), & x \in [ \, 0, \: 1 \, ].
		
	\end{array}$
\end{center}

\begin{figure}
\centering
    \includegraphics[width=0.5\linewidth]{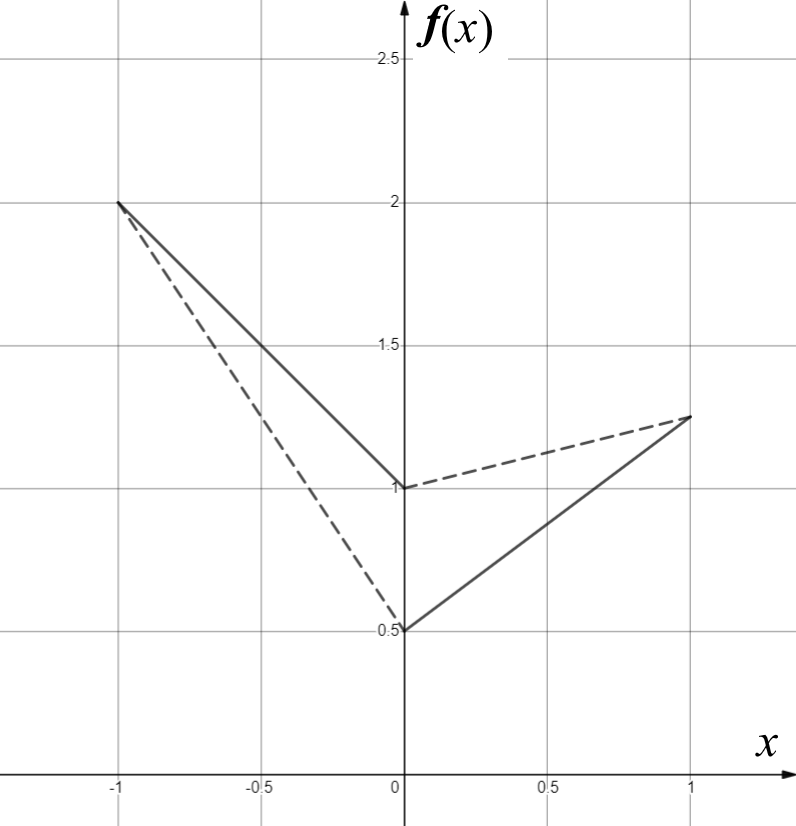}
    \caption{Сплошными линиями показаны исходные отрезки ломаной на левом и правом интервале. Верхняя пунктирная линия показывает достроение для $U_{\!\mbf{f}}(x)$ в центрально-ломаном базисе, нижняя --- для $L_{\!\mbf{f}}(x)$.}
    \label{fig:UL_example}
\end{figure}

\noindentТогда для центральной точки необходимо взять максимум (рис. \ref{fig:UL_example}):
\begin{center}
	
	$U_{\!\mbf{f}}(0) = \text{max} \Big\{ \, U_{\!\mbf{f}_{\!\!\textit{left}}}(0), \: U_{\!\mbf{f}_{\!\!\textit{right}}}(0) \, \Big\} = \text{max} \big\{ \, \overline{b}_{\textit{left}}, \: \overline{b}_{\textit{right}} \, \big\}.$
	
\end{center} 

Итак, получаем точки
\begin{center}

	$(-1, \: -\overline{k}_{\textit{left}} + \overline{b}_{\textit{left}}), \qquad (0, \: \text{max} \big\{ \, \overline{b}_{\textit{left}}, \: \overline{b}_{\textit{right}} \, \big\} ), \qquad (1, \: \overline{k}_{\textit{right}} + \overline{b}_{\textit{right}}).$
	
\end{center}

\noindentЧерез эти точки можно построить ломаную, представимую в центрально-ломаном базисе. Сначала заметим, что
\begin{center}

	$\overline{c} = \text{max} \big\{ \, \overline{b}_{\textit{left}}, \: \overline{b}_{\textit{right}} \, \big\}$.
	
\end{center}

Введём дополнительные обозначения:

\begin{center}
	$\begin{array}{rcl}
		
		\overline{k}_{\textit{left}} & = & \frac{1}{0 - (-1)} \cdot \Big( \text{max} \big\{ \, \overline{b}_{\textit{left}}, \: \overline{b}_{\textit{right}} \, \big\} - \big( -\overline{k}_{\textit{left}} + \overline{b}_{\textit{left}} \big) \Big) = \\
		
		& = & \text{max} \big\{ \, \overline{b}_{\textit{left}}, \: \overline{b}_{\textit{right}} \, \big\} - \big( -\overline{k}_{\textit{left}} + \overline{b}_{\textit{left}} \big), \\
		
		\overline{k}_{\textit{right}} & = & \frac{1}{1 - 0} \cdot \Big( \big( \overline{k}_{\textit{right}} + \overline{b} \big) - \text{max} \big\{ \, \overline{b}_{\textit{left}}, \: \overline{b}_{\textit{right}} \, \big\} \Big) = \\
		
		& = & \big( \overline{k}_{\textit{right}} + \overline{b} \big) - \text{max} \big\{ \, \overline{b}_{\textit{left}}, \: \overline{b}_{\textit{right}} \, \big\}.

	\end{array}$
\end{center}

\noindent$\overline{k}_{\textit{left}}$ --- наклон отрезка, соединяющего точки $\big( -1, \: U_{\!\mbf{f}}(-1) \big)$ и $\big( 0, \: U_{\!\mbf{f}}(0) \big)$, а $\overline{k}_{\textit{right}}$ --- наклон отрезка, соединяющего точки $\big( 0 , \: U_{\!\mbf{f}}(0) \big)$ и $\big( 1 , \: U_{\!\mbf{f}}(1) \big)$.

\noindentС учётом введённых обозначений получим следующую систему уравнений:

\begin{center}

	$\left\{\begin{array}{rrcl}
	
		x = -1: & - \overline{a} + \overline{b} + \overline{c} & = & - \overline{k}_{\textit{left}} + \overline{b}_{\textit{left}}, \\
	
		x = 1: & \overline{a} + \overline{b} + \overline{c} & = & \overline{k}_{\textit{right}} + \overline{b}_{\textit{right}}.
	
	\end{array}\right.$
	
\end{center}

\noindentИз этой системы находим, что
\begin{center}
	$\begin{array}{cc}
	
		\vspace{1mm}
		$\color{blue}\setlength{\fboxsep}{3mm}\fbox{\color{black}{$\overline{a} = \tfrac{1}{2} \cdot ( \overline{k}_{\textit{right}} + \overline{k}_{\textit{left}} + \overline{b}_{\textit{right}} - \overline{b}_{\textit{left}} )$}}$ \\
		
		$\color{blue}\setlength{\fboxsep}{3mm}\fbox{\color{black}{$\overline{b} = \tfrac{1}{2} \cdot ( \overline{k}_{\textit{right}} - \overline{k}_{\textit{left}} + \overline{b}_{\textit{right}} + \overline{b}_{\textit{left}} - 2 \cdot \overline{c} )$}}$
	
	\end{array}$

\end{center}

Аналогично рассмотрим построение $L_{\!\mbf{f}}(x)$ для:
\begin{center}

	$L_{\!\mbf{f}_{\!\!\textit{left}}}(x) = \underline{k}_{\textit{left}} \cdot x + \underline{b}_{\textit{left}}, \qquad L_{\!\mbf{f}_{\!\!\textit{right}}}(x) = \underline{k}_{\textit{right}} \cdot x + \underline{b}_{\textit{right}}.$
	
\end{center}
Для этого случая:
\begin{center}
	$\begin{array}{c}
	
		L_{\!\mbf{f}}(-1) = L_{\!\mbf{f}_{\!\!\textit{left}}}(-1) = -\underline{k}_{\textit{left}} + \underline{b}_{\textit{left}}, \\
		
		L_{\!\mbf{f}}(1) = L_{\!\mbf{f}_{\!\!\textit{right}}}(1) = \underline{k}_{\textit{right}} + \underline{b}_{\textit{right}}, \\
		
		L_{\!\mbf{f}}(0) = \text{min} \Big\{ \, L_{\!\mbf{f}_{\!\!\textit{left}}}(0), \: L_{\!\mbf{f}_{\!\!\textit{right}}}(0) \, \Big\} = \text{min} \big\{ \, \underline{b}_{\textit{left}}, \: \underline{b}_{\textit{right}} \, \big\}.
		
	\end{array}$
\end{center}

\noindentВведём обозначения:

\begin{center}
	$\begin{array}{rcl}
		
		\underline{k}_{\textit{left}} & = & \frac{1}{0 - (-1)} \cdot \Big( \text{min} \big\{ \, \overline{b}_{\textit{left}}, \: \underline{b}_{\textit{right}} \, \big\} - \big( -\underline{k}_{\textit{left}} + \underline{b}_{\textit{left}} \big) \Big) = \\
		
		& = & \text{min} \big\{ \, \underline{b}_{\textit{left}}, \: \underline{b}_{\textit{right}} \, \big\} - \big( -\underline{k}_{\textit{left}} + \underline{b}_{\textit{left}} \big), \\
		
		\underline{k}_{\textit{right}} & = & \frac{1}{1 - 0} \cdot \Big( \big( \underline{k}_{\textit{right}} + \underline{b} \big) - \text{min} \big\{ \, \underline{b}_{\textit{left}}, \: \underline{b}_{\textit{right}} \, \big\} \Big) = \\
		
		& = & \big( \underline{k}_{\textit{right}} + \underline{b} \big) - \text{min} \big\{ \, \underline{b}_{\textit{left}}, \: \underline{b}_{\textit{right}} \, \big\}.

	\end{array}$
\end{center}

\noindentСистема будет выглядеть следующим образом:

\begin{center}
	$\left\{\begin{array}{rrcl}
	
		x = -1: & - \underline{a} + \underline{b} + \underline{c} & = & - \underline{k}_{\textit{left}} + \underline{b}_{\textit{left}}, \\
	
		x = 1: & \underline{a} + \underline{b} + \underline{c} & = & \underline{k}_{\textit{right}} + \underline{b}_{\textit{right}}.
	
	\end{array}\right.$
\end{center}

\noindentРешением системы является: 
\begin{center}
	$\begin{array}{c}
	
		\vspace{1mm}
		$\color{blue}\setlength{\fboxsep}{3mm}\fbox{\color{black}{$\underline{a} = \tfrac{1}{2} \cdot (\underline{k}_{\textit{right}} + \underline{k}_{\textit{left}} + \underline{b}_{\textit{right}} - \underline{b}_{\textit{left}})$}}$ \\

		$\color{blue}\setlength{\fboxsep}{3mm}\fbox{\color{black}{$\underline{b} = \tfrac{1}{2} \cdot (\underline{k}_{\textit{right}} - \underline{k}_{\textit{left}} + \underline{b}_{\textit{right}} + \underline{b}_{\textit{left}} - 2 \cdot \underline{c})$}}$
		
	\end{array}$
\end{center}

Итак, если каким-либо образом были получены интервальные оценки на интервалах $[ \, -1, \: 0 \, ]$ и $[ \, 0, \: 1 \, ]$
\begin{center}
	$\begin{array}{c}

		\mbf{f}_{\!\!\textit{left}}(x) = \big[ \, \underline{k}_{\textit{left}} \cdot x + \underline{b}_{\textit{left}}, \: \overline{k}_{\textit{left}} \cdot x + \overline{b}_{\textit{left}} \, \big], \\
	
		\mbf{f}_{\!\!\textit{right}}(x) = \big[ \, \underline{k}_{\textit{right}} \cdot x + \underline{k}_{\textit{right}}, \: \overline{k}_{\textit{right}} \cdot x + \overline{b} \, \big],
	
	\end{array}$	
\end{center}

\noindentто

\begin{center}
	$\begin{array}{rcl}
	
		\mbf{f}_{\!\!\textit{left}} \, \bigcup \, \mbf{f}_{\!\!\textit{right}} & \subseteq & \Big[ \, \frac{1}{2} \cdot \big( \underline{k}_{\textit{right}} + \underline{k}_{\textit{left}} + \underline{b}_{\textit{right}} - \underline{b}_{\textit{left}} \big) \cdot x \: +\\

		& & \hspace{2cm} + \: \frac{1}{2} \cdot \big( \underline{k}_{\textit{right}} - \underline{k}_{\textit{left}} + \underline{b}_{\textit{right}} + \underline{b}_{\textit{left}} \: - \\

		& & \hspace{4cm} - \: 2 \cdot \text{min} \big\{ \, \underline{b}_{\textit{left}}, \: \underline{b}_{\textit{right}} \, \big\} \Big) \cdot |x| \: + \\
		
		& & \hspace{6cm} + \: \text{min} \big\{ \, \underline{b}_{\textit{left}}, \: \underline{b}_{\textit{right}} \, \big\}, \\

		& & \hspace{4mm} \frac{1}{2} \cdot \big( \overline{k}_{\textit{right}} + \overline{k}_{\textit{left}} + \overline{b}_{\textit{right}} - \overline{b}_{\textit{left}} \big) \cdot x \: + \\

		& & \hspace{2cm} + \: \frac{1}{2} \cdot \big( \overline{k}_{\textit{right}} - \overline{k}_{\textit{left}} + \overline{b}_{\textit{right}} + \overline{b}_{\textit{left}} \: - \\

		& & \hspace{4cm} - \: 2 \cdot \text{max} \big\{ \, \overline{b}_{\textit{left}}, \: \overline{b}_{\textit{right}} \big\} \big) \cdot |x| \: + \\
		
		& & \hspace{6cm} + \: \text{max} \big\{ \, \overline{b}_{\textit{left}} , \: \overline{b}_{\textit{right}} \, \big\} \, \Big].

	\end{array}$
\end{center}

Воспользуемся введёнными инструментами для функций

\begin{center}
	$\begin{array}{rcl}
	
		f_{1}(x) & = & x ^ {2}, \\
		
		f_{2}(x) & = & x \cdot |x|.
	
	\end{array}$
\end{center}

\noindentРассмотрим функцию $f_{1}(x)$. На интервале $[ \, -1, \: 0 \, ]$ получим интервальную оценку
\begin{center}

	$\mbf{f}_{\!\!\textit{left}}(x) = \big[ \, -x - \frac{1}{4}, \: -x \, \big],$
	
\end{center}
на интервале $[ \, 0, \: 1 \, ]$
\begin{center}

	$\mbf{f}_{\!\!\textit{right}}(x) = \big[ \, x, \: x + \frac{1}{4} \, \big].$
	
\end{center}

\noindentПосле <<склейки>> результатов получим
\begin{center}
	$\mbf{f}(x) = \big[ \, |x| - \frac{1}{4}, \: |x| \, \big].$
\end{center}

Рассмотрим теперь функцию $f_{2}(x)$. Она является нечётной, значит, согласно \cite{Alternance}, будем приближать её функцией $g(x)$, где $g_{0} = 0$.

\noindentНа интервале $[ \, -1 , \: 0 \, ]$ получим интервальную оценку
\begin{center}

	$\mbf{f}_{\!\!\textit{left}}(x) = \Big[ \, 2 \cdot \big(  \sqrt{2} - 1 \big) \cdot x - \big( 3 - 2 \cdot \sqrt{2} \big), \: 2 \cdot \big( \sqrt{2} - 1 \big) \cdot x + \big( 3 - 2 \cdot \sqrt{2} \big) \, \Big],$
	
\end{center}
на интервале $[ \, 0, \: 1 \, ]$
\begin{center}

	$\mbf{f}_{\!\!\textit{right}}(x) = \Big[ \, 2 \cdot \big( \sqrt{2} - 1 \big) \cdot x - \big( 3 - 2 \cdot \sqrt{2} \big), \: 2 \cdot \big( \sqrt{2} - 1 \big) \cdot x + \big( 3 - 2 \cdot \sqrt{2} \big) \, \Big].$
	
\end{center}

\noindentПосле <<склейки>> результатов получим 
\begin{center}

	$\mbf{f}(x) = \Big[ \, 2 \cdot \big( \sqrt{2} - 1 \big) \cdot x - \big( 3 - 2 \cdot \sqrt{2} \big), \: 2 \cdot \big( \sqrt{2} - 1 \big) \cdot x + \big( 3 - 2 \cdot \sqrt{2} \big) \, \Big].$
	
\end{center}

\subsubsection[Другой взгляд на приближение функций \\ от нескольких переменных]{Другой взгляд на приближение функций\\ от нескольких переменных}

Ранее было показано, как с помощью инструмента чебышёвского альтернанса и <<склейки>>, мы можем построить граничные функционалы $L(x_{1}, \dots, x_{n})$ и $U(x_{1}, \dots, x_{n})$, которые будут огибать совместное множество значений диаграммы зависимости снизу и сверху соответственно.

Рассмотрим построение функционалов границ $L(x_{1}, \dots, x_{n})$ и $U(x_{1}, \dots, x_{n})$ для интервальных оценок функций нескольких переменных. Автором работы было выявлено два основных подхода:

\begin{enumerate}
	\item \textit{Аналитический подход.}
	
	Данный подход подразумевает использование аналитических выкладок для получения явного вида $L(x_{1}, \dots, x_{n})$ и $U(x_{1}, \dots, x_{n})$, которые представимы в выбранном функциональном базисе.
	
	\item \textit{Подход редукции размерности.}
	
	Данный подход подразумевает сведение задачи нахождения интервальной оценки функции $f(x_{1}, \dots, x_{n})$ к задаче нахождения интервальной оценки функции одной переменной $x_{i}$.
	
	Этого можно достигнуть несколькими способами. Например, можно использовать классическую интервальную арифметику, то есть представить $f(x_{1}, \dots, x_{n})$ в виде 
\begin{center}
	$\begin{array}{rcl}
	
		f(x_{1}, \dots, x_{n}) & = & \sum_{i = 1}^{n} \frac{1}{n} \cdot f(x_{1}, \dots, x_{n}) \subseteq \\
		
		& \subseteq & \sum_{i = 1}^{n} f_{i}(\mbf{x}_{1}, \dots, \mbf{x}_{i - 1}, x_{i}, \mbf{x}_{i + 1}, \dots, \mbf{x}_{n}) \subseteq \\		
		
		& \subseteq & \sum_{i = 1}^{n} \mbf{f}_{\!\!i}(\mbf{x}_{1}, \dots, \mbf{x}_{i - 1}, x_{i}, \mbf{x}_{i + 1}, \dots, \mbf{x}_{n}) \subseteq \\
		
		& \subseteq & \sum_{i = 1}^{n} \Big[ \, L_{\!\mbf{f}_{\!\!i}}(\mbf{x}_{1}, \dots, \mbf{x}_{i - 1}, x_{i}, \mbf{x}_{i + 1}, \dots, \mbf{x}_{n}), \\
		
		& & \hspace{20mm} U_{\!\mbf{f}_{\!\!i}}(\mbf{x}_{1}, \dots, \mbf{x}_{i - 1}, x_{i}, \mbf{x}_{i + 1}, \dots, \mbf{x}_{n}) \, \Big ] \subseteq \\
		
		& \subseteq & \Big[ \, \sum_{i = 1}^{n} L_{\!\mbf{f}_{\!\!i}}(\mbf{x}_{1}, \dots, \mbf{x}_{i - 1}, x_{i}, \mbf{x}_{i + 1}, \dots, \mbf{x}_{n}), \\
		
		& & \hspace{20mm} \sum_{i = 1}^{n} U_{\!\mbf{f}_{\!\!i}}(\mbf{x}_{1}, \dots, \mbf{x}_{i - 1}, x_{i}, \mbf{x}_{i + 1}, \dots, x_{n}) \, \Big],
		
		\end{array}$
		
		где $f_{i}(x_{1}, \dots, x_{n}) = \frac{1}{n} \cdot f(x_{1}, \dots, x_{n}).$
\end{center}

Также возможен другой случай. Пусть $A$ --- множество аргументов функции $f(x_{1}, \dots, x_{n})$. Пусть $f(x_{1}, \dots, x_{n})$ можно представить в виде 
\begin{center}

	$f(x_{1}, \dots, x_{n}) = \sum_{i = 1}^{n} f_{i}(x_{k} \in A_{i})$, где 
	
	$A_{i} \bigcap A_{j} = \varnothing, \qquad i \neq j,$
	
	$\bigcup A_{i} = A, \qquad i = 1, \dots, n.$

\end{center}
	 Тогда интервальную оценку $f(x_{1}, \dots, x_{n})$ можно найти как
\begin{center}

	$\begin{array}{rcl}
	
		f(x_{1}, \dots, x_{n}) & = & \sum_{i = 1}^{n} f_{i} \Big( \big\{ \, x_{k} \in A_{i} \, \big| \, k = 1, \dots, \overline{\overline{A_{i}}} \, \big\} \Big) = \\
		
		& = & \sum_{i = 1}^{n} \sum_{j = 1}^{ \overline{\overline{A_{i}}} } f_{i} \Big( \big\{ \, x_{k} \in A_{i} \, \big| \, k = 1, \dots, \overline{\overline{A_{i}}} \, \big\} \Big) \, / \, \overline{\overline{A_{i}}} \subseteq \\
		
		& \subseteq & \sum_{i = 1}^{n} \sum_{j = 1}^{ \overline{\overline{A_{i}}} } \mbf{f}_{\!\!i} \Big( x_{j} \in A_{i}, \\
		
		& & \big\{ \mbf{x}_{k} \, \big| \, x_{k} \in A_{i}, \: k = 1, \dots, \overline{\overline{A_{i}}} , k \neq j \, \big\} \Big) \, / \, \overline{\overline{A_{i}}} = \\
		
		& = & \sum_{i = 1}^{n} \sum_{j = 1}^{ \overline{\overline{A_{i}}} } \mbf{f}_{\!\!i, j} \big( x_{j} \in A_{i} \big) \subseteq \\
		
		& \subseteq & \sum_{i = 1}^{n} \sum_{j = 1}^{ \overline{\overline{A_{i}}} } \Big[ \, L_{\!\mbf{f}_{\!\!i, j}}(x_{j} \in A_{i}), \: U_{\!\mbf{f}_{\!\!i, j}}(x_{j} \in A_{i}) \, \Big] \subseteq \\
		
		& \subseteq & \Big[ \, \sum_{i = 1}^{n} \sum_{i = 1}^{ \overline{\overline{A_{i}}} } L_{\!\mbf{f}_{\!\!i, j}}(x_{j} \in A_{i}), \\
		
		& & \hspace{4cm} \sum_{i = 1}^{n} \sum_{j = 1}^{ \overline{\overline{A_{i}}} } U_{\!\mbf{f}_{\!\!i, j}}(x_{j} \in A_{i}) \, \Big],
	
	\end{array}$
	
	$\begin{array}{l}
	
		\text{где } \mbf{f}_{\!\!i, j} \big( x_{j} \in A_{i} \big) = \mbf{f}_{\!\!i} \Big( x_{j} \in A_{i}, \\
	
		\hspace{4cm} \big\{ \mbf{x}_{k} \, \big| \, x_{k} \in A_{i}, \: k = 1, \dots, \overline{\overline{A_{i}}} , k \neq j \, \big\} \Big) \, / \, \overline{\overline{A_{i}}}.
	
	\end{array}$

\end{center}

\end{enumerate} 

Продемонстрируем данные подходы для получения интервальных оценок следующих членов, которые встречаются в итоговом выражении для операции умножения между функционалами границ
\begin{center}

	$x_{i} \: \cdot \: x_{j}, \qquad |x_{i}| \cdot x_{j}, \qquad i \neq j$.
	
\end{center}

\indentРассмотрим функцию
\begin{center}

	$f(x_{i}, \: x_{j}) = x_{i} \, \cdot \, x_{j}.$
	
\end{center}

\noindentИз неравенств

\begin{center}

	$\big( x_{i} - x_{j} \big) ^ {2} \geq 0, \qquad \big( x_{i} + x_{j} \big)^ {2} \geq 0,$

\end{center}
\noindentполучаем, что
\begin{center}

	$- \frac{1}{2} \cdot \big( x_{i} ^ {2} + x_{j} ^ {2} \big) \leq x_{i} \cdot x_{j},$
	
	$\frac{1}{2} \cdot \big( x_{i} ^ {2} + x_{j} ^ {2} \big) \geq x_{i} \cdot x_{j}.$
	
\end{center}

\noindentТаким образом,
\begin{center}
	$\begin{array}{rcl}
	
		\mbf{f}(x_{i}, \: x_{j}) & = & \big[ \, -\frac{1}{2} \cdot x_{i}^{2} - \frac{1}{2} \cdot x_{j}^{2}, \: \frac{1}{2} \cdot x_{i}^{2} + \frac{1}{2} \cdot x_{j}^{2} \, \big] \subseteq \\
		
		& \subseteq & \big[ \, -\frac{1}{2} \cdot |x_{i}| -\frac{1}{2} \cdot |x_{j}|, \: \frac{1}{2} \cdot |x_{i}| + \frac{1}{2} \cdot |x_{j}| \, \big].
		
	\end{array}$
\end{center}

Теперь применим подход \textit{редукции размерности}. Перепишем функцию $f(x_{i}, \dots, x_{j})$ следующим образом:
\begin{center}

	$\begin{array}{rcl}
		
		f(x_{i}, \: x_{j}) & = & \frac{1}{2} \cdot x_{i} \cdot x_{j} + \frac{1}{2} \cdot x_{i} \cdot x_{j} = \\
		
		& = & \frac{1}{2} \cdot f_{i}(x_{i}, \: x_{j}) + \frac{1}{2} \cdot f_{j}(x_{i}, \: x_{j}),
		
	\end{array}$
	
	\vspace{2mm}
	
	где $f_{i}(x_{i}, \: x_{j}) = x_{i} \cdot x_{j}, \qquad f_{j}(x_{i}, \: x_{j}) = x_{i} \cdot x_{j}.$
	
\end{center}

\noindentНайдём интервальное расширение функции $f_{i}(x_{i}, \: x_{j})$ считая, что данная функция является зависимой только от $x_{i}$:
\begin{center}

	$\mbf{f}_{\!\!i}(x_{i}, \: \mbf{x}_{j}) = x_{i} \cdot [ \, -1, \: 1 \, ].$
	
\end{center}

\noindentАналогичные действия проделаем для функции $f_{j}(x_{i}, \: x_{j})$. Найдём её интервальное расширение относительно переменной $x_{j}$:
\begin{center}

	$\mbf{f}_{\!\!j}(\mbf{x}_{i}, \: x_{j}) = [ \, -1, \: 1 \, ] \cdot x_{j}.$
	
\end{center}

\begin{figure}
\centering
    \includegraphics[width=0.5\linewidth]{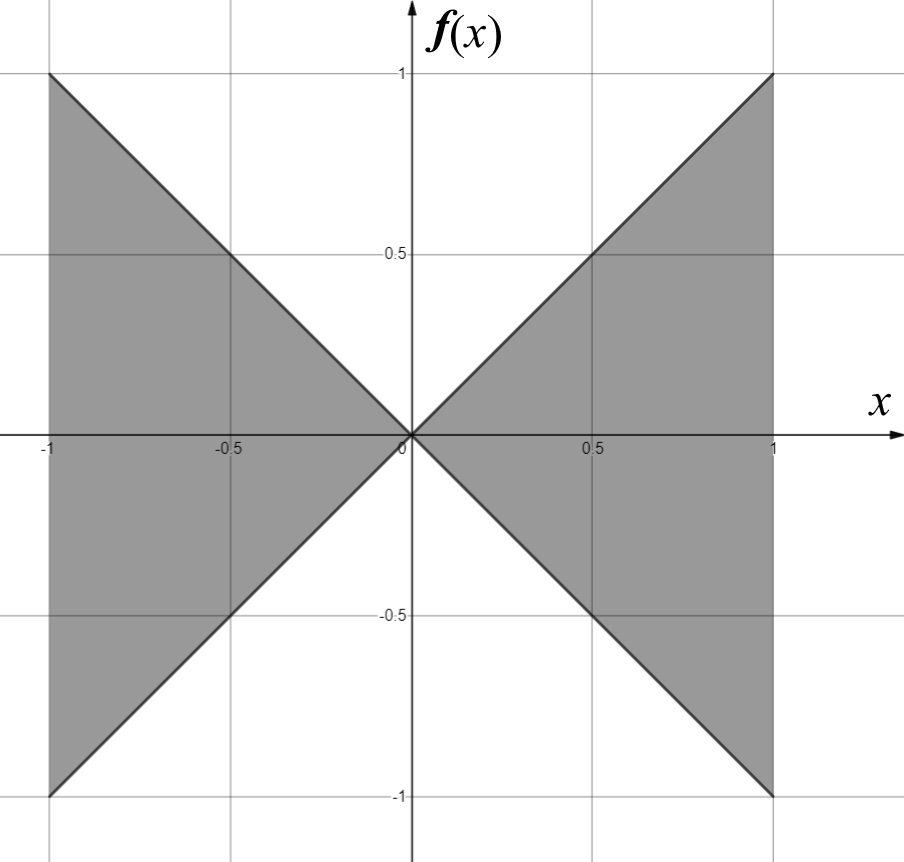}
    \caption{Множество совместных значений упорядоченной пары $\big(x, \mbf{f}(x) \big)$, где $\mbf{f}(x) = [ \, -1, \: 1 \, ] \cdot x, \: x \in [ \, -1, \: 1 \, ].$}
\label{fig:abs_example}
\end{figure}

\noindentМножество совместных значений упорядоченной пары $\big(x_{i}, \mbf{f}_{\!\!i}(x_{i}, \: \mbf{x}_{j}) \big)$ (рис. \ref{fig:abs_example}) можно точно описать в центрально-ломаном базисе
\begin{center}

	$\mbf{f}_{\!\!i}(x_{i}, \: \mbf{x}_{j}) = [\, -|x_{i}| , \: |x_{i}| \, ].$
	
\end{center}

\noindentТочно также при рассмотрении геометрической интерпретации интервального расширения $f_{j}(\mbf{x}_{i}, \: x_{j})$ получаем
\begin{center}

	$\mbf{f}_{\!\!j}(\mbf{x}_{i},\:  x_{j}) = [ \, -|x_{j}|, \: |x_{j}| \, ].$
	
\end{center}

Возвращаясь к первоначальному представлению функции, имеем:
\begin{center}
	$\begin{array}{rcl}
	
		f(x_{i}, \: x_{j}) & = & \frac{1}{2} \cdot f_{i}(x_{i}, \: x_{j}) + \frac{1}{2} \cdot f_{j}(x_{i}, \: x_{j}) \subseteq \\
		
		& \subseteq & \frac{1}{2} \cdot \mbf{f}_{\!\!i}(x_{i}, \: x_{j}) + \frac{1}{2} \cdot \mbf{f}_{\!\!j}(x_{i}, \: x_{j}) = \\
		
		& = & \frac{1}{2} \cdot \big[ \, -|x_{i}|, \: |x_{i}| \big] + \frac{1}{2} \cdot [ \, -|x_{j}|, \: |x_{j}|] = \\
		
		& = & \big[ \, -\frac{1}{2} \cdot |x_{i}| -\frac{1}{2} \cdot |x_{j}|, \: \frac{1}{2} \cdot |x_{i}| + \frac{1}{2} \cdot |x_{j}| \, \big].
		
	\end{array}$
\end{center}

Теперь перейдём к рассмотрению функции
\begin{center}

	$f(x_{i}, \: x_{j}) = |x_{i}| \cdot x_{j}.$
	
\end{center}

\noindentПрименим \textit{аналитический подход}. Используя неравенства
\begin{center}
	
	$\big( |x_{i}| - x_{j} \big)^{2} \geq 0, \qquad \big( |x_{i}| + x_{j} \big) ^ {2} \geq 0,$
	
\end{center}
\noindentполучим, что 
\begin{center}

	$- \frac{1}{2} \cdot \big( |x_{i}|^{2} + x_{j}^{2} \big) = - \frac{1}{2} \cdot \big( x_{i}^{2} + x_{j}^{2} \big) \leq |x_{i}| \cdot x_{j},$
	
	$\frac{1}{2} \cdot \big( |x_{i}|^{2} + x_{j}^{2} \big) = \frac{1}{2} \cdot \big( x_{i} ^ {2} + x_{j} ^ {2} \big) \geq |x_{i}| \cdot x_{j}.$
	
\end{center}

\noindentТаким образом,
\begin{center}
	$\begin{array}{rcl}
	
		\mbf{f}(x_{i}, \: x_{j}) & = & \big[ \, -\frac{1}{2} \cdot x_{i}^{2} -\frac{1}{2} \cdot x_{j}^{2}, \: \frac{1}{2} \cdot x_{i}^{2} + \frac{1}{2} \cdot x_{j}^{2} \, \big] \subseteq \\
		
		& \subseteq & \big[ \, -\frac{1}{2} \cdot |x_{i}| -\frac{1}{2} \cdot |x_{j}|, \: \frac{1}{2} \cdot |x_{i}| + \frac{1}{2} \cdot |x_{j}| \, \big].		

	\end{array}$
\end{center}

Теперь применим \textit{подход редукции размерности}. Перепишем функцию $f(x_{i}, \: x_{j})$ в следующем виде:
\begin{center}
	$\begin{array}{rcl}
	
		f(x_{i}, \: x_{j}) & = & |x_{i}| \cdot x_{j} = \frac{1}{2} \cdot |x_{i}| \cdot x_{j} + \frac{1}{2} \cdot |x_{i}| \cdot x_{j} = \\
		
		& = & \frac{1}{2} \cdot f_{i}(x_{i}, \: x_{j}) + \frac{1}{2} \cdot f_{j}(x_{i}, \: x_{j}),
		
	\end{array}$
	
	\vspace{2mm}	
	
	где $f_{i}(x_{i}, \: x_{j}) = |x_{i}| \cdot x_{j}, \qquad f_{j}(x_{i}, \: x_{j}) = |x_{i}| \cdot x_{j}.$
\end{center}

\noindentНайдём интервальное расширение функции $f_{i}(x_{i}, \: x_{j})$, считая, что она зависит только от $x_{i}$:
\begin{center}

	$\mbf{f}_{\!\!i}(x_{i}, \: \mbf{x}_{j}) = |x_{i}| \cdot [ \, -1, \: 1 \, ]$
	
\end{center}

\noindentМножество совместных значений упорядоченной пары $\big( x_{i}, \: \mbf{f}_{\!\!i}(x_{i}, \: \mbf{x}_{j}) \big)$ (рис. \ref{fig:abs_example}) можно точно описать в центрально-ломаном базисе как

\begin{center}

	$\mbf{f}_{\!\!i}(x_{i}, \: \mbf{x}_{j}) = \big[ \, -|x_{i}|, \: |x_{i}| \, \big].$
	
\end{center}

\noindentТеперь найдём интервальное расширение $f_{j}(x_{i}, \: x_{j})$ относительно $x_{j}$:

\begin{center}

	$\mbf{f}_{\!\!j}(\mbf{x}_{i}, \: x_{j}) = [ \, 0 , \: 1 \, ] \cdot x_{j}.$
	
\end{center}

\begin{figure}
\centering
    \includegraphics[width=0.5\linewidth]{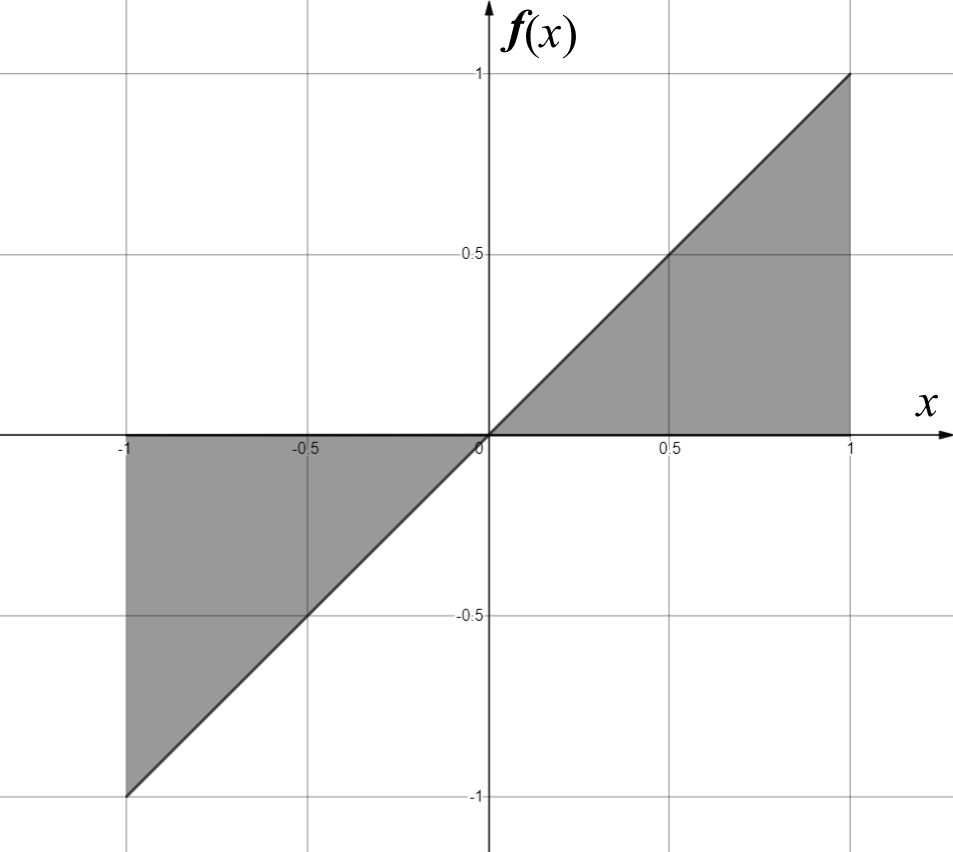}
    \caption{Множество совместных значений упорядоченной пары $\big( x, \: \mbf{f}(x) \big)$, где $\mbf{f}(x) = [ \, 0, \: 1 \, ] \cdot x, \: x \in [ \, -1, \:  1 \, ].$}
\label{fig:abs2_example}
\end{figure}

\noindentМножество совместных значений упорядоченной пары $\big( x_{j}, \: \mbf{f}_{\!\!j}(\mbf{x}_{i}, \: x_{j}) \big)$ (рис. \ref{fig:abs2_example}) можно точно описать в центрально-ломаном базисе как
	
	\begin{center}
	
		$\mbf{f}_{\!\!j}(\mbf{x}_{i}, \: x_{j}) = \big[ \, \frac{1}{2} \cdot x - \frac{1}{2} \cdot |x|, \: \frac{1}{2} \cdot x + \frac{1}{2} \cdot |x| \, \big].$
		
		\end{center}

Возвращаясь к первоначальному представлению функции, получим:
\begin{center}
	$\begin{array}{rcl}
	
		f(x_{i}, \: x_{j}) & = & \frac{1}{2} \cdot f_{i}(x_{i}, \: x_{j}) + \frac{1}{2} \cdot f_{j}(x_{i}, \: x_{j}) \subseteq \\
		
		& \subseteq & \frac{1}{2} \cdot \mbf{f}_{\!\!i}(x_{i}, \: x_{j}) + \frac{1}{2} \cdot \mbf{f}_{\!\!j}(x_{i}, \: x_{j}) \subseteq \\
		
		& \subseteq & \frac{1}{2} \cdot \mbf{f}_{\!\!i}(x_{i}, \: \mbf{x}_{j}) + \frac{1}{2} \cdot \mbf{f}_{\!\!j}(\mbf{x}_{i}, \: x_{j}) \subseteq \\
	
		& \subseteq & \frac{1}{2} \cdot \big[ \, -|x_{i}|, \: |x_{i}| \, \big] + \frac{1}{2} \cdot \big[ \, \frac{1}{2} \cdot x_{j} - \frac{1}{2} \cdot |x_{j}|, \: \frac{1}{2} \cdot x_{j} + \frac{1}{2} \cdot |x_{j}| \, \big] = \\
	
		& = & \big[ \, -\frac{1}{2} \cdot x_{i} + \frac{1}{4} \cdot x_{j} - \frac{1}{4} \cdot |x_{j}|, \: \frac{1}{2} \cdot x_{i} + \frac{1}{4} \cdot x_{j} + \frac{1}{4} \cdot |x_{j}| \, \big]
	
	\end{array}$
\end{center}

Сравнивая результаты этих двух подходов для данной функции получаем, что подход редукции размерности оказался лучше, поскольку напрямую задействует вид базиса для описания множества совместных значений диаграммы зависимости.

Наконец, рассмотрим функцию $f(x_{i}, \: x_{j}) = |x_{i}| \cdot |x_{j}|$. Сначала применим \textit{аналитический подход}. Используя неравенства
\begin{center}
	
		$\big( |x_{i}| - |x_{j}| \big) ^ {2} \geq 0 , \qquad \big( |x_{i}| + |x_{j}| \big) ^ {2} \geq 0,$
		
\end{center}
\noindentполучаем, что
\begin{center}

	$- \frac{1}{2} \cdot \big( |x_{i}|^{2} + |x_{j}|^{2} \big) = -\frac{1}{2} \cdot \big( x_{i}^{2} + x_{j}^{2} \big) \leq |x_{i}| \cdot |x_{j}|,$
	
	$ \frac{1}{2} \cdot \big( |x_{i}|^{2} + |x_{j}|^{2} \big) = \frac{1}{2} \cdot \big( x_{i}^{2} + x_{j}^{2} \big) \geq |x_{i}| \cdot |x_{j}|.$
	
\end{center} 

\noindentТаким образом,
\begin{center}
	$\begin{array}{rcl}
	
		\mbf{f}(x_{i}, \: x_{j}) & = & \big[ \, -\frac{1}{2} \cdot x_{i}^{2} - \frac{1}{2} \cdot x_{j}^{2}, \: \frac{1}{2} \cdot x_{i}^{2} + \frac{1}{2} \cdot x_{j}^{2} \, \big] \subseteq \\
	
		& \subseteq & \big[ \, -\frac{1}{2} \cdot |x_{i}| -\frac{1}{2} \cdot |x_{j}|, \: \frac{1}{2} \cdot |x_{i}| + \frac{1}{2} \cdot |x_{j}| \, \big].
	
	\end{array}$
\end{center}

Теперь будем использовать \textit{подход редукции размерности}. Представим функцию $f(x_{i}, \: x_{j})$ в виде
\begin{center}

	$\begin{array}{rcl}
	
		f(x_{i}, \: x_{j}) & = & \frac{1}{2} \cdot |x_{i}| \cdot |x_{j}| + \frac{1}{2} \cdot |x_{i}| \cdot |x_{j}| = \\
		
		& = & \frac{1}{2} \cdot f_{i}(x_{i}, \: x_{j}) + \frac{1}{2} \cdot f_{j}(x_{i}, \: x_{j}),	
	
	\end{array}$
	
	\vspace{2mm}	
	
	где $f_{i}(x_{i}, \: x_{j}) = |x_{i}| \cdot |x_{j}|, \qquad f_{j}(x_{i}, \: x_{j}) = |x_{i}| \cdot |x_{j}|.$
	
\end{center}

\noindentНайдём интервальное расширение функции $f_{i}(x_{i}, \: x_{j})$, считая что она зависит только от $x_{i}$:
\begin{center}

	$\mbf{f}_{\!\!i}(x_{i}, \: \mbf{x}_{j}) = [ \, 0, \: 1 \, ] \cdot x_{i}.$
	
\end{center}

\begin{figure}
\centering
    \includegraphics[width=0.5\linewidth]{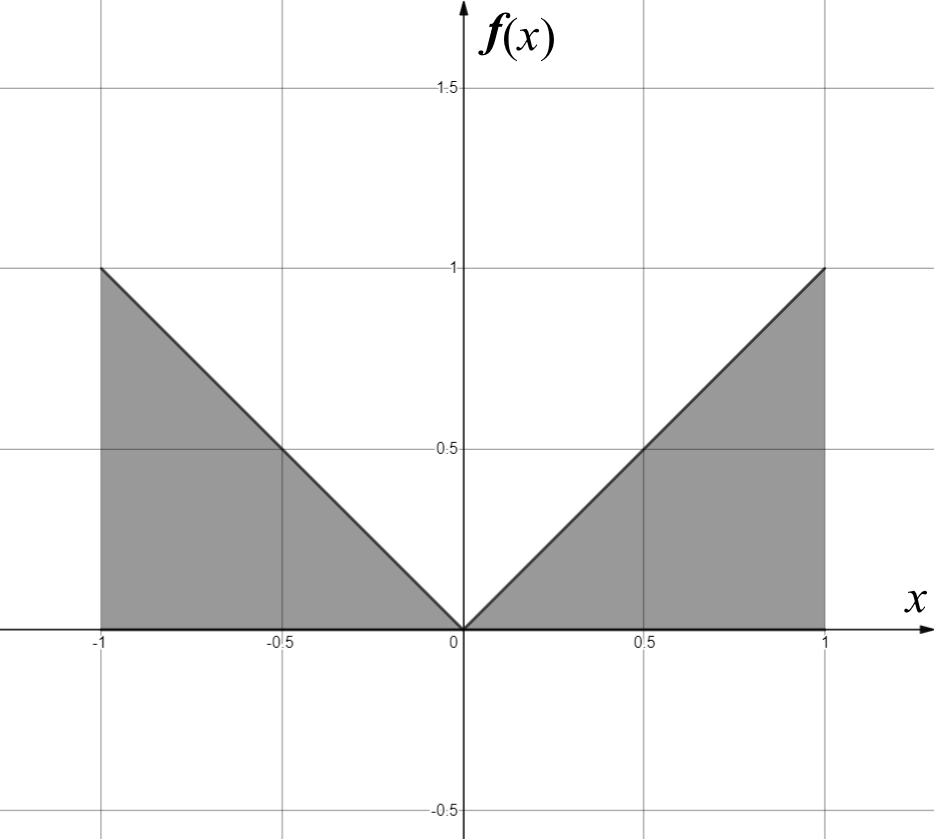}
    \caption{Множество совместных значений упорядоченной пары $\big( x, \: \mbf{f}(x) \big)$, где $\mbf{f}(x) = [ \, 0, \: 1 \, ] \cdot |x|, \: [ \, -1, \: 1 \, ].$}
\label{fig:abs3_example}
\end{figure}

\noindentМножество совместных значений упорядоченной пары $\big( x_{i}, \: \mbf{f}_{\!\!i}(x_{i}, \: \mbf{x}_{j}) \big)$ (рис. \ref{fig:abs3_example}) можно точно описать в центрально-ломаном базисе, как
\begin{center}

		$\mbf{f}_{\!\!i}(x_{i}, \: \mbf{x}_{j}) = \big[ \, 0, \: |x_{i}| \, \big].$
		
\end{center}

\noindentАналогично найдём интервальное расширение функции $f_{j}(x_{i}, \: x_{j})$ относительно $x_{j}$:
\begin{center}

	$f_{j}(\mbf{x}_{i}, \: x_{j}) = [ \, 0, \: 1 \, ] \cdot x_{j}.$

\end{center}

\noindentМножество совместных значений упорядоченной пары $\big( x_{j}, \: \mbf{f}_{\!\!j}(\mbf{x}_{i}, \: x_{j}) \big)$ (рис. \ref{fig:abs3_example}) можно точно описать в центрально-ломаном базисе как
\begin{center}

	$\mbf{f}_{\!\!j}(\mbf{x}_{i}, \: x_{j}) = \big[ \, 0, \: |x_{j}| \, \big].$
	
\end{center}

\noindentВозвращаясь к первоначальному представлению функции, получим
\begin{center}
	$\begin{array}{rcl}
	
	f(x_{i}, \: x_{j}) & = & \frac{1}{2} \cdot f_{i}(x_{i}, \: x_{j}) + \frac{1}{2} \cdot f_{j}(x_{i}, \: x_{j}) \subseteq \\
	
	& \subseteq & \frac{1}{2} \cdot \big[ \, 0, \: |x_{i}| \, \big] + \frac{1}{2} \cdot \big[ \, 0, \: |x_{j}| \, \big] = \big[ \, 0, \: \frac{1}{2} \cdot |x_{i}| + \frac{1}{2} \cdot |x_{j}| \, \big].
	
	\end{array}$
\end{center}

Сравнивая результаты двух подходов, снова видим преимущество подхода редукции размерности перед аналитическим.

\subsubsection{Умножение. Продолжение построения}

В предыдущих пунктах были найдены следующие интервальные оценки:

\begin{center}
	$\begin{array}{rcl}
		
		x_{i} ^ {2} & \subseteq & \big[ \, |x_{i}| - \frac{1}{4}, \: |x_{i}| \, \big], \\
		
		x_{i} \cdot |x_{i}| & \subseteq & \Big[ \, 2 \cdot \big( \sqrt{2} - 1 \big) \cdot x_{i} - \big( 3 - 2 \cdot \sqrt{2} \big), \\
		
		& & \hspace{5cm} 2 \cdot \big( \sqrt{2} - 1 \big) \cdot x_{i} + \big( 3 - 2 \cdot \sqrt{2} \big) \, \Big], \\
		
		x_{i} \cdot x_{j} & \subseteq & \big[ \, -\frac{1}{2} \cdot |x_{i}| -\frac{1}{2} \cdot |x_{j}|, \: \frac{1}{2} \cdot |x_{i}| + \frac{1}{2} \cdot |x_{j}| \, \big], \\
		
		|x_{i}| \cdot |x_{j}| & \subseteq & \big[ \, 0, \: \frac{1}{2} \cdot |x_{i}| + \frac{1}{2} \cdot |x_{j}| \, \big], \\	
		
		x_{i} \cdot |x_{j}| & \subseteq & \big[ \, -\frac{1}{2} \cdot x_{i} + \frac{1}{4} \cdot x_{j} + \frac{1}{4} \cdot |x_{j}|, \: \frac{1}{2} \cdot x_{i} + \frac{1}{4} \cdot x_{j} + \frac{1}{4} \cdot |x_{j}| \, \big], \\
		
	\end{array}$
	
	\vspace{3mm}	
	
	$i \neq j.$
\end{center}

Вернёмся к выражению операции умножения в центрально-ломаном базисе и используем полученные оценки (далее опустим аргументы функционалов границ):

\begin{center}
	$\begin{array}{rcl}
		
		L_{\mbf{a}} \cdot L_{\mbf{b}} & = & \sum_{i = 1}^{n} \Big( ( \underline{a}_{a, i} \cdot \underline{c}_{b} + \underline{a}_{b, i} \cdot \underline{c}_{a} ) \cdot x_{i} \: + \\
		
		& & \hspace{6mm} + \: ( \underline{b}_{a, i} \cdot \underline{c}_{b} + \underline{b}_{b, i} \cdot \underline{c}_{a} ) \cdot |x_{i}| \: + \\
		
		& &  \hspace{12mm} + \: ( \underline{a}_{a, i} \cdot \underline{a}_{b, i} + \underline{b}_{a, i} \cdot \underline{b}_{b, i} ) \cdot x_{i}^{2} \Big) \: + \\
		
		& + & \sum_{i = 1}^{n} \sum_{j = 1, j \neq i} \Big( ( \underline{b}_{a, i} \cdot \underline{b}_{b, j} ) \cdot |x_{i}| \cdot |x_{j}| \: + \\
		
		& & \hspace{6mm} + \: (\underline{a}_{a, i} \cdot \underline{a}_{b, j}) \cdot x_{i} \cdot x_{j} \: + \\
		
		& & \hspace{12mm} + \: (\underline{a}_{a, i} \cdot \underline{b}_{b, j} + \underline{a}_{b, i} \cdot \underline{b}_{a, j}) \cdot |x_{i}| \cdot x_{j} \Big) \subseteq \\
		
		& \subseteq & \sum_{i = 1}^{n} \Big( ( \underline{a}_{a, i} \cdot \underline{c}_{b} + \underline{a}_{b, i} \cdot \underline{c}_{a}) \cdot \big[ \, x_{i}, \: x_{i} \, \big] \: + \\
		
		& & \hspace{6mm} + \: (\underline{b}_{a, i} \underline{c}_{b} + \underline{b}_{b, i} \cdot \underline{c}_{a}) \cdot \big[ \, |x_{i}|, \: |x_{i}| \, \big] \: + \\
		
		& & \hspace{12mm} + \: (\underline{a}_{a, i} \cdot \underline{a}_{b, i} + \underline{b}_{a, i} \cdot \underline{b}_{b, i}) \cdot \big[ \, |x_{i}| - \frac{1}{4}, \: |x_{i}| \, \big] \Big) \: + \\
		
		& + & \sum_{i = 1}^{n} \sum_{j = 1, j \neq i}^{n} \Big( ( \underline{b}_{a, i} \cdot \underline{b}_{b, j}) \cdot \big[ \, 0, \: \frac{1}{2} \cdot |x_{i}| + \frac{1}{2} \cdot |x_{j}| \, \big] \: + \\
		
		& & \hspace{6mm} + \: (\underline{a}_{a, i} \cdot \underline{a}_{b, j}) \cdot \big[ \, -\frac{1}{2} \cdot |x_{i}| -\frac{1}{2} \cdot |x_{j}|, \: \frac{1}{2} \cdot |x_{i}| + \frac{1}{2} \cdot |x_{j}| \, \big] \: + \\
		
		& & \hspace{12mm} + \: (\underline{a}_{a, i} \cdot \underline{b}_{b, j} + \underline{a}_{b, i} \cdot \underline{b}_{a, j}) \cdot \big[ \, 0, \: \frac{1}{2} \cdot |x_{i}| + \frac{1}{2} \cdot |x_{j}| \, \big] \Big).
		
	\end{array}$
\end{center}

\noindentАналогичным образом вводятся интервальные оценки для
\begin{center}

	$L_{\mbf{a}}(x_{1}, \dots, x_{n}) \cdot U_{\mbf{b}}(x_{1}, \dots, x_{n})$,
	
	$U_{\mbf{a}}(x_{1}, \dots, x_{n}) \cdot L_{\mbf{b}}(x_{1}, \dots, x_{n})$,
	
	$U_{\mbf{a}}(x_{1}, \dots, x_{n}) \cdot U_{\mbf{b}}(x_{1}, \dots, x_{n})$.
	
\end{center}

Таким образом, была построена интервальная оценка для операции умножения в центрально-ломаном базисе.

\subsubsection{Построение процедур $F_{\downarrow}$ и $F_{\uparrow}$}

Отметим, что в силу свойств функционалов границ $L(x_{1}, \dots, x_{n})$ и $U(x_{1}, \dots, x_{n})$:
\begin{center}

	$L(x_{1}, \dots, x_{n}) \leq U(x_{1}, \dots, x_{n})$
	
	$x_{i} \in [ \, -1, \: 1 \, ], \qquad i = 1, \dots, n,$

\end{center}
верхняя граница будет получаться объединением верхних границ аргументов процедуры $F_{\uparrow}$, а нижняя граница --- нижних границ аргументов процедуры $F_{\downarrow}$.

Рассмотрим две какие-нибудь границы интервалов $\mathbb{F}\mathbb{R}$
\begin{center}

	$B_{1}(x_{1}, \dots, x_{n}) = \sum_{i = 1}^{n} a_{1, i} \cdot x_{i} + \sum_{i = 1}^{n} b_{1, i} \cdot |x_{i}| + c_{1}$,
	
	$B_{2}(x_{1}, \dots, x_{n}) = \sum_{i = 1}^{n} a_{2, i} \cdot x_{i} + \sum_{i = 1}^{n} b_{2, i} \cdot |x_{i}| + c_{2}$.

\end{center}

\noindentДля этих границ введём процедуру $F_{\uparrow} \big\{ \, B_{1}, \: B_{2} \, \big\}$ \big(далее опустим аргументы при $B_{1}(x_{1}, \dots, x_{n})$ и $B_{2}(x_{1}, \dots, x_{n})$\big):

\begin{center}
	$\begin{array}{rcl}
	
		F_{\uparrow} \big\{ \, B_{1}, \: B_{2} \, \big\} & = & F_{\uparrow} \big\{ \, \sum_{i = 1}^{n} a_{1, i} \cdot x_{i} + \sum_{i = 1}^{n} b_{1, i} \cdot |x_{i}| + c_{1}, \\
		
		& & \hspace{3cm} \sum_{i = 1}^{n} a_{2, i} \cdot x_{i} + \sum_{i = 1}^{n} b_{2, i} \cdot |x_{i}| + c_{2} \, \big\} \leq \\
		
		& \leq & \sum_{i = 1}^{n} F_{\uparrow} \big\{ \, a_{1, i} \cdot x_{i} + b_{1, i} \cdot |x_{i}|, \: a_{2, i} \cdot x_{i} + b_{2, i} \cdot |x_{i}| \, \big\} \: + \\
		
		& & \hspace{3cm} + \: F_{\uparrow} \big\{ \, c_{1}, \: c_{2} \, \big\} = \\
		
		& = & \sum_{i = 1}^{n} F_{\uparrow} \big\{ \, a_{1, i} \cdot x_{i} + b_{1, i} \cdot |x_{i}|, \: a_{2, i} \cdot x_{i} + b_{2, i} \cdot |x_{i}| \, \big\} \: + \\
	
		& & \hspace{3cm} + \: \text{max} \{ \, c_{1}, \: c_{2} \, \}. \\	
	
	\end{array}$

\end{center}

\noindentАналогичное выражение можно получить для процедуры $F_{\downarrow} \left\{ \, B_{1}, \: B_{2} \, \right\}$:
\begin{center}
	$\begin{array}{rcl}
	
		F_{\downarrow} \big\{ \, B_{1}, \: B_{2} \, \big\} & = & \sum_{i = 1}^{n} F_{\downarrow} \big\{ \, a_{1, i} \cdot x_{i} + b_{1, i} \cdot |x_{i}|, \: a_{2, i} \cdot x_{i} + b_{2, i} \cdot |x_{i}| \, \big\} \: + \\
		
		& & \hspace{3cm} + \: \text{min}\{ \, c_{1}, \: c_{2} \, \}
		
	\end{array}$

\end{center}

Теперь рассмотрим результат выражения
\begin{center}

	$F_{\uparrow} \big\{ \, a_{1} \cdot |x| + b_{1} \cdot x, \: a_{2} \cdot |x| + b_{2} \cdot |x| \, \big\}$.
	
\end{center}

\noindentПоскольку выбранный базис содержит вхождения линейной функции и функции модуля, то для построения достаточно рассмотреть три точки ломаной 
\begin{center}

	$x = -1, \qquad x = 0, \qquad x = 1$.
	
\end{center}

\noindentПо этим трём точкам можно построить верхнюю границу интервала $\mathbb{F}\mathbb{R}$. Аналогично можно построить
\begin{center}

	$F_{\downarrow} \big\{ \, a_{1} \cdot |x| + b_{1} \cdot x, \: a_{2} \cdot |x| + b_{2} \cdot |x| \, \big\}$.
	
\end{center}

\subsubsection{Деление}

По определению, данному ранее (далее опустим аргументы при функционалах границ и интервалах)

\begin{center}
    $\begin{array}{rcl}
    
    \mbf{a} \, / \, \mbf{b} & = & \Big[ \, F_{\downarrow} \big\{ \, L_{\mbf{a}} \, / \, L_{\mbf{b}}, \: L_{\mbf{a}} \, / \, U_{\mbf{b}}, \: U_{\mbf{a}} \, , \, L_{\mbf{b}}, \: U_{\mbf{a}} \, / \, U_{\mbf{b}} \, \big\}, \\
    
    & & \hspace{2cm} F_{\uparrow} \big\{ \, L_{\mbf{a}} \, / \, L_{\mbf{b}}, \: L_{\mbf{a}} \, / \, U_{\mbf{b}}, \: U_{\mbf{a}} \, / \, L_{\mbf{b}}, \: U_{\mbf{a}} \, / \, U_{\mbf{b}} \, \big\} \, \Big].
    
    \end{array}$ 
\end{center} 

Введём процедуры $F_{\downarrow}$ и $F_{\uparrow}$ для операции деления. Пусть в рассмотрении имеется функционал $H(x_{1}, \dots, x_{n})$, который является отношением двух функционалов
\begin{center}
	
	$H(x_{1}, \dots, x_{n}) = H_{n}(x_{1}, \dots, x_{n}) \, /  \, H_{d}(x_{1}, \dots, x_{n})$,
	
	$0 \notin \mbf{H}_{d}(x_{1}, \dots, x_{n})$.
	
\end{center}

Применим подход \textit{редукции размерности}. Перепишем функционал $H(x_{1}, \dots, x_{n})$ в виде
\begin{center}
	$\begin{array}{rcl}
	
		H(x_{1}, \dots, x_{n}) & = & H_{n}(x_{1}, \dots, x_{n}) \, / \, H_{d}(x_{1}, \dots, x_{n}) = \\
	
		& = & \sum_{i = 1}^{n} H_{i}(x_{i}) \, / \, H_{d}(\mbf{x}_{1}, \dots, \mbf{x}_{i - 1}, x_{i}, \mbf{x}_{i + 1}, \dots, \mbf{x}_{n}) \subseteq \\
	
		& \subseteq & \sum_{i = 1}^{n} \mbf{H}_{i}(x_{i}) \, / \, \mbf{H}_{d}(\mbf{x}_{1}, \dots, \mbf{x}_{i - 1}, x_{i}, \mbf{x}_{i + 1}, \dots, \mbf{x}_{n}) \subseteq \\
	
		& \subseteq & \bigg[ \, \sum_{i = 1}^{n} F_{\downarrow} \Big\{ L_{\!\mbf{H}_{i}}(x_{i}) \, / \, L_{\!\mbf{H}_{d}}(x_{i}), \: L_{\!\mbf{H}_{i}}(x_{i}) \, / \, U_{\!\mbf{H}_{d}}(x_{i}),  \\
	
		& & \hspace{3cm} U_{\!\mbf{H}_{i}}(x_{i}) \, / \, L_{\!\mbf{H}_{d}}(x_{i}), U_{\!\mbf{H}_{i}}(x_{i}) \, / \, U_{\!\mbf{H}_{d}}(x_{i}) \, \Big\}, \\
	
		& & \sum_{i = 1}^{n} F_{\uparrow} \Big\{ L_{\!\mbf{H}_{i}}(x_{i}) \, / \, L_{\!\mbf{H}_{d}}(x_{i}), \: L_{\!\mbf{H}_{i}}(x_{i}) \, / \, U_{\!\mbf{H}_{d}}(x_{i}), \\
	
		& & \hspace{3cm} U_{\!\mbf{H}_{i}}(x_{i}) \, / \, L_{\!\mbf{H}_{d}}(x_{i}), \: U_{\!\mbf{H}_{i}}(x_{i}) \, / \, U_{\!\mbf{H}_{d}}(x_{i}) \, \Big\} \bigg],

	\end{array}$
	
	\vspace{3mm}
	
	где $H_{i}(x_{i}) = a_{i} \cdot x_{i} + b_{i} \cdot |x_{i}| + \frac{c}{n}.$
\end{center}

\noindentТеперь для каждого $H_{i}$ найдём интервальную оценку на интервалах $[ \, -1, \: 0 \, ]$ и $[ \, 0, \: 1 \, ]$ с помощью аппарата чебышёвского альтернанса. Это возможно, так как на каждом из этих интервалов функция представляет собой выпуклую или вогнутую функцию, поскольку:
\begin{center}
	
	$\Big( \frac{a \, \cdot \, x \, + \, b}{c \, \cdot \, x \, + \, d} \Big)_{\hspace{-1mm} xx}^{\hspace{-1mm} ''} = \frac{2 \, \cdot \, c \, \cdot \, (b \, \cdot \, c \, - \, a \, \cdot \, d \,)}{(c \, \cdot \, x \, + \, d) ^ {3}},$ и если $0 \notin c \cdot \mbf{x} + d$, то 
	
	\vspace{2mm}
	
	$\Big( \frac{a \, \cdot \, x \, + \, b}{c \, \cdot \, x \, + \, d} \Big)_{\hspace{-1mm} xx}^{\hspace{-1mm} ''} \geq 0$ или $\Big( \frac{a \, \cdot \, x \, + \, b}{c \, \cdot \, x \, + \, d} \Big)_{\hspace{-1mm} xx}^{\hspace{-1mm} ''} \leq 0.$
	
\end{center}
По условию:
\begin{center}

	$0 \notin \mbf{H}_{d}(x_{1}, \dots, x_{n}),$
	
\end{center}
это значит, каждый аргумент процедур $F_{\downarrow}$ и $F_{\uparrow}$ будет вогнутой или выпуклой функцией на указанных интервалах.

После нахождения интервальных оценок на интервалах $[ \, -1, \: 0 \, ]$ и $[ \, 0, \: 1 \, ]$ проделаем процедуру <<склейки>>,  чтобы получить содержащий их интервал, представимый в центрально-ломаном базисе. Так можно найти интервальную оценку для результата операции деления между интервалами $\mathbb{F}\mathbb{R}$.

\clearpage
\section [Численные эксперименты] {Численные эксперименты}

\subsection[Характеристики уклонения\\ для сравнения интервальных арифметик]{Характеристики уклонения\\ для сравнения интервальных арифметик}

Можно рассматривать разные числовые характеристики, которые позволят сравнивать интервальные оценки, получаемые в различных интервальных арифметиках. В данной работе с этой целью будут рассматриваться две основные величины:
\begin{enumerate}

    \item Чебышёвское уклонение.
    
    \item Интегральное уклонение.
    
\end{enumerate}

\subsubsection{Чебышёвское уклонение}

Пусть имеется функция $f(x_{1}, \dots, x_{n})$, у которой найдена интервальная оценка
\begin{center}

	$\mbf{f}(x_{1}, \dots, x_{n}) = \big[ \, L_{\!\mbf{f}}(x_{1}, \dots, x_{n}), \: U_{\!\mbf{f}}(x_{1}, \dots, x_{n}) \, \big]$.
	
\end{center} 

Тогда \textit{чебышёвским уклонением} данной интервальной оценки будем называть величину
\begin{center}
	$\begin{array}{rcl}
	
    	T \big( \mbf{f}(x_{1}, \dots, x_{n}) \big) & = & \text{max}_{x_{1}, \dots, x_{n}} \big\{ \, U_{\!\mbf{f}}(x_{1}, \dots, x_{n}) - f(x_{1}, \dots, x_{n}), \\
    
    	& & \hspace{3cm} f(x_{1}, \dots, x_{n}) - L_{\!\mbf{f}}(x_{1}, \dots, x_{n}) \, \big\}.
    
	\end{array}$
\end{center}

Данная величина является оценкой сверху для погрешности описания функции базисом арифметики в чебышёвской норме
\begin{center}

	$\big\| f(x_{1}, \dots, x_{n}) \big\|_{\infty} = \text{max}_{x_{1}, \dots, x_{n}} \, | f(x_{1}, \dots, x_{n}) |.$
	
\end{center}

Чебышёвское уклонение не используется далее в работе, поскольку его точное вычисление затруднено знанием явного вида оцениваемой функции, но при этом доступна следующая оценка для этой величины:
\begin{center}

	$\begin{array}{rcl}
	
		T \big( \mbf{f}(x_{1}, \dots, x_{n}) \big) & \leq & \text{max}_{x_{1}, \dots, x_{n}} \big\{ U_{\!\mbf{f}}(x_{1}, \dots, x_{n}) - L_{\!\mbf{f}}(x_{1}, \dots, x_{n}) \big\} = \\

		& & = \| U_{\!\mbf{f}}(x_{1}, \dots, x_{n}) - L_{\!\mbf{f}}(x_{1}, \dots, x_{n}) \|_{\infty},
		
	\end{array}$
	
	\vspace{2mm}	
	
	это верно в силу $U_{\!\mbf{f}}(x_{1}, \dots, x_{n}) \geq f(x_{1}, \dots, x_{n}),$
	
	$L_{\!\mbf{f}}(x_{1}, \dots, x_{n}) \leq f(x_{1}, \dots, x_{n}).$

\end{center}

\subsubsection{Интегральное уклонение}

Пусть имеется функция $f(x_{1}, \dots, x_{n})$. Пусть также имеется интервальная оценка этой функции 
\begin{center}

	$\mbf{f}(x_{1}, \dots, x_{n}) = \Big[ \, L_{\mbf{f}}(x_{1}, \dots, x_{n}), \: U_{\mbf{f}}(x_{1}, \dots, x_{n}) \, \Big]$.

\end{center}

Тогда \textit{интегральным уклонением} данной интервальной оценки будем называть величину
\begin{center}
	$\begin{array}{rcl}
	
		I \big( \mbf{f}(x_{1}, \dots, x_{n}) \big) & = & \int_{-1}^{+1} \dots \int_{-1}^{+1} \Big( U_{\mbf{f}}(x_{1}, \dots, x_{n}) \: - \\
		
      & & \hspace{3cm} - \: f(x_{1}, \dots, x_{n}) \Big) \: dx_{1} \dots dx_{n} \: + \\
      
      & + & \int_{-1}^{+1} \dots \int_{-1}^{+1} \Big( f(x_{1}, \dots, x_{n}) \: - \\
      
      & & \hspace{3cm} - \: L_{\mbf{f}}(x_{1}, \dots, x_{n}) \Big) \: dx_{1} \dots dx_{n} = \\ 
      
      & = & \int_{-1}^{+1} \dots \int_{-1}^{+1} \Big( U_{\mbf{f}}(x_{1}, \dots, x_{n}) \: - \\
      
      & & \hspace{3cm} - \: L_{\mbf{f}}(x_{1}, \dots, x_{n}) \Big) \: dx_{1} \dots dx_{n}.
      
	\end{array}$
\end{center}

Данная характеристика является оценкой сверху для погрешности описания функции базисом арифметики в норме $L_{1}$
\begin{center}

	$\big\| f(x_{1}, \dots, x_{n}) \big\|_{1} = \int_{\mbf{x}_{1}} \dots \int_{\mbf{x}_{n}} \big| f(x_{1}, \dots, x_{n}) \big| \: dx_{1} \dots dx_{n}.$

\end{center}

Иными словами, интегральное уклонение показывает, какую по площади <<лишнюю>> область мы добавляем при описании функции выбранным в арифметике базисом. 

Достоинством данной характеристики является то, что её вычисление возможно без знания функции, у которой ищется интервальное расширение.

\subsection{Сравнение интервальных арифметик \\по интегральной характеристике уклонения}

Сравним классическую, функционально-базисную интервальные и аффинную арифметики по интегральному уклонению.

Далее приведены таблицы с результатами численных экспериментов (таб. \ref{tab1}, \ref{tab2}): в первом столбце представлены вычисляемые интервальные расширения выражений, во втором столбце --- значения интегрального уклонения для соответствующей арифметики.

После таблиц приведены пояснительные иллюстрации к некоторым примерам.

\begin{table}
\centering
\caption{\label{tab1}Сравнение интегрального уклонения $I$ для разных видов арифметик}
\vspace{3mm}
	\begin{small}
		\begin{tabular}{ c || c || c || c }
		
			 & Классическая & Аффинная & Функционально-граничная \\
		
			$f(x_{1}, \dots, x_{n})$ & интервальная & арифметика & интервальная \\
	
			 & арифметика &  & арифметика \\ \hline \hline \rule[-1mm]{0mm}{6mm}
	
			$1$ & $0$ & $0$ & $0$ \\
	
			$x$ & $4$ & $0$ & $0$ \\
	
			$x ^ {2}$ & $4$ & $4$ & $0.5$ \\
		
			$x ^ {3}$ & $4$ & $4$ & $0.81129150101524$ \\
	
			$x ^ {4}$ & $2$ & $4$ & $0.712626265847085$ \\
	
			$x ^ {5}$ & $4$ & $4$ & $0.96166573212244$ \\
			
			$x ^ {6}$ & $4$ & $4$ & $0.87830923447167$ \\
			
			$x ^ {7}$ & $4$ & $4$ & $1.037794268789845$ \\
	
			$x ^ {8}$ & $4$ & $4$ & $0.97285398143004$ \\
			
			$x ^ {9}$ & $4$ & $4$ & $1.08576770451591$ \\
			
			$x ^ {10}$ & $4$ & $4$ & $1.03154556274759$ \\ \hline \hline \rule[-1mm]{0mm}{6mm}
	
			$x \cdot [ \, -1, \: 0 \, ]$ & $4$ & $2$ & $1$ \\
			
			$x \cdot [ \, 0, \: 1 \, ]$ & $4$ & $2$ & $1$ \\	
			
			$x \cdot [ \, -1, \: -0.5 \, ]$ & $4$ & $1$ & $0.5$ \\
			
			$x \cdot [ \, -0.5, \: 0 \, ]$ & $2$ & $1$ & $0.5$ \\
			
			$x \cdot [ \, 0, \: 0.5 \, ]$ & $2$ & $1$ & $0.5$ \\
			
			$x \cdot [ \, 0.5, \: 1 \, ]$ & $4$ & $1$ & $0.5$ \\ 
	
			$x \cdot [ \, -1, \: 1 \, ]$ & $4$ & $4$ & $2$ \\ \hline \hline \rule[-1mm]{0mm}{6mm}
			
			$x \cdot y$ & $8$ & $8$ & $2$ \\
			
			$y \cdot x$ & $8$ & $8$ & $2$ \\ \hline \hline \rule[-1mm]{0mm}{6mm}
			
			$x ^ {2} \cdot y$ & $8$ & $8$ & $1.75$ \\
			
			$x \cdot y \cdot x$ & $8$ & $8$ & $2.51471862576143$ \\
			
			$y \cdot x ^ {2}$ & $8$ & $8$ & $2.51471862576143$ \\ \hline \hline \rule[-1mm]{0mm}{6mm}
			
			$x ^ {3} \cdot y$ & $8$ & $8$ & $2$ \\
			
			$x ^ {2} \cdot y \cdot x$ & $8$ & $8$ & $2.51471862576143$ \\
			
			$x \cdot y \cdot x^{2}$ & $8$ & $8$ & $2.68376618407357$ \\
			
			$y \cdot x^{3}$ & $8$ & $8$ & $2.68376618407357$ \\ \hline \hline \rule[-1mm]{0mm}{6mm}
			
			$x ^ {4} \cdot y$ & $8$ & $8$ & $1.68658008588991$ \\
			
			$x ^ {3} \cdot y \cdot x$ & $8$ & $8$ & $2.42640687119286$ \\
			
			$x ^ {2} \cdot y \cdot x ^ {2}$ & $8$ & $8$ & $2.68376618407357$ \\
			
			$x \cdot y \cdot x ^ {3}$ & $8$ & $8$ & $2.78344186487968$ \\
			
			$y \cdot x ^ {4}$ & $8$ & $8$ & $2.78344186487968$ \\ \hline \hline \rule[-1mm]{0mm}{6mm}
			
			$x ^ {5} \cdot y$ & $8$ & $8$ & $2$ \\
			
			$x ^ {4} \cdot y \cdot x$ & $8$ & $8$ & $2.50367965644036$ \\
			
			$x ^ {3} \cdot y \cdot x ^ {2}$ & $8$ & $8$ & $2.65476220843933$ \\
			
			$x ^ {2} \cdot y \cdot x ^ {3}$ & $8$ & $8$ & $2.78344186487968$ \\
			
			$x \cdot y \cdot x ^ {5}$ & $8$ & $8$ &$2.83068004995108$ \\
			
			$y \cdot x ^ {5}$ & $8$ & $8$ & $2.83068004995108$

		\end{tabular}
	\end{small}
\end{table}

\begin{table}
\centering
\caption{\label{tab2}Сравнение интегральной характеристики $I$ для разных видов интервальных арифметик. Операция деления.}
\vspace{3mm}
	\begin{small}
		\begin{tabular}{ c || c || c }
		
			 & Классическая & Функционально-граничная \\
		
			$f(x_{1}, \dots, x_{n})$ & \hspace{5mm} интервальная \hspace{5mm} & интервальная \\
	
			 & арифметика & арифметика \\ \hline \hline \rule[-1mm]{0mm}{6mm}
	
			$(x + 2) \, / \, (x + 2)$ & $5.333333333333333$ & $0.000000000000002$ \\
	
			$1 \, / \, (x + 2)$ & $1.333333333333333$ & $0.154335453948487$ \\
	
			$1 \, / \, (x + 3)$ & $0.5$ & $0.0309600827457410$ \\
		
			$1 \, / \, (x + 4)$ & $0.266666666666667$ & $0.0111669633764997$ \\
	
			$1 \, / \, (x + 5)$ & $0.166666666666667$ & $0.0052557816242133$ \\ \hline \hline \rule[-1mm]{0mm}{6mm}
	
			$1 \, / \, (x + [ \, 2, \: 3 \, ])$ & $1.5$ & $0.473886926876691$ \\
			
			$1 \, / \, (x + [ \, 3, \: 4 \, ])$ & $0.6$ & $0.20859862653166603$ \\
			
			$1 \, / \, (x + [ \, 4, \: 5 \, ])$ & $0.333333333333333$ & $0.118389836375957$ \\
	
			$1 \, / \, (x + [ \, 5, \: 10 \, ])$ & $0.318181818181818$ & $0.21638703599257883$ \\ \hline \hline \rule[-1mm]{0mm}{6mm}
	
			$x \, / \, [ \, 1, \: 2 \, ]$ & $4$ & $0.25$ \\
			
			$x \, / \, [ \, 2, \: 3 \, ]$ & $2$ & $0.083333333333334$ \\	
			
			$x \, / \, [ \, 3, \: 4 \, ]$ & $1.333333333333333$ & $0.0416666666666671$ \\
			
			$x \, / \, [ \, 4, \: 5 \, ]$ & $1$ & $0.0250000000000002$ \\
			
			$x \, / \, [ \, 5, \: 10 \, ]$ & $0.8$ & $0.050000000000001$ \\ \hline \hline \rule[-1mm]{0mm}{6mm}
			
			$x \, / \, (y + 2)$ & $8$ & $0.666666666666671$ \\
			
			$y \, / \, (x + 2)$ & $8$ & $0.666666666666671$ \\ \hline \hline \rule[-1mm]{0mm}{6mm}
			
			$x ^ {2} \, / \, (y + 2)$ & $8$ & $1.35416666666667$ \\
			
			$x \, / \, (y + 2) \cdot x$ & $8$ & $1.67647908384096$ \\
			
			$1 \, / \, (y + 2) \cdot x ^ {2}$ & $8$ & $1.74235739103973$ \\ \hline \hline \rule[-1mm]{0mm}{6mm}
			
			$x ^ {3} \, / \, (y + 2)$ & $8$ & $1.63645985547802$ \\
			
			$x ^ {2} \, / \, (y + 2) \cdot x$ & $8$ & $1.97449798090104$  \\
			
			$x \, / \, (y + 2) \cdot x^{2}$ & $8$ & $1.93790283299493$ \\
			
			$1 \, / \, (y + 2) \cdot x^{3}$ & $8$ & $2.03222087986149$ \\ \hline \hline \rule[-1mm]{0mm}{6mm}
			
			$x ^ {4} \, / \, (y + 2)$ & $8$ & $1.44449780742915$ \\
			
			$x ^ {3} \, / \, (y + 2) \cdot x$ & $8$ & $1.7585136523794$ \\
			
			$x ^ {2} \, / \, (y + 2) \cdot x ^ {2}$ & $8$ & $1.78234721543857$ \\
			
			$x \, / \, (y + 2) \cdot x ^ {3}$ & $8$ & $1.81505575983538$ \\
			
			$1 \, / \, (y + 2) \cdot x ^ {4}$ & $8$ & $1.98006435610048$ \\ \hline \hline \rule[-1mm]{0mm}{6mm}
			
			$x ^ {5} \, / \, (y + 2)$ & $8$ & $1.72751020178193$ \\
			
			$x ^ {4} \, / \, (y + 2) \cdot x$ & $8$ & $2.04570366923249$ \\
			
			$x ^ {3} \, / \, (y + 2) \cdot x ^ {2}$ & $8$ & $2.06395221667211$ \\
			
			$x ^ {2} \, / \, (y + 2) \cdot x ^ {3}$ & $8$ & $2.08978012172187$ \\
			
			$x \, / \, (y + 2) \cdot x ^ {5}$ & $8$ & $2.09916240309717$ \\
			
			$1 \, / \, (y + 2) \cdot x ^ {5}$ & $8$ & $2.20737939401541$

		\end{tabular}
	\end{small}
\end{table}

\begin{figure}
	\centering
    \includegraphics[width=0.47\linewidth]{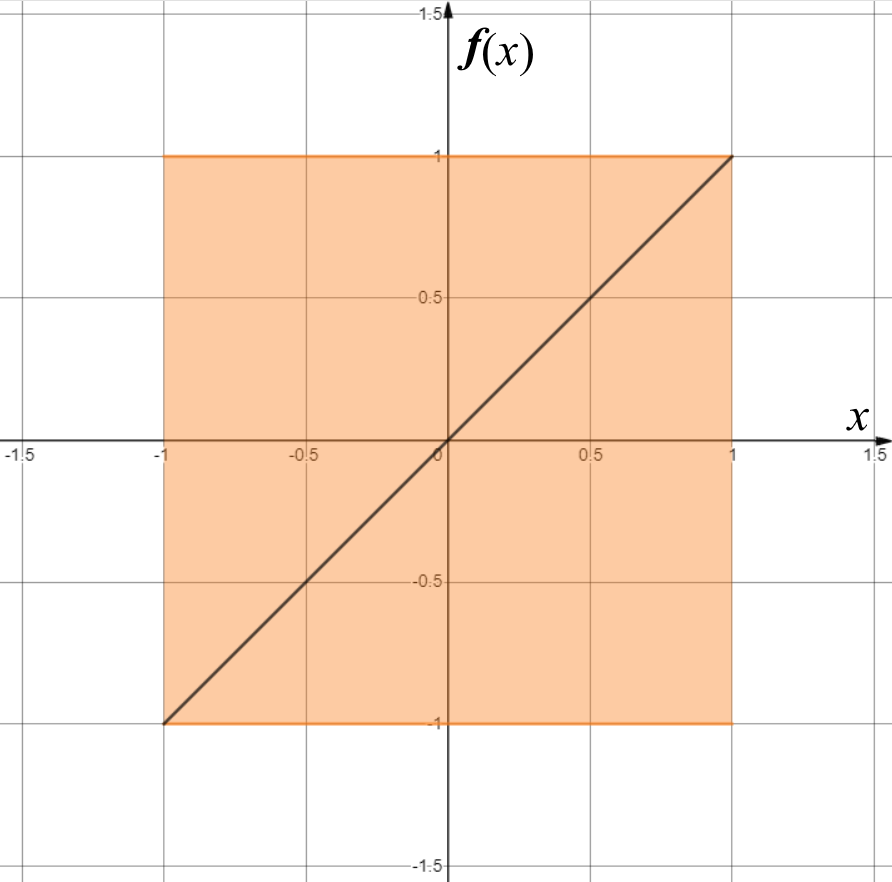}
    \hspace{0.04\linewidth}
    \includegraphics[width=0.47\linewidth]{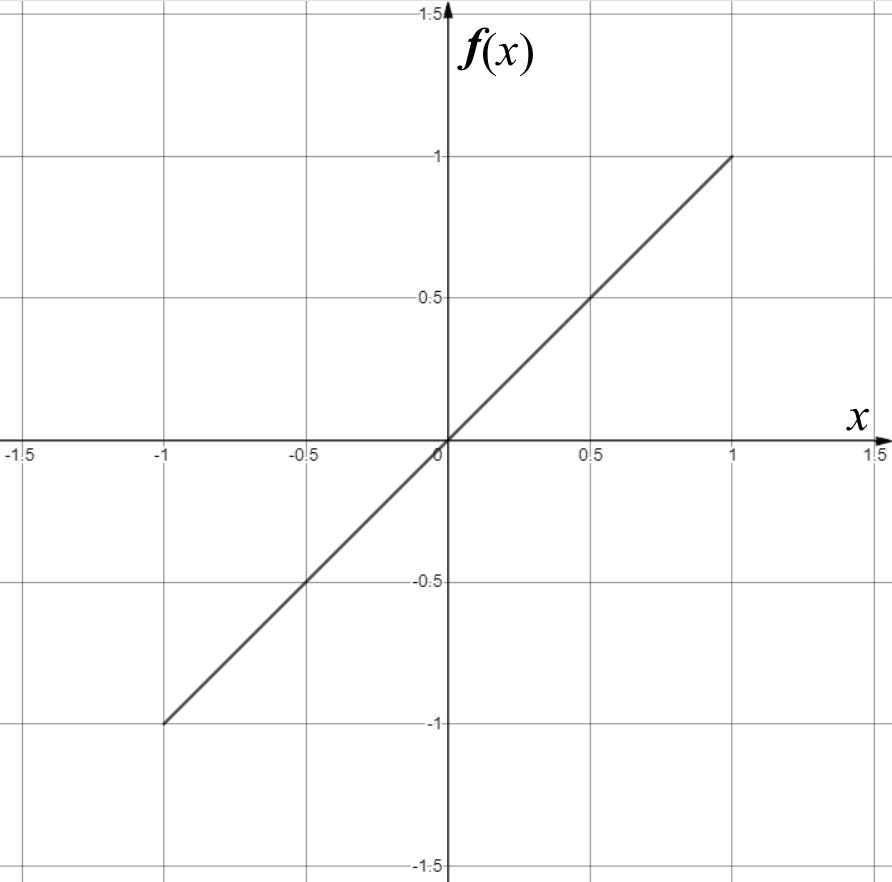}
    \caption{Множество совместных значений упорядоченной пары $\big( x, \: f(x) \big)$, где $f(x) = x$ (чёрная линия). На левом рисунке показано его приближение в классической интервальной арифметике (оранжевая область), а также в аффинной арифметике (совпадает с множеством совместных значений). На правом --- в функционально-граничной интервальной арифметике (совпадает с множеством совместных значений).}
	\label{fig:x_difference}
\end{figure}

\begin{figure}
	\centering
    \includegraphics[width=0.47\linewidth]{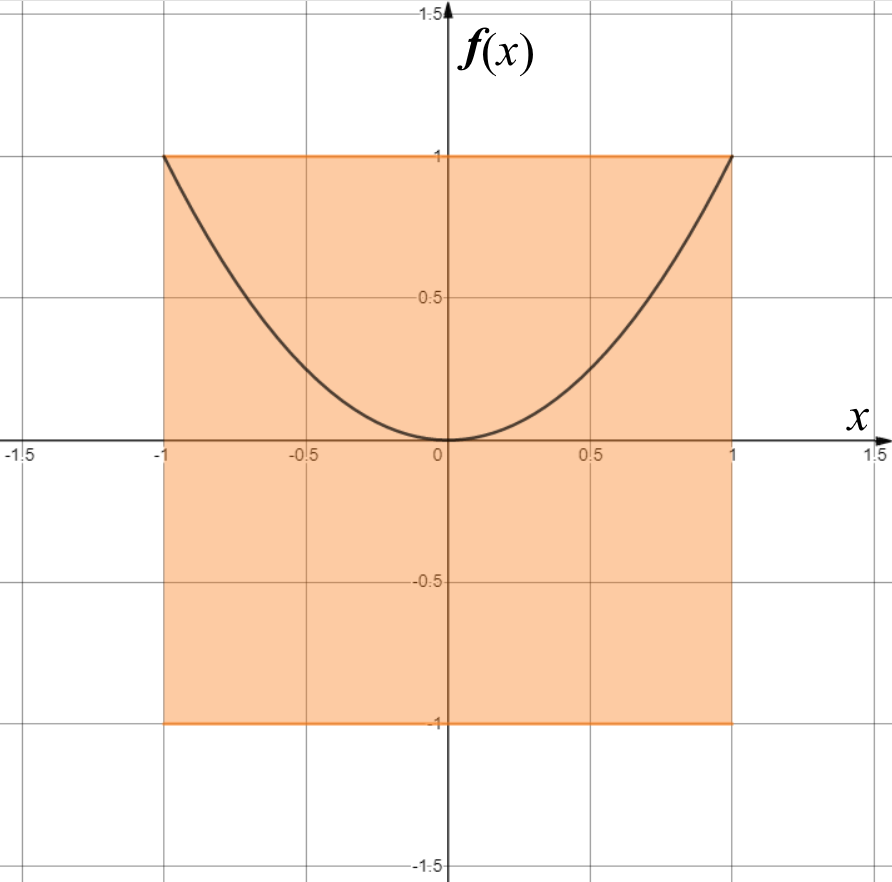}
    \hspace{0.04\linewidth}
    \includegraphics[width=0.47\linewidth]{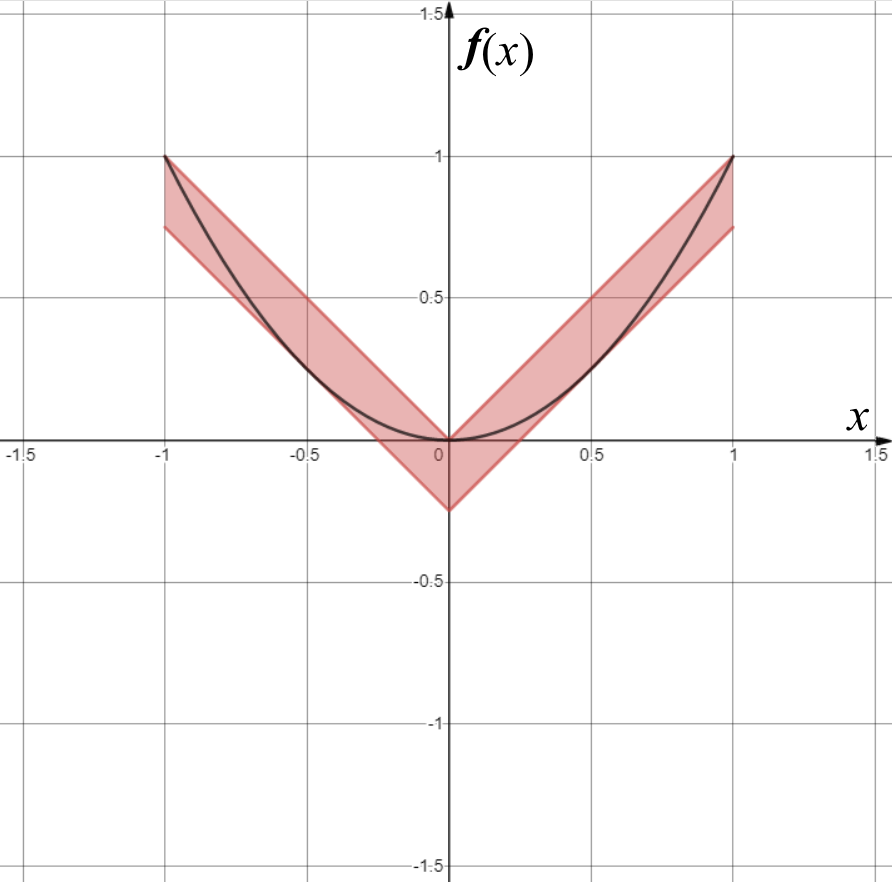}
    \caption{Множество совместных значений упорядоченной пары $\big( x, \: f(x) \big)$, где $f(x) = x^{2}$ (чёрная линия). На левом рисунке показано его приближение в классической интервальной и аффинной арифметиках (оранжевая область), на правом --- в функционально-граничной интервальной арифметике (красная область).}
	\label{fig:x2_difference}
\end{figure}

\begin{figure}
	\centering
    \includegraphics[width=0.47\linewidth]{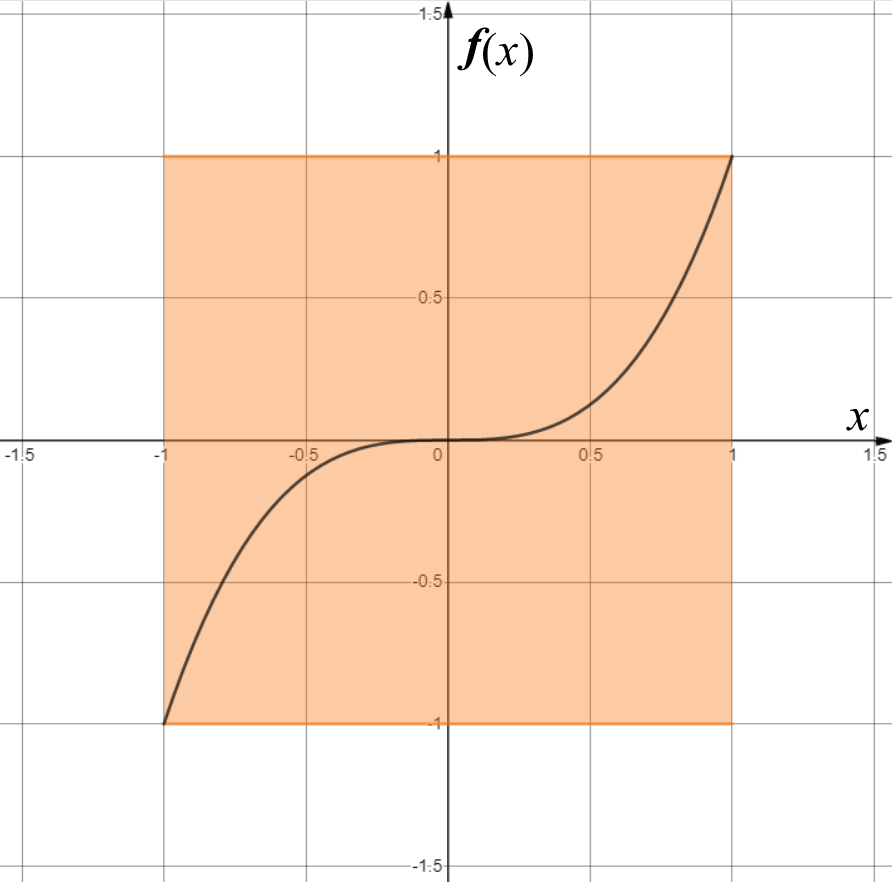}
    \hspace{0.04\linewidth}
    \includegraphics[width=0.47\linewidth]{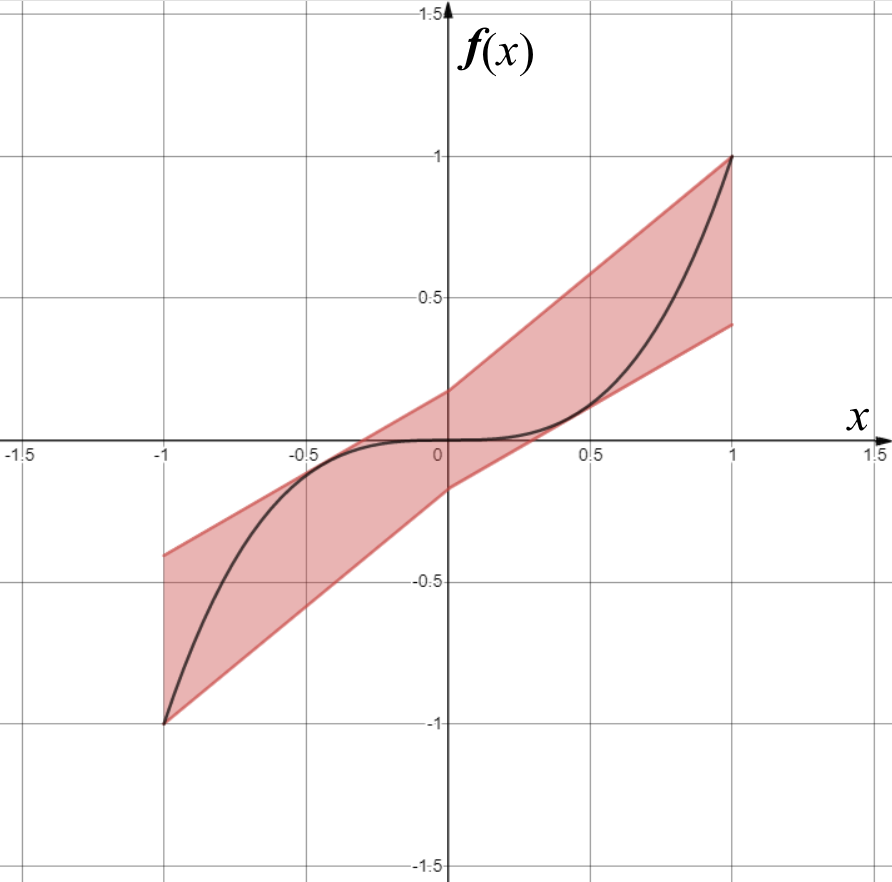}
    \caption{Множество совместных значений упорядоченной пары $\big( x, \: f(x) \big)$, где $f(x) = x^{3}$ (чёрная линия). На левом рисунке показано его приближение в классической интервальной и аффинной арифметиках (оранжевая область), на правом --- в функционально-граничной интервальной арифметике (красная область).}
	\label{fig:x3_difference}
\end{figure}

\begin{figure}
	\centering
    \includegraphics[width=0.47\linewidth]{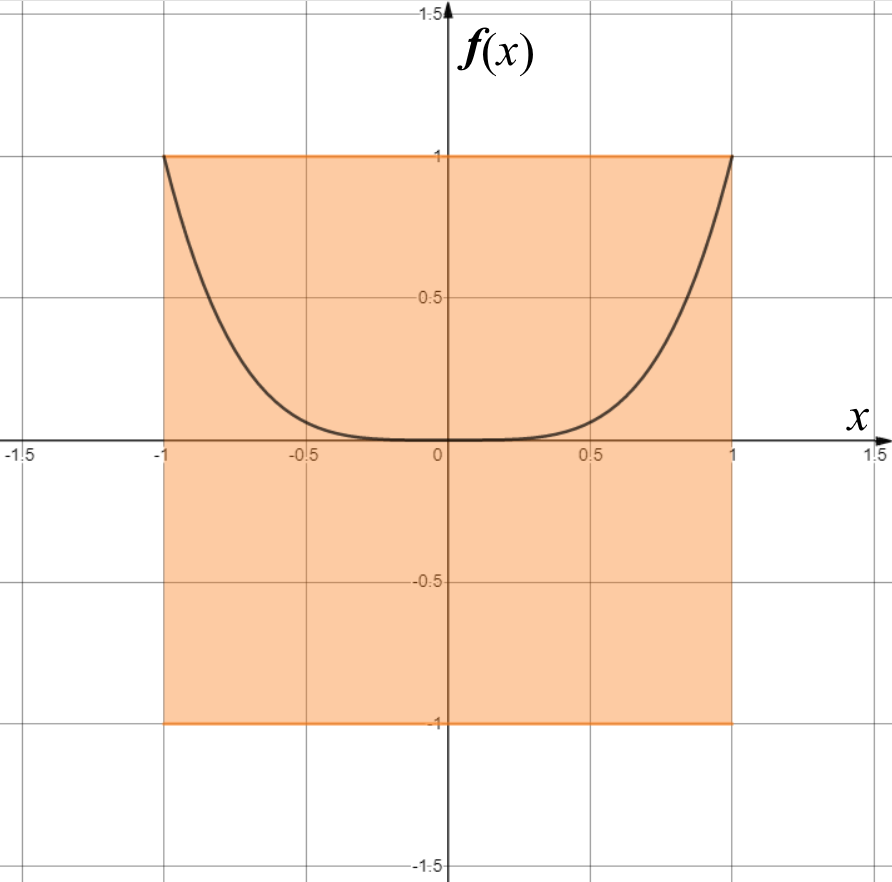}
    \hspace{0.04\linewidth}
    \includegraphics[width=0.47\linewidth]{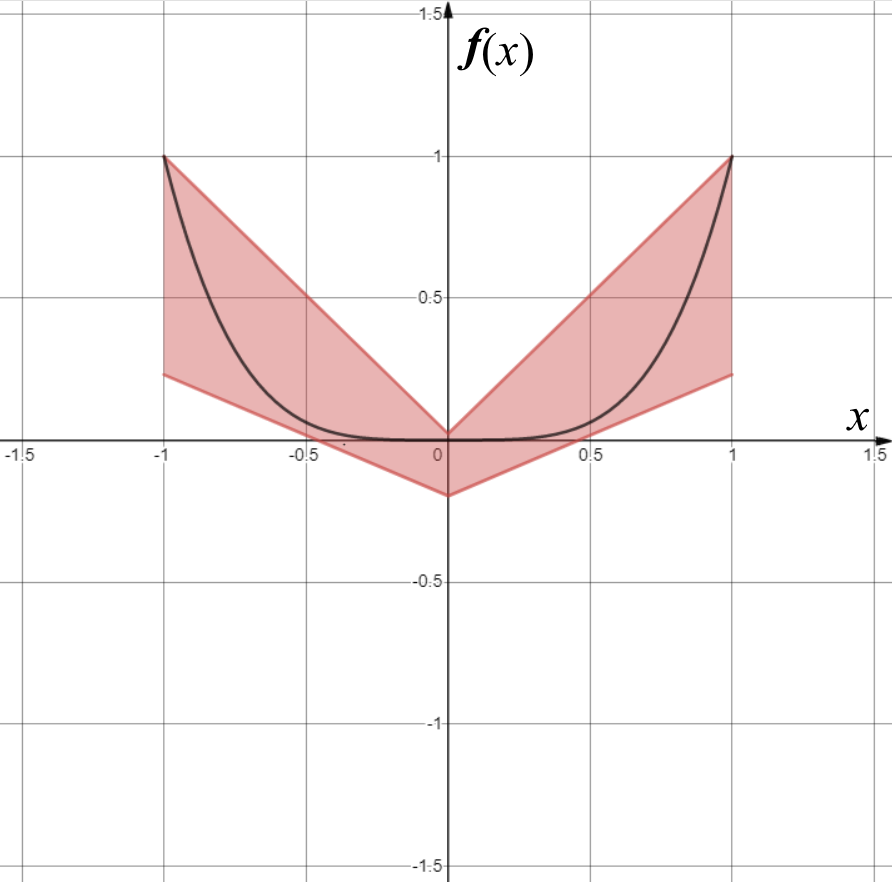}
    \caption{Множество совместных значений упорядоченной пары $\big( x, \: f(x) \big)$, где $f(x) = x^{4}$ (чёрная линия). На левом рисунке показано его приближение в классической интервальной и аффинной арифметике (оранжевая область), на правом --- в функционально-граничной интервальной арифметике (красная область).}
	\label{fig:x4_difference}
\end{figure}

\begin{figure}
	\centering
    \includegraphics[width=0.47\linewidth]{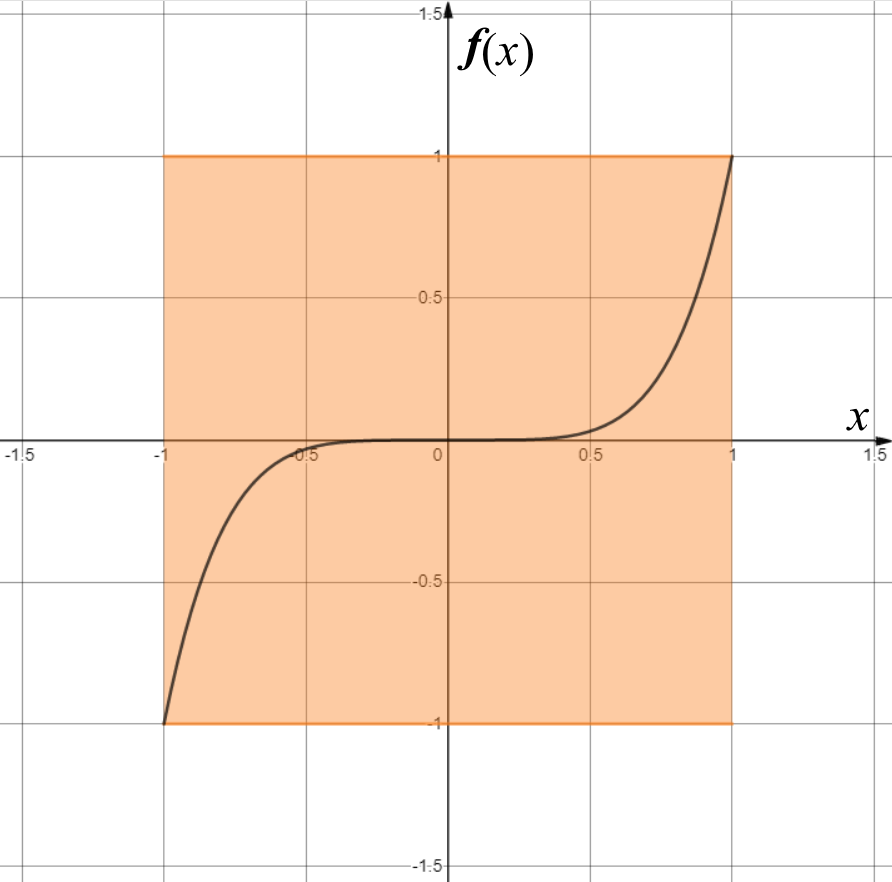}
    \hspace{0.04\linewidth}
    \includegraphics[width=0.47\linewidth]{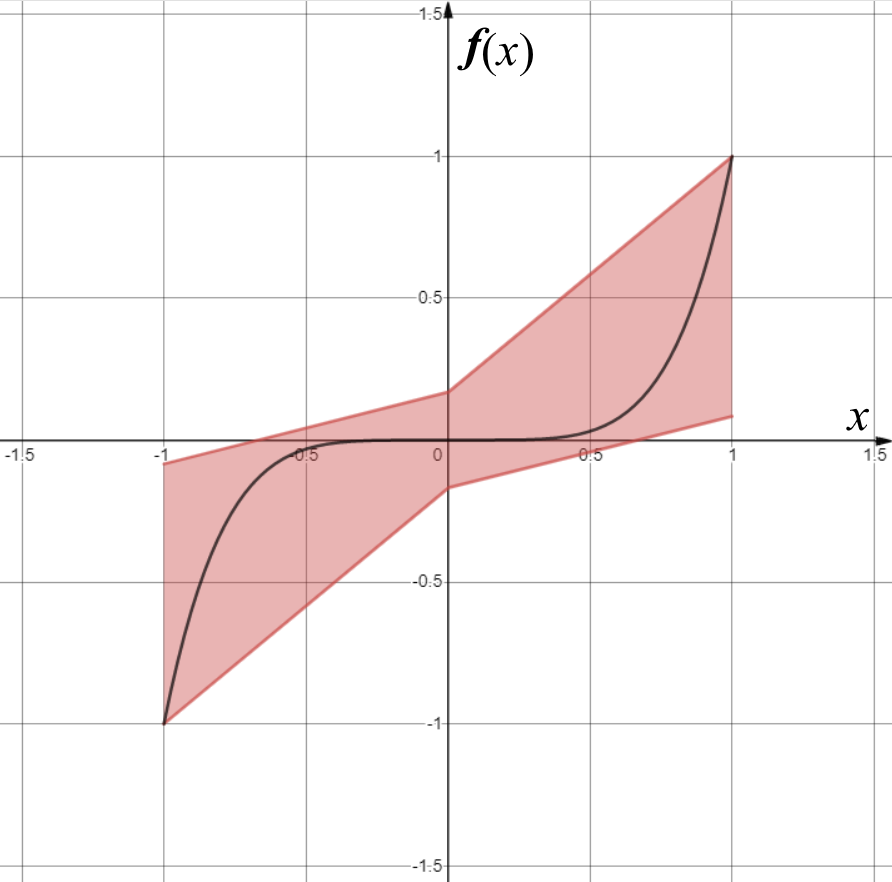}
    \caption{Множество совместных значений упорядоченной пары $\big( x, \: f(x) \big)$, где $f(x) = x^{5}$ (чёрная линия). На левом рисунке показано его приближение в классической интервальной и аффинной арифметиках (оранжевая область), на правом --- в функционально-граничной интервальной арифметике (красная область).}
	\label{fig:x5_difference}
\end{figure}

\begin{figure}
	\centering
    \includegraphics[width=0.47\linewidth]{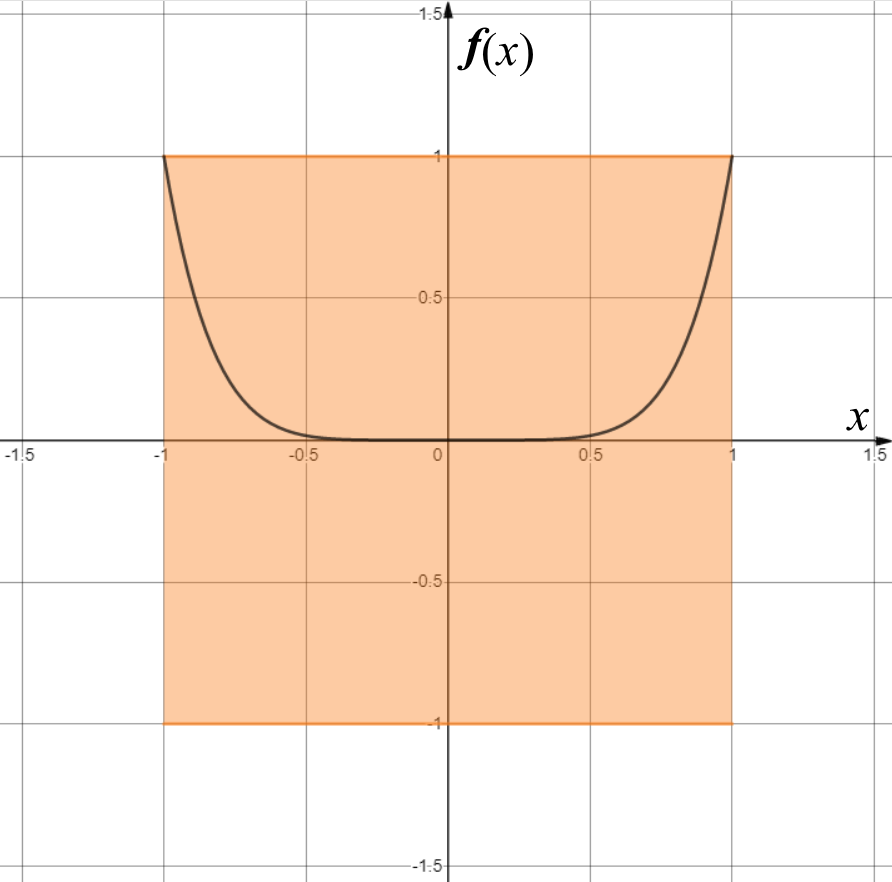}
    \hspace{0.04\linewidth}
    \includegraphics[width=0.47\linewidth]{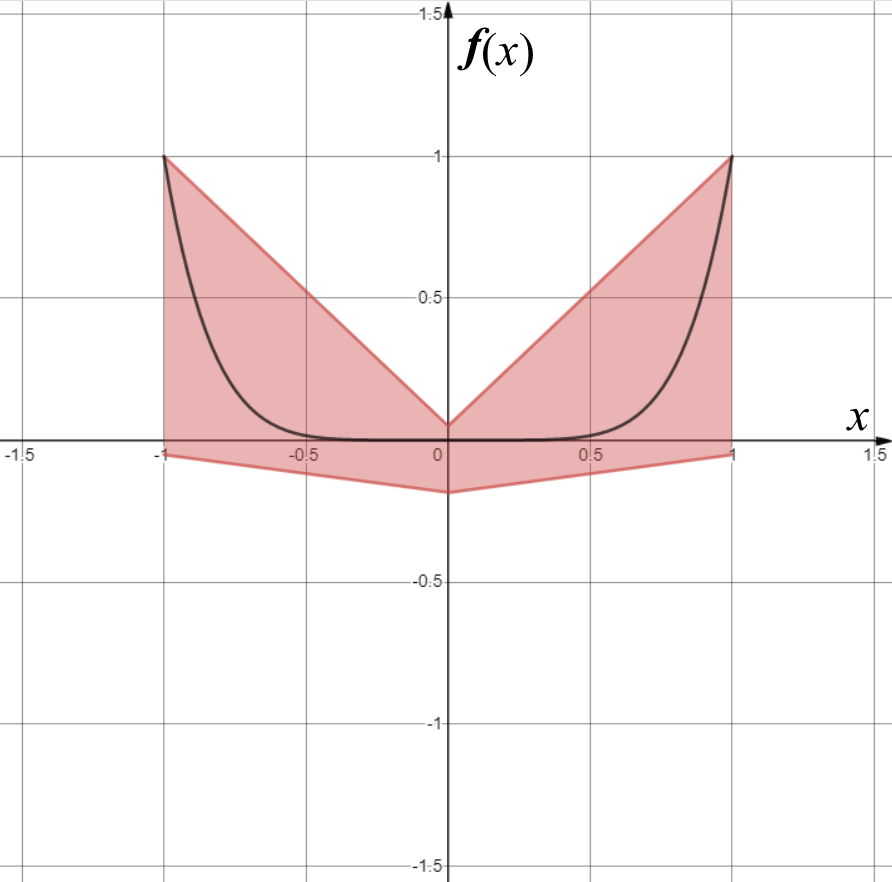}
    \caption{Множество совместных значений упорядоченной пары $\big( x, \: f(x) \big)$, где $f(x) = x^{6}$ (чёрная линия). На левом рисунке показано его приближение в классической интервальной и аффинной арифметиках (оранжевая область), на правом --- в функционально-граничной интервальной арифметике(красная область).}
	\label{fig:x6_difference}
\end{figure}

\begin{figure}
	\centering
    \includegraphics[width=0.47\linewidth]{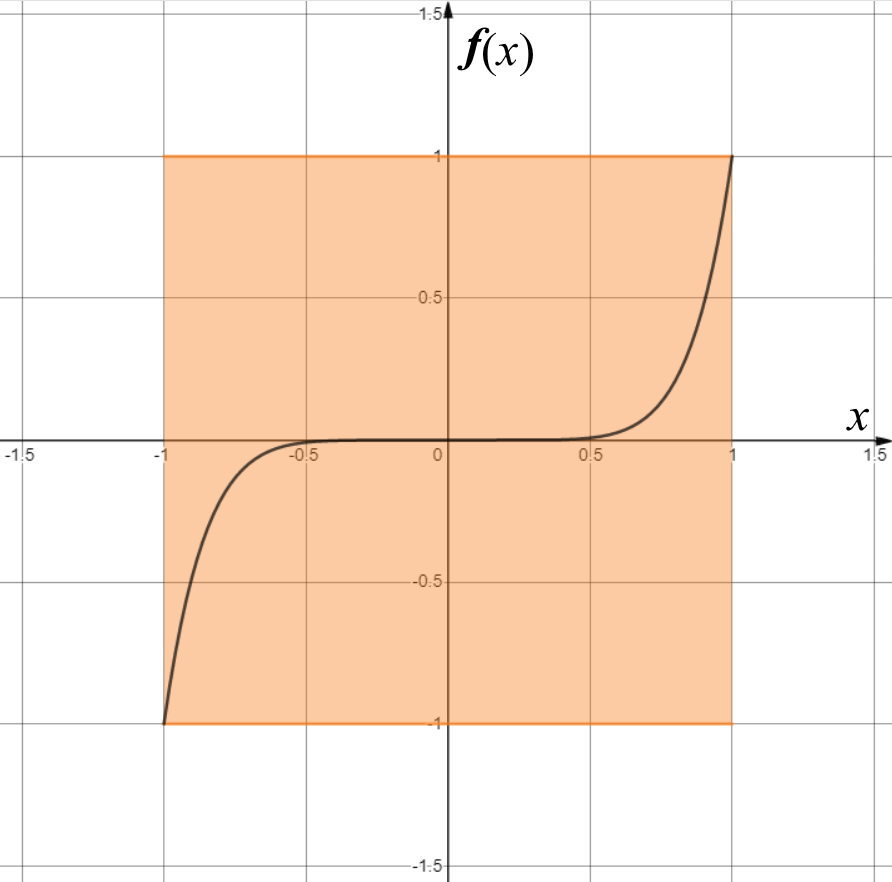}
    \hspace{0.04\linewidth}
    \includegraphics[width=0.47\linewidth]{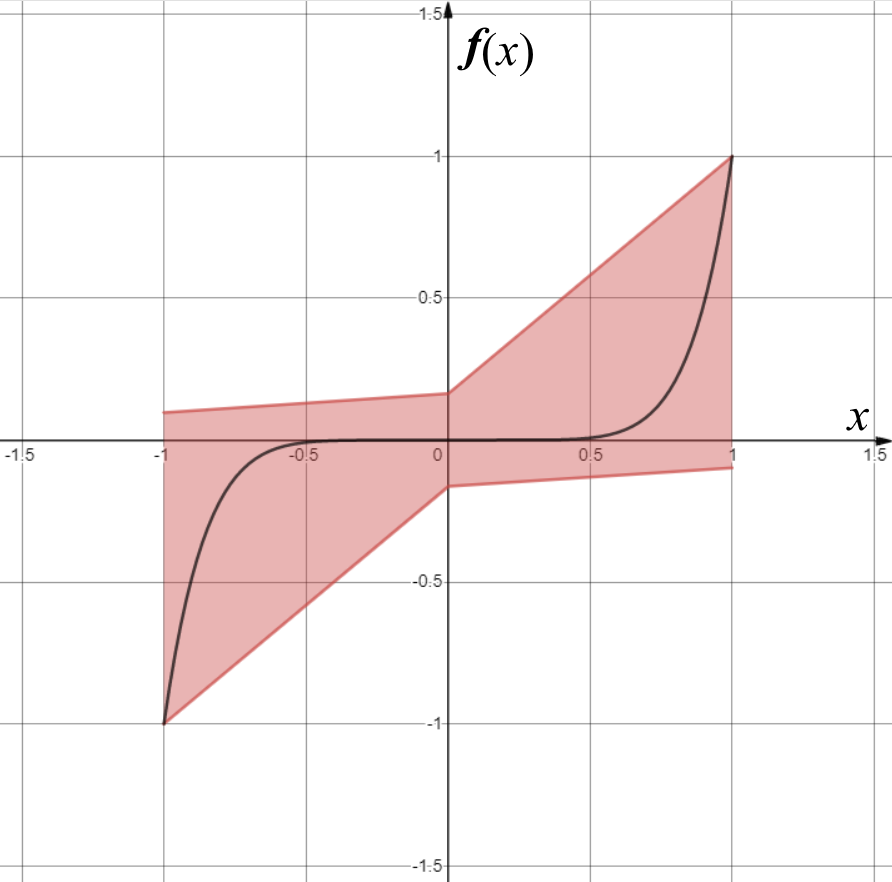}
    \caption{Множество совместных значений упорядоченной пары $\big( x, \: f(x) \big)$, где $f(x) = x^{7}$ (чёрная линия). На левом рисунке показано его приближение в классической интервальной и аффинной арифметиках (оранжевая область), на правом --- в функционально-граничной интервальной арифметике (красная область).}
	\label{fig:x7_difference}
\end{figure}

\begin{figure}
	\centering
    \includegraphics[width=0.47\linewidth]{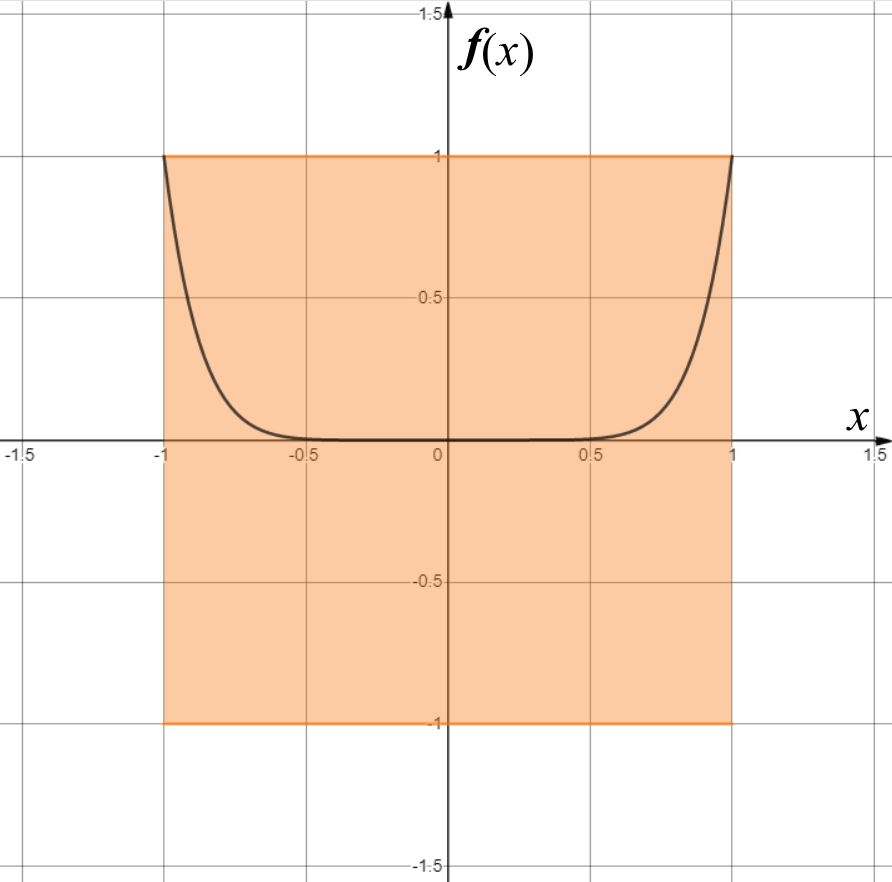}
    \hspace{0.04\linewidth}
    \includegraphics[width=0.47\linewidth]{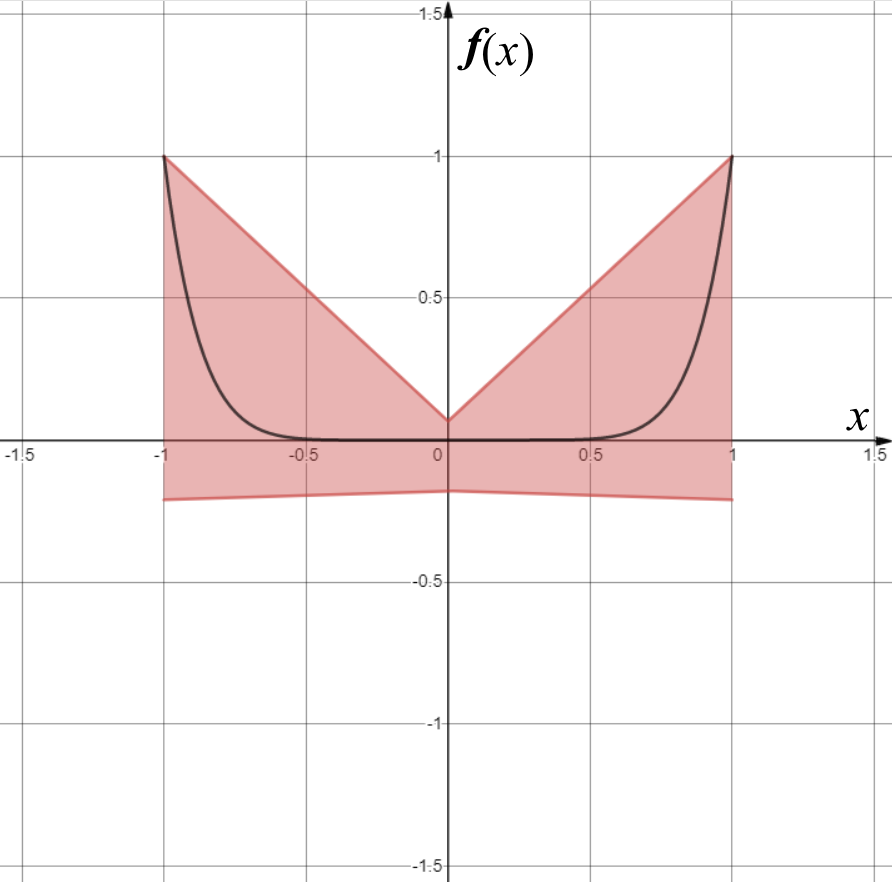}
    \caption{Множество совместных значений упорядоченнйо пары $\big( x, \: f(x) \big)$, где $f(x) = x^{8}$ (чёрная линия). На левом рисунке показано его приближение в классической интервальной и аффинной арифметиках (оранжевая область), на правом --- в функционально-граничной интервальной арифметике (красная область).}
	\label{fig:x8_difference}
\end{figure}

\begin{figure}
	\centering
    \includegraphics[width=0.47\linewidth]{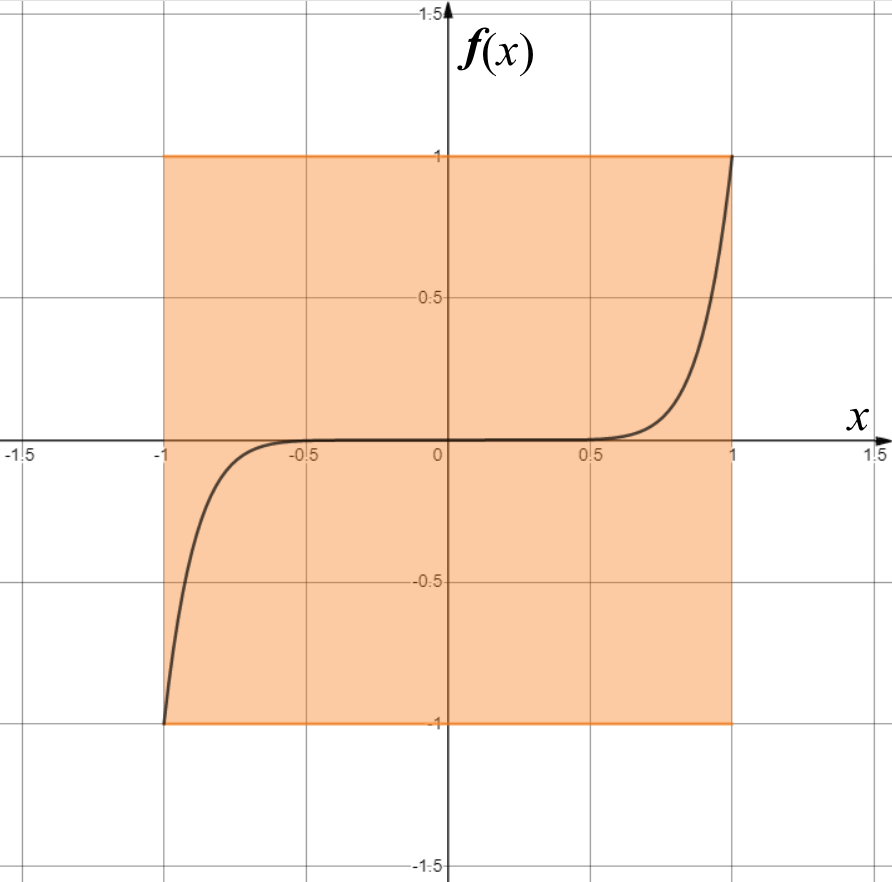}
    \hspace{0.04\linewidth}
    \includegraphics[width=0.47\linewidth]{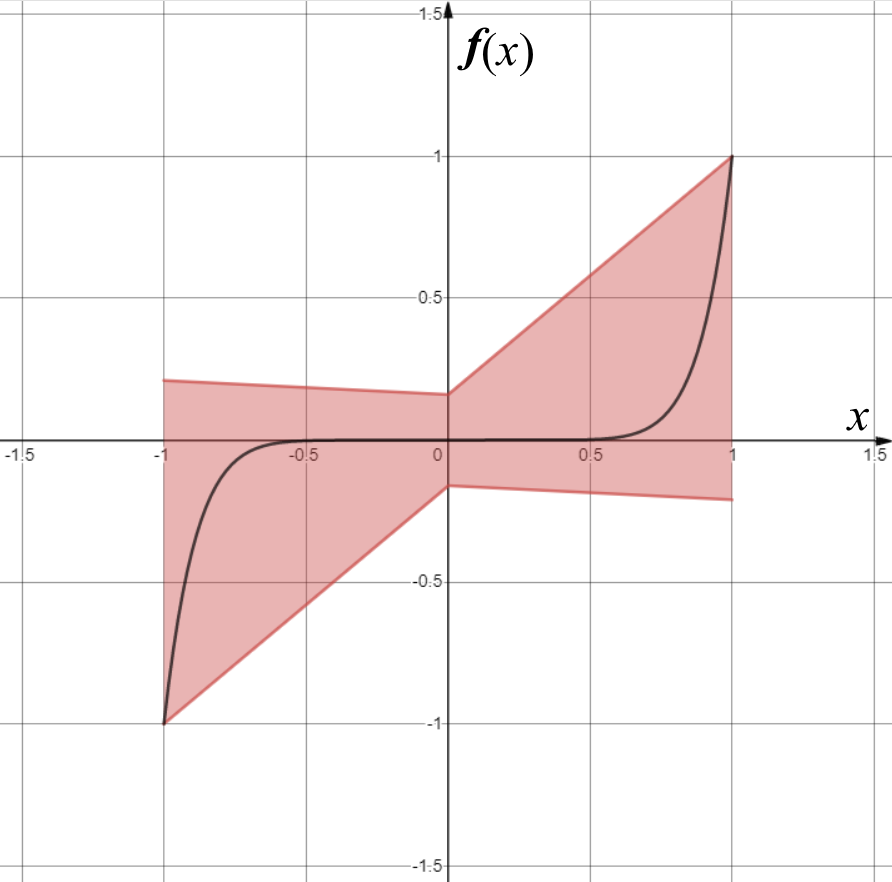}
    \caption{Совместное множество значений упорядоченной пары $\big( x, \: f(x) \big)$, где $f(x) = x^{9}$ (чёрная линия). На левом рисунке показано его приближение в классической интервальной и аффинной арифметиках (оранжевая область), на правом --- в функционально-граничной интервальной арифметике (красная область).}
	\label{fig:x9_difference}
\end{figure}

\begin{figure}
	\centering
    \includegraphics[width=0.47\linewidth]{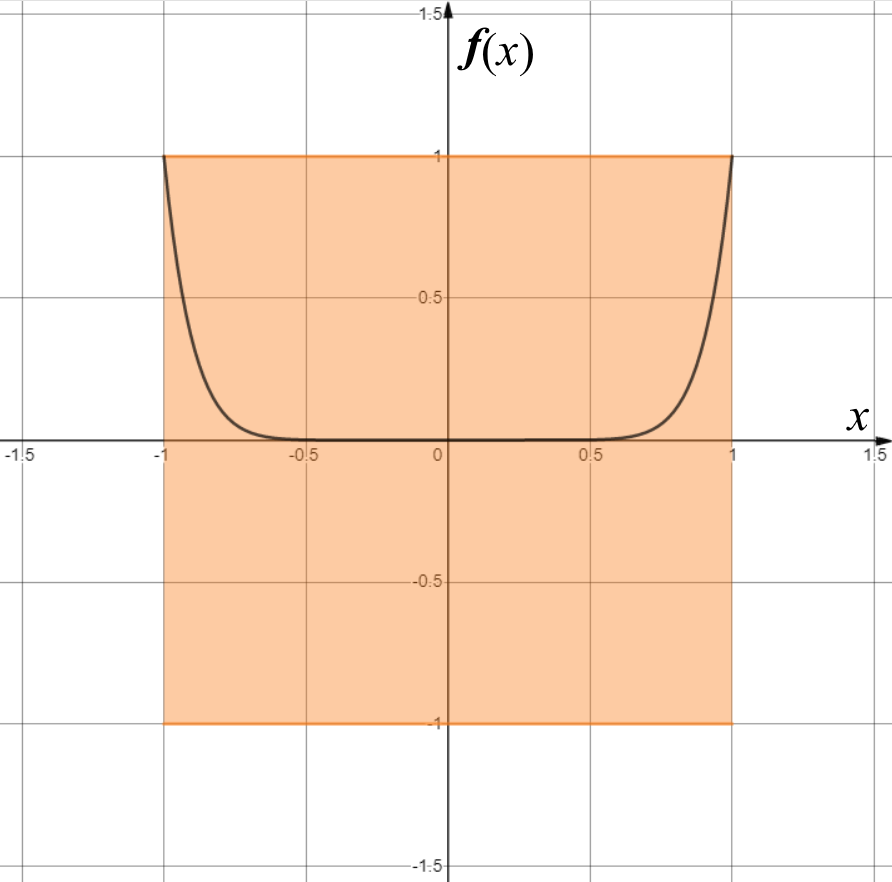}
    \hspace{0.04\linewidth}
    \includegraphics[width=0.47\linewidth]{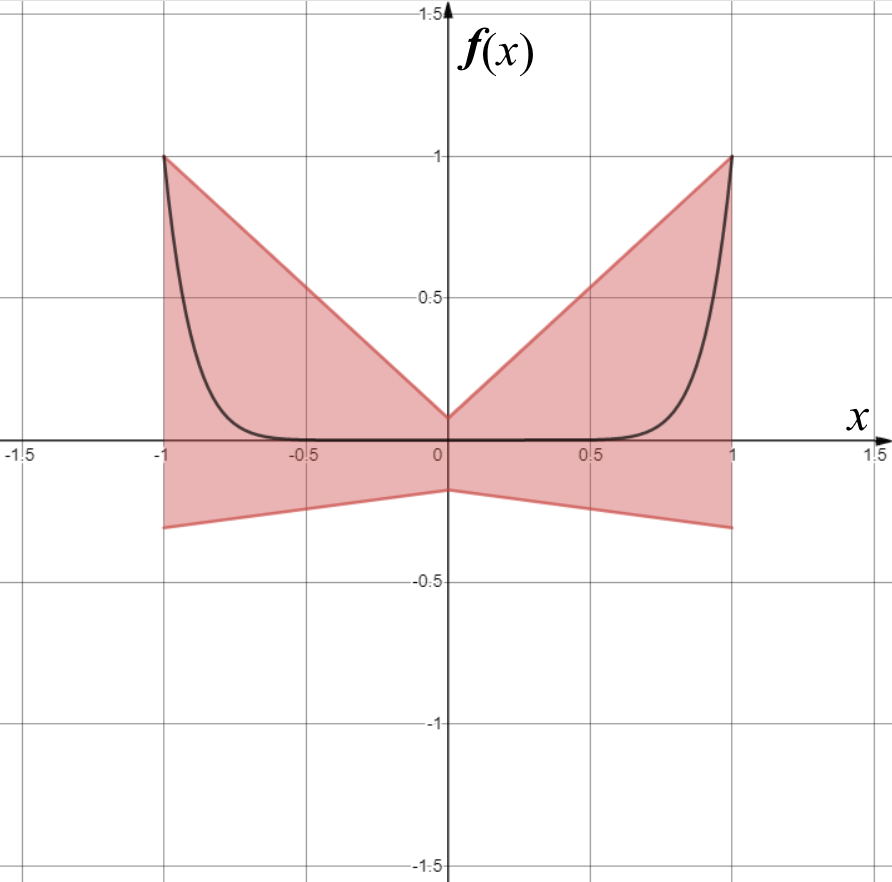}
    \caption{Множество совместных значений упорядоченной пары $\big( x, \: f(x) \big)$, где $f(x) = x^{10}$ (чёрная линия). На левом рисунке показано его приближение в классической интервальной и аффинной арифметиках (оранжевая область), на правом --- в функционально-граничной интервальной арифметике (красная область).}
	\label{fig:x10_difference}
\end{figure}

\begin{figure}
	\centering
    \includegraphics[width=0.47\linewidth]{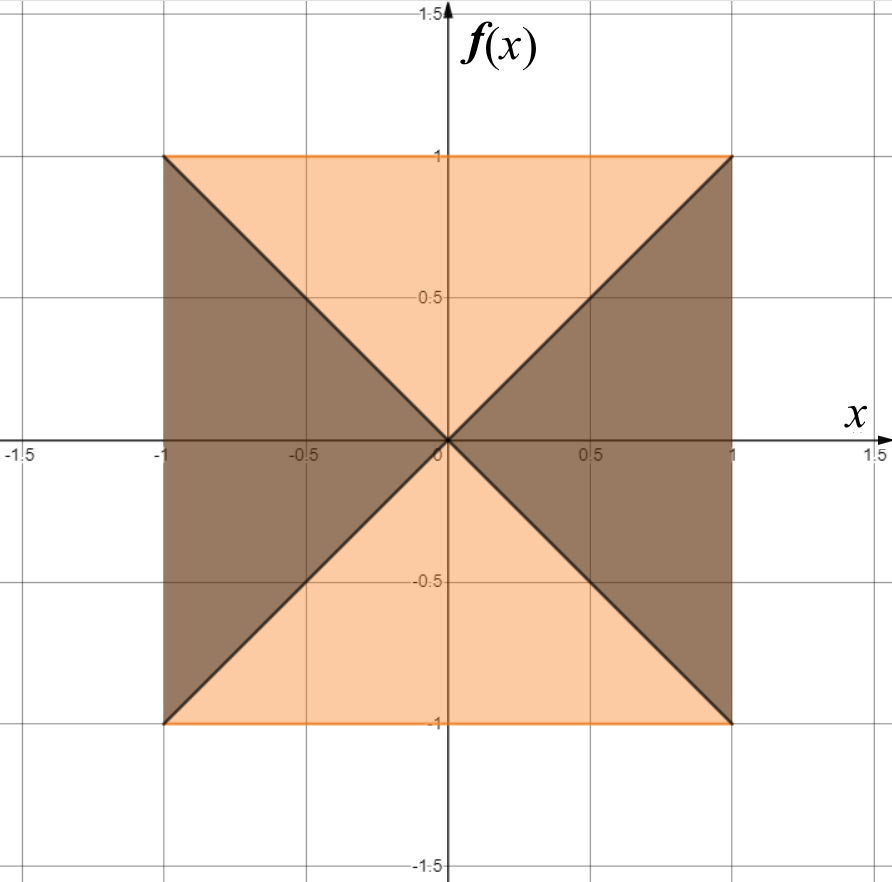}
    \hspace{0.04\linewidth}
    \includegraphics[width=0.47\linewidth]{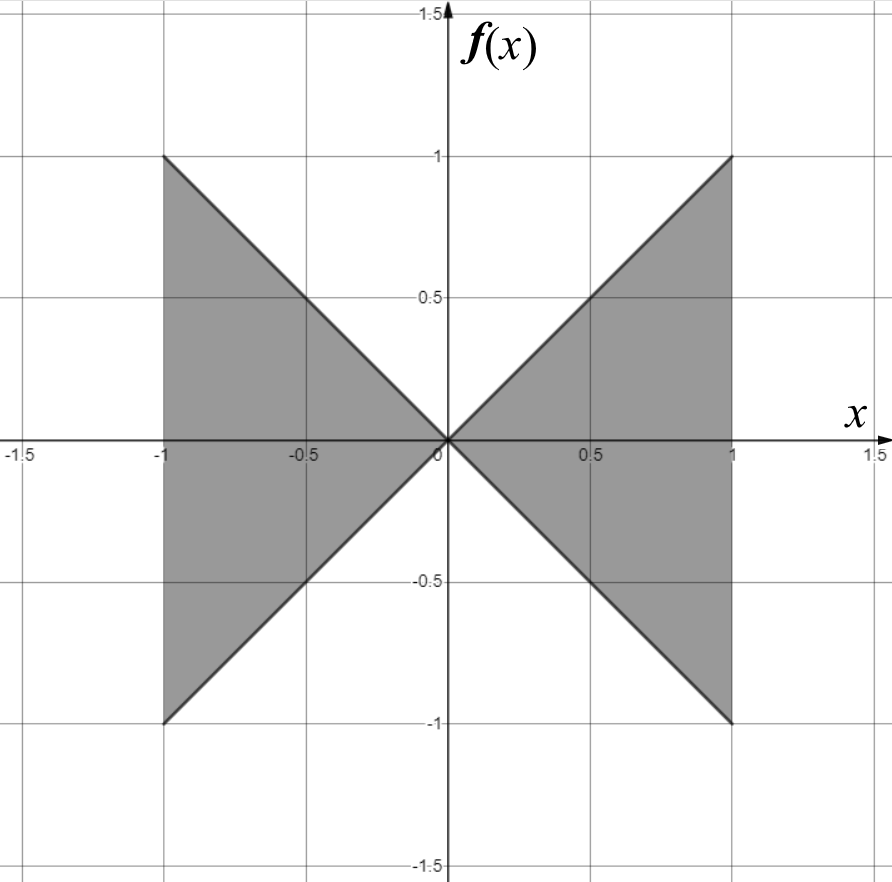}
    \caption{Множество совместных значений упорядоченной пары $\big( x, \: \mbf{f}(x) \big)$, где $\mbf{f}(x) = [ \, -1, \: 1 \, ] \cdot x$ (серая область). На левом рисунке показано его приближение в классической интервальной и аффинной арифметике (оранжевая область), на правом --- в функционально-граничной интервальной арифметике (совпадает с множеством совместных значений).}
	\label{fig:[-1,1]x_difference}
\end{figure}

\begin{figure}
	\centering
    \includegraphics[width=0.47\linewidth]{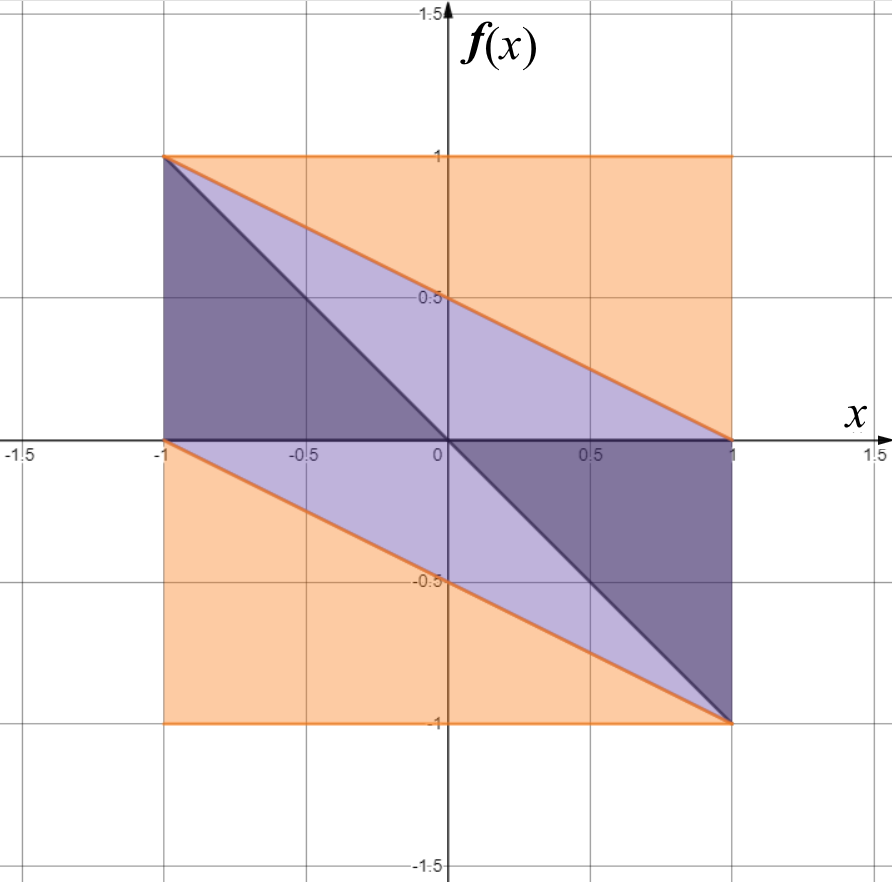}
    \hspace{0.04\linewidth}
    \includegraphics[width=0.47\linewidth]{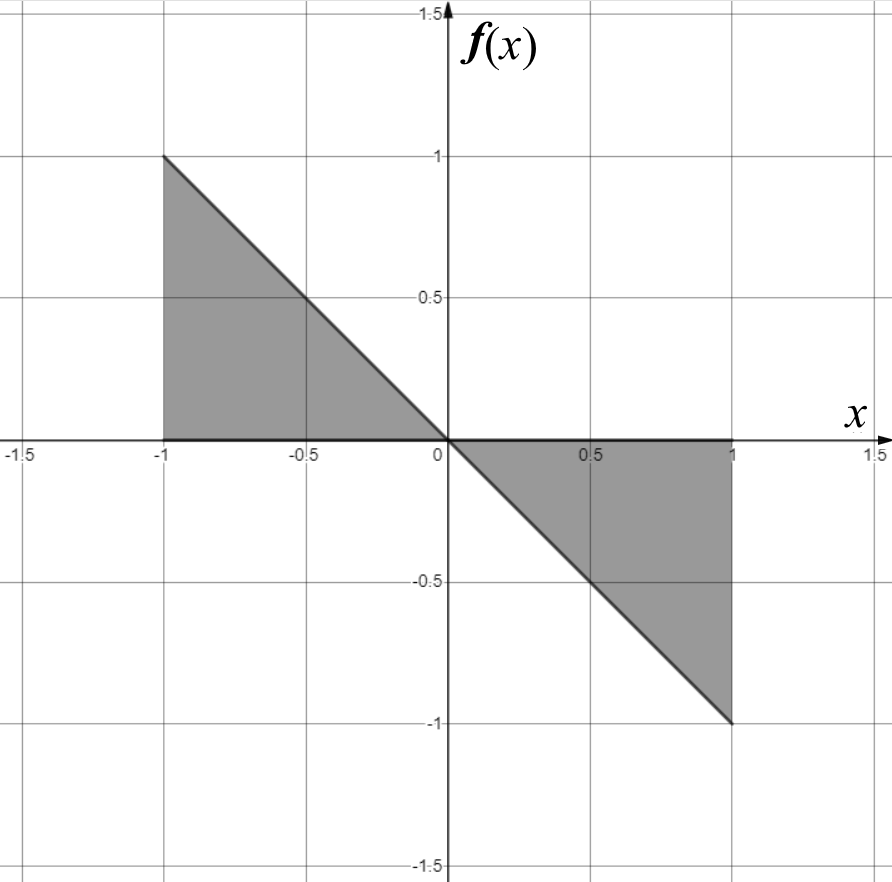}
    \caption{Множество совместных значений упорядоченной пары $\big( x, \mbf{f}(x) \big)$, где $\mbf{f}(x) = [ \, -1, \: 0 \, ] \cdot x$ (серая область). На левом рисунке она показана в классической интервальной арифметике (оранжевая область), а также в аффинной арифметике (фиолетовая область). На правом рисунке показана диаграмма зависимости в функционально-граничной интервальной арифметике (совпадает с множеством совместных значений).}
	\label{fig:[-1,0]x_difference}
\end{figure}

\begin{figure}
	\centering
    \includegraphics[width=0.47\linewidth]{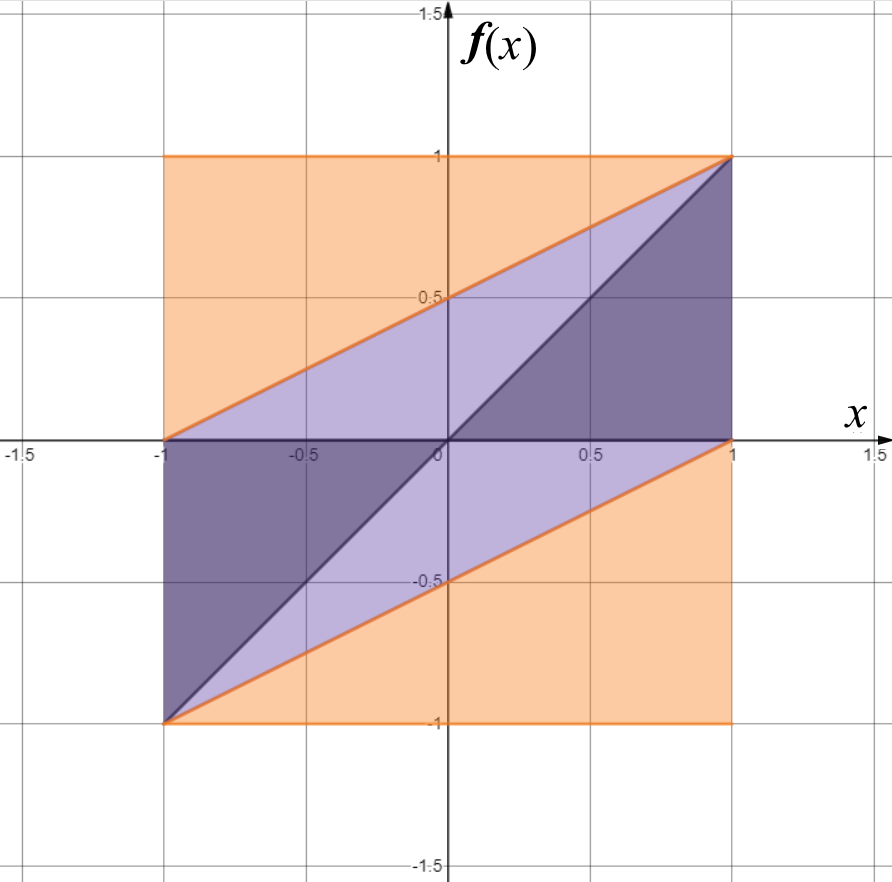}
    \hspace{0.04\linewidth}
    \includegraphics[width=0.47\linewidth]{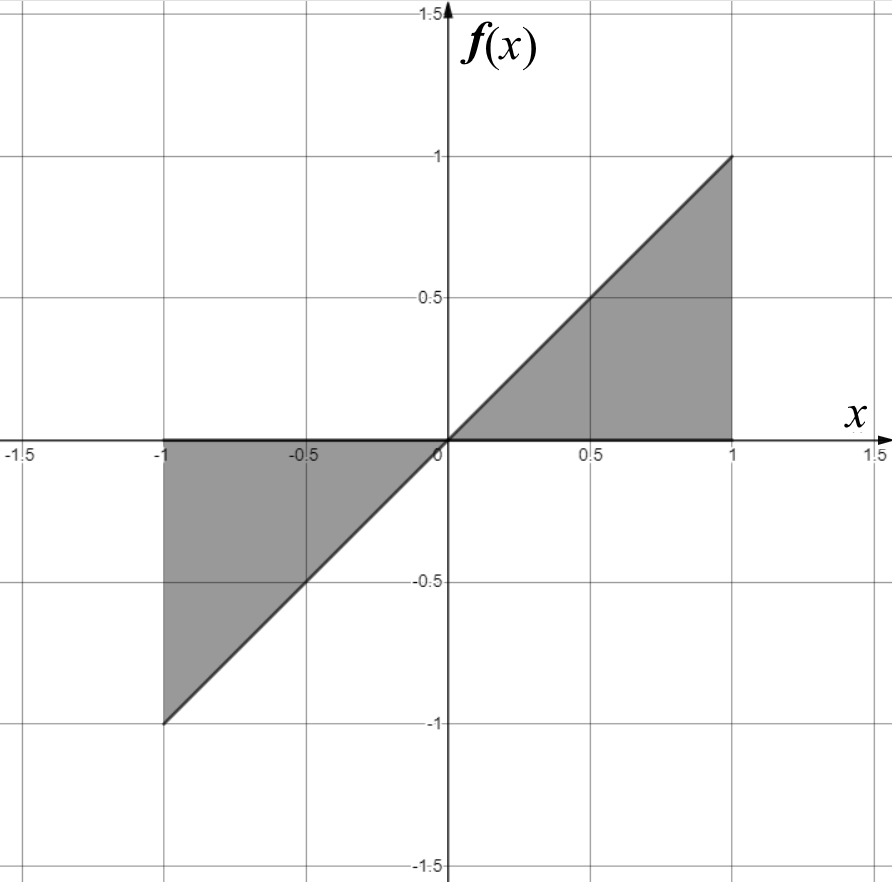}
    \caption{Множество совместных значений упорядоченной пары $\big( x, \: \mbf{f}(x) \big)$, где $\mbf{f}(x) = [ \, 0, \: 1 \,] \cdot x$ (серая область). На левом рисунке показано его приближение в классической интервальной арифметике (оранжевая область), а также в аффинной арифметике (фиолетовая область). На правом рисунке --- в функционально-граничной интервальной арифметике (совпадает с множеством совместных значений).}
	\label{fig:[0,1]x_difference}
\end{figure}

\begin{figure}
	\centering
    \includegraphics[width=0.47\linewidth]{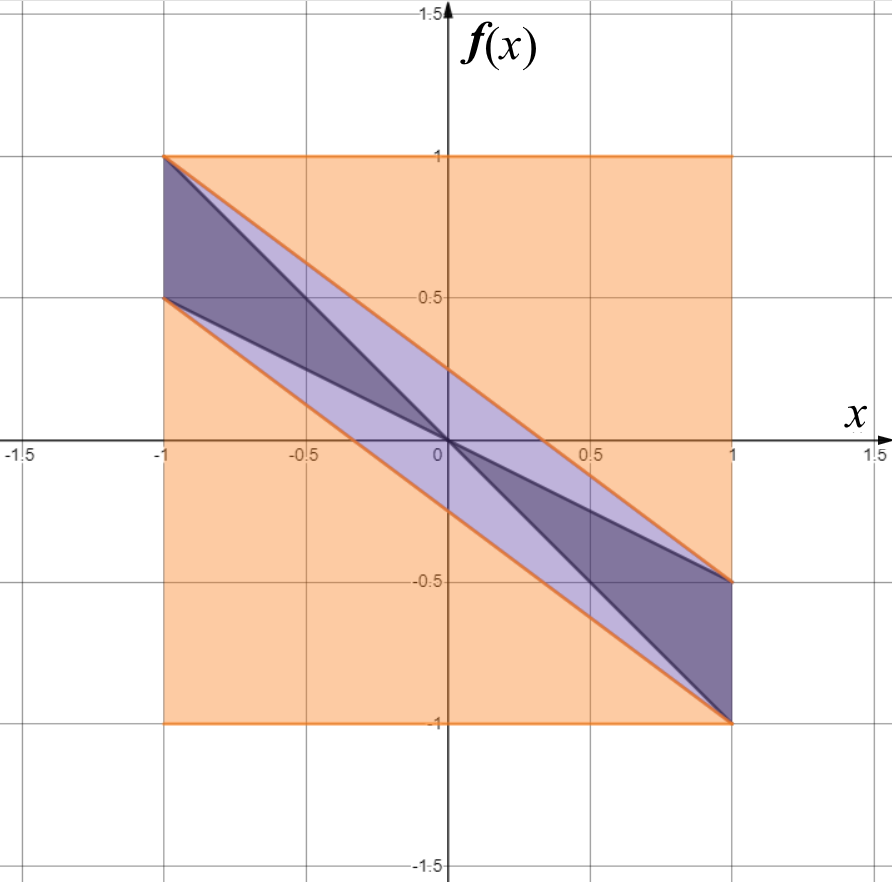}
    \hspace{0.04\linewidth}
    \includegraphics[width=0.47\linewidth]{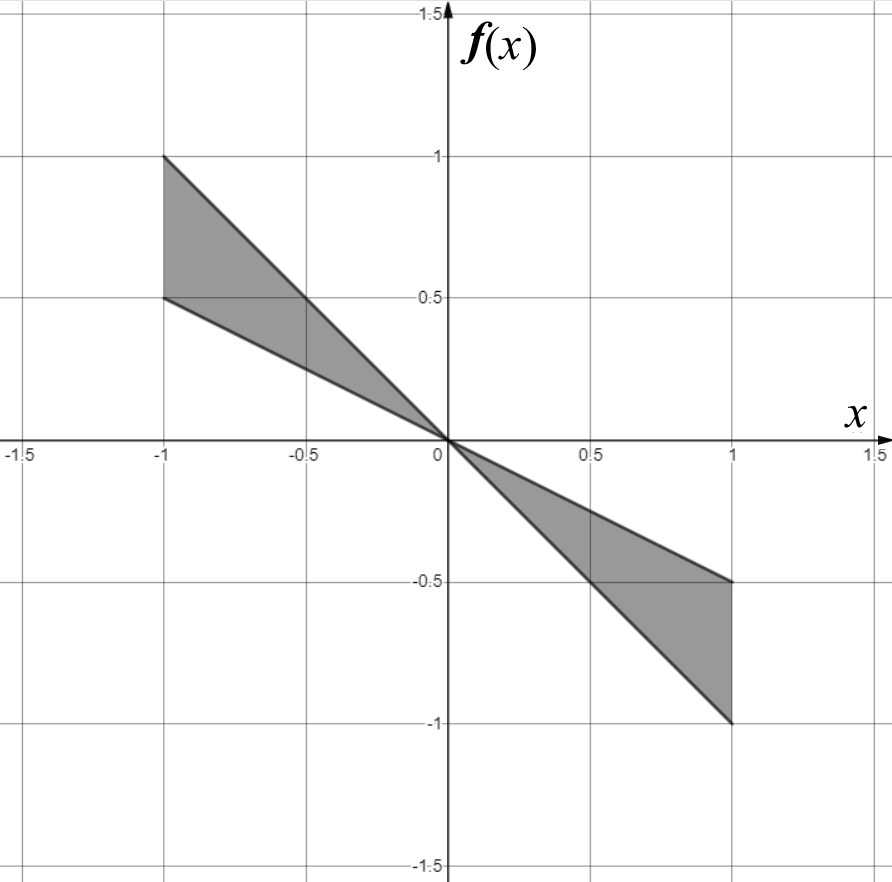}
    \caption{Множество совместных значений упорядоченной пары $\big( x, \: \mbf{f}(x) \big)$, где $\mbf{f}(x) = \big[ \, -1, \: -\frac{1}{2} \, \big] \cdot x$ (серая область). На левом рисунке показано его приближение в классической интервальной арифметике (оранжевая область), а также в аффинной арифметике (фиолетовая область). На правом рисунке показана диаграмма зависимости в функционально-граничной интервальной арифметике (совпадает с множеством совместных значений).}
	\label{fig:[-1,-0.5]x_difference}
\end{figure}

\begin{figure}
	\centering
    \includegraphics[width=0.47\linewidth]{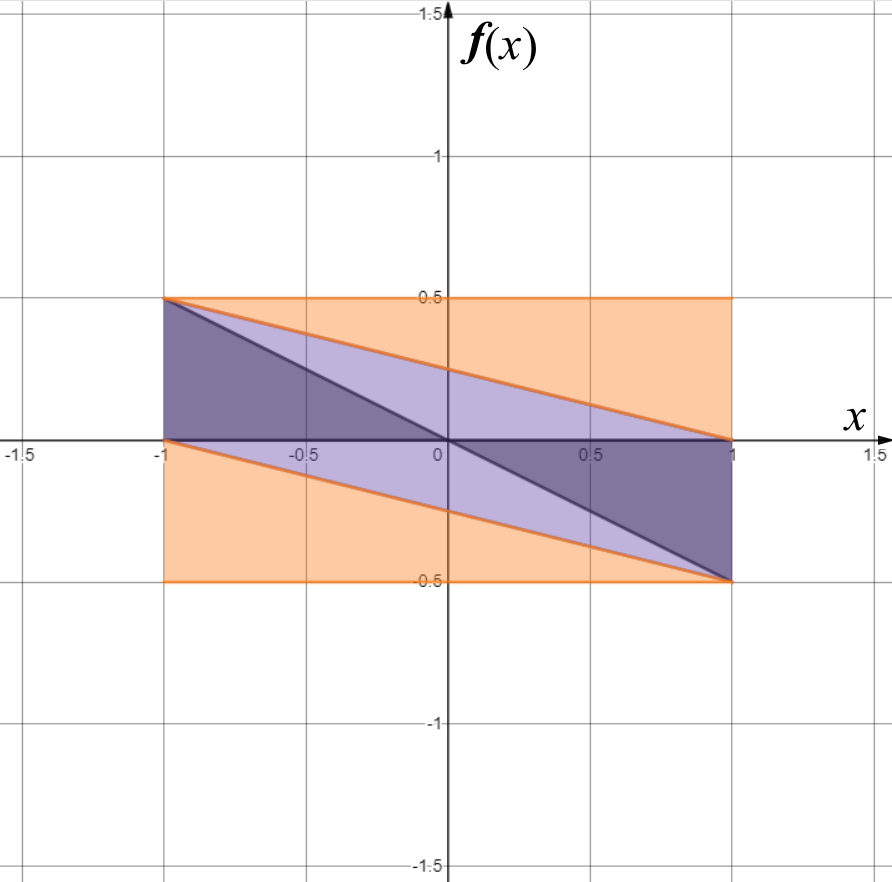}
    \hspace{0.04\linewidth}
    \includegraphics[width=0.47\linewidth]{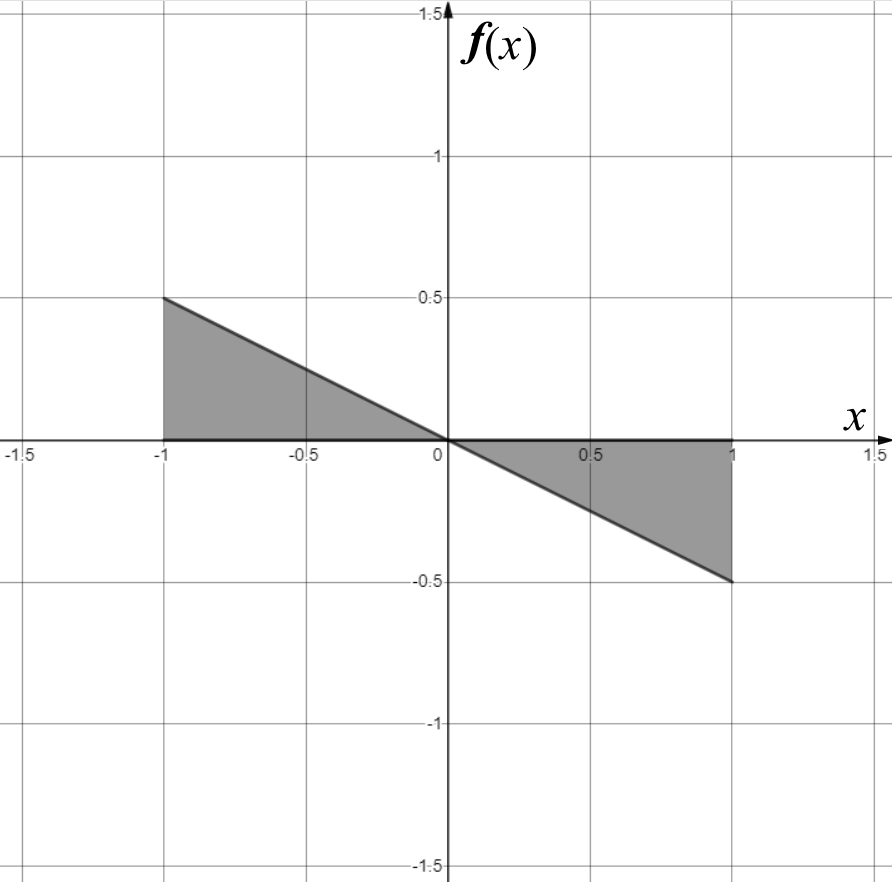}
    \caption{Множество совместных значений упорядоченной пары $\big( x, \: \mbf{f}(x) \big)$, где $\mbf{f}(x) = \big[ \, -\frac{1}{2}, \: 0 \, \big] \cdot x$ (серая область). На левом рисунке показано его приближение в классической интервальной арифметике (оранжевая область), а также в аффинной арифметике (фиолетовая область). На правом рисунке --- в функционально-граничной интервальной арифметике (совпадает с множеством совместных значений).}
	\label{fig:[-0.5,0]x_difference}
\end{figure}

\begin{figure}
	\centering
    \includegraphics[width=0.47\linewidth]{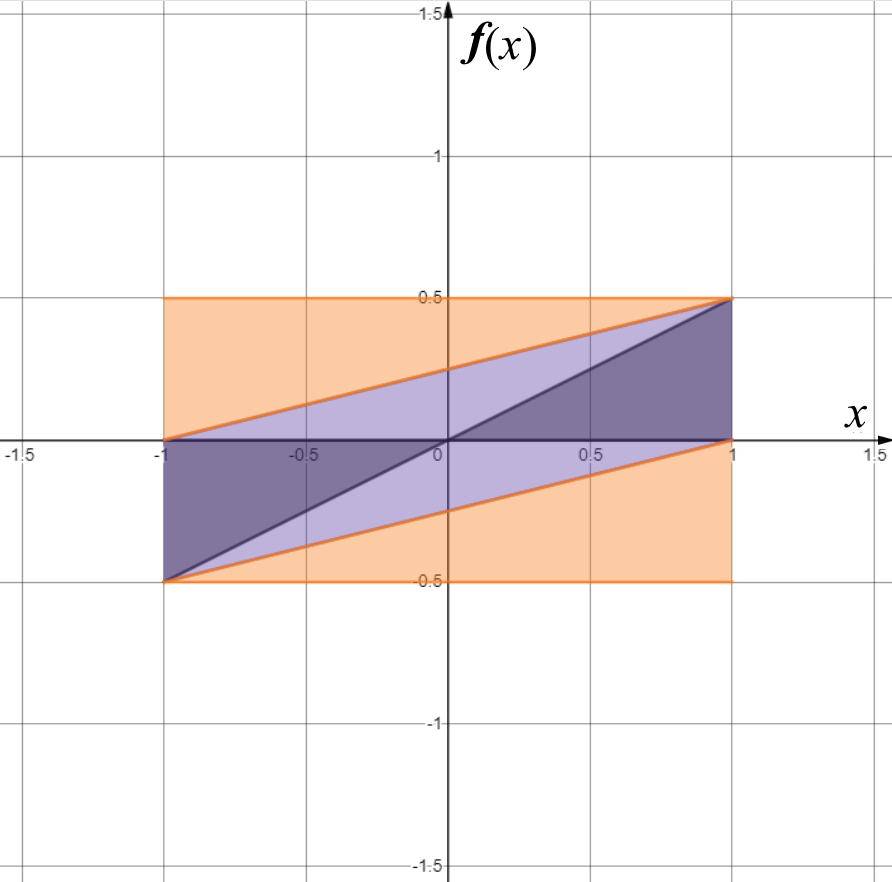}
    \hspace{0.04\linewidth}
    \includegraphics[width=0.47\linewidth]{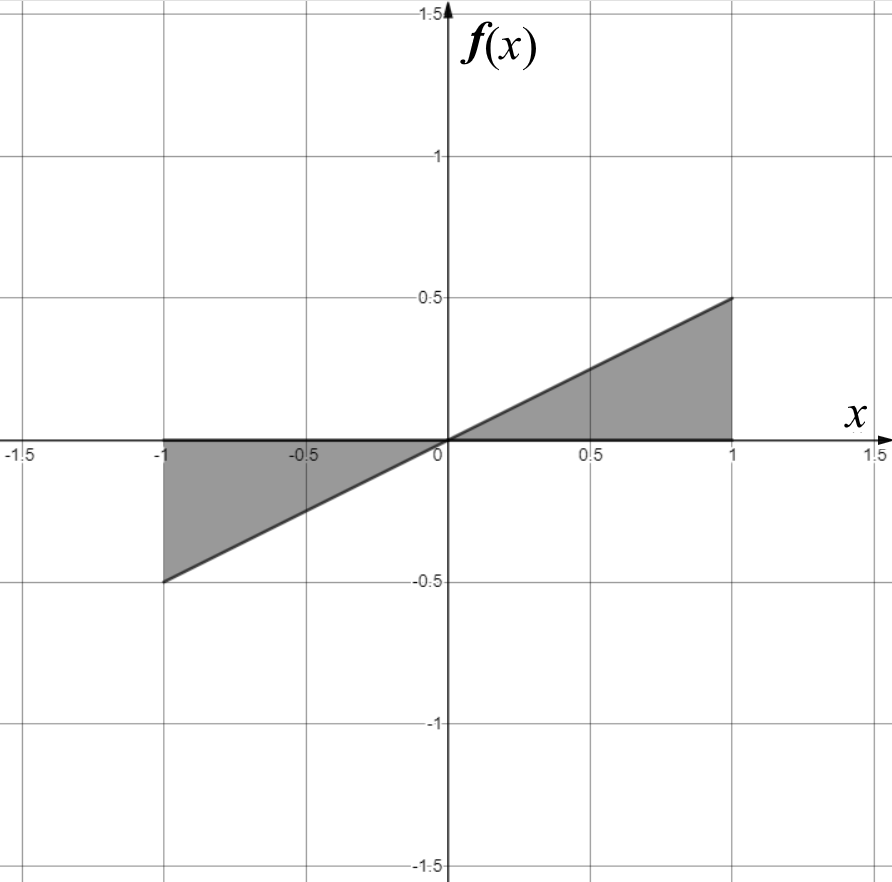}
    \caption{Множество совместных значений упорядоченной пары $\big( x, \: \mbf{f}(x) \big)$, где $\mbf{f}(x) = \big[ \, 0, \: \frac{1}{2} \, \big] \cdot x$ (серая область). На левом рисунке показано его приближение в классической интервальной арифметике (оранжевая область), а также в аффинной арифметике (фиолетовая область). На правом рисунке --- в функционально-граничной интервальной арифметике (совпадает с множеством совместных значений).}
	\label{fig:[0,0.5]x_difference}
\end{figure}

\begin{figure}
	\centering
    \includegraphics[width=0.47\linewidth]{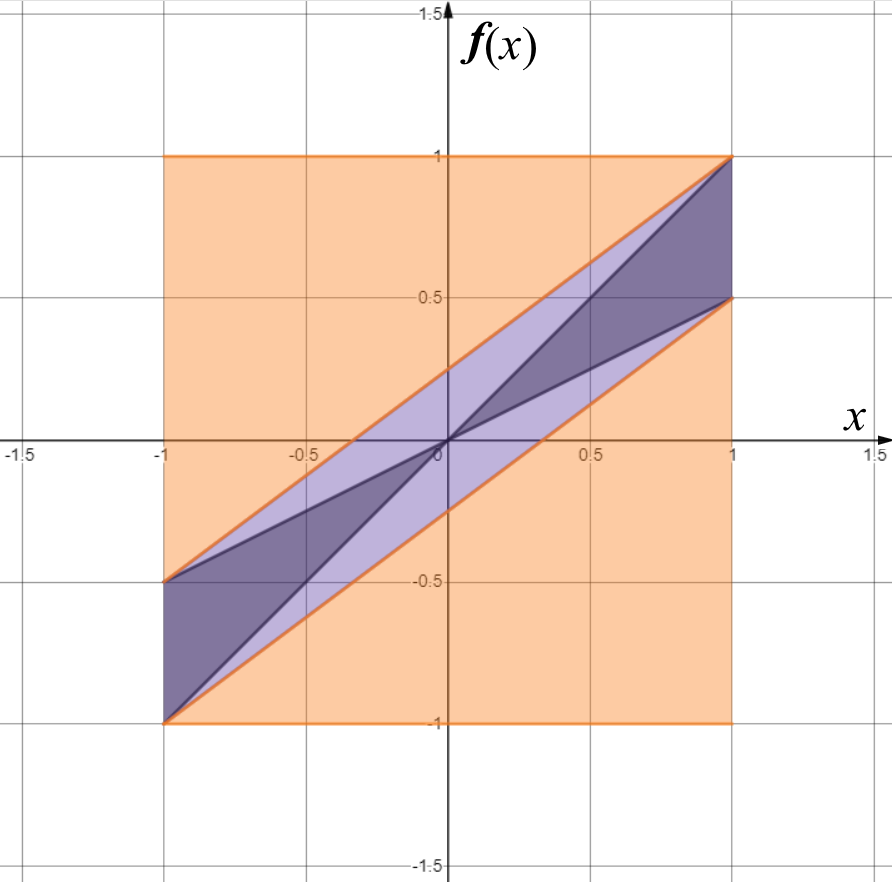}
    \hspace{0.04\linewidth}
    \includegraphics[width=0.47\linewidth]{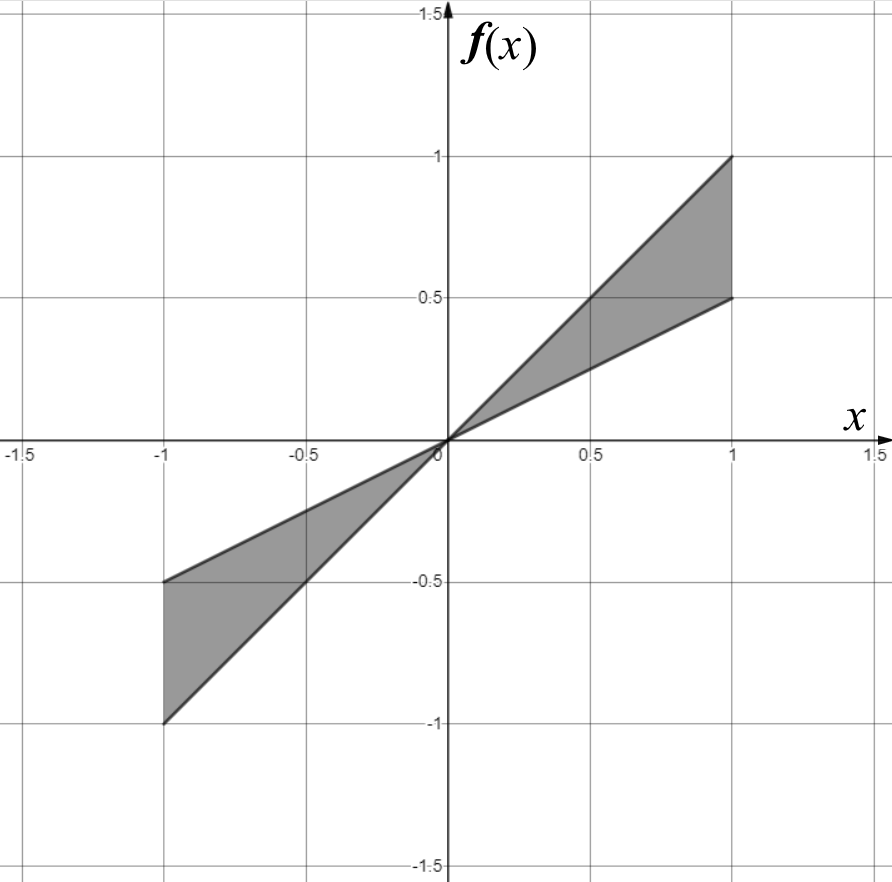}
    \caption{Множество совместных значений упорядоченной пары $\big( x, \: \mbf{f}(x) \big)$, где $f(x) = \big[ \, \frac{1}{2}, \: 1 \, \big] \cdot x$ (серая область). На левом рисунке показано его приближение в классической интервальной арифметике (оранжевая область), а также в аффинной арифметике (фиолетовая область). На правом рисунке --- в функционально-граничной интервальной арифметике (совпадает с множеством совместных значений).}
	\label{fig:[0.5,1]x_difference}
\end{figure}

\begin{figure}
	\centering
    \includegraphics[width=0.47\linewidth]{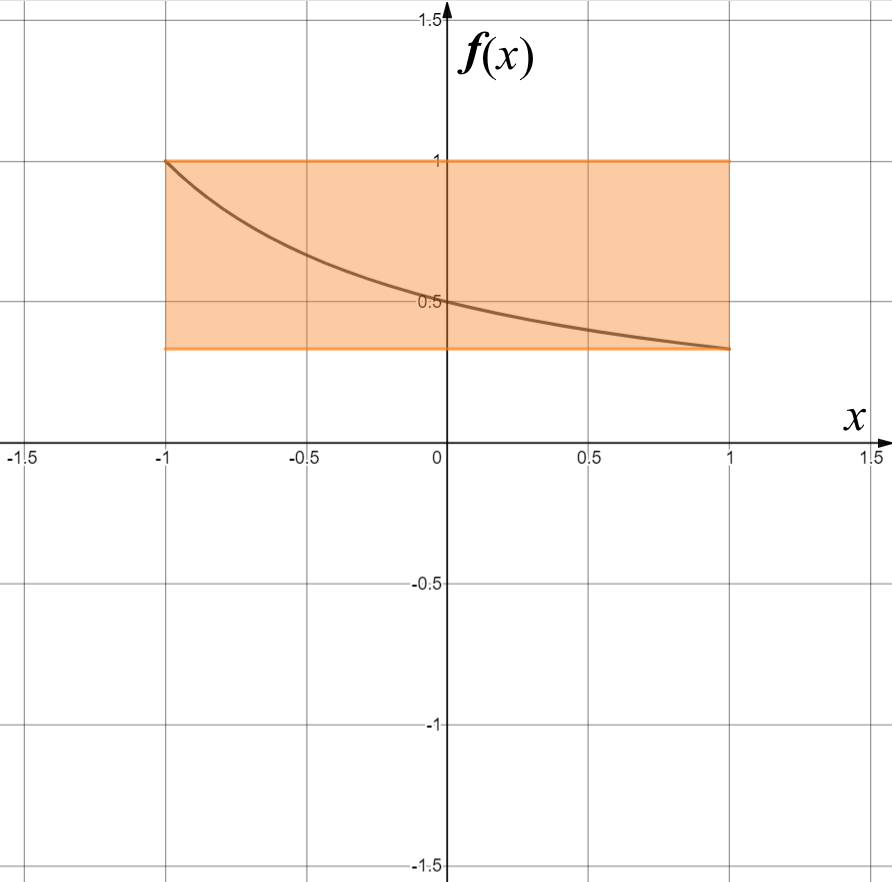}
    \hspace{0.04\linewidth}
    \includegraphics[width=0.47\linewidth]{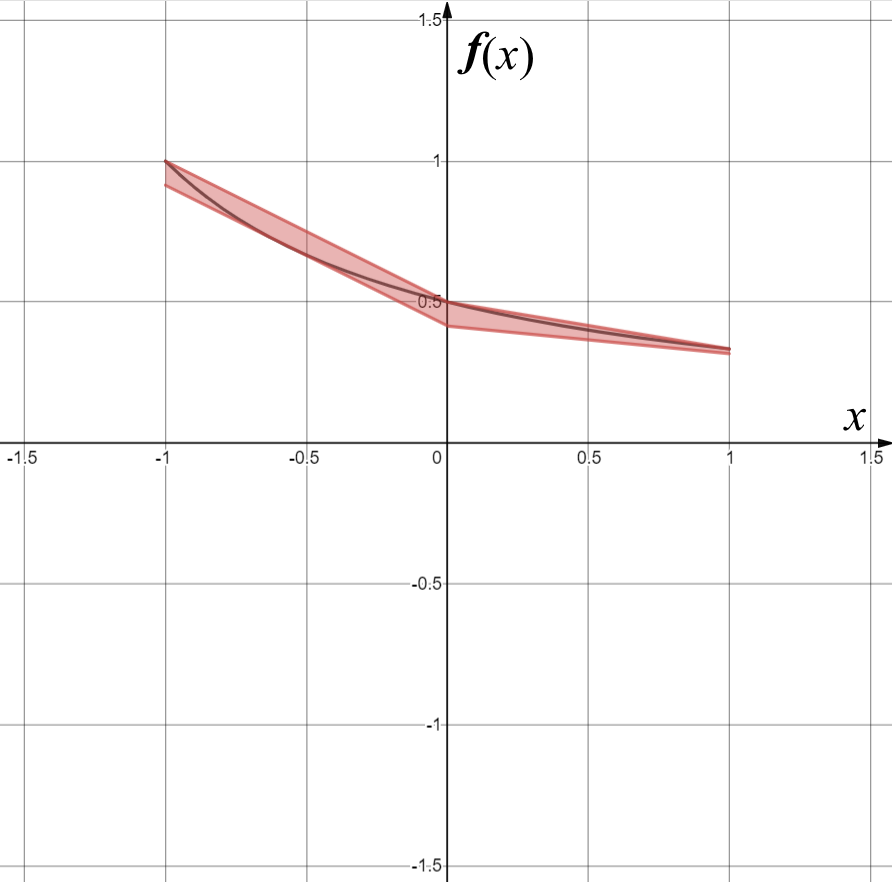}
    \caption{Множество совместных значений упорядоченной пары $\big( x, \: f(x) \big)$, задаваемой выражением $f(x) = \frac{1}{x \, + \, 2}$ (чёрная линия). На левом рисунке показано его приближение в классической интервальной арифметике (оранжевая область). На правом рисунке --- в функционально-граничной интервальной арифметике (красная область).}
	\label{fig:1div(x+2)_difference}
\end{figure}

\begin{figure}
	\centering
    \includegraphics[width=0.47\linewidth]{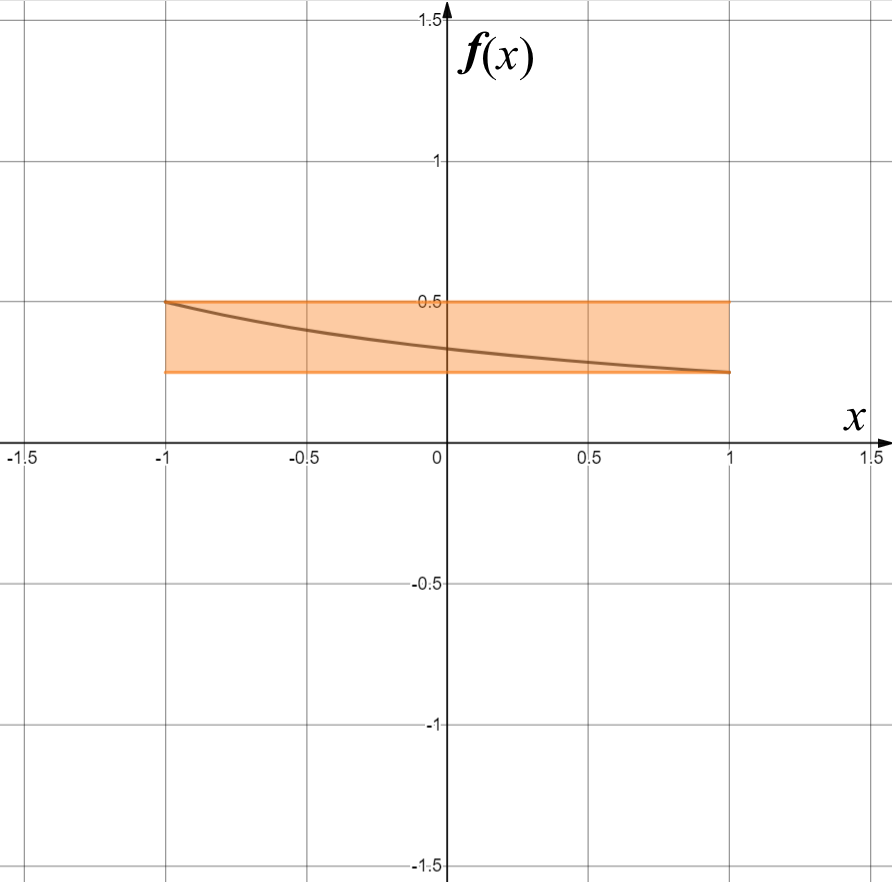}
    \hspace{0.04\linewidth}
    \includegraphics[width=0.47\linewidth]{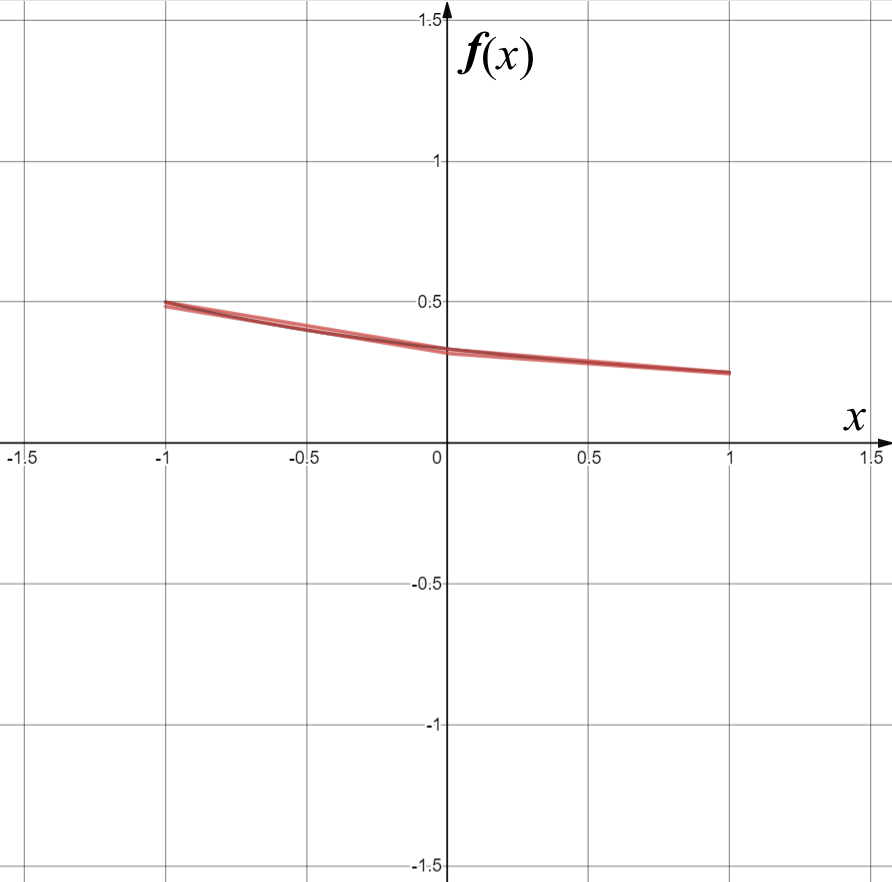}
    \caption{Множество совместных значений упорядоченной пары $\big( x, \: f(x) \big)$, где $f(x) = \frac{1}{x \, + \, 3}$ (чёрная линия). На левом рисунке показано его приближение в классической интервальной арифметике (оранжевая область). На правом рисунке --- в функционально-граничной интервальной арифметике (красная область).}
	\label{fig:1div(x+3)_difference}
\end{figure}

\begin{figure}
	\centering
    \includegraphics[width=0.47\linewidth]{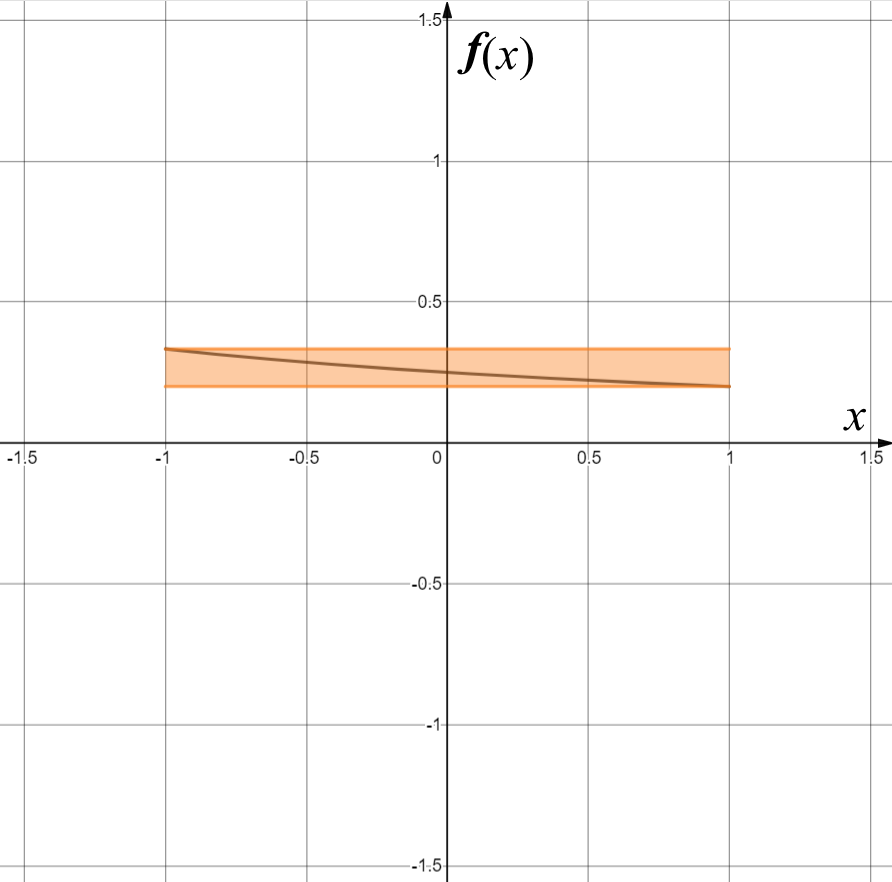}
    \hspace{0.04\linewidth}
    \includegraphics[width=0.47\linewidth]{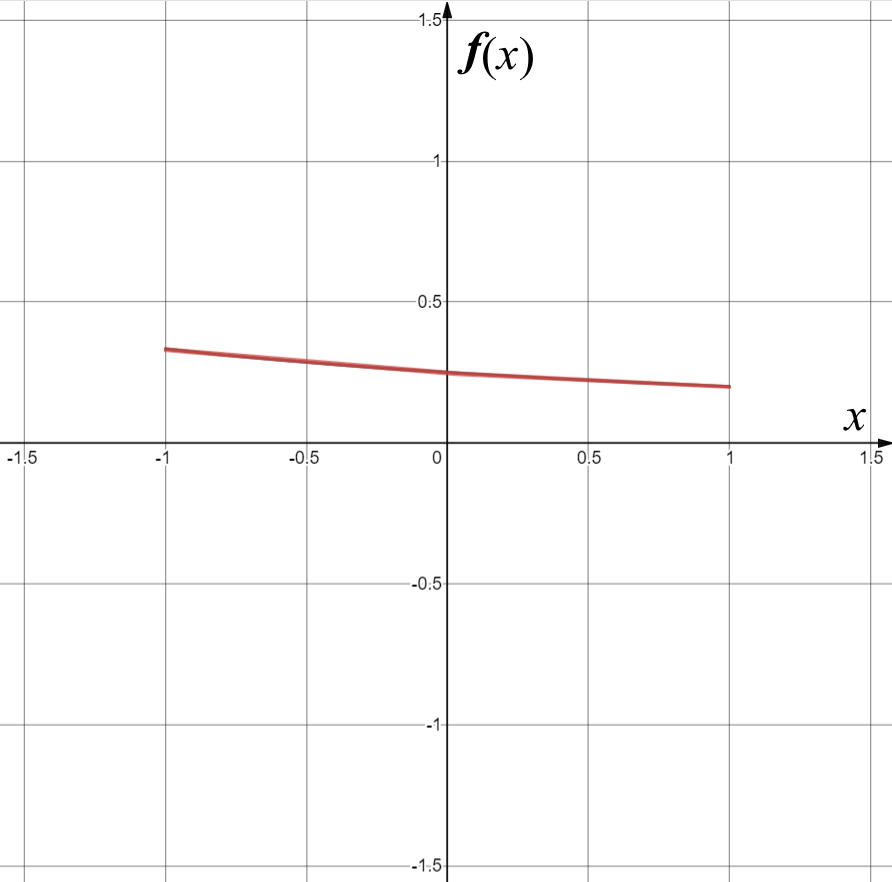}
    \caption{Множество совместных значений упорядоченной пары $\big( x, \: f(x) \big)$, где $f(x) = \frac{1}{x \, + \, 4}$ (чёрная линия). На левом рисунке показано его приближение в классической интервальной арифметике (оранжевая область). На правом рисунке --- в функционально-граничной интервальной арифметике (красная область).}
	\label{fig:1div(x+4)_difference}
\end{figure}

\begin{figure}
	\centering
    \includegraphics[width=0.47\linewidth]{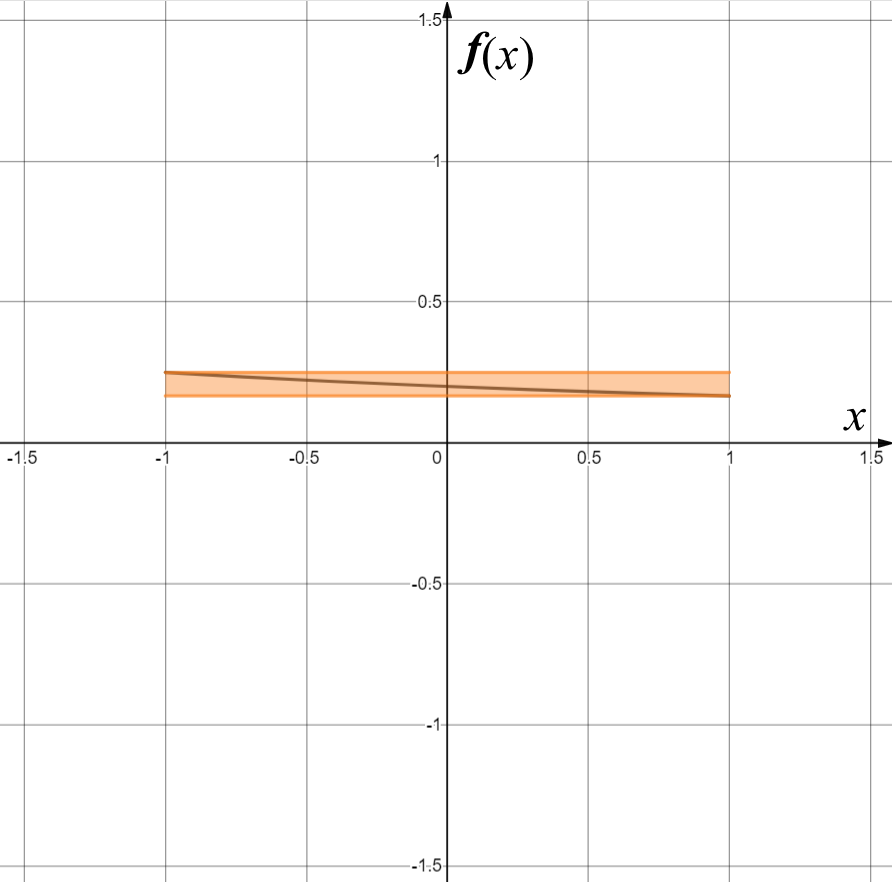}
    \hspace{0.04\linewidth}
    \includegraphics[width=0.47\linewidth]{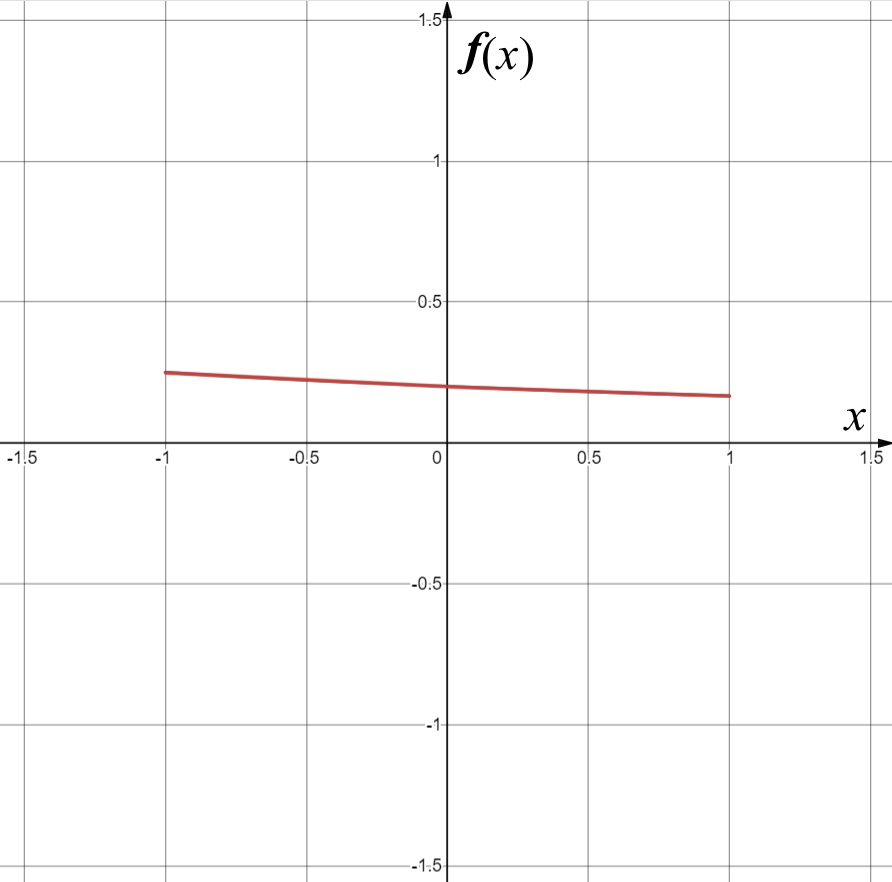}
    \caption{Множество совместных значений упорядоченной пары $\big( x, \: f(x) \big)$, где $f(x) = \frac{1}{x \, + \, 5}$ (чёрная линия). На левом рисунке показано его приближение в классической интервальной арифметике (оранжевая область). На правом рисунке --- в функционально-граничной интервальной арифметике (красная область).}
	\label{fig:1div(x+5)_difference}
\end{figure}

\begin{figure}
	\centering
    \includegraphics[width=0.47\linewidth]{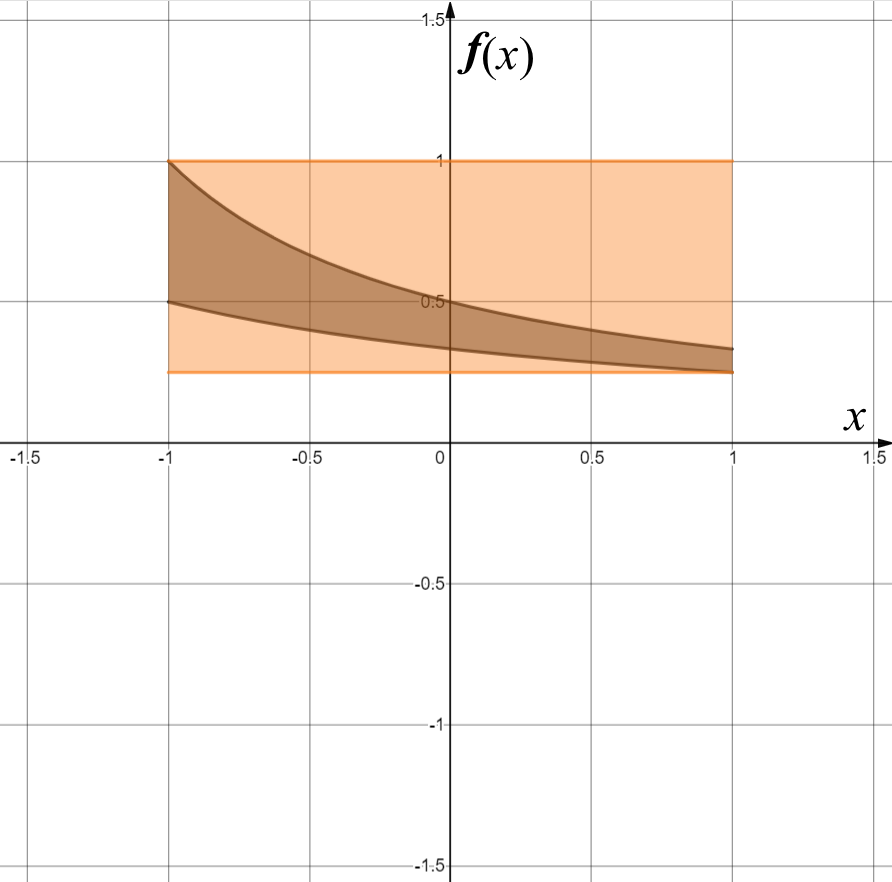}
    \hspace{0.04\linewidth}
    \includegraphics[width=0.47\linewidth]{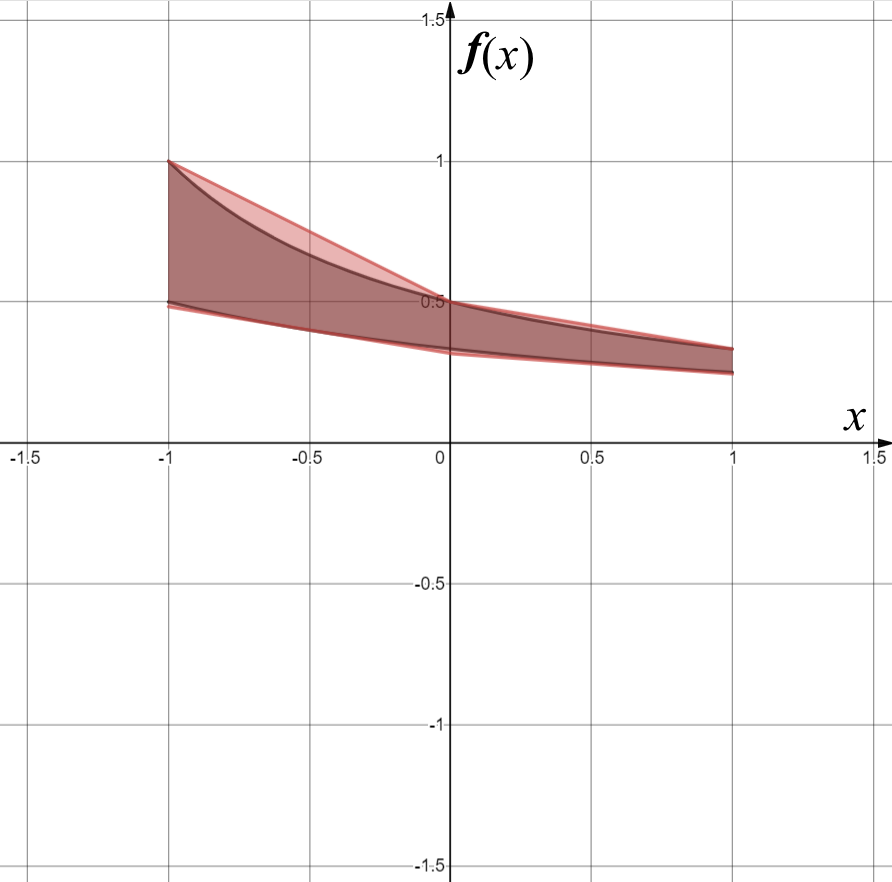}
    \caption{Множество совместных значений упорядоченной пары $\big( x, \: \mbf{f}(x) \big)$, где $\mbf{f}(x) = \frac{1}{x \, + \, [ \, 2, \: 3 \, ]}$ (серая область). На левом рисунке показано его приближение в классической интервальной арифметике (оранжевая область). На правом рисунке --- в функционально-граничной интервальной арифметике (красная область).}
	\label{fig:1div(x+[2,3])_difference}
\end{figure}

\begin{figure}
	\centering
    \includegraphics[width=0.47\linewidth]{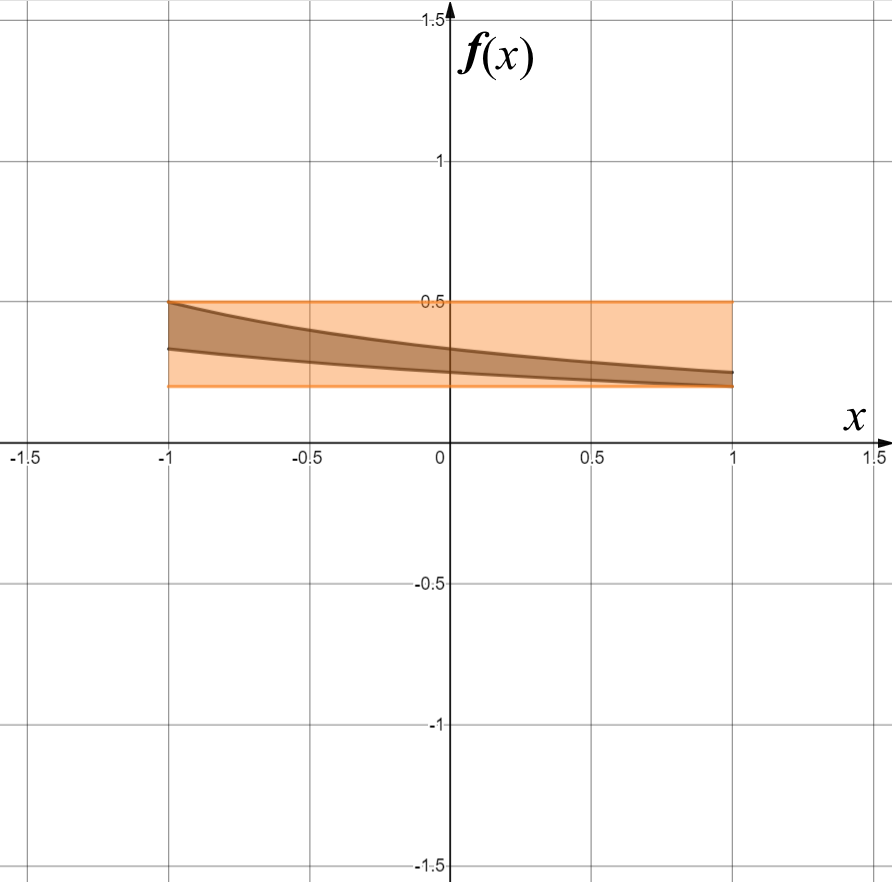}
    \hspace{0.04\linewidth}
    \includegraphics[width=0.47\linewidth]{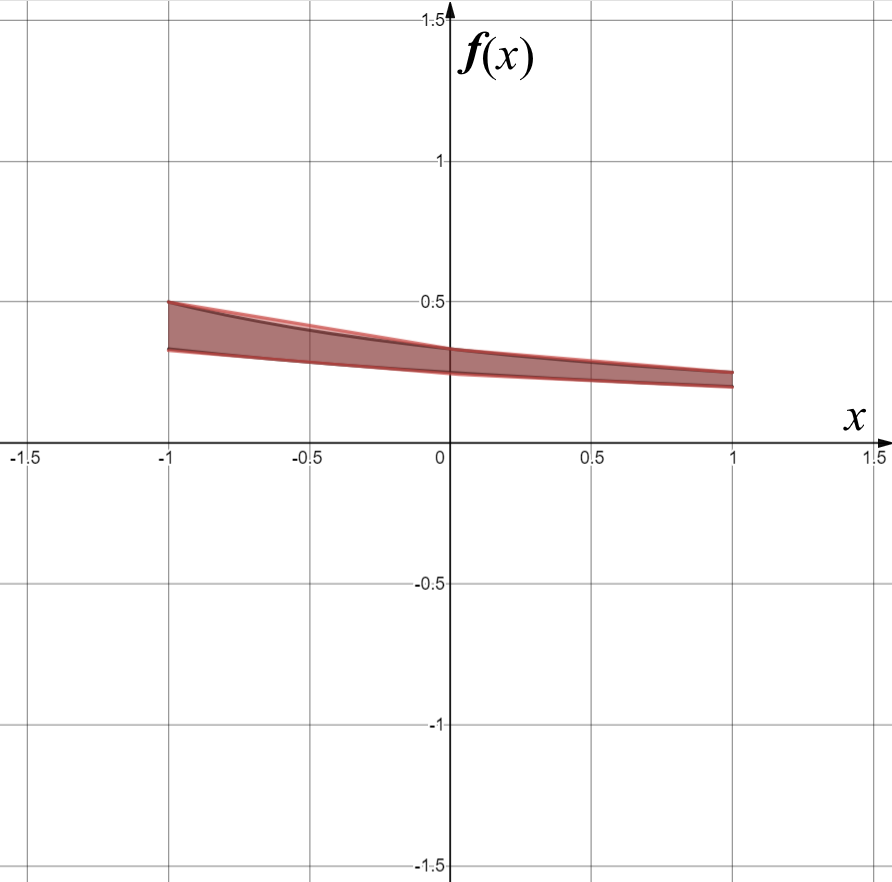}
    \caption{Множество совместных значений упорядоченной пары $\big( x, \: \mbf{f}(x) \big)$, где $\mbf{f}(x) = \frac{1}{x \, + \, [ \, 3, \: 4 \, ]}$ (серая область). На левом рисунке показано его приближение в классической интервальной арифметике (оранжевая область). На правом рисунке --- в функционально-граничной интервальной арифметике (красная область).}
	\label{fig:1div(x+[3,4])_difference}
\end{figure}

\begin{figure}
	\centering
    \includegraphics[width=0.47\linewidth]{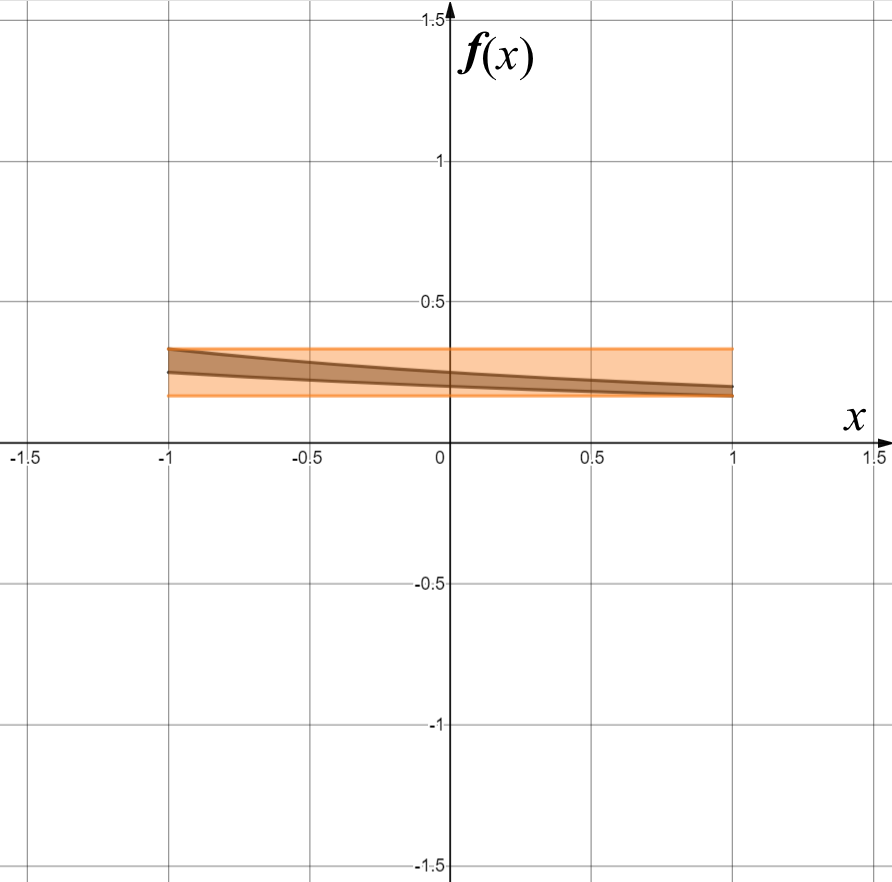}
    \hspace{0.04\linewidth}
    \includegraphics[width=0.47\linewidth]{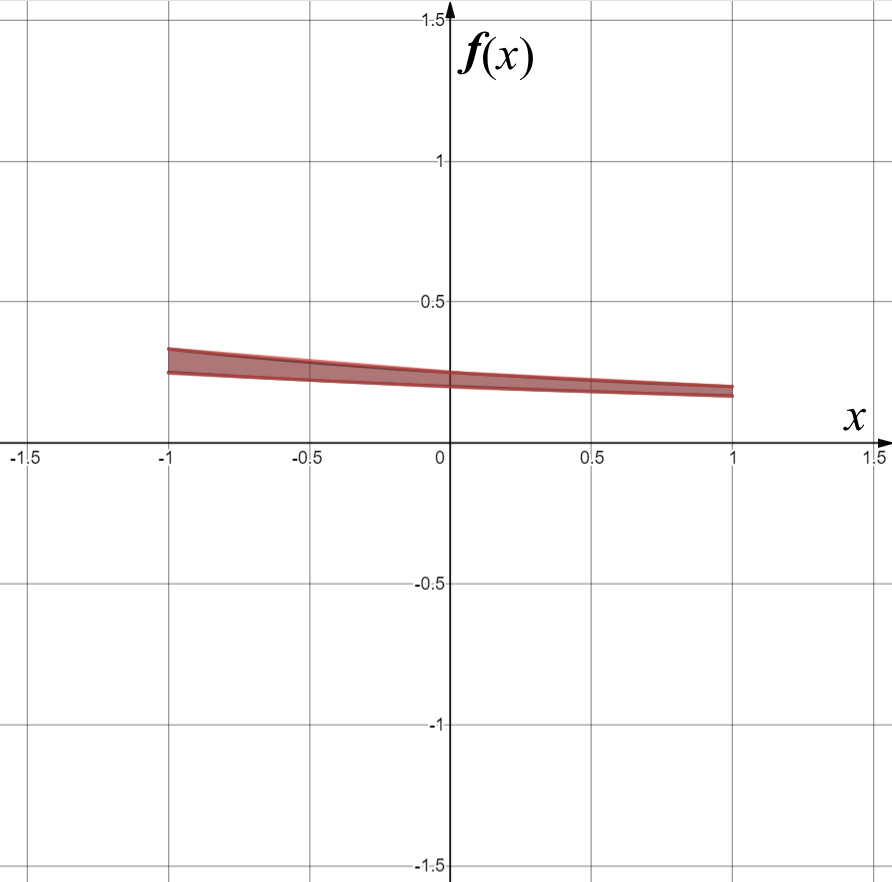}
    \caption{Множество совместных значений упорядоченной пары $\big( x, \: \mbf{f}(x) \big)$, где $\mbf{f}(x) = \frac{1}{x \, + \, [ \, 4, \: 5 \, ]}$ (серая область). На левом рисунке показано его приближение в классической интервальной арифметике (оранжевая область). На правом рисунке --- в функционально-граничной интервальной арифметике (красная область).}
	\label{fig:1div(x+[4,5])_difference}
\end{figure}

\begin{figure}
	\centering
    \includegraphics[width=0.47\linewidth]{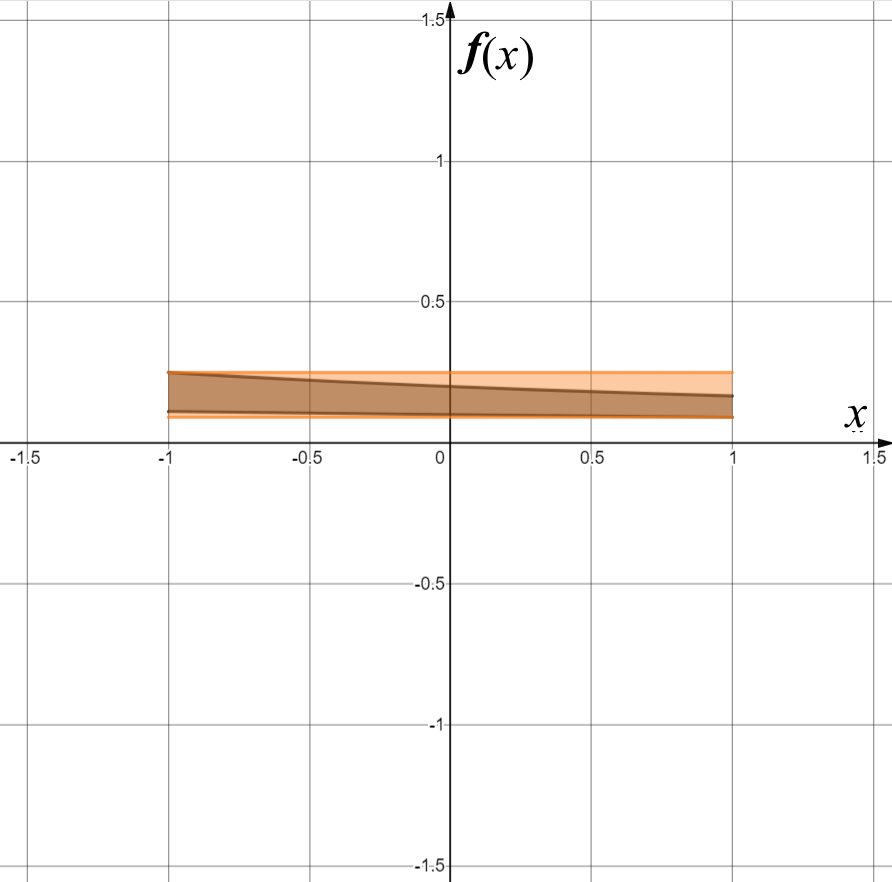}
    \hspace{0.04\linewidth}
    \includegraphics[width=0.47\linewidth]{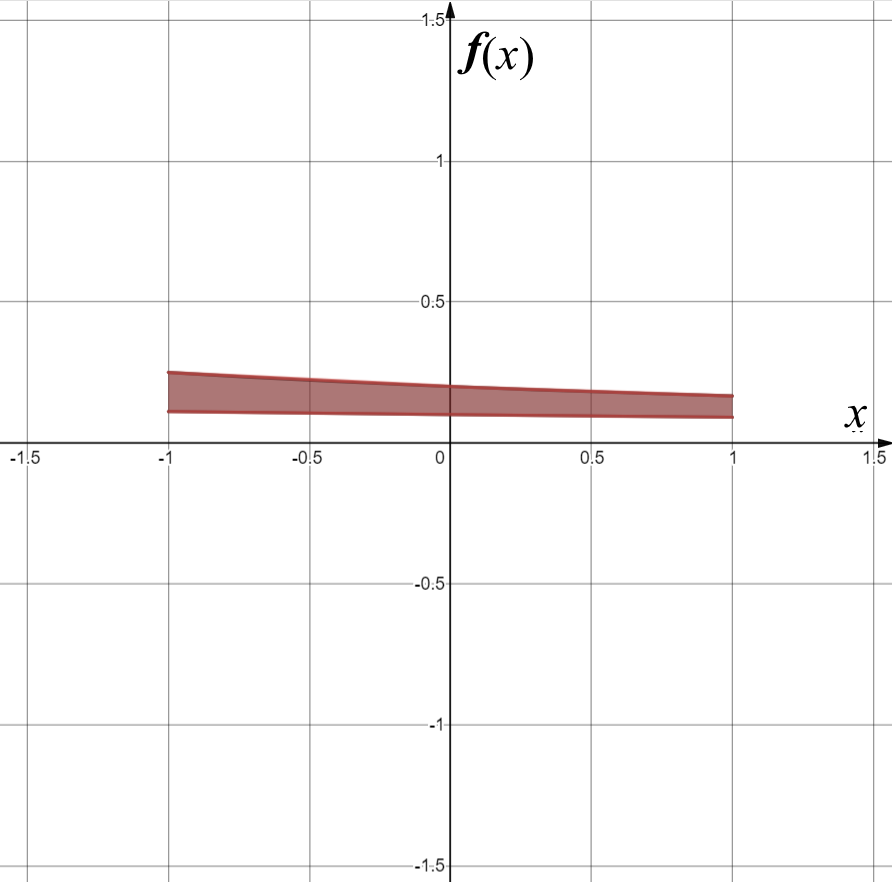}
    \caption{Множество совместных значений упорядоченной пары $\big( x, \: \mbf{f}(x) \big)$, где $\mbf{f}(x) = \frac{1}{x \, + \, [ \, 5, \: 10 \, ]}$ (серая область). На левом рисунке показано его приближение в классической интервальной арифметике (оранжевая область). На правом рисунке --- в функционально-граничной интервальной арифметике (красная область).}
	\label{fig:1div(x+[5,10])_difference}
\end{figure}

\clearpage

\subsection{Решение ИСЛАУ со связями с применением \\ функционально-граничной арифметики}

Рассмотрим ИСЛАУ, которая рассматривалась в работе Ахмерова \cite{AchmerovZonotopesQuestion}. Она решалась методом Гаусса при помощи интервально-аффинной арифметики. Асимптотика трудоёмкости данного метода равняется $O(n^{5})$, где $n$ --- размер матрицы системы.

\begin{center}

	$\left ( \begin{array}{ccc}
	
		[ \, 0.7, \: 1.3 \,] & [ \, -0.3, \: 0.3] & [ \, -0.3, \: 0.3 \, ] \\
			
		[ \, -0.3, \: 0.3 \, ] & [ \, 0.7, \: 1.3 \, ] & [ \, -0.3, \: 0.3 \, ] \\
			
		[ \, -0.3, \: 0.3 \, ] & [ \, -0.3, \: 0.3 \, ] & [ \, 0.7, \: 1.3 \, ]	
	
	\end{array} \right ) \cdot \left ( \begin{array}{c}
	
		x_{1} \\
		
		x_{2} \\
		
		x_{3}	
	
	\end{array} \right )
	=
	\left ( \begin{array}{c}
	
		[ \, -14, \: -7 \, ] \\	
		
		[ \, 9, \: 12 \, ] \\
		
		[ \, -3, \: 3 \, ]
	
	\end{array} \right )$

\end{center}

Далее приведены таблицы результатов решений в разных арифметиках. Отметим, что трудоёмкость метода Гаусса для решения ИСЛАУ в функционально-граничной интервальной арифметике равна $O(n^{7})$. Рассмотрим случай симметричной матрицы (таб. \ref{tab:slau_sym}), когда:
\begin{center}

	$\mbf{a}_{i j} = \mbf{a}_{j i},$
	
\end{center}
а также случай кососимметричной матрицы (таб. \ref{tab:not_sym_slau}), когда:
\begin{center}

	$\mbf{a}_{i j} = - \mbf{a}_{j i}$.
	
\end{center}

Проанализируем полученные результаты. Классическая интервальная арифметика даёт самые грубые оценки среди приведённых, поскольку не способна учитывать связи между элементами ИСЛАУ. Использование аффинной арифметики и функционально-граничной интервальной арифметики за счёт наличия такой специфики позволяет получить для оценок более узкие интервалы. 

А так как функционально-граничная интервальная арифметика позволяет эффективнее обрабатывать квадратичные члены и описывает границы интервалов независимо, то она показывает лучшую оценку бруса неизвестных среди приведённых за счёт увеличения трудоёмкости алгоритма.

\begin{table}[ht]
	\caption{Оценки объединённого множества решений ИСЛАУ для случая симметричной матрицы.}
	\label{tab:slau_sym}
	\vspace{2mm}
	\begin{tabular}{r||c||c||c}

		& Классическая & Интервально- & Функционально-граничная \\
		
		$\mbf{x}_{i}$ & интервальная & аффиная & интервальная \\
		
		& арифметика & арифметика & арифметика \\ \hline \hline \rule[-1mm]{0mm}{6mm}
		
		$\mbf{x}_{1}$ & $[ \, -101.00, \: 71.00 \, ]$ & $[ \, -101.00, \: 64.80 \, ]$ & $[ \, -101.00, \: 54.21 \, ]$ \\
		
		$\mbf{x}_{2}$ & $[ \, -62.25, \: 99.00 \, ]$ & $[ \, -56.06, \: 99.00 \, ]$ & $[ \, -41.38, \: 99.00 \, ]$ \\
		
		$\mbf{x}_{3}$ & $[ \, -90.00, \: 90.00 \, ]$ & $[ \, -90.00, \: 90.00 \, ]$ & $[ \, -90.00, \: 90.00 \, ]$		

	\end{tabular}

\end{table}

\begin{table}[ht]
	\caption{Оценки объединённого множества решений ИСЛАУ для случая кососимметричной матрицы.}
	\label{tab:not_sym_slau}
	\vspace{2mm}
	\begin{tabular}{r||c||c||c}

		& Классическая & Интервально- & Функционально-граничная \\
		
		$\mbf{x}_{i}$ & интервальная & аффиная & интервальная \\
		
		& арифметика & арифметика & арифметика \\ \hline \hline \rule[-1mm]{0mm}{6mm}
		
		$\mbf{x}_{1}$ & $[ \, -101.00, \: 71.00 \, ]$ & $[ \, -46.58, \: 21.44 \, ]$ & $[ \, -35.71, \: 17.50 \, ]$ \\
		
		$\mbf{x}_{2}$ & $[ \, -62.25, \: 99.00 \, ]$ & $[ \, -14.98, \: 42.03 \, ]$ & $[ \, -16.07, \: 29.37 \, ]$ \\
		
		$\mbf{x}_{3}$ & $[ \, -90.00, \: 90.00 \, ]$ & $[ \, -31.33, \: 31.33 \, ]$ & $[ \, 	-23.13, \: 23.13 \, ]$		

	\end{tabular}

\end{table}

\clearpage
\section[Детали компьютерной реализации]{Детали компьютерной реализации}

Важным вопросом являются детали реализации интервальных арифметик на ЭВМ. Особенностью вычислений с помощью чисел с плавающей точкой является то, что промежуточные и итоговые результаты округляются из-за невозможности точного представления любого вещественного числа непрерывной числовой оси в дискретной памяти машины.

Первая специфика состоит в том, что нам необходима возможность читать и изменять регистр процессора для округления чисел с плавающей точкой. Эта деталь создаёт много трудностей при создании вычислительных пакетов.

Например, в последнее время наиболее популярными языками становятся кроссплатформенные и машинонезависимые. Это влечёт за собой то, что программист не имеет встроенных в язык средств низкоуровневого контроля своих программ.

Автор работы решил данную проблему следующим образом: была написана динамическая низкоуровневая <<.dll>> библиотека, содержащая все необходимые инструменты для использования на языке высокого уровня.

В данной библиотеке были реализованы следующие сущности:
\begin{enumerate}

    \item Контроллер режима округления
    
    \item Тестирующий модуль классической интервальной арифметики
    
    \item Классический интервал
    
    \item Псевдослучайный генератор чисел
    
\end{enumerate}

\subsection{Контроллер режима округления}

При разработке <<низкоуровневой>> библиотеки для работы с направленным округлением, следует иметь ввиду, что это весьма сужает список доступных языков программирования для её разработки. Данное обстоятельство связано с тем, что в угоду кросс-платформенности и машинонезависимости разработчики отказываются от предоставления доступа к низкоуровневым системным операциям. Поэтому часто доступ к таким операциям, как смена направления округления вещественных чисел закрыт.

Даже при использовании низкоуровневого языка необходимо специально настраивать компилятор для генерации кода, способного правильно работать с направленным округлением. Это требование необходимо, поскольку используя различные оптимизации и ускорения при генерации кода и вычислении выражений с плавающей точкой, компилятор способен искажать  итоговый результат.

Разработка пакета, представленного в работе, велась в интегрированной среде разработки (\texttt{IDE}) \texttt{Microsoft Visual Studio 2017} на языке \texttt{C++} стандарта 2011 года.

Для настройки компилятора необходимо включить специальную опцию --- \texttt{/fp:strict}. Она предназначена для установки специального режима работы компилятора с вещественными числами. Далее идёт перевод некоторых положений относительно этой опции из официальной документации \texttt{Microsoft} \cite{MicrosoftDocumentationCpp}.

По умолчанию компилятор использует режим \texttt{/fp:precise}. В режиме \texttt{/fp:precise} компилятор сохраняет порядок вычислений и режим округления вычисляемых выражений с плавающей точкой, когда он оптимизирует и создает объектный код. Компилятор производит округление вычислений в 4 случаях:

\begin{enumerate}
    \item Присваивания
    \item Приведение типов
    \item Передача аргумента с плавающей точкой в функцию
    \item Возврат результата с плавающей точкой из функции
\end{enumerate}

Промежуточные вычисления могут вычисляться с машинной точностью. Приведение типов может быть использовано для явного округления вычислений. В данном режиме компилятор не выполняет никаких математических преобразований над вычисляемыми выражениями, кроме как в тех случаях, когда гарантируется побитовая идентичность результатов. Выражения, содержащие специальные значения ($+\infty$, $-\infty$, NaN, $-0.0$) вычисляются в соответствии со спецификацией IEEE-754.

Компилятор генерирует код по умолчанию, предполагая, что среда с плавающей точкой не изменится во время выполнения. То есть предполагается, что код не перехватывает исключения с плавающей запятой, считывает или записывает регистры состояния с плавающей запятой и не изменяет режимы округления.

Включенный режим \texttt{/fp:fast} позволяет компилятору переупорядочивать, комбинировать, или каким-либо образом упрощать вычисляемые выражения с плавающей точкой для повышения быстродействия или уменьшения размера кода. Компилятор может пропустить округление в тех 4 случаях, перечисленных выше для режима \texttt{/fp:precise}. Из-за данного поведения результаты, полученные при этом режиме и других (\texttt{/fp:precise}, \texttt{/fp:strict}) могут заметно отличаться. Определенное поведение при работе со специальными значениями ($+\infty$, $-\infty$, NaN, $-0.0$) не указывается.

В данном режиме компилятор также предполагает инвариантность среды во время исполнения программы.

Компилятор в режиме \texttt{/fp:strict} ведёт себя похожим образом, как и в режиме \texttt{/fp:precise}. То есть он сохраняет порядок вычислений выражений, также округляет результаты вычислений в указанных 4 случаях, а также обрабатывает специальные значения вещественных переменных согласно стандарту IEEE-754. Основное отличие состоит в том, что в данном режиме программа может безопасно получить доступ к среде вычислений чисел с плавающей точкой, а также изменять ее. В данном режиме не делаются усечения вещественных чисел во время вычислений.

Но режим \texttt{/fp:strict} более трудоёмок, чем режим \texttt{/fp:precise}, так как компилятор вставляет дополнительный код для перехвата исключений при операциях с вещественными числами, а также код, получающий информацию о среде вычислений, или код, изменяющий эту среду. Если код не использует данные возможности, то стоит использовать режим \texttt{/fp:precise}.

В данном пакете, благодаря тому, что программно отделяется высокоуровневый и низкоуровневый контроль, задача контроллера режима округления состоит в том, чтобы до момента создания первого интервала для оперирования режим округления был изменён на направленное округление вверх.

\subsection{Тестирующий модуль интервальной арифметики}

При работе с направленным округлением сложно даже в режиме отладки следить за тем, чтобы округление чисел с плавающей точкой происходило верно. Поэтому было решено написать тестирующий модуль, который можно запускать для непосредственной проверки правильности округления. 

Данный модуль производит действия с интервалами по базовому принципу (\ref{eq:BaseIntervalPrincipe}). При этом, при поиске минимума, режим округления переключается на округление вниз, а при поиске максимума --- в режим округления вверх.

\subsection{Классический интервал}

Интервал классической интервальной арифметики является центральным объектом разрабатываемой низкоуровневой библиотеки. 

\subsubsection{Поля интервального объекта}

Интервал представляет в виде двух вещественных чисел-концов. В программной реализации для хранения этих чисел будет использоваться тип \texttt{double}.

Левый (или нижний) конец интервала интервального объекта назовём <<\texttt{LeftBound}>>, а правый (или верхний) --- <<\texttt{RightBound}>>.

\subsubsection{Создание интервального объекта}

Создать интервальный объект также можно будет по двум вещественным числам.

Опишем алгоритм создания интервального объекта по двум вещественным числам --- концам интервала.

\begin{center}

	\underline{\textit{Алгоритм} <<Constructor>>}

\end{center}

\noindent \underline{\textit{Вход:}}

\comment{левый конец создаваемого интервала}

\texttt{double leftBound;}

\comment{правый конец создаваемого интервала}

\texttt{double rightBound;} 

\noindent \underline{\textit{Псевдокод:}}

\comment{присваивание левого конца интервала}

\texttt{Result.LeftBound} $\leftarrow$ \texttt{leftBound;}

\comment{присваивание правого конца интервала}

\texttt{Result.RightBound} $\leftarrow$ \texttt{rightBound;}

\noindent \underline{\textit{Выход:}}

\comment{созданный интервальный объект $[ \, \texttt{leftBound}, \: \texttt{rightBound} \, ]$}

\texttt{Interval Result;}

Будет удобно также уметь создавать интервальные объекты из одного вещественного числа, опишем данный алгоритм.

\begin{center}

	\underline{\textit{Алгоритм} <<ConstructorByNumber>>}

\end{center}

\noindent \underline{\textit{Вход:}}

\comment{число для создания интервала}

\texttt{double number;} 

\noindent \underline{\textit{Псевдокод:}}

\texttt{Result $\leftarrow$ Constructor(number, number);}

\noindent \underline{\textit{Выход:}}

\comment{созданный интервальный объект $[ \, \text{number}, \: \text{number} \, ]$}

\texttt{Interval Result;}

\subsubsection{Обработка интервальных объектов}

В пункте ??? была приведена выжимка из документации языка \texttt{C++} для настройки работы компилятора с числами с плавающей точкой. Воспользуемся этим для реализации простейших операций для обработки интервальных объектов.

\subsubsection*{Простейшая програмнная реализация}
\addcontentsline{toc}{section}{Простейшая программная реализация}

Согласно базовому принципе интервального анализа, необходимо сохранить внутри интервала все результаты операций. Это можно достигнуть переключением регистра округления процессора.
По умолчанию, режим округления процессора стоит в режиме "до ближайшего". 

Поэтому, чтобы не потерять результат, при нахождении левого конца интервала необходимо переключить процессор в режим округления "округление вниз до ближайшего", а при нахождении правого конца интервала --- в режим "округление вверх до ближайшего".

Проиллюстрируем это описанием алгоритмов для сложения и вычитания интервальных объектов, согласно определению операций ???.

\begin{center}

	\underline{\textit{Алгоритм} <<Operator$+$>>}
	
\end{center}

\noindent \underline{\textit{Вход:}}

\comment{первый интервальный объект для сложения}

\texttt{Interval firstInterval;}

\comment{второй интервальный объект для сложения}

\texttt{Interval secondInterval;} 

\noindent \underline{\textit{Псевдокод:}}

\texttt{LastRoundRegistryState} $\leftarrow$

\hspace{1cm} \texttt{[Текущий режим округления процессора];}

\texttt{Перевод процессора в режим <<округление вниз до ближайшего>>;}

\texttt{ResultLeftBound} $\leftarrow$ \texttt{firstInterval.LeftBound} $+$

\hspace{1cm} \texttt{secondInterval.LeftBound};

\texttt{Перевод процессора в режим <<округление вверх до ближайшего>>}

\texttt{ResultRightBound} $\leftarrow$ \texttt{firstInterval.RightBound} $+$

\hspace{1cm} \texttt{secondInterval.RightBound};

\texttt{Перевод процесора в режим <<LastRoundRegistryState>>}

\texttt{Result} $\leftarrow$ \texttt{Constructor(ResultLeftBound, ResultRightBound)};

\noindent \underline{\textit{Выход:}}

\comment{результат сложения интервальных объектов}

\comment{firstInterval и secondInterval}

\texttt{Interval Result;}

\begin{center}

	\underline{\textit{Алгоритм} \texttt{Operator$-$}}
	
\end{center}

\noindent \underline{\textit{Вход:}}

\comment{уменьшаемый интервальный объект}

\texttt{Interval firstInterval;}

\comment{вычитаемый интервальный объект}

\texttt{Interval secondInterval;}

\noindent \underline{\textit{Псевдокод:}}

\texttt{LastRoundRegistryState} $\leftarrow$

\hspace{1cm} \texttt{[Текущий режим округления процессора];}

\texttt{Перевод процессора в режим <<округление вниз до ближайшего>>;}

\texttt{ResultLeftBound} $\leftarrow$

\hspace{1cm} \texttt{firstInterval.LeftBound - secondInterval.RightBound;}

\texttt{Перевод процессора в режим <<округление вверх до ближайшего>>;}

\texttt{ResultRightBound} $\leftarrow$

\hspace{1cm} \texttt{firstInterval.RightBound - secondInterval.LeftBound;}

\texttt{Перевод процессора в режим LastRoundRegistryState;}

\texttt{Result} $\leftarrow$ \texttt{Constructor(ResultLeftBound, ResultRightBound);}

\noindent \underline{\textit{Выход:}}

\comment{результат вычитания интервальных объектов}

\comment{firstInterval и secondInterval}

\texttt{Interval Result;}

\subsubsection*{Продвинутая програмнная реализация}
\addcontentsline{toc}{section}{Продвинутая программная реализация}

Пользуясь тем свойством, что $\text{max} \big\{ f(x) \big\} = - \text{min} \big\{ -f(x) \big\}$, можно убрать одну из инструкций смены режима округления процессора. Приведём далее усовершенствованные алгоритмы, по умолчанию считая, что для работы будем использовать режим <<кругление вверх до ближайшего>>.

\begin{center}

	\underline{\textit{Алгоритм} <<Operator$+$>>}
	
\end{center}

\noindent \underline{\textit{Вход:}}

\comment{первый интервальный объект для сложения}

\texttt{Interval firstInterval;}

\comment{второй интервальный объект для сложения}

\texttt{Interval secondInterval;} 

\noindent \underline{\textit{Псевдокод:}}

\texttt{LastRoundRegistryState} $\leftarrow$

\hspace{1cm} \texttt{[Текущий режим округления процессора];}

\texttt{Перевод процессора в режим <<округление вверх до ближайшего>>;}

\texttt{ResultLeftBound} $\leftarrow -(-\texttt{firstInterval.LeftBound} -$

\hspace{1cm} $\texttt{secondInterval.LeftBound});$

\texttt{ResultRightBound} $\leftarrow$ \texttt{firstInterval.RightBound} $+$

\hspace{1cm} \texttt{secondInterval.RightBound};

\texttt{Перевод процесора в режим <<LastRoundRegistryState>>}

\texttt{Result} $\leftarrow$ \texttt{Constructor(ResultLeftBound, ResultRightBound)};

\noindent \underline{\textit{Выход:}}

\comment{результат сложения интервальных объектов}

\comment{firstInterval и secondInterval}

\texttt{Interval Result;}

\begin{center}

	\underline{\textit{Алгоритм} \texttt{Operator$-$}}
	
\end{center}

\noindent \underline{\textit{Вход:}}

\comment{уменьшаемый интервальный объект}

\texttt{Interval firstInterval;}

\comment{вычитаемый интервальный объект}

\texttt{Interval secondInterval;}

\noindent \underline{\textit{Псевдокод:}}

\texttt{LastRoundRegistryState} $\leftarrow$

\hspace{1cm} \texttt{[Текущий режим округления процессора];}

\texttt{Перевод процессора в режим <<округление вверх до ближайшего>>;}

\texttt{ResultLeftBound} $\leftarrow$

\hspace{1cm} $-(-\texttt{firstInterval.LeftBound + secondInterval.RightBound);}$

\texttt{ResultRightBound} $\leftarrow$

\hspace{1cm} \texttt{firstInterval.RightBound - secondInterval.LeftBound;}

\texttt{Перевод процессора в режим LastRoundRegistryState;}

\texttt{Result} $\leftarrow$ \texttt{Constructor(ResultLeftBound, ResultRightBound);}

\noindent \underline{\textit{Выход:}}

\comment{результат вычитания интервальных объектов}

\comment{firstInterval и secondInterval}

\texttt{Interval Result;}

Тогда естественным решением

Для обработки были разработаны следующие процедуры:
\begin{enumerate}
	\item Сложение интервалов
	\item Вычитание интервалов
	\item Умножение интервалов
	\item Деление интервалов
	\item Возведение интервала в квадрат
	\item Извлечение квадратного корня из интервала
	\item Алгебраический минус из интервала
	\item Взятие модуля интервала
	\item Нахождение ширины интервала
	\item Нахождение радиуса интервала
	\item Нахождение середины интервала
	\item Нахождение объединения интервалов
	\item Нахождение, пересекаются ли интервалы
	\item Нахождение пересечения интервалов
	\item Взятие левого конца интервала
	\item Взятие правого конца интервала
	\item Нахождение, является ли интервал нуль-содержащим
\end{enumerate}

\subsection{Псевдослучайный генератор чисел}

Для различных целей в пакете реализован псевдослучайный генератор чисел. Он имеет следующие функции:

\begin{enumerate}
	\item Генерация вещественного числа, равномерно распределенного в отрезке.
	\item Генерация целого числа, равномерно распределенного в отрезке.
	\item Генерация вещественного числа, имеющего нормальное распределение с заданными параметрами $\alpha$, $\mu$.
\end{enumerate}

Был использован язык \texttt{C\#,} на котором был написан <<класс-обертка>>. Он импортирует все функции из низкоуровневой библиотеки и отвечает за исполнение <<неуправляемого>> (\texttt{unmanaged}) кода. Данный инструментарий обеспечивает возможности для быстрого и безопасного исполнения пакета. Он позволяет решить следующие проблемы:

\begin{enumerate}
    \item Автоматическое слежение среды исполнения управляемого кода за выделенной в куче памятью с помощью <<сборщика мусора>> (\texttt{Garbage Collector}).
    \item Слежение за режимом округления из неуправляемого кода, что позволяет разделить логику работы с вещественными числами
    \item Использование \texttt{LINQ}-запросов, которые могут использовать удобный встроенный инструментарий, например встроенные функции сортировки и т. д.
    \item Использование \texttt{Windows}-специфичных инструментов для создания графических файлов
    \item Встроенная потоковая параллелизация
    \item Программное разделение уровней абстракций низкого и высокого уровня
\end{enumerate}

\clearpage
\section* {Заключение}
\addcontentsline{toc}{section}{Заключение}

\addcontentsline{toc}{subsection}{Выводы}
\subsection* {Выводы}

В ходе написания работы

\begin{enumerate}

	\item Была выдвинута концепция для создания новой функционально-граничной арифметики, которая улучшает свойства классической интервальной и аффинной африметик.

	\item Были введены арифметические операции между интервалами $\mathbb{F}\mathbb{R}$.
	
	\item Была введена интегральная характеристика для оценки качества приближения диаграмм зависимости.	
	
	\item Была написана библиотека для работы с направленным округлением.
	
	\item Была написана программа оперирования классическими и функционально-граничными интервалами.
	
	\item Проведено сравнение классической, функционально-граничной интервальных и аффинной арифметик по интегральной характеристике в выражения с одной  двумя переменными. По итогам сравнения показано превосходство функционально-граничной интервальной арифметики.
	
	\item Функционально-граничная арифметика была применена для нахождения оценок объединённых множеств решений различных ИСЛАУ со связями. В результате было показано превосходство функционально-граничной интервальной арифметики над классической и аффинной.

\end{enumerate}

\addcontentsline{toc}{subsection}{Дальнейшие планы работы}
\subsection* {Дальнейшие планы работы}

Написанная работа открыла перспективы для дальнейшего изучения и развития функционально-граничной арфиметики:

\begin{enumerate}

	\item Получение аппроксимации квадратичных членов при операции умножения за время работы $O(n)$, где $n$ --- количество интервальных параметров.
	
	\item Внедрение в базис квадратичных членов $x_{i} ^ {2}$.

	\item Получение априорных теоретических оценок для интегральной характеристики.
	
	\item Повышение эффективности для операции деления.
	
	\item Оптимизация процедуры $F_{\downarrow}$ и $F_{\uparrow}$ для более эффективного приближения объединения интервалов $\mathbb{F}\mathbb{R}$.
	
	\item Улучшение архитектуры программы и повышение её быстродействия.
	
	\item Введение других операций, таких как: взятие модуля, возведение в квадрат, взятие квадратного корня, и другие.

\end{enumerate}


\end{doublespace}
\end{document}